\documentclass[12pt,a4paper,fleqn]{report}
\usepackage[a4paper,left=3cm,right=2cm,top=3cm,bottom=2cm]{geometry}
\usepackage[active]{srcltx}
\usepackage[brazilian]{babel}
\usepackage[utf8]{inputenc}
\usepackage[T1]{fontenc}
\usepackage{indentfirst}
\usepackage{verbatim}
\usepackage{setspace} 
\usepackage{amsfonts}
\usepackage{amssymb, amsthm}
\usepackage{amsmath}
\usepackage{amssymb}
\usepackage{array}
\usepackage{enumerate}
\usepackage{enumitem}
\usepackage{color,graphicx}
\usepackage{graphics}
\usepackage{float}
\usepackage[backref=page]{hyperref}
\usepackage{bookmark}
\usepackage[nottoc]{tocbibind} 
\usepackage[all]{xy}
\usepackage{tikz}
\usepackage{faktor}
\usepackage{nicefrac}
\usepackage{tasks}
\usepackage{multicol}
\usepackage{mathtools}
\usepackage{makeidx}
\usepackage{tablefootnote}
\usepackage[square,numbers]{natbib}
\usepackage[small]{caption} 
\newcommand{\nocontentsline}[3]{}
\newcommand{\tocless}[2]{\bgroup\let\addcontentsline=\nocontentsline#1{#2}\egroup}

\newcommand{\C}{\mathbb C}

\newcommand{\W}{\ensuremath{\mathbb W}}
\newcommand{\T}{\mathbb T}

\newcommand{\Hilb}{\ensuremath{Hilb}}

\makeatletter
\newcommand{\mathleft}{\@fleqntrue\@mathmargin0pt}
\newcommand{\mathcenter}{\@fleqnfalse}
\makeatother

\def \P#1{\mathbb{P}^#1} 
\def \bb#1{\ensuremath{\mathbb{#1}}} 
\def \p#1{\ensuremath{\mathbb P^{#1}}} 
\def \cl#1{\ensuremath{\mathcal{#1}}}

\newtheorem{thm}{Teorema}[chapter]
\def \bthm{\begin{thm}} \def \ethm{\end{thm}}

\newtheorem{cor}[thm]{Corol\'{a}rio}
\def \bcor{\begin{cor}} \def \ecor{\end{cor}}

\newtheorem{lem}[thm]{Lema}
\def \blem{\begin{lem}} \def \elem{\end{lem}}

\newtheorem{prop}[thm]{Proposi\c{c}\~{a}o}
\def\bprop{\begin{prop}}\def\eprop{\end{prop}}

\newtheorem{axi}[thm]{Axioma}
\def \baxi{\begin{axi}} \def \eaxi{\end{axi}}

\newtheorem{defn}[thm]{Defini\c{c}\~{a}o}
\def \bdefn{\begin{defn}} \def \edefn{\end{defn}}

\newtheorem{exe}[thm]{Exemplo}
\def \bexe{\begin{exe}} \def \eexe{\end{exe}} 

\newtheorem{obs}[thm]{Observa\c{c}\~{a}o}
\def \bobs{\begin{obs}}\def\eobs{\end{obs}}

\def\bdem{\begin{proof}}\def\edem{\end{proof}}

\def\beq{\begin{equation}}\def\eeq{\end{equation}}

\def \segre{\mbox{Segre}}

\makeindex
\def\bi{\begin{itemize}}\def\ei{\end{itemize}}
\def\ba#1{\begin{array}{#1}}\def\ea{\end{array}}
\def\na#1{\noalign{\vskip#1pt}}
\def\lb{\linebreak}
\def\ra{\ensuremath{\rightarrow}}
\def\lar{\ensuremath{\longrightarrow}}
\def\ov#1{\ensuremath{\overline{#1}}}

\def\mb#1{\ensuremath{\mathbf{#1}}}

\def\pol{po\/li\/nô\/mio}
\def\nin{\ensuremath{\not\in}}
\def\wed#1{\ensuremath{\stackrel{#1}\wedge}}
\makeatletter

\@addtoreset{equation}{chapter}

\makeatother

\begin{document}
\pagestyle{empty}

\begin{titlepage}


\begin{center}
{\LARGE Universidade Federal de Minas Gerais}

{\Large Instituto de Ciências Exatas}

{\large Departamento de Matemática}

\vspace{100pt}
{\LARGE Weversson Dalmaso Sellin}
\par
\vspace{100pt}
{\Huge Enumeração de hipersuperfícies com subesquemas singulares}
\par
\vfill
\textbf{{\large Belo Horizonte - MG}\\
{\large \the\year}}
\end{center}
\end{titlepage}




\thispagestyle{empty}

\begin{center}
{\LARGE Weversson Dalmaso Sellin}
\par
\vspace{100pt}
{\Huge Enumeração de hipersuperfícies com subesquemas singulares}
\end{center}
\par
\vspace{100pt}
\noindent\hspace*{0.48\linewidth}\begin{minipage}{0.52\linewidth}
{\large Tese apresentada ao Programa de Pós-graduação em Matemática do Instituto de Ciências Exatas da Universidade Federal de Minas Gerais, como parte dos requisitos para a obtenção do grau de doutor em Matemática.
\par 
\vspace{1em}
Área de concentração: Geometria Algébrica}
\end{minipage}

\par
\vspace{1em}
\noindent\hspace*{0.48\linewidth}\parbox{0.52\linewidth}{{\large Orientador: Israel Vainsencher}}

\par
\vfill
\begin{center}
\textbf{{\large Belo Horizonte - MG}\\
{\large \the\year}}
\end{center}

\newpage

\begin{center}\textbf{\Huge Dedicatória}\end{center}
\vspace{1cm}


\begin{figure}[ht]
\centering
\includegraphics[scale=0.6]{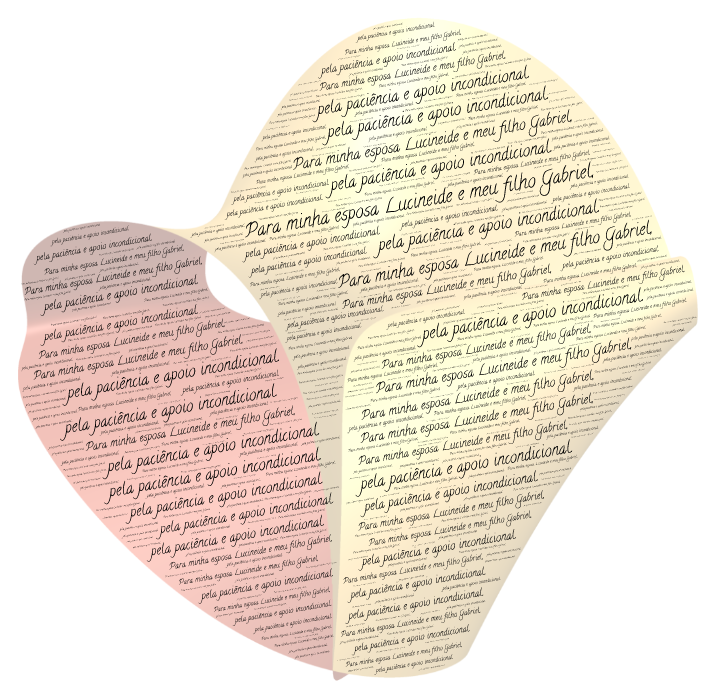}
\end{figure}

\newpage

\vspace*{0.4\textheight}
{\doublespacing
Eu vejo a barra do dia surgindo

Pedindo pra gente cantar

Eu tenho tanta alegria, adiada

Abafada, quem dera gritar

...

\begin{flushright}
Chico Buarque
\end{flushright}

\newpage


\begin{center}\textbf{\Huge Agradecimentos}\end{center}
\vspace{1cm}

{\onehalfspacing
 Já dizia o escritor Machado de Assis: "A gratidão de quem recebe um benefício é sempre menor que o prazer daquele de quem o faz". 
 
 No intuito de maximizar o prazer àqueles que algum benefício me proporcionou durante a minha trajetória é que  venho agradecê-los.
 
 Aos meus pais, pela vida e pelo amor incondicional  que nos momentos de sucesso pode até parecer irrelevante, mas que nas ocasiões de fracasso se torna um consolo e uma segurança que não se encontram em qualquer outro lugar.
 
 A minha esposa Lucineide e meu filho Gabriel, pelo apoio, paciência e por serem meu porto seguro nessa jornada.
 
 Aos amigos que ganhei no decorrer da vida e que tanto me enriqueceu ao revelarem cada vez mais de mim. Não vou nomeá-los para evitar  o risco de pecar por omissão.
 
 Aos mestres professores que tive durante a trajetória acadêmica e que tanto contribuíram para minha formação, não só como pesquisador na área de Matemática, mas acima de tudo, como um ser humano melhor. Em especial, ao meu orientador Israel Vainsencher, pois sem sua ajuda essa tese não se tornaria realidade.
 
 Aos colegas do DCEX/UFVJM pela confiança e apoio.
 
 Por fim, à Universidade Federal dos Vales do Jequitinhonha e Mucuri, pelo apoio aos meus estudos.}

\newpage

\vspace*{10pt}
\begin{center}\textbf{\Huge Resumo}\end{center}
\vspace{1cm}

\noindent{\onehalfspacing
Nesta tese estuda-se o lugar das hipersuperfícies
com singularidades não\break
isoladas. Especificamente, fixada uma subvariedade
irredutível fechada de um esquema de Hilbert, $\W
\subset \Hilb_{p_W(t)}(\P{n})$,
considera-se a subvariedade $\Sigma(\W, d) \subset
\mathbb{P}^{N_d} =
\mathbb{P}(H^0(\mathcal{O}_{\P{n}}(d)))$, formada pelas
hipersuperfícies de grau $d$ em $\P{n}$, que são
singulares em algum elemento (variável) $W\in
\W$. Sob a condição de
que um membro geral $W\in\W$ é suave, irredutível e de dimensão
positiva, mostra-se que o grau de $\Sigma(\W, d)$
se expressa por um polinômio $p^\W(d)$  para $d\gg0$.
O polinômio $p^\W(d)$ é calculado explicitamente
para certas famílias $\W$, distinguidas pela
existência de literatura com descrição
adequada.  Notadamente, estudam-se os casos

\medskip

$\begin{array}{lll}
  \bullet\, \W_{(k+1,n+1)}&:= &\{\,\P{k} - \mbox{
  linear em }\P{n}\,\} ,\\\na6
\bullet\,
  \W_m&:=&\{\mbox{  curvas planas de grau $m$ em
  \p3} \} ,\\\na6
\bullet\,  \W_{twc}&:=& \{\mbox{ 
curvas cúbicas reversas em \p3 }\} ,\\\na6
\bullet\,  \W_{rc}&:=&\{ \mbox{
superfícies cúbicas regradas em \p4 }\} ,\\\na6
\bullet\,
\W_{sg} &:=&\{ \mbox{ 3-variedades de Segre em \p5
  }
  \} \ \rm e\\\na6
  \bullet\,  \W_{eqc}&:=&\{ \mbox{
 curvas quárticas elípticas
   em \p3  }\} .
\end{array}$

\medskip

\noindent
O método consiste em descrever uma
dessingularização $\W'\rightarrow\W$ que\lb
parametriza uma família plana de subesquemas de 
\p{n} cujo membro geral é definido por um ideal da
forma $I_W^2$, quadrado do ideal de um membro
geral $W\in\W$.  A variedade $\W'$ é munida de
um subfibrado $\mathcal{E}_d$ do fibrado trivial
$\W'\times{}H^0(\mathcal{O}_{\P{n}}(d))$ para $d\gg0$,
cuja fibra sobre um membro geral $W'\in\W'$ é
formada pelas $F\in{}H^0(\mathcal{O}_{\P{n}}(d))$
com gradiente nulo ao longo de $W$. Além disso, o
mapa induzido na projetivização
$\mathbb{P}(\mathcal{E}_d) \rightarrow
\mathbb{P}^{N_d}$ tem imagem a \lb
variedade
$\Sigma({\W,d})$ e é genericamente injetivo para
$d\gg0$. Polinomialidade do grau de
$\Sigma({\W,d})$ é provada
usando Grothendieck-Riemann-Roch.
Nos casos acima destacados,
usa-se a fórmula de resíduos de Bott para o
cálculo explícito do grau de $\Sigma({\W,d})$.

\par
\vspace{1em}
\noindent\textbf{Palavras-chave:} Geometria
enumerativa, fórmula de resíduos de Bott,\lb
singularidades.}
\newpage

\begin{center}\textbf{\Huge Abstract}\end{center}

\vspace{1cm}

\noindent{\onehalfspacing
This thesis investigates the locus of
hypersurfaces with nonisolated singularities.
More precisely, given a closed, irreducibe subvariety
of a Hilbert scheme, $\W\subset \Hilb_{p_W(t)}(\P{n})$, 
we define a subvariety $\Sigma(\W, d) \subset
\mathbb{P}^{N_d} =
\mathbb{P}(H^0(\mathcal{O}_{\P{n}}(d)))$, formed
by the hypersurfaces of degree $d$ in $\P{n}$ 
which are singular along some (variable) member
$W\in \W$. Assuming that a general member
$W\in\W$ is smooth, irreducible and  positive 
dimensional, we show that the degree of $\Sigma(\W,d)$
is  expressed by a  polinomial  $p^\W(d)$  for all
$d\gg0$. The polynomial $p^\W(d)$ is made explicit 
for a few families $\W$, distinguished by
the existence of an adequate description in the
literature.  Notably, we study the cases

\medskip

$\begin{array}{lll}
  \bullet\, \W_{(k+1,n+1)}&:= &\{\,\P{k} - \mbox{
  linear in }\P{n}\,\} ,\\\na6
\bullet\,
  \W_m&:=&\{\mbox{  plane curves of degree $m$ in
  \p3} \} ,\\\na6
\bullet\,  \W_{twc}&:=& \{\mbox{
    twisted cubic curves in \p3 }\} ,\\\na6
\bullet\,  \W_{rc}&:=&\{ \mbox{  ruled cubic
  surfaces in  \p4 }\} ,\\\na6
\bullet\,
\W_{sg} &:=&\{ \mbox{  Segre cubic 3-folds in \p5
  }
  \} \ \rm and\\\na6
  \bullet\,  \W_{eqc}&:=&\{ \mbox{
    elliptic quartic curves in \p3 }\} .
\end{array}$

\medskip

\noindent
The method  consists in describing
a desingularization $\W'\rightarrow\W$ such that $\W'$
parameterizes a flat family of subschemes of \p n the
general member of which is defined by an ideal of
the form $I_W^2$, square of the ideal of a general
member $W\in\W$. The variety $\W'$ comes equipped, for
$d\gg0$, with a vector subbundle $\mathcal{E}_d$ of the
trivial bundle $\W'\times H^0(\mathcal{O}_{\P{n}}(d))$,
with  fiber over a general member $W'\in\W'$ formed
by the $F\in{}H^0(\mathcal{O}_{\P{n}}(d))$ such that
the gradient vanishes along $W$. Moreover, the map
$\mathbb{P}(\mathcal{E}_d)\rightarrow
\mathbb{P}^{N_d}$
induced in the projectivization has image the variety
$\Sigma({\W,d})$ and is generically injective for
$d\gg0$.  Polynomiality follows using
Grothendieck-Riemann-Roch. In the  cases
above displayed  Bott's localization at
fixed points is
employed to derive 
explicit formula for the degree of $\Sigma({\W,d})$.

\par
\vspace{1em}
\noindent\textbf{Keywords:}
Enumerative Geometry, Bott’s residues formula,
singularities.}

\newpage



\begin{center}\textbf{\Huge Notações etc}\end{center}
\vspace{1cm}
{\onehalfspacing

\noindent$K$:= corpo algebricamente fechado, em
geral  de característica zero.

\noindent$K[x_0,\dots, x_n]$:= anel de polinômios nas variáveis $x_0, \dots, x_n$.

\noindent$\mathcal{F}_d=
H^0(\mathcal{O}_{\P3}(d)):=$ fibrado vetorial
trivial com fibra o espaço de polinômios
homogêneos de grau $d$ nas variáveis $x_0, \dots, x_n$.

\noindent$S_i(\_):=$ i-ésima potência simétrica.

\noindent$\wed i \_ :=$ i-ésima potência exterior.

\noindent $\Hilb_{p_W(t)}(\P{n}):=$ o esquema de Hilbert de subesquemas de $\P{n}$ com polinômio de Hilbert fixado $p_W(t)$

\noindent $\mathbb{P}^{N_d} = \mathbb{P}(H^0(\mathcal{O}_{\P{n}}(d))):=$ o espaço de parâmetros de hipersuperfícies de grau $d$ em $\P{n}$

\noindent $\mathcal{N}_{Y/X}:=$ fibrado normal

\noindent $\mathcal{T}\mathbb{X}:=$ fibrado tangente

}

\newpage

\thispagestyle{empty}
\addtocontents{toc}{\protect\thispagestyle{empty}}
\addtocontents{toc}{\protect\sloppy}
\tableofcontents
\newpage
\pagestyle{plain} 

\onehalfspacing 
\chapter*{Introdução}
\addcontentsline{toc}{chapter}{Introdução}

O objetivo deste trabalho é estudar o lugar das
hipersuperfícies que apresentam certas
singularidades não isoladas.

Especificamente, seja $\Hilb_{p_\W(t)}(\P{n})$ o esquema de Hilbert de subesquemas de $\P{n}$ com polinômio de Hilbert fixado $p_\W(t)$ e seja $\W$ uma subvariedade irredutível fechada do $\Hilb_{p_\W(t)}(\P{n})$ cujo membro geral estamos admitindo ser suave. Denotamos por $\Sigma(\W{}, d) \subset \mathbb{P}^{N_d} = \mathbb{P}(H^0(\mathcal{O}_{\P{n}}(d)))$, o lugar das hipersuperfícies de grau $d$ em $\P{n}$ que são singulares em algum elemento (variável) $W\in \W$. O propósito aqui é determinar o grau de $\Sigma(\W{}, d)$ para certas famílias $\W$, distinguidas pela existência de literatura com descrição adequada. Notadamente, estudam-se os casos

\medskip

$\begin{array}{lll}
  \bullet\, \W_{(k+1,n+1)}&:= &\{\,\P{k} - \mbox{
  linear em }\P{n}\,\} ,\\\na6
\bullet\,
  \W_m&:=&\{\mbox{  curvas planas de grau $m$ em
  \p3} \} ,\\\na6
\bullet\,  \W_{twc}&:=& \{\mbox{ 
curvas cúbicas reversas em \p3 }\} ,\\\na6
\bullet\,  \W_{rc}&:=&\{ \mbox{
superfícies cúbicas regradas em \p4 }\} ,\\\na6
\bullet\,
\W_{sg} &:=&\{ \mbox{ 3-variedades de Segre em \p5
  }
  \} \ \rm e\\\na6
  \bullet\,  \W_{eqc}&:=&\{ \mbox{
 curvas quárticas elípticas
   em \p3  }\} .
\end{array}$

\medskip

Em cada caso, descreve-se uma
dessingularização $\W'\rightarrow\W$ que\lb
parametriza uma família plana de subesquemas de 
\p{n} cujo membro geral é definido por um ideal da
forma $I_W^2$, quadrado do ideal de um membro
geral $W\in\W$.  A variedade $\W'$ é munida de
um subfibrado $\mathcal{E}_d$ do fibrado trivial
$\W'\times{}H^0(\mathcal{O}_{\P{n}}(d))$ para $d\gg0$,
cuja fibra sobre um membro geral $W'\in\W'$ é
formada pelas $F\in{}H^0(\mathcal{O}_{\P{n}}(d))$
com gradiente nulo ao longo de $W$. Além disso, o
mapa induzido na projetivização
$\mathbb{P}(\mathcal{E}_d) \rightarrow
\mathbb{P}^{N_d}$ tem imagem a \lb
variedade
$\Sigma({\W,d})$ e é genericamente injetivo para
$d\gg0$ (vide Proposição \ref{prop_fibrado} \ e
Teorema \ref{teorema_injetividade_generica}). E
então, seguimos a estratégia de usar a fórmula de resíduos de Bott (vide Teorema \ref{formula_Bott}),
como explicado em Ellingsrud \& Str{\o}mme \cite{Ellingsrud_Stromme_Bott_96}, para o cálculo
do grau de $\Sigma({\W{},d})$. 

A utilização da fórmula de resíduos de Bott é
justificada pelo fato que nos problemas tratados
neste trabalho, em geral, precisamos efetuar a
explosão de alguma variedade $X$ ao longo de uma
subvariedade $Y$ (não singular) e considerar o
divisor excepcional, o qual é a projetivização do
fibrado normal,
$\mathbb{P}(\mathcal{N}_{Y/X})$. No entanto, na
maioria dos problemas considerados, conhecemos os
sucessivos fibrados normais somente fibra a fibra
e isso não é suficiente para fazermos cálculos
enumerativos usando a teoria de interseção
clássica. Nessa conjuntura é que a fórmula de
resíduos de Bott torna-se uma ferramenta
indispensável, pois requer apenas certos dados das
fibras do fibrado normal $\mathcal{N}_{Y/X}$ sobre
os pontos fixos por uma ação de um toro (vide
seção \ref{secao_bott} para maiores
detalhes). Observação análoga se aplica aos
fibrados $\cl E_d$ acima mencionados.

É conhecido que o grau de $\Sigma(\W{},d)$,
para $d\gg 0$, é um polinômio em $d$ cujo grau
$\deg_\W$  é
menor do que ou igual a $n\times \dim \W$ (vide
Proposição \ref{grau_polinomio}). No entanto, em
todos os casos efetivamente tratados, encontramos
$\deg_\W =(k+1)\times \dim(\W)$, onde $k =
\dim(W)$ para $W\in \W$. A validade para $\W$
qualquer resta conjectural. No caso de
singularidades isoladas Vainsencher
\cite{Vainsencher_2013} obtém fórmulas polinomiais
explícitas para o número de hipersuperfícies em
uma família $k$-dimensional adequada apresentando
$k \leq 6$ pontos duplos ordinários, cujo grau
obedece a fórmula acima. {Ainda no caso de
singularidades isoladas para curvas $\delta$-nodal, isto é, curvas que apresentam $\delta$ nós e nenhuma outra singularidade, temos o resultado conhecido como conjectura de Göttsche, a qual possui várias provas (cf. \cite{Kazarian_03}, \cite{kool2011}, \cite{liu2000} e \cite{tzeng2012}). A conjectura de Göttsche afirma que para um fibrado em retas $L$ suficientemente amplo sobre uma superfície $S$, o número de curvas $\delta$-nodal em um sistema linear geral $\delta$-dimensional é dado por um polinômio universal de grau $\delta$ em 4 variáveis, independentes de $S$ e $L$, $N_\delta (c_1(L)^2, c_1(L)\cdot c_1(S), c_1(S)^2, c_2(S))$. Em Rennemo
\cite{Rennemo_2017}, encontra-se uma generalização da conjectura de Göttsche, onde é permitido outros tipos de singularidades isoladas, além de nós, para o caso de curvas (cf. \cite{Rennemo_2017}, Prop. 7.2, pág. 4). O mesmo método utilizado para o caso de curvas com singularidas isoladas aplica-se para hipersuperfícies em dimensão arbitrária (cf. \cite{Rennemo_2017}, Prop. 7.8, pág. 4). Cabe ressaltar que o método utilizado na demonstração dos resultados em \cite{Rennemo_2017} não é construtivo e não apresenta um caminho aparente de transformar a prova em um algoritmo para obtenção de informação sobre o polinômio em questão.}

No que segue faremos uma breve descrição dos conteúdos dos respectivos capítulos e apêndices que formam o corpo desta tese.
No Capítulo \ref{cap_pl} estabelecemos os resultados gerais que são necessários ao desenvolvimento e entendimento dos demais capítulos. Especificamente, na seção \ref{pot_s_o} tratamos da relação entre potência simbólica e potência ordinária de um ideal  e como essa relação se reflete na condição de uma hipersuperfície ser singular em algum subesquema de $\P{n}$. Verifica-se que uma hipersuperfície $F$ de $\mathbb{P}^n$ é singular ao longo de $W\in\W$ (pelo menos nos casos considerados aqui) se, e somente se, ela contém a primeira vizinhança infinitesimal de $W$ em $\mathbb{P}^n$, isto é, o subesquema fechado de $\mathbb{P}^n$ tendo $(\mathcal{I}_W)^2$ como o feixe ideal. Na seção \ref{secao_fibrado_polinomio}\, é estabelecida a existência do subfibrado $\mathcal{E}_d$ do
fibrado trivial $\W{'}\times H^0(\p n,\cl O(d))$,  cuja
projetivização $\mathbb{P}(\mathcal{E}_d)$ produz
um mapa genericamente injetivo (ao menos para $d>>0$)
$$\mathbb{P}(\mathcal{E}_d)
\rightarrow \mathbb{P}^{N_d}$$ com imagem o
conjunto $\Sigma({\W{},d})$ (vide Proposição
\ref{prop_fibrado} e Teorema
\ref{teorema_injetividade_generica}).
Aqui $\W'\ra\W$ denota uma dessingularização,
composta de explosões com centros lisos e tal que
$\W'$ parametriza uma família plana de subesquemas
com membro geral dado pelo quadrado do ideal de um
membro geral de \W.
E nessa seção também é enunciado o resultado que garante o comportamento polinomial do grau de $\Sigma({\W{},d})$ em função de $d$ (vide Proposição \ref{grau_polinomio}). Na seção \ref{secao_bott} apresentamos a principal ferramenta para o cálculo do grau de $\Sigma({\W{},d})$, a saber a fórmula de resíduos de Bott (vide Teorema \ref{formula_Bott}).

No Capítulo \ref{cap_pklinear} tratamos do caso de hipersuperfícies em $\P{n}$ singulares ao longo de um $\P{k}$-linear, $k < n-1$. Aqui, o espaço de parâmetros consiste na grassmanniana $\mathbb{G}(k+1,n+1)$. E determinamos o grau da família de hipersuperfícies de grau $d$ singulares em algum membro de $\mathbb{G}(k+1,n+1)$, o qual vamos denotar por $\deg \Sigma(\W{}_{(k+1,n+1)},d)$, para os seguintes valores de $(k+1,n+1)$ explicitados nas fórmulas \ref{grau_reta}, \ref{grau_p2_p4}, \ref{grau_p2_p5} e \ref{grau_p3_p5}.

No Capítulo \ref{cap_cplane} consideramos o caso de curvas planas em $\P3$, onde obtemos resultados explícitos para o grau da família de superfícies de grau $d$ singulares  ao longo de alguma curva plana (variável) e com grau $m$ fixado, denotada por $\Sigma(\W_m{}, d)$. Explicitamente, obtemos fórmulas polinomiais que expressam o grau de $\Sigma(\W_m{}, d)$ em função de $d$ para os casos $m=2, 3$ (vide as fórmulas \ref{grau_conica_1} e \ref{grau_cubica}). Indicamos o Apêndice \ref{ap_curva_plana} para o
leitor  interessado em obter o $\deg
\Sigma(\W_m{}, d)$ para outros valores de $m$. 

No Capítulo \ref{cap_QD} tratamos dos casos de
hipersuperfícies em $\P{n}$ ($n=3,4,5$) singulares
ao longo do lugar de base
de uma rede de quádricas (variável) do tipo determinantal, isto é, gerada pelos menores
$2\times 2$ de uma matriz $3 \times 2$ de formas lineares. Especificamente, na seção \ref{sec_QD_twc} é considerada a família $\W_{twc}$ de cúbicas reversa em $\P3$. Uma cúbica reversa  
é uma curva racional, suave de
grau 3 em $\mathbb{P}^3$; seu polinômio
de Hilbert é $3t+1$. Piene \& Schlessinger
\cite{Piene_Schlessinger_1985} mostraram que a
componente $\W_{twc}$ de cúbicas reversa do esquema
de Hilbert $\Hilb_{3t+1}(\mathbb{P}^3)$ é uma
variedade projetiva suave de dimensão
12. Posteriormente, Ellingsrud, Piene \& Str{\o}mme
\cite{Ellingsrud_Piene_Stromme_87} provam que a
subvariedade $\bb X$ da Grassmanniana $G(3,10)$
formada pelas redes de quádricas do tipo
determinantal,  é uma variedade suave e que a componente
$\W_{twc}$ é a explosão de $\bb X$ ao longo da subvariedade de
redes com uma componente fixa. Ellingsrud \&
Str{\o}mme \cite{Ellingsrud_Stromme_89} também
mostraram que $\bb X$ é um quociente geométrico da
variedade de
matrizes de formas lineares $2\times 3$
semiestáveis. Essa descrição permite então o
cálculo do Anel de Chow de $\bb X$ e $\W_{twc}$. Vainsencher
\& Xavier \cite{Vainsencher_Xavier_02} apresentam
uma compactificação suave explícita de um espaço
de parâmetros para a família de cúbicas reversa
adequada para a aplicação da fórmula de resíduos
de Bott. Aqui utilizamos de início a 
descrição dos pontos fixos e respectivos tangentes explicitados em
Ellingsrud \& Str{\o}mme
\cite{Ellingsrud_Stromme_Bott_96}. Entretanto,
a família definida por $(\mathcal{I}_W)^2$ para
algum  $W\in \W_{twc}$ não é plana (flat).
Faz-se  necessário efetuar
uma explosão  $\W'_{twc}\rightarrow\W_{twc}$   ao
 longo de uma certa subvariedade suave
 para obtermos uma família plana, para a qual 
 obtemos de fato o fibrado $\mathcal{E}_d$ (vide subseção \ref{fibras_twc}). Na subseção \ref{espaco_parametro_twc} fazemos uma breve descrição do espaço de parâmetros $\W_{twc}$ baseada em Ellingsrud, Piene \& Str{\o}mme \cite{Ellingsrud_Piene_Stromme_87} e Vainsencher \cite{Vainsencher_1987}. Na subseção \ref{aplicando_bott_twc} damos a descrição dos pontos fixos oriundos da ação induzida do toro $\mathbb{T} = \mathbb{C}^*$, bem como as descrições do espaço tangente e das fibras de $\mathcal{E}_d$  sobre cada ponto fixo, necessárias à aplicação da fórmula de resíduos de Bott. Para o grau da família de superfícies singulares ao longo de alguma cúbica reversa (variável), denotada por $\Sigma({\W_{twc}{},d})$, obtemos a fórmula polinomial explicitada em \ref{grau_twc}.

Na seção \ref{sec_QD_rc} tomamos como espaço de
parâmetros a família de superfícies cúbicas
regradas em $\P4$, $\W_{rc}$, e obtemos uma
fórmula polinomial que expressa o grau da família
de hipersuperfícies de grau $d$ singulares em
algum membro $W\in \W_{rc}$ (variável), denotada
por $\Sigma(\W_{rc}{},d)$. Inicialmente, na
subseção \ref{espaco_parametros_rc}, fazemos um
resumo da construção do espaço de parâmetros
$\W_{rc}$ de acordo com Elingsrud \& Str{\o}mme
\cite{Ellingsrud_Stromme_89}. Posteriormente,
efetuamos a descrição dos pontos fixos para uma
escolha adequada da ação de um toro e em seguida
aplicamos a fórmula de resíduos de Bott no cálculo
do grau de $\Sigma(\W_{rc}{},d)$. Cabe ressaltar
que para a obtenção do  fibrado $\mathcal{E}_d$,
mostrou-se necessário efetuar uma explosão análoga ao caso de cúbicas reversa (vide subseção \ref{fibras_ed_rc}). Por fim, obtemos uma fórmula polinomial explícita que dá o grau de $\Sigma(\W_{rc}{},d)$ em função de $d$ (vide \ref{grau_cubica_regrada_P4}).

Na seção \ref{sec_QD_sg}, a família de 3-variedades de Segre em $\P5$ é tomada como espaço de parâmetros sobre o qual impomos que uma hipersuperfície de grau $d$ seja singular em algum membro variável da família. Denotamos por $\W_{sg}$ a família de 3-variedades de Segre  em $\P5$ e por $\Sigma(\W_{sg}{},d)$ a família de hipersuperfícies de grau $d$ em $\P5$ singulares em algum membro (variável) $W\in \W_{sg}$. Na subseção \ref{espaco_parametros_sg} é descrito o espaço de parâmetros $\W_{sg}$ e na subseção \ref{aplicando_bott_sg} são reunidas as informações necessárias à aplicação da fórmula de resíduos de Bott: descrição dos pontos fixos, contribuição numérica, etc,  para o cálculo do grau de $\Sigma(\W_{sg}{},d)$. Embora teoricamente o grau de $\Sigma(\W_{sg}{},d)$ seja dado por uma fórmula polinomial em $d$ (vide Proposição \ref{grau_polinomio}), por limitações computacionais não conseguimos calcular o $\deg \Sigma(\W_{sg}{},d)$ para valores suficientes de $d$, a saber, pelo menos $120+1 ( = 5 \times 24 + 1)$ valores distintos de $d$, para então interpolar os resultados e obter a fórmula polinomial em questão. Na Tabela \ref{graus_segre} explicitamos os valores obtidos para $\deg \Sigma(\W_{sg}{},d)$ com $4 \leq d \leq 28$.

No Capítulo \ref{cap_eqc} estudamos o caso da
família de quárticas elípticas. Cada curva quártica
elíptica é a interseção 
completa de um  feixe de superfícies
quádricas. Vainsencher \& Avritzer
\cite{Vainsencher_Avritzer_92} obtiveram uma
descrição explícita da componente $\W_{eqc}$ de
quárticas elípticas do esquema de Hilbert
$\Hilb_{4t}(\mathbb{P}^3)$, a qual tem dimensão
16, que é a dimensão da grassmanniana de feixes de
quádricas em \p3. Essa descrição foi utilizada em
Ellingsrud \& Str{\o}mme \cite{Ellingsrud_Stromme_Bott_96} para enumerar
curvas em certos Calabi-Yau 3-folds e em
Cukierman, Lopez \& Vainsencher \cite{Cukierman_Lopez_Vainsencher_14} para
enumerar superfícies contendo uma curva quártica
elíptica. Para o cálculo do grau de
$\Sigma(\W{}_{eqc},d)$, análogo ao caso de
cúbicas reversas, mostrou-se necessário efetuar
uma explosão adicional para obtermos efetivamente
um subfibrado $\mathcal{E}_d$
(vide subseção \ref{sub_fibra_eqc}). Na seção
\ref{espaco_parametros_eqc} é apresentada uma
descrição sumária do espaço de parâmetros $\W_{eqc}$ construído em   \cite{Vainsencher_Avritzer_92}. Na seção \ref{aplicando_bott_eqc} é efetuada a descrição dos pontos fixos  para a escolha adequada da ação de um toro, para posterior aplicação da fórmula de resíduos de Bott no cálculo do grau da família de superfícies de grau $d$ singulares ao longo de alguma quártica elíptica (variável), denotada por $\Sigma(\W{}_{eqc},d)$. O resultado obtido para o grau de $\Sigma(\W{}_{eqc},d)$ é uma fórmula polinômial em $d$ (vide \ref{grau_singular_eqc}).

Por fim, incluímos nos Apêndices \ref{ap_line_schubert2}, \dots, \ref{ap_quartica_eliptica}  os códigos em Macaulay2 \cite{Macaulay2}, Singular \cite{Singular} e Maple \cite{Maple_2015} utilizados nos cálculos enumerativos presentes nesta tese.

\chapter{Preliminares\label{cap_pl}}

Neste capítulo apresentamos alguns resultados de
caráter geral que usaremos no decorrer do
trabalho. Inicialmente tratamos da relação entre
potência simbólica e potência ordinária de um
ideal  e como essa relação se reflete na condição
de uma hipersuperfície ser singular em algum
subesquema de $\P{n}$. Especificamente, seja $\W$
uma subvariedade irredutível fechada do esquema de
Hilbert, $\Hilb_{P_{\W}(t)}(\mathbb{P}^n)$,
de subesquemas em $\mathbb{P}^n$ com polinômio de
Hilbert $P_\W(t)$ fixado, cujo membro geral estamos
admitindo ser suave. Consideremos a família
de subesquemas de $\mathbb{P}^n$ definido por
$(\mathcal{I}_W)^2$ para algum membro $W$ em $\W$,
onde $(\mathcal{I}_W)^2$ indica o quadrado do
feixe ideal de $W$. Esta nova família não é plana
em geral. Entretanto, existe um morfismo birracional
$\W'\rightarrow\W$, com $\W'$ liso e projetivo, que
planifica a família. Isto significa dizer que o
mapa racional de \W\ em $\Hilb(\p n)$ definido por $W\mapsto
W'$ onde o ideal $\cl I_{W'}=\cl I_W^2$ estende-se
a um morfismo de $\W'$ em $\Hilb$.
Na seção \ref{pot_s_o},
verifica-se que uma hipersuperfície $F$ de
$\mathbb{P}^n$ é singular ao longo de $W\in \W$
(pelo menos nos casos aqui considerados) se, e
somente se, ela contém a primeira vizinhança
infinitesimal de $W$ em $\mathbb{P}^n$, isto é, o
subesquema fechado de $\mathbb{P}^n$ com o feixe ideal
$(\mathcal{I}_W)^2$. Na seção
\ref{secao_fibrado_polinomio}\, obtém-se um
subfibrado vetorial $\mathcal{E}_d$ do
fibrado trivial $\W'{}\times
H^0(\mathcal{O}_{\P{n}}(d))$, cuja
projetivização  produz um mapa
genericamente injetivo 
(ao menos para $d\gg0$)
$$\mathbb{P}(\mathcal{E}_d)
\rightarrow \mathbb{P}^{N_d}$$ com imagem o
conjunto $\Sigma({\W{},d})$, este último indicando o lugar das hipersuperfícies de grau $d$ em $\P{n}$ que são singulares em algum elemento (variável) $W\in \W$, (vide Proposição \ref{prop_fibrado} e Teorema \ref{teorema_injetividade_generica}). E nessa seção também é enunciado o resultado que garante o comportamento polinomial do grau de $\Sigma({\W{},d})$ em função de $d$ (vide Proposição \ref{grau_polinomio}). Na seção \ref{secao_bott} apresentamos a principal ferramenta para o cálculo do grau de $\Sigma({\W{},d})$, a saber a fórmula de resíduos de Bott (vide Teorema \ref{formula_Bott}).

\section{Potência Simbólica $\times$ Potência Ordinária\label{pot_s_o}}

No que segue, $K$ é um corpo algebricamente fechado de característica zero e $S = K[x_1, \dots, x_n]$ é o anel de polinômios em $n$ variáveis. 

Seja $\mathfrak{p}$ um ideal primo de um anel
noetheriano $R$. A componente $\mathfrak{p}$-primária da $n$-ésima potência de $\mathfrak{p}$ é chamada a $n$-ésima potência simbólica de $\mathfrak{p}$, e denotada por $\mathfrak{p}^{(n)}$. Temos que $$\mathfrak{p}^{(n)} = Ker(R \rightarrow R_\mathfrak{p}\big/\mathfrak{p}^nR_\mathfrak{p})= \{r\in R;sr \in \mathfrak{p}^n \mbox{ para algum } s\in R\setminus \mathfrak{p}\}  .$$

Observe que $\mathfrak{p}^n \subset
\mathfrak{p}^{(n)}$, no entanto, a igualdade não é
sempre válida. Por exemplo, em \cite{2017arXiv170803010D} encontramos o ideal primo $$p = \left\langle x^4 -yz, y^2-xz, x^3y - z^2\right\rangle \subseteq K[x,y,z]$$  para o qual temos  
$\mathfrak{p}^{(2)} \neq \mathfrak{p}^2$. De fato, a potência ordinária 

{ \baselineskip 24 pt \parskip 2 pt 
\noindent
\centerline{$\mathfrak{p}^2= \langle y^4-2xy^2z+x^2z^2,x^3y^3-x^4yz-y^2z^2+xz^3,$}

\centerline{$x^4y^2-x^5z-y^3z+xyz^2,x^7z+x^2y^3z-3x^3yz^2+z^4,$}

\centerline{$x^7y-x^3y^2z-x^4z^2+yz^3,x^8-2x^4yz+y^2z^2\rangle,$}

\noindent ao passo que a componente $\mathfrak{p}$-primária é dada por 

\centerline{$\mathfrak{p}^{(2)} = \langle y^4-2xy^2z+x^2z^2,x^3y^3-x^4yz-y^2z^2+xz^3,$}

\centerline{$x^4y^2-x^5z-y^3z+xyz^2,x^7+x^2y^3-3x^3yz+z^3\rangle$.}
}

No caso geométrico, as potências simbólicas de $\mathfrak{p}$ têm uma descrição dada pelo Teorema \ref{ZariskiNagata}, devido a Zariski e Nagata. 

Seja $\mathfrak{p}\subset S$ um ideal primo e 
$X$ a variedade afim correspondente. Definimos $$\mathfrak{p}^{\langle n\rangle} = \{f\in S; f \mbox{ anula-se à ordem } \geq n \mbox{ em cada ponto de } X\}.$$

A condição que $f$ anula-se à ordem $n$ em cada ponto $x\in X$, significa que $f \in \mathfrak{m}_x^n$, onde $\mathfrak{m}_x$ é o ideal maximal do ponto $x$. De forma equivalente, significa que a expansão de Taylor de $f$ em torno de $x$ inicia com termos de ordem maior do que ou igual a $n$. Daí podemos escrever $$\mathfrak{p}^{\langle n\rangle} = \bigcap_{x\in X} \mathfrak{m}_x^n.$$

Além disso, se a característica de $K$ é zero, podemos descrever $\mathfrak{p}^{\langle n\rangle}$ como o conjunto de polinômios que anulam-se junto com suas derivadas parciais de ordem menor do que $n$ em todos os pontos de $X$.

\bthm[Zariski, Nagata] \label{ZariskiNagata} Seja $S  = K[x_1, \dots, x_n]$ e $\mathfrak{p}$ um ideal primo de $S$. Então $\mathfrak{p}^{\langle n\rangle} = \mathfrak{p}^{( n )}$, a $n$-ésima potência simbólica de $\mathfrak{p}$. \ethm 
\bdem Ver Eisenbud \cite{Eisenbud_95}, Teorema 3.14.\edem

Em \cite{Hochster_73}, M. Hochster estabelece vários critérios para igualdade entre a $n$-ésima potência $\mathfrak{p}^n$ de um primo $\mathfrak{p}$ em um anel Noetheriano $R$ e sua $n$-ésima potência simbólica $\mathfrak{p}^{( n )}$ para cada inteiro positivo $n$. Como aplicação desses critérios tem-se que $\mathfrak{p}^n = \mathfrak{p}^{( n )}$ no caso em que $\mathfrak{p}$ é um primo gerado por uma $R$-sequência (uma interseção completa). {Já em Dao et al.  \cite{2017arXiv170803010D} é mostrado que se $I$ é o ideal homogêneo saturado definindo um subesquema $X\subset \P{n}$ tal que $\textrm{codim}(X) = 2$, $X$ é aritmeticamente Cohen-Macaulay e uma interseção completa local então $I^{(m)} = I^m$ para todo $m\geq 1$. Como consequência, se $W \subset  \mathbb{P}^n$ é um elemento (suave) de uma das seguintes famílias
\bi
\item $\P{k}$-linear em $\P{n}$, $k<n-1$ (vide Capítulo \ref{cap_pklinear}),
\item curva plana suave de grau $m$ em $\mathbb{P}^3$ (vide  Capítulo \ref{cap_cplane}),
\item curva quártica elíptica suave em $\mathbb{P}^3$ (vide Capítulo  \ref{cap_eqc}),
\item cúbica reversa suave (vide seção  \ref{sec_QD_twc}), 
\item cúbica regrada suave em $\P4$ (vide seção  \ref{sec_QD_rc}), ou
\item uma 3-variedade de Segre em $\P5$ (vide seção \ref{sec_QD_sg})
\ei 
e $W{'}$ 
 a primeira vizinhança infinitesimal de $W$ em $\mathbb{P}^n$, isto é, o subesquema fechado de $\mathbb{P}^n$ tendo $(\mathcal{I}_W)^2$ como o feixe ideal, teremos que uma hipersuperfície $F$ de $\mathbb{P}^n$ é singular ao longo de $W$ se, e somente se, ela contém $W{'}$.}

 \section{Existe{m} um fibrado vetorial e uma fórmula polinomial\label{secao_fibrado_polinomio}}
 
 No decorrer deste trabalho nos deparamos com o
 cálculo do grau de certas subvarie\-dades de
 $\mathbb{P}^{N}$ que surgem como imagem de certos
 subfibrados projetivos de um fibrado trivial
 $X\times \mathbb{P}^{N}$. Para essa finalidade o próximo resultado é de grande utilidade.

\blem \label{lema_segre} Sejam $V$ um espaço
vetorial de dimensão $N+1$ e $\mathcal{V}$ um
subfibrado do fibrado trivial $X\times V$ sobre
uma variedade $X$ de dimensão $n$. Considere
$\mathbb{P}(\mathcal{V})
\subset X\times\p N
$ a projetivização do fibrado $\mathcal{V}$ e
$\mathbb{P}^N = \mathbb{P}(V)$. Então temos o
seguinte diagrama
induzido pelas projeções nos fatores:
$$ \xymatrix{ & \mathbb{P}(\mathcal{V})\ar[dl]_-{p_1} \ar[dr]^-{p_2} & \\
  X & & \mathbb{P}^N }
$$  
Denote por $M \subseteq \mathbb{P}^N$ a imagem de $p_2$ e suponha $p_2$ genericamente injetiva. Então $$\deg M = \int \segre(n, \mathcal{V})\cap [X].$$ onde $\segre(n, \_)$ indica a $n$-ésima classe de Segre. \elem

\bdem Ver Cuadrado \cite{Ferrer_2010}, Lema 1.22,
pag. 14, ou Fulton \cite{Fulton_1998}, Exemplo 8.3.14, pág. 143.\edem
 
\label{W'}
Seja $\W$ uma subvariedade irredutível fechada do
 esquema de Hilbert,\break
 $\Hilb_{P_{\W}(t)}(\mathbb{P}^n)$, de
 subesquemas em $\mathbb{P}^n$ com polinômio de
 Hilbert $P_\W(t)$, cujo membro geral estamos
 admitindo ser suave. E considere a subvariedade
 fechada $\W{'}$ em
 $\Hilb_{P_{\W'}(t)}(\mathbb{P}^n)$ obtida
 como o fecho da família de subesquemas de
 $\mathbb{P}^n$ definido por $(\mathcal{I}_W)^2$
 para um membro genérico $W$ em $\W$, onde
 $(\mathcal{I}_W)^2$ indica o quadrado do feixe
 ideal de $W$.
{Note a aparição do novo polinômio de Hilbert,
$P_{\W'}(t)$, graças a planitude genérica.}
 Em geral, a família de subesquemas de $\mathbb{P}^n$ definido por $(\mathcal{I}_W)^2$ para algum membro $W$ em $\W$ não é plana (vide Capítulo \ref{cap_eqc} e Seções \ref{sec_QD_twc}, \ref{sec_QD_rc} e \ref{sec_QD_sg}).
 
 Denote por $\mathbb{P}^{N_d}=\mathbb{P}(H^0(\mathcal{O}_{\P{n}}(d)))$ o espaço de parâmetros das hipersuperfícies de grau $d$ em $\P{n}$.  Pela seção \ref{pot_s_o}, exigir que uma hipersuperfície $F$ de grau $d$ seja singular ao longo de um membro geral $W\in \W$ é equivalente dizer que $F$ é um elemento de $H^0((\mathcal{I}_W)^2(d))$. 

 Escreveremos  $\mbox{Sing}(F)$ o lugar singular de $F$.
 \blem \label{lema_singF}
Suponha $\cl J_d:=(\mathcal{I}_W)^2(d)$ globalmente gerado.
Seja $F$ um elemento geral de $H^0(\cl J_d)$. Então, o lugar singular de $F$  é igual a $W$ conjuntistamente. 
\elem
\begin{proof}
A hipótese de
$\cl J_d$ ser globalmente gerado, implica
por Bertini (cf.\,\cite[10.9.2]{Hartshorne_1977}) que $\mbox{Sing}(F)\subseteq W$.\footnote{Obrigado a Angelo F. Lopez por esclarecer este ponto.}
A inclusão $W\subseteq\mbox{Sing}(F)$ é evidente.
\end{proof}

Em geral a decomposição primária do lugar
singular de $F$ revela a existência de pontos
imersos.
\bexe \em 
Considere o guarda-chuva de Whitney,  $F:= \langle
x_0^2x_3 - x_1^2x_2\rangle$, cujo lugar singular é
dado por $\mbox{Sing}(F) = \langle
2x_0x_3,-2x_1x_2,-x_1^2,x_0^2\rangle$. A
decomposição primária de $\mbox{Sing}(F)$ revela,
além da reta $l = \langle x_0,x_1\rangle$, os
pontos imersos $\langle x_0, x_1, x_2\rangle$ e
$\langle x_0, x_1, x_3\rangle$, os quais estão
sobre $l$. Cabe ressaltar que sobre estes pontos especiais, {``pinch points'', o cone tangente é um plano duplo}.
\eexe


Notação como em \S\ref{W'},
considere os mapas de projeção de $\W{'}\times \mathbb{P}^{N_d}$
  
\beq\label{diag_fund} 
~\hskip3cm
\xymatrix{ & \W{'} \times \mathbb{P}^{N_d} \ar[dl]_-{p_1}\ar[dr]^-{p_2} & \\
 \W{'} & & \mathbb{P}^{N_d}}\eeq  
 
 A regularidade de Castelnuovo-Mumford \cite{Mumford_1966} mostra que para todo $d\gg 0$, o subconjunto $\widetilde{\Sigma}(\W',d)$ de pares $(Z,F)\in \W'\times \mathbb{P}^{N_d}$ tal que $ Z\subset F$, é um fibrado projetivo sobre $\W'$ via a primeira projeção $p_1$.
Explicitamente, com base em Cukierman, Lopez \& Vainsencher \cite{Cukierman_Lopez_Vainsencher_14}, seja $\widetilde{Z}\subset \W'\times\P{n}$ o subesquema universal e da mesma forma $\widetilde{F}\subset \mathbb{P}^{N_d}\times \P{n}$ a hipersuperfície universal de grau $d$. E indique por $\widehat{Z}, \widehat{F}$ suas imagens inversas para $\W'\times\mathbb{P}^{N_d}\times \P{n}$. Daí segue que temos o seguinte diagrama de feixes sobre $Y:= \W'\times \mathbb{P}^{N_d}\times \P{n}$:

\centerline{
\xymatrix{ \mathcal{O}_Y \ar[r] \ar[dr]_-{\rho} & \mathcal{O}_Y(\widehat{F}) \ar[r]\ar@{->>}[d] & \mathcal{O}_{\widehat{F}}(\widetilde{F})\\
& \mathcal{O}_{\widehat{Z}}(\widehat{F})}
}

Por construção, a seta oblíqua $\rho$ anula-se em um ponto $(Z,F,x) \in \W'\times \mathbb{P}^{N_d}\times\P{n}$ se, e somente se, $x\in F\cap Z$. Assim, temos que $Z\subset F$ quando a condição anterior é válida para todo $x\in Z$. Desse modo, $\widetilde{\Sigma}(\W',d)$ é igual ao esquema de zeros de $\rho$ ao longo das fibras da projeção $p_{12}: \widehat{Z} \rightarrow \W'\times \mathbb{P}^{N_d}$. Relembrando Altman \& Kleiman \cite{Altman_Kleiman_1977}, (2.1) pag. 14, isto é igual aos zeros da seção adjunta do fibrado vetorial imagem direta ${p_{12}}_*(\mathcal{O}_{\widehat{Z}}(\widehat{F}))$. 

Denotamos por 

\centerline{\xymatrix{\W'\times\P{n} \ar[r]^-{q_2}\ar[d]_-{q_1} & \P{n}\\
\W'}
}
\noindent os mapas de projeção de $\W'\times \P{n}$. Como $\mathcal{O}(\widetilde{F}) = \mathcal{O}_{\mathbb{P}^{N_d}}(1)\otimes\mathcal{O}_{\P{n}}(d)$, pela fórmula de projeção temos que produzir uma seção de $\mathcal{O}_{\mathbb{P}^{N_d}}(1)\otimes\mathcal{D}_d$, onde $\mathcal{D}_d= {q_1}_*(\mathcal{O}_{\widetilde{Z}}(d))$. Por Castelnuovo-Mumford e teoria de mudança de base, existe um inteiro $d_0$ (= regularidade)  tal que $\mathcal{D}_d$ é um fibrado vetorial de posto $P_{\W'}(d)$ para todo $d\geq d_0$,  onde $P_{\W'}(t)$ indica o polinômio de Hilbert dos membros de $\W'$. De fato, $\mathcal{D}_d$ se encaixa na seguinte sequência de fibrados vetoriais sobre $\W'$:
\begin{equation}
{
\xymatrix{ 0 \ar[r] & {q_1}_*(\mathcal{I}_{\widehat{Z}}\otimes \mathcal{O}_{\P{n}}(d)) \ar@{=}[d]  \ar[r]& {q_1}_*(\mathcal{O}_{\P{n}}(d)) \ar[r] \ar@{=}[d]& {q_1}_*(\mathcal{O}_{\widetilde{Z}}(d)) \ar[r] \ar@{=}[d]& 0\\
& \mathcal{E}_d \ \ar@{>->}[r] & \mathcal{F}_d \ar@{->>}[r] & \mathcal{D}_d}
}\label{eq_fibrados}
\end{equation}
\noindent
onde $\mathcal{F}_d = H^0(\mathcal{O}_{\P{n}}(d))$ indica (fibrado trivial com fibra) o espaço de polinômios homogêneos de grau $d$.

Tomando projetivização e  a imagem inversa para $\W'\times \mathbb{P}^{N_d}$, obtemos:

\mathcenter{\beq\label{eq:fibrado}
  \vcenter{\vbox{
  \xymatrix{ & \mathcal{O}_{\mathbb{P}^{N_d}}(-1)
\vphantom{_{|_{|_|}}}
    \ar@{>->}[d] \ar[dr]^-{\overline{\rho}} & \\
\mathcal{E}_d \ \ar@{>->}[r] & \mathcal{F}_d \ar@{->>}[r] & \mathcal{D}_d}}}
  \eeq}

\noindent
Por construção, $\overline{\rho}$ anula-se
precisamente sobre $\widetilde{\Sigma}(\W',d)$. E
isso nos diz que
\\\centerline{$\widetilde{\Sigma}(\W',d)=\mathbb{P}(\mathcal{E}_d)$.}

Como $\mbox{codim}_{\W'\times\mathbb{P}^{N_d}}\widetilde{\Sigma}(\W',d) = P_{\W'}(d)$, o que coincide com o posto de $\mathcal{D}_d$, segue que $\widetilde{\Sigma}(\W',d)$ representa a classe de Chern top dimensional de $\mathcal{O}_{\P{n}}(1)\otimes \mathcal{D}_d$, (cf. Fulton \cite{Fulton_1998}, 3.2.16, p. 61). 

Agora, denotamos por $\Sigma(\W{},d)$ o
subconjunto de $\mathbb{P}^{N_d}$ formado pelas\lb
hipersuperfícies que contêm algum membro de
$\W'$.  Para um ponto geral
 de $\W'$, isto é equivalente a dizer que a superfície é singular em algum membro geral de $\W$. Assim, com a notação de (\ref{diag_fund}), temos 
 \beq
 \ \ \ \
\hskip 1.95 cm
 \Sigma(\W{},d) = p_2(\widetilde{\Sigma}(\W',d)). \eeq
E da injetividade genérica do  mapa $p_2$ (vide Teorema \ref{teorema_injetividade_generica} ) 
\beq\label{mapa_principal}\ \ \ 
\begin{tabular}{rrl}
 $\W'\times \mathbb{P}^{N_d} \supset \widetilde{\Sigma}(\W{},d)$ & $\xrightarrow{p_2}$ &  $\Sigma(\W,d) \subset \mathbb{P}^{N_d}$\\ 
 $(Z,F)$ & $\mapsto$ & $F$ \\ 
  \end{tabular}  
\eeq
e do Lema \ref{lema_segre}, obtemos que 
\beq\label{grau_geral}\ \ \ \hskip 1.5 cm
\deg \Sigma(\W{},d) = \int \segre(\omega, \mathcal{E}_d)\cap [\W{'}],\eeq
onde $\omega := \dim \W{'} = \dim \W$.

Para referências futuras registramos as argumentações acima na Proposição \ref{prop_fibrado}.

\bprop\label{prop_fibrado}
Seja $\W\subset
\Hilb_{P_{\W}(t)}(\mathbb{P}^n)$
 uma subvariedade fechada irredutível,
cujo membro geral corresponde a uma subvariedade
lisa de \p n.   Notação como acima, 
existe uma dessingularização $\W'\ra\W$ 
tal que
\\$\mb{(i)}$
para $d\gg 0$ existe um subfibrado vetorial
$\mathcal{E}_d\subset\W'\times\cl F_d$
cuja fibra sobre $W'\in\W'$ geral é o subespaço
$H^0((\mathcal{I}_W)^2(d))\subset\cl F_d$ formado
pelos \pol s homogêneos de grau $d$ com gradiente
nulo em $W$;
\\$\mb{(ii)}$
o mapa 
 $\widetilde{\Sigma}(\W{},d):=\mathbb{P}(\mathcal{E}_d)
\lar\p N=\p{}(\cl F_d)$ induzido por projeção é
genericamente injetivo e
 $$deg \Sigma(\W{},d) = \int \segre(\omega, \mathcal{E}_d)\cap [\W{'}],$$ onde $\omega = \dim \W{'} = \dim \W$.
\eprop

A seguir, pagamos nossa dívida quanto a
injetividade genérica do mapa em
(\ref{mapa_principal}). Para isso, fazemos
inicialmente um lema técnico (Lema
\ref{familia_plana}) necessário para garantir a
existência do subfibrado $\cl E_d$ acima mencionado.

\blem \label{familia_plana}
Seja $\W\subset
\Hilb_{P_{\W}(t)}(\mathbb{P}^n)$
uma subvariedade fechada irredutível.
A família formada pelos subesquemas de
$\P{n}$ definidos por $(\mathcal{I}_W)^2$ para
$W\in \W$ suave é plana.  
\elem
\begin{proof}
A afirmação decorre do seguinte fato bem conhecido.
  Seja $B$ uma $A$-álgebra e $I$ um ideal de $B$
  tal que  $ C=B/I$ é plano sobre $A$ e $I/I^2$
  é localmente livre
  como $C$-módulo. Então a vizinhança infinitesimal
  $B/I^2$ é plana sobre $A$.
  De fato, temos a sequência exata 
  $$I/I^2\subset B/I^2\twoheadrightarrow B/I,$$
onde os extremos são planos sobre $A$.  Por fim,
tomamos $A\ra B\twoheadrightarrow C$ como anéis de coordenadas
de abertos afins de $\W,\W\times\p n\supset Z=$
subesquema universal.
\end{proof}

Enfim, vamos à injetividade genérica.

\bthm\label{teorema_injetividade_generica}  Seja
$\W \subset \Hilb_{P_{W}(t)}(\mathbb{P}^n)$,
subvariedade fechada irredutível cujo membro geral
é irredutível, liso, e de dimens\~ao $\leq n-2$.
Seja $\ov\W'\subset \Hilb_{P_{W'}(t)}(\mathbb{P}^n)$
o fecho da imagem do mapa racional
que associa a cada $W\in\W$
liso, o subesquema $W'$ com ideal $\cl I_{W'}=(\cl
I_W)^2$. 
Então para $d\gg 0$ o mapa 
\begin{center}\begin{tabular}{rrl}
 $\ov\W{'}\times \mathbb{P}^{N_d} \supset \widetilde{\Sigma}(\W{},d)$ & $\xrightarrow{p_2}$ &  $\Sigma(\W{},d) \subset \mathbb{P}^{N_d}$\\ 
 $(Z\subset{}F)$ & $\mapsto$ & $F$ \\ 
  \end{tabular}\end{center}
  é genericamente injetivo.  
\ethm
  
\begin{proof}

Seja $\W'_0$ um aberto de $\ov\W'$ formado por
subesquemas $W'$ com ideal da forma
$\cl I_{W'}=(\cl I_W)^2$ com $W\in\W$ liso.
Fixe $W'\in \W{'}_0$. Seja $F$ uma hipersuperfície
geral de grau $d$ contendo $W'$. Isto significa que
$F$ é geral em
$H^0(\mathcal{I}_{W'}(d))=H^0((\mathcal{I}_W)^2(d))$
com $W=W'_{red}$ suave, irredutível. Pelo Lema
\ref{lema_singF} segue que $\mbox{Sing}(F) = W$
(conjuntistamente).
Seja $Z'\in \W{'}_0$  tal que $Z'\subset F$. Por
construção de $\W{'}_0$, temos $\cl I_{Z'}=
(\cl I_{Z})^2$ para  $Z=Z'_{red}
\in\W
$ suave, irredutível.
Em particular, o espaço tangente da hipersuperfície $F$ em cada ponto
de $Z'$ é igual ao espaço tangente de \p n. 
Daí vem $Z'\subseteq \mbox{Sing}(F) = W$
e, $(Z')_{red} = W$, seguindo $Z'=W'$.
Isso mostra que 
$$ \widetilde{\Sigma}(\W{},d)_{|\W_0'}
 \xrightarrow{\ p_2 \ }
\Sigma(\W{},d)$$
é injetivo. \end{proof}

\noindent{
\bobs
Podemos mudar de base 
passando a uma dessingularização
$\W'\lar\ov\W'$.
\eobs} 

Em Cukierman, Lopez \& Vainsencher \cite{Cukierman_Lopez_Vainsencher_14} ou em Vainsencher \cite{Vainsencher_2015},\break encontra-se um resultado que estabelece o comportamento polinomial do grau em (\ref{grau_geral}) em função do grau $d$ da hipersuperfície em $\P{n}$, o qual registramos na Proposição \ref{grau_polinomio} para referência futura.
\bprop\label{grau_polinomio}
Nas hipóteses anteriores e para $\W \subset \Hilb_{P_{\W}(t)}(\mathbb{P}^n)$ fixado, temos que  $\deg(\Sigma(\W{},d))$ é um polinômio em $d$ de grau $\leq n \dim(\W)$, para todo $d\gg 0$.
\eprop
 \begin{proof}Seja $w = \dim \W=\dim \W'$.
 Como $\cl E_d$ e $\cl D_d$ são imagem direta de um feixe sobre $\W'\times \P{n}$ (vide \ref{eq_fibrados}) e $\segre(w,\cl E_d) = \mbox{Chern}(w, \cl D_d)$, podemos aplicar Grothendieck-Riemann-Roch (c.f. \cite[p. 286]{Fulton_1998}) para expressar o caracter de Chern de $\cl D_d$ como 
 $$ch(\cl D_d) = ch(q_1!(\cl O_{\widetilde{Z}}(d)) = q_1{_*}(ch(\cl O_{\widetilde{Z}})\cdot ch(\cl O_{\P n}(d))\cdot \mbox{todd}(\P n)).$$
 Note que  o lado direito é um polinômio em $d$ de grau $\leq n$. Por outro lado, a classe de Chern $c_w$ é um polinômio de grau $w$ sobre os coeficientes do caracter de Chern e isto nos leva a deduzir que $\mbox{Chern}(w,\cl D_d) = \segre(w,\cl E_d)$ é um polinômio em $d$ de grau $\leq n w$.
 \end{proof}
 
 No entanto, em todos os casos tratados neste trabalho, encontramos como grau do polinômio na Proposição \ref{grau_polinomio} precisamente 
$(k+1)\times \dim (\W)$, onde $k$ indica a
 dimensão de um membro de $\W$. A validade para
 $\W$ arbitrário resta conjectural. 
 
\section{Fórmula de resíduos de Bott\label{secao_bott}}

Nesta seção apresentamos a principal ferramenta para tratar dos casos considerados neste trabalho, a saber a fórmula de resíduos de Bott. Os resultados aqui apresentados são baseados nas referências Ellinsgrud \& Str{\o}mme \cite{Ellingsrud_Stromme_Bott_96} e Meireles \& Vainsencher \cite{Meireles_Vainsencher_2001}.  

No que segue, $X$ será uma variedade completa
suave de dimensão $n$ munida de uma ação algébrica
do toro $\mathbb{T}=\C^\star$. O lugar dos pontos fixos indicamos por $X^\T$. 

Sejam $E_1, \dots, E_s$ fibrados vetoriais $\mathbb{T}$-equivariantes sobre $X$ e $$p(x_1^1, \dots, x_s^1, \dots, x_1^n, \dots, x_s^n)$$ um polinômio homogêneo ponderado de grau $n$ nas variáveis $x_j^i$, onde $x_j^i$ tem grau $i$. Denotemos por $p(E_1,\dots, E_s)$ o polinômio nas classes de Chern de $E_1, \dots, E_s$, obtido pela substituição $x_j^i = c_i(E_j)$.

A fórmula de integração de Bott (cf. \cite{Meireles_Vainsencher_2001}, Corolário 2.4.7, pág. 33) calcula o grau do zero ciclo $p(E_1, \dots, E_s)\cap [X]$ em termos das restrições dos fibrados $E_i's$ ao lugar dos pontos fixos, $X^\T \subset X$.

Escrevemos, abreviadamente, 
$$p(E) = p(E_1,\dots, E_s) \ \mbox{ e } \ p^\T(E) = p(E_{1T}, \dots, E_{sT})$$
para o polinômio correspondente às classes de Chern $\mathbb{T}$-equivariantes dos fibrados $E_1, \dots, E_s$. Note que 
$$p(E)\cap [X] = i^*(p^\T(E)\cap [X]_\T),$$
onde $i: X^\T  \hookrightarrow X$ indica o mapa de inclusão.

\bthm[Fórmula de resíduos de Bott]\label{formula_Bott} Sejam $E_1, \dots, E_s$, fibrados vetoriais $\mathbb{T}$-equivariantes sobre uma variedade suave e completa $X$. Então
\beq \int_X (p(E)\cap [X]) = \sum\limits_{F\subset X^\T} \pi_{F_*}\left(\dfrac{P^\T(E_{|F}\cap [F]_\T}{c_{d_F}(\mathcal{N}_{F/X})}\right),
\eeq
onde $d_F$ denota a codimensão da componente $F$ em $X$.
\ethm
\begin{proof}
Vide \cite{Meireles_Vainsencher_2001}, Teorema 2.6, pág. 35.
\end{proof}

Dada uma componente conexa (= irredutível) $F \subset X^\T$ do lugar dos pontos fixos, temos que as classes de Chern $\mathbb{T}$-equivariantes $c_k^\T(E_{|_F})$ e $c_{d_F}^\T(\mathcal{N}_{F/X})$ podem ser calculadas no anel de Chow equivariante $A_*^\T(F)$ em termos dos caracteres que comparecem na decomposição de $E_{|_F}$ e $\mathcal{N}_{F/X}$ em auto-subfibrados e das classes de Chern destes últimos.

No caso em que $X^\T$ é um conjunto finito de pontos, ou seja, os pontos fixos da ação de $\mathbb{T}$ sobre $X$ são isolados, as classes 
$$c_k^\T(E_{|_F}) \mbox{ e } c_{d_F}^\T(\mathcal{N}_{F/X}) = c_{\dim(X)}^\T(\mathcal{T}X) $$
podem ser descritas puramente em termos dos caracteres associados aos auto-fibrados. Precisamente, feita a decomposição $E_{|_F} = \oplus_\chi E_{|_F}^\chi$ em auto-espaços, as classes de Chern equivariantes de cada somando é dada por
\beq\label{pesos_bott} c_k^\T(E_{|_F}) = \binom{r}{k}\chi^k, \qquad r = \mbox{ posto de } E_{|_F}^\chi. \eeq

Em particular, {temos que a classe de Chern top dimensional $c_{max}^\T(E_{|_F})$, no caso em que o $\textrm{posto de }E=\dim (X)$ (por exemplo $E =\mathcal{T}_X$) },  é
representada no anel de Chow equivariante do ponto
fixo $F$ pelo produto de todos os caracteres que
aparecem na decomposição da fibra $E_{|_F}$ em
auto-espaços, com as respectivas
multiplicidades.

Nos casos considerados neste trabalho, tipicamente tomamos $\mathbb{T}$ como o grupo multiplicativo $\mathbb{C}^*$ agindo sobre $X$ de tal forma que o conjunto de pontos fixos $X^\T$ seja finito. Basicamente, a ação é tomada de maneira que cada caracter seja dado por $\chi_i = t^{a_i}$, $a_i \in \mathbb{Z}$ e assim o operador induzido no anel de Chow equivariante é $a_i\cdot t$, onde $t$ agora indica a classe hiperplana. 

No lado direito da fórmula de resíduos de Bott (\ref{formula_Bott}), o numerador $p^\T(E_{|_F})$ consiste em um polinômio homogêneo de grau $n = \dim X$ nas variáveis que são os caracteres que ocorrem na decomposição em auto-subfibrados. Tipicamente, suponha que o polinômio em questão contenha um termo igual a $c_1^{n-2}\cdot c_2$, enquanto que, digamos, $E_{|_F} = 2\chi_1+\chi_2$. Temos assim, $c_1^\T(E_{|_F}) = 2\chi_1+\chi_2$ e, analogamente, $c_2^\T(E_{|_F}) = \chi_{1}^2+2\chi_1\cdot\chi_2$. Daí segue que o operador de grau $n$ é dado por $(2\chi_1+\chi_2)^{n-2}\cdot(\chi_1^2+2\chi_1\chi_2)$, e o referido termo ganha a forma final $(2a_1+a_2)^{n-2}\cdot(a_1^2+2a_1\cdot a_2)\cdot t^n$ sobre o  anel de Chow equivariante. Portanto, o numerador e o denominador no lado direito de (\ref{formula_Bott}) são múltiplos inteiros de $t^n$. Ao simplificarmos, obtemos um número racional e, consequentemente, o lado direito da fórmula de resíduos de Bott (\ref{formula_Bott}) consiste em uma soma finita de números racionais obtidos a partir dos pesos como descrito em (\ref{pesos_bott}).

Mais precisamente, denote por $\tau_1(E,F),
\dots, \tau_r(E,F)$ os pesos que ocorrem na
decomposição de $E_{|_F}$ em auto-subfibrados, e
para cada inteiro $k\geq 0$, seja $\sigma_k(E,F)$
a $k$-ésima função simétrica elementar desses
pesos. Temos então os seguintes corolários.

\bcor
Nas notações acima, cada classe de Chern equivariante $c_k^\T(E_{|_F})$ é representada no anel de Chow equivariante do ponto fixo $F$ por $\sigma_k(E,F)$.
\ecor

\bcor
A classe de Chern equivariante top dimensional do fibrado tangente de $X$ é dada no anel de Chow equivariante de um ponto fixo $F$ pelo produto dos pesos que ocorrem na decomposição da respectiva fibra.
\ecor

\chapter{Hipersuperfícies singulares ao longo de um \texorpdfstring{$\P{k}$}{Pk}-linear\label{cap_pklinear}}
Neste capítulo tratamos do caso de
hipersuperfícies em $\P{n}$ singulares ao longo de
um $\P{k}$-linear, $k < n-1$. Aqui, o espaço de
parâmetros consiste na grassmanniana
$\mathbb{G}(k+1,n+1)$, cuja dimensão é dada por
$
(k+1)\times (n-k)$. O nosso interesse é determinar o grau da família de hipersuperfícies de grau $d$ singulares em algum membro de $\mathbb{G}(k+1,n+1)$, o qual vamos denotar por $\deg \Sigma(\W{}_{(k+1,n+1)},d)$.

Inicialmente, tratamos em detalhe do caso
específico $k=1$ e $n=3$ (superfícies em \p3 singulares ao longo de uma reta variável) explicitando a aplicação da fórmula de resíduos de Bott. Posteriormente, listamos os resultados obtidos para outros valores de $k$ e $n$. 

\section{Impondo uma reta singular}

Considere a família de retas em $\P 3$ parametrizada pela grassmanniana de retas $\W_{(2,4)} = \mathbb{G}(2,4)$. Estamos interessados em determinar o grau de $\Sigma(\W_{(2,4)}{},d) \subset \mathbb{P}^{N_d}$, a família de superfícies de grau $d$ que contém algum membro $W \in \W_{(2,4)}$ em seu lugar singular.

Pela Proposição \ref{prop_fibrado} temos que para $d\gg 0$ existe um subfibrado $\mathcal{E}_d$ do fibrado trivial $\W_{(2,4)}\times \mathbb{P}^{N_d}$, cuja projetivização $\mathbb{P}(\mathcal{E}_d)=\widetilde{\Sigma}(\W_{(2,4)}{},d)$ e temos um mapa $p_2$ obtido pela segunda projeção, o qual é genericamente injetivo (vide Teorema \ref{teorema_injetividade_generica}) e com imagem o conjunto $\Sigma(\W_{(2,4)}{},d)$. 
$$ \xymatrix{  \widetilde{\Sigma}(\W_{(2,4)}{},d)
  &=&
  \mathbb{P}(\mathcal{E}_d)
\ar[dll]^-{p_1} \ar[drr]_-{p_2}
  &\subset& \W_{(2,4)}{} \times \mathbb{P}^{N_d}  \\
            \W_{(2,4)}{}  & &&&\Sigma(\W_{(2,4)}{},d) \subset\mathbb{P}^{N_d} }$$  
%
%
Assim, pelo Lema \ref{lema_segre} temos que 
\beq\label{deg_line} \deg \Sigma(\W_{(2,4)}{},d) = \int \segre(4,\mathcal{E}_d)\cap [\W'_{(2,4)}{}]. \eeq

Neste caso, vamos explicitar o fibrado
$\mathcal{E}_d$. Para isso, consideremos a
sequência tautológica da grassmanniana de retas em
$\P3$ 
\mathcenter{
\begin{equation}
\mathcal{S} \rightarrowtail \mathcal{F}_1 \twoheadrightarrow \mathcal{Q},
\end{equation}}
onde $\mathcal{F}_1$ denota o fibrado trivial com fibra $\langle x_0, x_1, x_2, x_3\rangle$, o espaço de formas lineares de grau um. Cada fibra de $\mathcal{S}$ é um subespaço bidimensional das equações de uma reta e a correspondente fibra de $\mathcal{Q}$ é o espaço bidimensional de coordenadas homogêneas sobre a reta.

Exigir que uma superfície de grau $d$ em $\P3$
contenha uma reta em seu lugar singular é
equivalente requerer que a (equação da) superfície pertença ao quadrado do ideal que define a reta. Daí segue que obtemos um mapa de fibrados $$ S_2(\mathcal{S}) \otimes \mathcal{F}_{d-2} \twoheadrightarrow \mathcal{E}_d \subset \mathcal{F}_d,$$ o qual possui um núcleo dado por $\mathcal{S}\otimes \stackrel2\wedge\mathcal{S}\otimes \mathcal{F}_{d-3}$. Portanto, obtemos uma sequência exata 
\begin{equation}
\mathcal{S}\otimes \wed 2\mathcal{S}\otimes \mathcal{F}_{d-3} \rightarrowtail S_2(\mathcal{S}) \otimes \mathcal{F}_{d-2} \twoheadrightarrow \mathcal{E}_d \subset \mathcal{F}_d.
\end{equation}

Usando o pacote Schubert2 para Macaulay2 \cite{Macaulay2} encontramos 
\begin{small}\beq\label{grau_reta} \deg \Sigma(\W_{(2,4)}{},d) =
\dfrac{1}{64}d(d-1)(27d^6-117d^5+269d^4-375d^3+312d^2-132d+48)
\eeq\end{small}
para o grau da subvariedade
$\Sigma(\W_{(2,4)}{},d) \subset \mathbb{P}^{N_d}$
de superfícies em $\P3$ de grau $d$ singulares ao
longo de alguma reta (variável). O leitor pode
consultar o Apêndice \ref{ap_line_schubert2} para
obter o código utilizado no cálculo anterior. O
valor para $d=2$ é 10, como esperado para o grau
da imagem de $\p{}(\cl F_1)\times\p{}(\cl
F_1)\ra\p{}(\cl F_2)$, formada pelos pares de planos.

\subsection{Impondo uma reta singular via Bott\label{reta_via_Bott}}
Nesta seção, repetimos o cálculo anterior para obter o $\deg\Sigma(\W_{(2,4)}{},d)$ com o intuito de elucidar a aplicação da fórmula de resíduos de Bott, como explicado em Ellingsrud \& Str{\o}mme\cite{Ellingsrud_Stromme_Bott_96} e Meurer \cite{Meurer_96}, e que será utilizada nos demais casos tratados neste trabalho.

Considere o toro $\mathbb{T}=\mathbb{C}^*$ agindo
diagonalmente sobre $(\mathbb{C}^4)^\vee$ via
\\\centerline{$t\circ x_i := t^{w_i}x_i$,} com pesos apropriados, digamos:
\beq\label{pesoscurvasplanas} w_0 = 4, w_1 = 11, w_2=17, w_3=32.\eeq 
A ação de $\mathbb{T}$ sobre $(\mathbb{C}^4)^\vee$ induz uma acão natural sobre $\P3$, de tal forma que os fibrados vetoriais $O_{\P3}(d)$ e $\mathcal{T}\P3$ são equivariantes. Assim, relembrando a sequência de Euler vemos que na fibra, digamos sobre $\langle e_0\rangle \in \P3$, o termo do meio é $\langle x_0 \rangle \otimes \langle e_0, \dots, e_3\rangle$. Consequentemente, o espaço tangente $\mathcal{T}\P3_{\langle e_0\rangle} = \langle x_0\rangle \otimes \langle e_1, e_2, e_3\rangle$ como $\mathbb{T}$-espaços, o qual decompõe-se em $\mathbb{T}$-espaços unidimensional com pesos $u_i = w_0 - w_i$, $i=1,2,3$. Nessas condições as classes de Chern $\mathbb{T}$-equivariantes de $\mathcal{T}\P3$ no ponto fixo $f = \langle e_0\rangle$ são apenas as funções simétricas elementares: $c_1^f = u_1+u_2+u_3$, $c_2^f = u_1u_2+u_1u_3+u_2u_3$ e $c_3^f = u_1u_2u_3$.

De modo análogo, obtemos uma ação natural induzida sobre $\W_{(2,4)} = \mathbb{G}(2,4)$ tal que os fibrados vetoriais sobre $\W_{(2,4)}$ sejam $\mathbb{T}$-equivariantes. No caso da Grassmanniana de retas em $\P3$ temos seis pontos fixos, correspondentes aos eixos coordenados $\langle x_0,x_1\rangle,\langle x_0,x_2\rangle \dots, \langle x_2,x_3\rangle$. De (\ref{deg_line}), computamos 
\beq\label{line_bott} \deg \Sigma(\W_{(2,4)}{},d) = \int_{\W'_{(2,4)}{}} \segre(4,\mathcal{E}_d) = \sum\limits_f \dfrac{c^f_4(-\mathcal{E}_d)}{c^f_4(\mathcal{T}\W'_{(2,4)})}\eeq
localizando nos seis pontos fixos. A última igualdade em (\ref{line_bott}) é devido à fórmula de resíduos de Bott (vide Teorema \ref{formula_Bott}) e,  $c_4^f(.)$ indica a classe de Chern $\mathbb{T}$-equivariante no ponto fixo $f$. 

{
  \bobs
Considere 
uma sequência de fibrados vetoriais $$0\rightarrow \mathcal{S}\rightarrow \mathcal{T} \rightarrow \mathcal{Q} \rightarrow 0,$$ onde $\mathcal{T}$ indica um fibrado trivial. Pela fórmula de Whitney temos a seguinte igualdade envolvendo a classe de Chern total dos fibrados $\mathcal{S}$ e $\mathcal{Q}$: $c(\mathcal{S})c(\mathcal{Q})=1$, o que implica $c(\mathcal{Q}) = c(\mathcal{S})^{-1}$. Além disso, temos que $c(\mathcal{Q}) = c(-\mathcal{S})$ e $c(\mathcal{S})=\segre(\mathcal{S})^{-1}$, donde segue que $\segre(\mathcal{S})=c(-\mathcal{S})$ e isto justifica a ocorrência do fibrado $\mathcal{-E}_d$ em (\ref{line_bott}). \eobs}

O denominador em (\ref{line_bott}), digamos para
$l = \langle x_0, x_1\rangle$, é obtido da
seguinte forma: primeiro encontramos a fibra do tangente
$$\mathcal{T}_l \W_{(2,4)} = Hom(S_l,Q_l) = \langle x_0,x_1\rangle^\vee \otimes \langle x_2, x_3 \rangle = \frac{x_2}{x_0}+\frac{x_3}{x_0}+\frac{x_2}{x_1}+\frac{x_3}{x_1},$$
onde $\frac{x_i}{x_j}$ indica o $\mathbb{T}$-espaço com peso $w_i - w_j$.

Dessa forma, obtemos $c^l_4(\mathcal{T}\W_{(2,4)}) = (w_2-w_0)(w_3-w_0)(w_2-w_1)(w_3-w_1)$, que com a escolha dos pesos em (\ref{pesoscurvasplanas}) nos fornece o valor 45864.

Da mesma forma, o numerador requer a decomposição
em pesos do fibrado $\mathcal{E}_d$.  A fibra
sobre a reta $l = \langle x_0,x_1\rangle$ consiste
nas superfícies $F$ de grau $d$ que contêm $l$ em
seu lugar singular, ou de forma equivalente que
$F$ pertence ao quadrado do ideal da reta
$l$. Portanto, $F$ é da forma $F = Ax_0^2+Bx_0x_1
+ Cx_1^2$, com $A\in K[x_0, \dots, x_3]$ e $B,C\in
K[x_1,x_2,x_3]$ polinômios de grau
$d-2$. Efetuando uma simples contagem de monômios
vemos que o posto do fibrado $\mathcal{E}_d$ é
dado por $rk (\mathcal{E}_d) = \binom{d+1}{3}+
2\binom{d}{2}$.  Ou ainda, como a família de
subesquemas de $\P3$ definida pelo quadrado do
ideal da reta $l$ é plana, com polinômio de
Hilbert constante $p_{\W'_{(2,4)}}(t) = 3t+1$,
segue que $rk(\mathcal{E}_d) = \binom{d+3}{3} -
p_{\W'_{(2,4)}}(d) = \binom{d+3}{3} - (3d+1)$.
Daí segue que  $\dim \Sigma(\W_{(2,4)}{},d) = 4 + \binom{d+3}{3} - (3d+1) -1.$

Para fixar as ideias tomamos $d=3$, o que nos dá
$\dim \Sigma(\W_{(2,4)}{},3) = 13$ e
$rk(\mathcal{E}_3) = 10$. Já a fibra de
$\mathcal{E}_3$ sobre $l$, que consiste nas $F\in H^0(\mathcal{O}_{\P3}(3))$ com gradiente nulo ao longo de $l$,  tem a decomposição em
pesos dada por:
$$\mathcal{E}_{3_l} = x_0^3+x_0^2x_1+
x_0^2x_2+x_0^2x_3 +
x_0x_1^2+x_0x_1x_2+x_0x_1x_3+x_1^3+x_1^2x_2+x_1^2x_3
$$
e, consequentemente, $-\mathcal{E}_{3_l}$ tem a decomposição em pesos dada por:
$$ -\mathcal{E}_{3_l} =
x_0x_2^2+x_1x_2^2+x_2^3+x_0x_2x_3+x_1x_2x_3+x_2^2x_3+x_0x_3^2+x_1
x_3^2+x_2x_3^2+x_3^3,
$$
onde $x_i^\alpha x_j^\beta x_k^\gamma$ indica o $\mathbb{T}$-espaço com peso $\alpha w_i + \beta w_j+\gamma w_k$.

A correspondente contribuição numérica  referente
a $c^l_4(-\mathcal{E}_3)$, utilizando os pesos em
(\ref{pesoscurvasplanas}), é 3217978137. O ponto
fixo escolhido produz a fração $3217978137/45864$.
A contribuição total dos seis pontos fixos fornece 

\noindent$$\ba c \frac{3217978137}{45864} - \frac{2152229961}{17640} + \frac{774359841}{28665} + \frac{1227942219}{28665} - \frac{392711889}{17640} + \frac{218302833}{45864} = 504. \ea$$

Este é o grau da subvariedade de $\mathbb{P}^{19}$
consistindo nas superfícies de grau 3 em $\P3$ que
são singulares ao longo de alguma reta. Daí segue,
por exemplo, que existem 504 "guarda-chuvas"\,  de
Whitney (Figura \ref{umbrella_whitney}) passando
por  13 pontos gerais
(cf. Coray \& Vainsencher \cite{Coray1986}).
\begin{figure}[ht]
\centering\ifpdf
\includegraphics[scale=0.2]{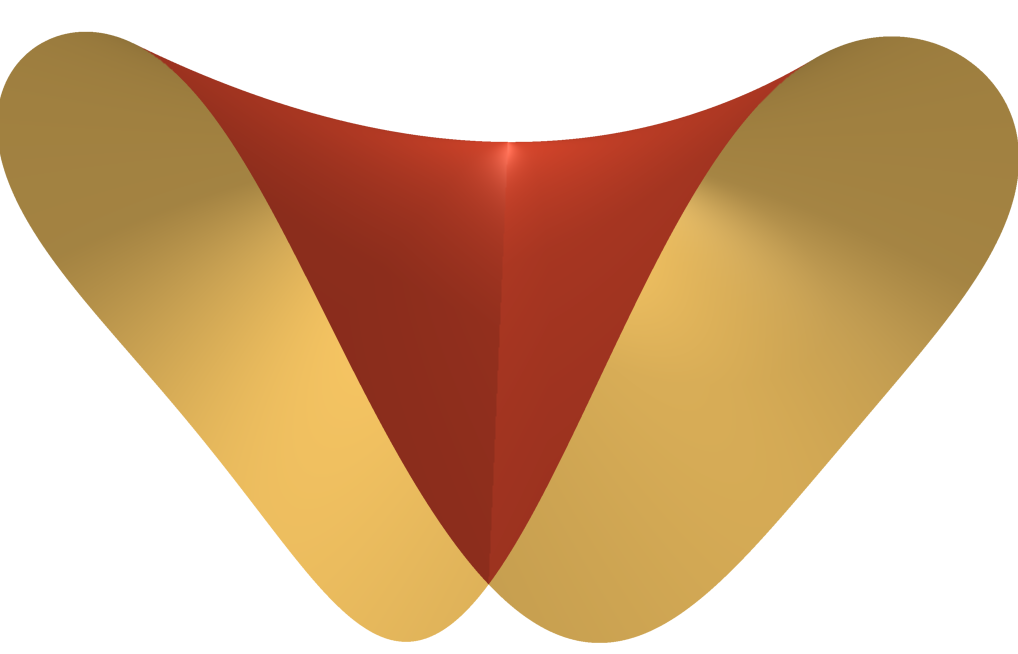}
\caption{$x^2 = y^2z$}\fi
\label{umbrella_whitney}
\end{figure}

Como $\deg \Sigma(\W_{(2,4)}{},d)$ é um polinômio
em $d$ de grau menor do que ou igual $12 ( =
3\times 4)$ (Vide Proposição \ref{grau_polinomio})
é suficiente encontrarmos o grau para 13 valores
de $d$ e interpolar os resultados. No Apêndice
\ref{ap_pk_linear} o leitor pode consultar um
código para executar os cálculos acima usando o
software Maple \cite{Maple_2015}. Encontramos
a mesma expressão antecipada em (\ref{grau_reta})
para o grau da subvariedade $\Sigma(\W_{(2,4)}{},d) \subset \mathbb{P}^{N_d}$ de superfícies em $\P3$ de grau $d$ singulares ao longo de alguma reta (variável). Observe que o grau em (\ref{grau_reta}) é $(1+1)\times\dim(\W_{(2,4)})$.
\section{Impondo um \texorpdfstring{$\P{k}$}{}-linear como lugar singular de hipersuperfícies}
Análogo ao caso de retas em $\P3$, consideremos
hipersuperfícies singulares ao longo de um
$\P{k}$-linear em $\mathbb{P}^n$, $1 < k <
n-1$. Aqui, o espaço de parâmetros consiste na
grassmanniana $\mathbb{G}(k+1,n+1)$. O nosso interesse é determinar o grau da família de hipersuperfícies de grau $d$ singulares em algum membro de $\mathbb{G}(k+1,n+1)$, o qual vamos denotar por $\deg \Sigma(\W{}_{(k+1,n+1)},d)$. No Apêndice \ref{ap_pk_linear} o leitor encontra os códigos utilizados para efetuar os cálculos necessários, utilizando a fórmula de resíduos de Bott,  para obtenção do $\deg \Sigma(\W{}_{(k+1,n+1)},d)$. Abaixo listamos os resultados obtidos para $(k,n) = (2,4), (2,5) \mbox{ e } (3,5)$.
\beq\label{grau_p2_p4}
\ba {c}
\deg \Sigma(\W{}_{(3,5)},d)= \dfrac{1}{82944}d(d-1)(d+1)(d+2)(9d^{14}-18d^{13}\\\na9
-63d^{12}+396d^{11}-405d^{10}-1530d^9+5328d^8-4176d^7-9414d^6\\\na9
+27208d^5-24347d^4-4696d^3+36572d^2-32544d+14400)
\ea
\eeq 
\beq\label{grau_p2_p5}\ba{c}
\deg \Sigma(\W{}_{(3,6)},d)=\dfrac{1}{4251528000}d(d-1)(d+1)(d+2)\\\na9
(12800d^{23}-25600d^{22}-224000d^{21}+966400d^{20}+520800d^{19}\\\na9
-10632000d^{18}+18128000d^{17}+35186000d^{16}-170677265d^{15}\\\na9
+145358830d^{14}+449576760d^{13}-1292773830d^{12}+778144037d^{11}\\\na9
+2164141556d^{10}-5208921230d^9+3728975455d^8+3332483181d^7\\\na9
-10452711042d^6+10781927010d^5-2523245175d^4-7609562253d^3\\\na9
+11511503406d^2-8323547040d+3637418400)
\ea\eeq
\beq\label{grau_p3_p5}\ba{c}
\deg \Sigma(\W{}_{(4,6)},d)=\dfrac{1}{54358179840000}d(d-1)(d+2)(d+1)(1125d^{28}\\\na9
+15750d^{27}+86625d^{26}+168750d^{25}-187875d^{24}-38250d^{23}\\\na9
+8824725d^{22}+23473350d^{21}-32467725d^{20}-128183670d^{19}\\\na9
+426415635d^{18}+1377078570d^{17}-2137554049d^{16}\\\na9
-7117020302d^{15}+15925316455d^{14}+37514746370d^{13}\\\na9
-82840806388d^{12}-125157483544d^{11}+422227932240d^{10}\\\na9
+287672117600d^{9}-1529648949952d^{8}+207120164224d^{7}\\\na9
+4517312266240d^{6}-3047085731840d^{5}-6253154779136d^{4}\\\na9
+11893749153792d^{3}+2911913902080d^{2}\\\na9
-8455245004800{d}+2378170368000)
\ea
\eeq

Observe que novamente os graus em
(\ref{grau_p2_p4}), (\ref{grau_p2_p5}) e
\ref{grau_p3_p5}), são iguais a
\lb
$(k+1)\times \dim \mathbb{G}(k+1,n+1)$.

\chapter{Superfícies singulares ao longo de curvas planas\label{cap_cplane}}

A família, $\W_m$, de curvas planas de grau $m>1$ em
$\mathbb{P}^3$ é
parametrizada por um $\mathbb{P}^{N_m}$-fibrado
sobre $\check{\mathbb{P}}^3$ $$\mathcal{X}_m \rightarrow
\check{\mathbb{P}}^3,$$ onde $N_m =
\binom{m+2}{2}-1$. De fato, tomemos a sequência tautológica de $\check{\mathbb{P}}^3$, o espaço projetivo dual:
\beq \mathcal{S}=\mathcal{O}_{\check{\mathbb{P}}^3}(-1) \rightarrowtail \mathcal{F}_1 \twoheadrightarrow Q.\eeq 
Temos que $\W_m = \mathbb{P}(S_m(Q))$. E consequentemente, $\dim \W_m = N_m+3$. A sequência tautológica sobre $\W_m = \mathbb{P}(S_m(Q))$ é dada por 

\centerline{
\xymatrix{\mathcal{O}_{\W_m}(-1)\ \ar@{^{(}->}[r]
  \ar[dr] & \ S_m(\mathcal{Q})\ar@{->>}[d]\\
& \mathcal{O}_{\check{\mathcal{Q}}}(m)}
}

Considere a curva universal $C_m$ junto com o mapa de projeção $\pi$ 

\centerline{\xymatrix{C_m \subset
    \mathbb{P}(\check{\mathcal{Q}})\times_{\check{\mathbb{P}}^3}
    \overbrace{\mathbb{P}(S_m(\mathcal{Q}))}^{\W_m}
    \ \ar@{^{(}->}[r]\ar[d]_-{\pi}& \mathbb{P}^3\times \W_m\\
\W_m }}

Daí segue a sequência exata

\centerline{
\xymatrix{ S_{d-m}(\mathcal{Q}) \otimes
  \mathcal{O}_{\W_m}(-1) \ \ar@{>->}[r] & S_d(\mathcal{Q}) \ar@{->>}[r] & \pi_*\mathcal{O}_{C_m}(d)}.
}

E por outro lado, temos

\centerline{
\xymatrix{ \mathcal{O}_{\P3}(-1) = \mathcal{O}_{\check{\mathcal{Q}}}(-1)\quad  \ar@{>->}[r] \ar[rrd] & \check{\mathcal{Q}}\quad   \ar@{>->}[r] & \mathbb{C}^4  \ar@{->>}[d]\\
& & \mathcal{O}_{\check{\mathbb{P}}^3}(1)}
}

e 

\centerline{
\xymatrix{ & & & \mathcal{O}_{\W_m}(-1) \otimes \mathcal{O}_{\check{\mathcal{Q}}}(d-k)\ar[d]\\
& & & \mathcal{O}_{\mathbb{P}(\check{\mathcal{Q}})|_{\W_m}}(d)\ar[d]\\
\left(\faktor{\mathcal{I}_{C_m}}{\mathcal{I}_{C_m}^2}\right)\ \ar@{>->}[r] & \mathcal{O}_{\P3}(d) \ar[r] & \left(\faktor{\mathcal{O}}{\mathcal{I}_{C_m}^2}\right) \otimes \mathcal{O}_{\P3}(d)\ar@{->>}[r] & \mathcal{O}_{C_m}(d)}
}

Agora, da sequência conormal torcida por $\mathcal{O}_{\mathcal{C}_m}(d)$, a saber

$$\left(\!\!\ba{*5c}
\left(\mathcal{O}_{\check{\mathbb{P}}^3}(-1)\otimes
\mathcal{O}_{\check{\mathcal{Q}}}(-1)\right)_{|_{C_m}}
&&&&\\
||&&&&\\
\check{\mathcal{N}}_{\mathbb{P}(\check{\mathcal{Q}})|_{\P3}}|_{C_m}
&\!\!\!\!\!\!\!\!\!\!\!\!\!\!
\xymatrix@C2.5pc{\ar@{^{(}->}[r] &}
&\!\!
\check{\mathcal{N}}_{C_m}|_{\P3_{\W_m}}
&\!\!\!
\xymatrix@C2.5pc{\ar@{->>}[r] &}
&
\check{\mathcal{N}}_{C_m}\mathbb{P}(\check{\mathcal{Q}})_{\W_m}
\ea\!\!
\right)\otimes \mathcal{O}_{C_m}(d)
$$
%
segue que 
\beq\deg \Sigma(\W_m{},d) = \int chern(N_m+3, \mathcal{E}_d),\eeq
onde 

$\ba c \mathcal{E}_d =  (\mathcal{O}_{W_m}(-1)\otimes (S_{d-m}(\mathcal{Q})-S_{d-2m}(\mathcal{Q})\otimes \mathcal{O}_{W_m}(-1))\oplus \\\na9
\mathcal{O}_{\check{\mathbb{P}}^3}(-1)\otimes(S_{d-1}(\mathcal{Q})-S_{d-1-m}(\mathcal{Q})\otimes \mathcal{O}_{W_m}(-1)))\oplus\\\na9 
(S_d(\mathcal{Q})-S_{d-m}(\mathcal{Q})\otimes\mathcal{O}_{W_m}(-1)).\ea$

No Apêndice \ref{ap_curva_m_schubert2} encontram-se os códigos para o cálculo do $\deg\Sigma(\W_m{},d)$, grau da família de superfícies de grau $d$ singulares ao longo de uma curva plana de grau $m$, utilizando o pacote Schubert2 do Macaulay2  \cite{Macaulay2}.  Por exemplo, para $m=2$ encontramos 
\begin{small}\begin{equation}\label{grau_conica_1}
\ba c \deg\Sigma(\W_2{},d) = \dfrac{1}{2580480}(d-2)(150903d^{15}-3809754d^{14}\\\na9 
+44834472d^{13}-317080224d^{12}+1422290970d^{11}-3579080844d^{10}\\\na9
-455933988d^9+47928493544d^8 -237841700217d^7+712127741206d^6\\\na9
-1498533401372d^5+2287674925704d^4-2504345972608d^3\\\na9 
+1873638158208d^2-859900216320d+182801203200)\ea
\end{equation}\end{small}

\noindent para o grau da subvariedade
$\Sigma(\W_2{},d)\subset \P{{N_d}}$ de superfícies
de grau $d$ singulares ao longo de alguma cônica
(variável). Observe que o grau do polinômio em (\ref{grau_conica_1}) é igual a $(1+1)\times\dim(\W_2)$.

O leitor interessado em obter $\deg\Sigma(\W_m{},d) $ para outros valores de $m$, basta alterar seu valor no código do Apêndice  \ref{ap_curva_m_schubert2}.

\section{Impondo uma curva plana singular via Bott}

Nesta seção, seguimos o mesmo princípio aplicado na seção \ref{reta_via_Bott} para obtermos o $\deg\Sigma(\W_m{},d)$ via a fórmula de resíduos de Bott.  

Tome, digamos, $\langle x_0, C_m\rangle \in
\mathcal{X}_m$. Ponha $\overline{\mathcal{F}_i} =
\mathcal{F}_i\big/\langle x_0\mathcal{F}_{i-1}\rangle$
e considere o subespaço vetorial de \pol s
homogêneos de grau $d$,
\beq\label{esp_superf_singulares_curvas_planas} x_0^2\mathcal{F}_{d-2} + x_0C_m \overline{\mathcal{F}}_{d-m-1} + C_m^2\overline{\mathcal{F}}_{d-2m} \subset \mathcal{F}_d.\eeq

Seja $F_d$ uma superfície em $\P3$ de grau $d$
singular ao longo da curva plana $\langle
x_0,C_m\rangle$, ou seja, $F_d$ pertence ao
quadrado do ideal dessa curva e assim é um
elemento do espaço vetorial definido em (\ref{esp_superf_singulares_curvas_planas}). Pela Proposição \ref{grau_polinomio}, para cada $d\gg 0$, obtemos um subfibrado $\mathcal{E}_d
\subset \mathcal{X}_m \times \cl{F}_d$. Através de uma simples contagem de parâmetros em (\ref{esp_superf_singulares_curvas_planas}) vemos que o posto do fibrado $\mathcal{E}_d$ é dado por
\beq rk(\mathcal{E}_d) =  \binom{d+1}{3}+\binom{d-m+1}{2}+\binom{d-2m+2}{2}.\eeq
%

Além disso, a projetivização $\mathbb{P}(\mathcal{E}_d)$ produz um
mapa genericamente injetivo (ver Teorema
\ref{teorema_injetividade_generica})
$$\mathbb{P}(\mathcal{E}_d)
\longrightarrow \mathbb{P}^{N_d}
$$ com imagem, $\Sigma(\W_m{},d)\subset
\mathbb{P}^{N_d}$,  a subvariedade
 formada pelas superfícies em $\P3$ de grau
$d$  singulares ao longo de algum membro $W\in \W_m$. Daí segue que a dimensão de $\Sigma(\W_m{},d)$ é dada por 
\begin{small}
\beq \dim
\Sigma(\W_m{},d) = N_m + 3 +rk(\cl E_d)
-1.\eeq
Pelo Lema \ref{lema_segre} temos que o grau de $\Sigma(\W_m{},d)$ é dado por
\beq\label{grau_sup_sing_m} \deg \Sigma(\W_m{},d) = \displaystyle\int_{\W'_m{}} \segre(N_m+3,\mathcal{E}_d).\eeq\end{small}

Para o cálculo do grau em (\ref{grau_sup_sing_m}) utilizamos a fórmula de resíduos de Bott
\beq\label{grau_sup_sing_m_bott} \deg \Sigma(\W_m{},d) = \displaystyle\int_{\W'_m{}} \segre(N_m+3,\mathcal{E}_d) = \sum\limits_f \frac{c^f_{N_m+3}(-\mathcal{E}_d)}{c^f_{N_m+3}(\mathcal{T}\W'_m)}.\eeq
 Por exemplo, vamos analisar as contas para o caso de cônicas, $m=2$. Assim, $\W_2$, a família de cônicas em $\P3$, é parametrizada por um $\P5$-fibrado sobre $\check{\mathbb{P}}^3$
$$\mathcal{X}_2 \rightarrow
\check{\mathbb{P}}^3,$$
o que nos dá $\dim \W_ 2 = 5+ 3 = 8$. Aqui o
polinômio de Hilbert $p_{\W{'}_2}(t)=6t-3$
(quadrado da cônica) e, consequentemente, o $rk(\mathcal{E}_d) = \binom{d+3}{3}-(6d-3) = \binom{d+1}{3}+\binom{d-1}{2}+\binom{d-2}{2}$. O mapa $\mathbb{P}(\mathcal{E}_d) \rightarrow \mathbb{P}^{N_d}$ é genericamente injetivo para $d\geq 4$. 

De modo análogo ao que fizemos para o cálculo do $\deg\Sigma(\W_{(2,4)}{},d)$, tomamos o toro $\mathbb{T}=\mathbb{C}^*$ agindo diagonalmente sobre $(\mathbb{C}^4)^\vee$ via $t\circ x_i := t^{w_i}x_i$, com pesos apropriados, digamos:
\beq \label{pesos_curvas_planas_grau_m} w_0=11, w_1 = 17, w_2 = 32, w_3 = 55. \eeq
Tais pesos são escolhidos com a preocupação de não
anular denominadores na aplicação da fórmula de Bott.
A ação de $\mathbb{T}$ induz uma ação natural sobre $\check{\mathbb{P}}^3$ de tal forma que os fibrados sobre este último sejam $\mathbb{T}$-equivariantes e com pontos fixos isolados. No caso de $\W_2$, os pontos fixos são dados por $\langle x_i, x_jx_k\rangle, i\nin\{j,k\}$, o que nos fornece um total de 24 pontos fixos. Em geral, dado qualquer $m>1$, o total de pontos fixos é $4\times\binom{m+2}{2}$. A representação da fibra do denominador em (\ref{grau_sup_sing_m_bott}), por exemplo no ponto fixo $f = \langle x_0, x_1^2\rangle$, é dada por
\beq\ba {lcl} \mathcal{T}_f \W_2 &=& \langle x_1,x_2,x_3\rangle \otimes \langle x_0\rangle^\vee + \langle x_1x_2, x_1x_3, x_2^2, x_2x_3, x_3^2\rangle \otimes \langle x_1^2\rangle^\vee \\\na9
 &=& \dfrac{x_1}{x_0}+\dfrac{x_2}{x_0}+\dfrac{x_3}{x_0}+\dfrac{x_1x_2}{x_1^2}+\dfrac{x_1x_3}{x_1^2}+\dfrac{x_2^2}{x_1^2}+\dfrac{x_2x_3}{x_1^2}+\dfrac{x_3^2}{x_1^2},\ea\eeq

\noindent onde, digamos $\dfrac{x_2x_3}{x_1^2}$, indica o $\mathbb{T}$-espaço com peso $w_2+w_3-2w_1$.

Dessa forma, obtemos
\\\centerline{$
\ba r
c_8^f(\mathcal{T}\W_2)=(w_1-w_0)(w_2-w_0)(w_3-w_0)
(w_2-w_1)\\
(w_3-w_1)(2w_2-2w_1)(w_2+w_3-2w_1)(2w_3-2w_1),
\ea$}
cuja contribuição numérica com a escolha dos pesos em (\ref{pesos_curvas_planas_grau_m}) nos fornece  $381864067200$. 

Resta-nos determinar a contribuição numérica do numerador referente ao ponto fixo $f = \langle x_0, x_1^2\rangle$. Fixamos, por exemplo $d=4$.  Neste caso, o posto do fibrado $\mathcal{E}_4$ é igual a 14, e a fibra sobre o ponto fixo $f = \langle x_0, x_1^2\rangle$ tem a decomposição em pesos dada por
\begin{equation}
\begin{array}{lcl}
\mathcal{E}_{4_f} &=& x_0^4+x_0^3x_1+x_0^2x_1^2+x_0x_1^3+x_1^4+x_0^3x_2+x_0^2x_1x_2+ x_0x_1^2x_2  \\ 
 & & +x_0^2x_2^2+x_0^3x_3+ x_0^2x_1x_3+x_0x_1^2x_3+x_0^2x_2x_3+x_0^2x_3^2
\end{array} 
\end{equation}
e, consequentemente, a decomposição em pesos de
$-\mathcal{E}_{4_f}$ é dada pelos monômios complementares,
\begin{equation}
\begin{array}{lcl}
-\mathcal{E}_{4_f} &=& x_1^3x_2+x_0x_1x_2^2+x_1^2x_2^2+x_0x_2^3+x_1x_2^3+x_2^4+x_1^3x_3\\
&&+x_0x_1x_2x_3+x_1^2x_2x_3+x_0x_2^2x_3+x_1x_2^2x_3+x_2^3x_3\\
&&+x_0x_1x_3^2+x_1^2x_3^2+x_0x_2x_3^2+x_1x_2x_3^2+x_2^2x_3^2+x_0x_3^3\\
&&+x_1x_3^3+x_2x_3^3+x_3^4.
\end{array} 
\end{equation}
A correspondente contribuição numérica referente a
$c_8^f(-\mathcal{E}_4)$ é dada por
$26219809342420614792105$. O ponto fixo $f =
\langle x_0, x_1^2\rangle$ produz a fração
$26219809342420614792105/381864067200$.
A contribuição total dos 24 pontos fixos é

\begin{small}
$$\ba c -\frac{8294254683787989619313}{1033322564736}
+\frac{11480507816638689194433}{88599162960}
-\frac{14346107840478803887953}{22095279360}\\\\
+\frac{16213827152205024969825}{21115987200}
-\frac{20531505056648637457965}{66503883600}
+\frac{26219809342420614792105}{381864067200}\\\\
+\frac{5140092469753813952921}{938689424640}
-\frac{7188709550595989072715}{108659350800}
+\frac{9956072891927111859681}{6396001920}\\
-\frac{10279852364804860281657}{6382787040}
+\frac{14535210787862565542283}{102263348880}
-\frac{20933013157315975732425}{759177619200}\\\\
-\frac{1535076342887470619756}{6643313994240}
+\frac{2810740530320214286200}{138095496000}
-\frac{3114225431344314975636}{102336030720}\\\\
+\frac{5567673174379196709972}{48207882240}
-\frac{6248256141284973091968}{35759465280}
+\frac{7025744180090961243900}{100989504000}\\\\
+\frac{550879843165884905497}{549474710400}
-\frac{735813438864604392697}{29436145200}
+\frac{830277686763712168477}{13736867760}\\\\
-\frac{1000713635855169897697}{11213769600}
+\frac{1137602715408000198097}{15699277440}
-\frac{1297105409181088576537}{65936965248}\\\\ = 151420.\ea $$
\end{small}

Este é o grau da subvariedade $\Sigma(\W_2{},4) \subset \P{{N_4}}$ consistindo nas superfícies de grau 4 que são singulares em alguma cônica (variável). Na figura \ref{surfaceconic} vemos o exemplo de uma quártica com uma cônica em seu lugar singular.
\begin{figure}[!ht]
\centering\ifpdf
\includegraphics[scale=0.13]{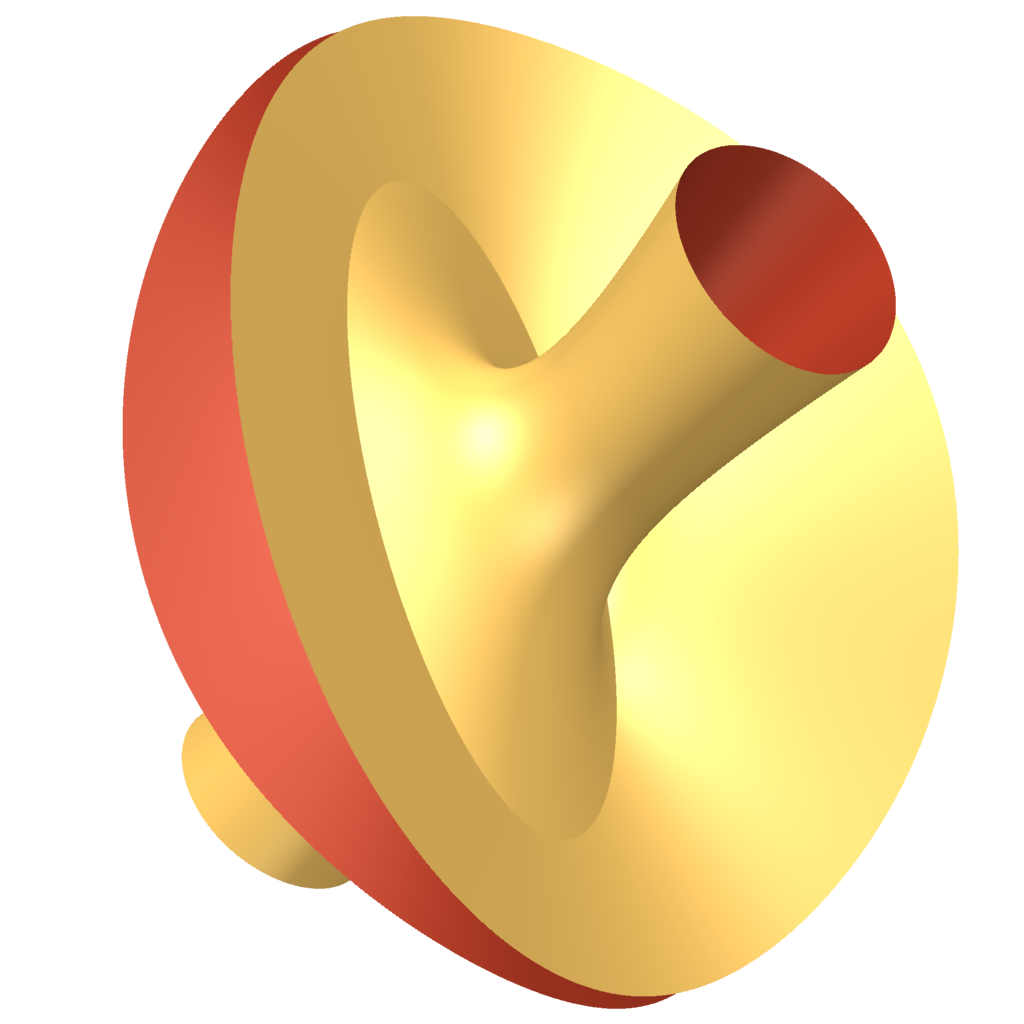}\fi
\begin{small}\caption{$-x^4-8x^2y^2-2xy^3+1/8y^4+6x^3z+x^2yz-3xy^2z+1/4y^3z-6x^2z^2-3xyz^2+3/8y^2z^2-xz^3+1/4yz^3+1/8z^4-2x^3+x^2y-2xy^2+1/4y^3+4x^2z-3xyz+1/2y^2z-xz^2+1/2yz^2+1/4z^3+x^2+8xy-7/8y^2+4xz-3/4yz-7/8z^2-y-z+2 = 0$}\label{surfaceconic}\end{small}
\end{figure}

Além disso, vide Proposição \ref{grau_polinomio},
sabemos que $\deg\Sigma(\W_2{},d)$ é um polinômio
em $d$ de grau menor do que ou igual a $24 (
=3\times 8)$. Isto   nos diz que é suficiente obtermos os graus para 25 valores distintos de $d$ e interpolar os resultados. No Apêndice \ref{ap_curva_plana} encontra-se à disposição do leitor um código para efetuar os cálculos acima com a utilização do software Macaulay2 \cite{Macaulay2}. Encontramos, como era de esperar, o mesmo resultado de (\ref{grau_conica_1}).

O código do Apêndice \ref{ap_curva_plana} utilizado nas contas
anteriores é facilmente adaptado para o cálculo
do $\deg \Sigma(\W_m{},d)$ para qualquer outro valor
de $m>1$. No entanto, ao aumentar o valor de $m$ são
necessárias mais iterações para obter o polinômio
interpolador que dá o grau de $\Sigma(\W_m{},d)$ e
isto gera um custo computacional significativo. Por exemplo, para o caso $m=3$ temos que $\dim \W_3 = 12$ e daí segue que é necessário o cálculo do $\deg\Sigma(\W_3{},d)$ para  $37 (= 3\times 12 +1)$ valores distintos de $d$ e depois efetuar a interpolação polinomial, o que fornece o resultado:
{\small
\begin{equation}\label{grau_cubica}
\ba c \deg\Sigma(\W_3{},d) = \dfrac{1}{32699842560}(13286025d^{24}\\\na9
-1038081420d^{23}+39146062158d^{22}-946074434976d^{21}\\\na9
+16407919974303d^{20}-216603408547548d^{19}+2251372103607528d^{18}\\\na9
-18776305509313968d^{17}+126579622223230407d^{16}\\\na9
-686155959955971780d^{15}+2911999863446866566d^{14}\\\na9
-8886007643094113376d^{13}+12799827743693355329d^{12}\\\na9
+50456388588134712812d^{11}-483658040042985949724d^{10}\\\na9
+2229927488252098274992d^9-7358275057877141245584d^8\\\na9
+18804143410678335462720d^7-38007885859704936084800d^6\\\na9
+60658830486712279959808d^5-75133955486596446561280d^4+\\\na9
69793667761693681135616d^3-45744106516543857328128d^2\\\na9
+18819557445986636267520d-3636764182567924531200)
\ea
\end{equation}}

Novamente, observe que o grau em (\ref{grau_cubica}) é igual a $(1+1)\times\dim \W_3$.

\chapter{Hipersuperfícies singulares
ao longo de redes de quádricas
 do tipo determinantal\label{cap_QD}}
Neste capítulo tratamos do caso de hipersuperfícies em $\P{n}$ ($n=3,4,5$) singulares ao longo de uma rede quádricas (variável) do tipo determinantal, isto é, gerada pelos menores
$2\times 2$ de uma matriz $3 \times 2$ de formas lineares. Especificamente, consideramos as famílias $\W_{twc}$ de cúbicas reversas em $\P3$, $\W_{rc}$ de cúbicas regradas em $\P 4$ e $\W_{sg}$ de 3-variedades de Segre em $\P5$, para as quais aplicamos a fórmula de resíduos de Bott no cálculo do grau da família de hipersuperfícies de grau $d$ que contém algum membro (variável) $W$ da família ($\W_{twc}, \W_{rg}$ ou $\W_{sg}$) em seu lugar singular.

\section{Superfícies singulares ao longo de cúbicas reversa\label{sec_QD_twc}}

Uma cúbica reversa  
é uma curva racional, suave de
grau 3 em $\mathbb{P}^3$.  Seu polinômio
de Hilbert é $3t+1$. Piene \& Schlessinger
\cite{Piene_Schlessinger_1985} mostraram que a
componente $\W_{twc}$ de cúbicas reversa do esquema
de Hilbert $\Hilb_{3t+1}(\mathbb{P}^3)$ é uma
variedade projetiva suave de dimensão
12. Posteriormente, Ellingsrud, Piene \& Str{\o}mme
\cite{Ellingsrud_Piene_Stromme_87} provam que a
subvariedade $\mathbb{X}$ da Grassmanniana $G(3,\mathcal{F}_2)$
formada pelas redes de quádricas do tipo
determinantal, isto é, gerada pelos menores
$2\times 2$ de uma matriz $3 \times 2$ de formas
lineares, é uma variedade suave e que a componente
$\W_{twc}$ é a explosão de $\mathbb{X}$ ao longo da subvariedade de
redes com uma componente fixa, fato este também
verificado por Vainsencher
\cite{Vainsencher_1987}. Ellingsrud \& 
Str{\o}mme \cite{Ellingsrud_Stromme_89} 
mostraram que $\mathbb{X}$ é um quociente geométrico da
variedade de
matrizes de formas lineares $2\times 3$
semiestáveis. Essa descrição permite então o
cálculo do Anel de Chow de $\mathbb{X}$ e $\W_{twc}$. Vainsencher
\& Xavier \cite{Vainsencher_Xavier_02} apresentam
uma compactificação suave explícita de um espaço
de parâmetros para a família de cúbicas reversa
adequada para a aplicação da fórmula de Bott. No presente trabalho utilizamos a 
descrição dos pontos fixos e respectivos tangentes explicitados em Ellingsrud \& Str{\o}mme \cite{Ellingsrud_Stromme_95}, \cite{Ellingsrud_Stromme_Bott_96},  para aplicarmos a fórmula de resíduos de Bott no cálculo do grau de $\Sigma(\W_{twc}{},d) \subset \mathbb{P}^{N_d}$, a família de superfícies de grau $d$ que contém algum membro $W \in \W_{twc}$ em seu lugar singular. Cabe ressaltar que o número de pontos fixos aumenta significativamente por conta da explosão adicional (vide Seção \ref{fibras_twc}) necessária para a planificação da família  formada pelos subesquemas de $\P3$ definido por $\mathcal{I}^2_W$ para algum $W \in \W_{twc}$.

Mantendo a notação coerente com a seção \ref{secao_fibrado_polinomio}, segue da Proposição \ref{prop_fibrado} e Teorema \ref{teorema_injetividade_generica} que para $d\gg 0$ existe um fibrado vetorial $\mathcal{E}_d$ tal que $\widetilde{\Sigma}(\W_{twc}{},d) = \mathbb{P}(\mathcal{E}_d)$ e o grau 
\beq\label{grau_twc_1} \deg \Sigma(\W_{twc}{},d) = \int \segre(w, \mathcal{E}_d)\cap [\W'_{twc}{}],\eeq
onde $w = dim(\W_{twc}) = 12$.

E pela fórmula de resíduos de Bott o grau em (\ref{grau_twc_1}) é obtido por
\beq\label{deg_twc_bott} \deg \Sigma(\W_{twc}{},d) = \int_{\W'_{twc}{}} \segre(12,\mathcal{E}_d) = \sum\limits_f \dfrac{c^f_{12}(-\mathcal{E}_d)}{c^f_{12}(\mathcal{T}\W'_{twc})}.\eeq 

\subsection{O espaço de parâmetros\label{espaco_parametro_twc}}

No que segue fazemos um resumo da construção de $\W_{twc}$, diagrama (\ref{diagrama_twc}), como descrito em Ellingsrud, Piene \& Str{\o}mme \cite{Ellingsrud_Piene_Stromme_87} e Vainsencher \cite{Vainsencher_1987}. Ponha $\mathbb{G}(3,\mathcal{F}_2)$ a grassmanniana de redes de quádricas em $\P3$ e $\mathbb{X}$ a subvariedade de $\mathbb{G}(3,\mathcal{F}_2)$ formada pelas redes de quádricas do tipo determinantal, isto é, gerada pelos menores $2\times 2$ de uma matriz $3 \times 2$ de formas lineares. Temos que $\W_{twc} = \widetilde{\mathbb{X}}$,
\beq \label{diagrama_twc}
\begin{gathered}
\xymatrix{
\W_{twc} &=& \widetilde{\mathbb{X}} \ar[d] &\supset& \widetilde{\mathbb{E}} \ar[d] \\
\mathbb{G}(3,\mathcal{F}_2) & \supset & \mathbb{X} &\supset & Z
}
\end{gathered}
\eeq
onde $\widetilde{\mathbb{X}}$ é a explosão de $\mathbb{X}$ ao longo da subvariedade $Z$, esta última consistindo nas redes de quádricas determinantais com um plano fixado e um ponto imerso nesse plano, isto é, são da forma $L_0\cdot \langle L_0, L_1, L_2\rangle$ com $L_0=0$ a equação do plano e $L_0=L_1=L_2=0$ a equação do ponto imerso,
\begin{center}
$Z = \left\{
\resizebox{!}{1.1cm}{
\begin{tikzpicture}
\draw (0,0) -- (1,1);
\draw (0,0) -- (2,0);
\draw (1,1) -- (3,1);
\draw (2,0) -- (3,1);
\draw (1.5, 0.5) node[black]{$\star$};
\end{tikzpicture}}\right\}$.
\end{center}

\subsection{Aplicando Bott\label{aplicando_bott_twc}}
Para a aplicação da fórmula de resíduos de Bott é necessário obter a contribuição numérica dos pontos fixos para uma escolha adequada da ação de um toro, tipicamente $\mathbb{T}=\mathbb{C}^*$ agindo diagonalmente sobre $(\mathbb{C}^4)^\vee$ via $t\circ x_i := t^{w_i}x_i$, com pesos apropriados, digamos:
\beq \label{pesos_twc} w_0=11, w_1 = 17, w_2 = 32, w_3 = 55. \eeq 

A ação de $T$ induz ações naturais sobre $\mathbb{G}(3,\mathcal{F}_2)$, $\mathbb{X}\subset\mathbb{G}(3,\mathcal{F}_2)$ e  $\widetilde{\mathbb{X}}$ de tal forma que os pontos fixos resultam isolados. 

No que segue fazemos uma descrição dos pontos fixos dessa ação, bem como a contribuição numérica do numerador e denominador em (\ref{deg_twc_bott}).

\subsubsection{Pontos fixos em \texorpdfstring{$\mathbb{X}$}{}}
De acordo com Ellingsrud \& Stromme \cite{Ellingsrud_Stromme_Bott_96}, os pontos fixos (isolados) da ação de  $\mathbb{T}$ induzida sobre $\mathbb{X}$ são projetivamente equivalentes a um dos seguintes tipos:
\begin{multicols}{2}
\begin{enumerate}[label=(\arabic*)]
\item $\langle x_0x_1, x_1x_2, x_2x_3\rangle$\label{pf_1}
\item $\langle x_0x_1, x_1x_2, x_0x_2\rangle$ \label{pf_2}
\item $\langle x_0x_1, x_2^2, x_0x_2\rangle$ \label{pf_3}
\item $\langle x_0^2, x_0x_1, x_1^2\rangle$ \label{pf_4}
\item $\langle x_0^2, x_0x_1, x_0x_2\rangle$ \label{pf_5}
\end{enumerate}
\end{multicols}
\bobs Existem vários pontos fixos sobre cada classe de isomorfismo acima. De fato, verifica-se que permutando as variáveis de um dado $\mathbb{P}^3$ o número de pontos fixos do tipo \ref{pf_1} - \ref{pf_5} são, respectivamente, 12, 4, 24, 6 e 12 \eobs

Os pontos fixos que estão sobre o centro de explosão $Z$ (vide diagrama (\ref{diagrama_twc})) são os do tipo \ref{pf_5}. Deste modo, o centro de explosão $Z$ contém 12 pontos fixos. Os demais 46 pontos fixos estão fora do centro de explosão $Z$ e, sobre estes a contribuição é calculada sobre $\mathbb{X}$, ou seja, os divisores excepcionais não contribuem nestes pontos fixos. Portanto, temos contribuições imediatas para (\ref{deg_twc_bott}), ou seja, dado $f$ em \ref{pf_1} - \ref{pf_4} temos:
$$f \Rightarrow \dfrac{c_{12}^f(-\mathcal{E}_d)}{c_{12}^f(\mathcal{T}\mathbb{X})},$$
onde o denominador $c_{12}^f(\mathcal{T}\mathbb{X})$ é o produto dos pesos da representação de $\mathbb{T}_f \mathbb{X}$ para cada um dos pontos fixos $f$ em \ref{pf_1} - \ref{pf_4}. De acordo com Ellingsrud \& Stromme \cite{Ellingsrud_Stromme_Bott_96}, Proposição 3.10, pág. 11, a representação de $\mathcal{T}_f \mathbb{X}$ para cada $f$ do tipo \ref{pf_1} - \ref{pf_4} é dada por:
\beq\label{rep_1_4_twc} f \Rightarrow\mathcal{T}_f \mathbb{X} = Hom(F,E)\otimes(x_0+x_1+x_2+x_3) - End(E) - End(F) + 1 \eeq

onde as representações $E$ e $F$ são explicitadas na Tabela \ref{table_E_F_1_4}:
\begin{table}[h!]
\centering
\begin{tabular}{|c|c|c|c|}
\hline 
Tipo & $\mathcal{I}_x$ & E & F \\ 
\hline 
(1) & $\langle x_0x_1, x_1x_2, x_2x_3\rangle$ & $x_0x_1+ x_1x_2+ x_2x_3$ & $x_0x_1x_2+x_1x_2x_3$ \\ 
\hline 
(1) & $\langle x_0x_1, x_1x_2, x_0x_2\rangle$ & $x_0x_1+ x_1x_2+ x_0x_2$ & $2x_0x_1x_2$ \\ 
\hline 
(3) & $\langle x_0x_1, x_2^2, x_0x_2\rangle$ & $x_0x_1+ x_2^2+ x_0x_2$ & $x_0x_1x_2+x_0x_2^2$ \\ 
\hline 
(4) & $\langle x_0^2, x_0x_1, x_1^2\rangle$ & $x_0^2+ x_0x_1+ x_1^2$ & $x_0x_1^2+x_0^2x_1$ \\ 
\hline 
\end{tabular} 
\caption{Representações de E e F}\label{table_E_F_1_4}
\end{table}

Intrinsecamente, $E$ e $F$ são descritos como $E = H^0(\P3,\cl I_x(2)) \subseteq H^0(\cl O_{\P3}(2))$, \break $F = \mbox{Ker}(E\otimes H^0(\cl O_{\P3}(1)) \stackrel{mult}{\longrightarrow}  H^0(\cl O_{\P3})))$, onde $\cl I_x$ indica o feixe ideal de uma cúbica reversa.
\bobs 
Estamos cometendo um abuso de notação identificando $x_i$ com o caracter $\lambda_i$ sobre $\mathbb{T}$ tal que para qualquer $t\in \mathbb{T}$ temos $t\cdot x_i = \lambda_i(t)x_i$. Assim, nas representações de $E$ e $F$ deveríamos trocar $x_i$ por $\lambda_i$, mas para efeitos computacionais (vide Apêndice \ref{ap_twc}) essa troca não causa confusão.
\eobs

\bobs Outro ponto que merece destaque é o fato de que as mesmas representações dos espaços tangentes descritas na fórmula \ref{rep_1_4_twc} podem ser obtidas utilizando os códigos  no Apêndice \ref{codigo-macaulay2_tangente_P4}, onde basta alterarmos a variável $n$ para 3.
\eobs

Tomando como exemplo o ponto fixo $f = \langle x_0x_1, x_1x_2, x_2x_3\rangle$ temos que a representação em (\ref{rep_1_4_twc}) exprime-se da seguinte forma:
\begin{small}
$$\ba {ccl}
\mathcal{T}_f\mathbb{X} &=& Hom(x_0x_1x_2+x_1x_2x_3, x_0x_1+ x_1x_2+ x_2x_3)\otimes(x_0+x_1+x_2+x_3)\\\na9 
&-& End(x_0x_1+ x_1x_2+ x_2x_3) - End(x_0x_1x_2+x_1x_2x_3)+1\\\na9
&=& (x_0x_1x_2+x_1x_2x_3)^\vee\otimes(x_0x_1+ x_1x_2+ x_2x_3)\otimes(x_0+x_1+x_2+x_3)\\\na9
&-& (x_0x_1+ x_1x_2+ x_2x_3)^\vee \otimes(x_0x_1+ x_1x_2+ x_2x_3) \\\na9
&-& (x_0x_1x_2+x_1x_2x_3)^\vee\otimes(x_0x_1x_2+x_1x_2x_3)+1\\\na9
&=& \frac{x_0}{x_2}+\frac{x_0}{x_3}+\frac{x_3}{x_1}+\frac{x_3^2}{x_0x_1}+\frac{x_3}{x_0}+\frac{x_0^2}{x_2x_3}+\frac{x_1}{x_0}+\frac{x_0}{x_1}+\frac{x_2}{x_1}+\frac{x_1}{x_2}+\frac{x_3}{x_2}+\frac{x_2}{x_3}
\ea$$
\end{small}
onde, por exemplo, $\dfrac{x_i}{x_j}$ indica o $\mathbb{T}$-espaço com peso $w_i-w_j$. A contribuição numérica do denominador em (\ref{cont_D_tilde_twc}) é obtida como o produto dos pesos da representação acima.

\bobs A contribuição numérica do numerador em  (\ref{deg_twc_bott}) sobre cada ponto fixo será descrita na seção \ref{fibras_twc}.  \eobs
\subsubsection{Pontos fixos em \texorpdfstring{$\widetilde{\mathbb{X}}$}{}}

De acordo com Ellingsrud \& Str{\o}mme \cite{Ellingsrud_Stromme_Bott_96} os pontos fixos do tipo \ref{pf_5}, $f = \langle x_0^2, x_0x_1, x_0x_2\rangle$, dão origem, a menos de equivalência projetiva, aos pontos fixos no divisor excepcional $\widetilde{\mathbb{E}}$ listados a seguir:
\begin{multicols}{2}
\begin{enumerate}[label=(5.\arabic*)]
\item $\langle x_0^2, x_0x_1, x_0x_2, x_1x_2x_3\rangle$\label{pf_5_1}
\item $\langle x_0^2, x_0x_1, x_0x_2, x_1x_2^2\rangle$
\item $\langle x_0^2, x_0x_1, x_0x_2, x_2^2x_3\rangle$
\item $\langle x_0^2, x_0x_1, x_0x_2, x_2^3\rangle$\label{pf_5_4}
\end{enumerate}
\end{multicols}
\bobs Verifica-se que permutando as variáveis de um dado $\P3$ o número de pontos fixos dos tipos \ref{pf_5_1} - \ref{pf_5_4} são, respectivamente, 12, 24, 24, 24.\eobs

Para os 84 pontos fixos provenientes de \ref{pf_5_1} - \ref{pf_5_4} a contribuição em (\ref{deg_twc_bott}) é calculada sobre $\widetilde{\mathbb{X}}$ e, nestes casos, temos que considerar a contribuição dada pelo divisor excepcional $\widetilde{\mathbb{E}}$. Sobre cada um dos pontos fixos $f'$ do tipo \ref{pf_5_1} - \ref{pf_5_4} temos a contribuição
\beq \label{cont_D_tilde_twc}f' \Rightarrow\dfrac{c_{12}^{f'}(-\mathcal{E}_d)}{c_{12}^{f'}(\mathcal{T}_{f'}\widetilde{\mathbb{X}})}.\eeq

De acordo com Ellingsrud \& Str{\o}mme \cite{Ellingsrud_Stromme_Bott_96}, Proposição 3.11, pág. 12, a representação de $\mathcal{T}_{f'}\widetilde{\mathbb{X}}$ para cada um dos pontos fixos do tipo \ref{pf_5_1} - \ref{pf_5_4} é dada por
\beq\label{rep_5_twc}\mathcal{T}_{f'}\widetilde{\mathbb{X}} = A + \mu^{-1}(B-\mu)+(x_0x_1x_2)^{-1}\mu, \eeq
onde $\mu$ é o caracter do gerador cúbico minimal, i.e, $x_1x_2x_3, x_1x_2^2, x_2^2x_3$ e $x_2^3$, respectivamente, e $A=x_0^{-1}(x_1+x_2+x_3)+x_3(x_1^{-1}+x_2^{-1})$, $B = x_1^3+x_1^2x_2+x_1^2x_3+x_1x_2^2+x_1x_2x_3+x_2^3+x_2^2x_3$. 

Por exemplo, para o ponto fixo $f' = \langle x_0^2, x_0x_1, x_0x_2, x_1x_2x_3\rangle$ a contribuição do denominador em (\ref{cont_D_tilde_twc}) é obtida como o produto dos pesos da representação de $\mathcal{T}_{f'}\widetilde{\mathbb{X}}$, onde:
\beq\label{tang_twc_blowup}\mathcal{T}_{f'}\widetilde{\mathbb{X}} = \frac{x_1}{x_0}+\frac{x_2}{x_0}+\frac{2x_3}{x_0}+\frac{x_3}{x_1}+\frac{x_3}{x_2}+\frac{x_1^2}{x_2x_3}+\frac{x_1}{x_3}+\frac{x_1}{x_2}+\frac{x_2}{x_3}+\frac{x_2^2}{x_1x_3}+\frac{x_2}{x_1}.\eeq

Uma outra forma de analisarmos as contribuições dos pontos fixos no divisor excepcional $\widetilde{\mathbb{E}}\subset \widetilde{\mathbb{X}}$ da explosão de $\mathbb{X}$ ao longo de $Z$ é lembrarmos  que este é a projetivização $\bb P (\cl N_{Z/\mathbb{X}})$ do fibrado normal $\cl N_{Z/\mathbb{X}}$. Este último é o quociente $\cl T_\bb X/\cl T_Z$. Assim, sobre cada ponto fixo $f\in Z$, a fibra $\widetilde{\bb E}_f$ é o espaço projetivo $\bb P (\cl N_{_f Z/\bb X})$.

\bobs \label{qt_pt_fixos_div_excep} Caso a decomposição do espaço normal $\cl N_{_f Z/\bb X}$ apresente somente caracteres distintos, então em $\bb P (\cl N_{_f Z/\bb X})$ haverá somente um número finito de pontos fixos sobre o ponto fixo $f\in Z$, a saber, a dimensão de $\cl N_{_f Z/\bb X}$.\eobs
Agora, podemos calcular os espaços normais sobre cada um dos 12 pontos fixos do tipo \ref{pf_5}, por exemplo, sobre $f = \langle x_0^2, x_0x_1, x_0x_2\rangle$ temos: 

$$\begin{array}{lcl}\cl N_{_f Z/\bb X} &=& \left(\frac{x_1^2}{x_0x_2}+2\frac{x_1}{x_0}+2\frac{x_2}{x_0}+\frac{x_2^2}{x_0x_1}+\frac{x_1x_3}{x_1x_2}+\frac{x_1x_3}{x_0x_2}+\frac{x_3}{x_1}+2\frac{x_3}{x_0}+\frac{x_2x_3}{x_0x_1}\right) \\
&-&\left(\frac{x_1}{x_0}+\frac{x_2}{x_0}+\frac{x_3}{x_0}+\frac{x_3}{x_1}+\frac{x_3}{x_2}\right)\\
&=& \frac{x_1^2}{x_0x_2}+\frac{x_1}{x_0}+\frac{x_2}{x_0}+\frac{x_2^2}{x_0x_1}+\frac{x_1x_3}{x_0x_2}+\frac{x_3}{x_0}+\frac{x_2x_3}{x_0x_1}
\end{array}$$

De modo análogo, efetuando o cálculo sobre os demais 11 pontos fixos, verificamos que a $\bb T$-ação induzida sobre $\widetilde{\bb X}$ tem um número finito de pontos fixos, ao todo $84 = ( 7 \times 12)$. 

Agora, precisamos estudar o espaço tangente $\cl T_{f'} \widetilde{\bb X}$ em um ponto fixo $f'$ do divisor excepcional, na fibra sobre $f\in Z$. Este espaço é dado pela decomposição 
\beq \label{decomposicao_tangente}\cl T_{f'} \widetilde{\bb X} =  \cl L_{f'}\oplus \cl T_fZ\oplus \cl T _{[\cl L_{f'}]}\bb P (\cl N_{_f Z/\bb X}),\eeq 
onde $\cl L_{f'}$ é a reta representada pelo ponto $f'$ no espaço projetivo $\bb P (\cl N_{_f Z/\bb X})$. Assim, por exemplo, sobre o ponto fixo $f = \langle x_0^2, x_0x_1, x_0x_2\rangle \in Z$ para o qual já obtemos a expressão 
$$N_{_f Z/\bb X} = \frac{x_1^2}{x_0x_2}+\frac{x_1}{x_0}+\frac{x_2}{x_0}+\frac{x_2^2}{x_0x_1}+\frac{x_1x_3}{x_0x_2}+\frac{x_3}{x_0}+\frac{x_2x_3}{x_0x_1}$$
segue que temos 7 pontos fixos, um para cada auto-espaço que aparece na decomposição acima. Tomando para $f'$ o ponto correspondente ao auto-espaço com caracter $\frac{x_3}{x_0}$, obtemos
$$\begin{array}{lcl}\cl T_{f'} \widetilde{\bb X} &=& \left(\frac{x_3}{x_0}\right)+\left(\frac{x_1}{x_0}+\frac{x_2}{x_0}+\frac{x_3}{x_0}+\frac{x_3}{x_1}+\frac{x_3}{x_2}\right)\\
&+&\left(\frac{x_1^2}{x_0x_2}+\frac{x_1}{x_0}+\frac{x_2}{x_0}+\frac{x_2^2}{x_0x_1}+\frac{x_1x_3}{x_0x_2}+\frac{x_2x_3}{x_0x_1}\right)\otimes \left(\frac{x_3}{x_0}\right)^\vee\\
&=& \frac{x_1}{x_0}+\frac{x_2}{x_0}+2\frac{x_3}{x_0}+\frac{x_3}{x_1}+\frac{x_3}{x_2}+\frac{x_1^2}{x_2x_3}+\frac{x_1}{x_3}+\frac{x_1}{x_2}+\frac{x_2}{x_3}+\frac{x_2^2}{x_1x_3}+\frac{x_2}{x_1}\end{array}.$$

Note que a representação acima é a mesma que aparece em \ref{tang_twc_blowup}. 

Por outro lado, em Vainsencher \cite{Vainsencher_1987} temos a descrição geométrica do divisor excepcional $\widetilde{\bb E}_f$ como o espaço projetivo de curvas cúbicas no plano $x_0=0$ singulares no ponto $(0:0:0:1)$, ou seja, $\widetilde{\bb E}_f$ é a projetivização do espaço gerado pelas formas cúbicas $x_1^3, x_1^2x_2, x_1x_2^2, x_2^3, x_1^2x_3, x_1x_2x_3, x_2^2x_3$. Por comparação direta da auto-decomposição dos espaços tangentes obtemos $$x_1^3\leftrightarrow \frac{x_1^2}{x_0x_2}, \, x_1^2x_2 \leftrightarrow \frac{x_1}{x_0}, \, x_2^3 \leftrightarrow \frac{x_2^2}{x_0x_1}, \, x_1^2x_3 \leftrightarrow \frac{x_1x_3}{x_0x_2}, \, x_1x_2x_3 \leftrightarrow \frac{x_3}{x_0}, \, x_2^2x_3 \leftrightarrow \frac{x_2x_3}{x_0x_1}.$$

De modo análogo, procedemos a análise das contribuições sobre cada um dos 84 pontos fixos provenientes do tipo \ref{pf_5} em $\widetilde{\bb X}$, bem como sobre os demais 46 pontos dos tipos \ref{pf_1} - \ref{pf_4}. O cálculo efetivo das contribuições é realizado utilizando os códigos escritos no Macaulay2, \cite{Macaulay2}, listados no Apêndice \ref{ap_twc}.
\subsubsection{Fibras de \texorpdfstring{$\mathcal{E}_d$}{}\label{fibras_twc}}
Precisamos agora determinar a contribuição numérica do numerador em (\ref{deg_twc_bott}), o que exige o conhecimento da representação da fibra de $\mathcal{E}_d$ sobre cada um dos pontos fixos. Temos que a família formada pelos subesquemas de $\P3$ definido por $\mathcal{I}^2_W$ para algum $W\in \W_{twc}$ não é plana. De fato, o ponto fixo do tipo \ref{pf_4} $\langle x_0^2, x_0x_1, x_1^2\rangle$ é um membro de $\W_{twc}$ (polinômio de Hilbert $p_{\W_{twc}}(t) = 3t+1$), mas seu quadrado tem polinômio de Hilbert "ruim", a saber $p_{\W{}_{twc}}=10t-10$, onde o esperado é $9t-7$. O locus dos pontos fixos do tipo \ref{pf_4}, o qual vamos denotar por $\mathbb{G}$, é a variedade contida em $\W_{twc}$ obtida como imagem da Grassmanniana de retas em $\P3$, $\mathbb{G}(2,4)$, via o mapa $\rho: \mathbb{G}(2,4) \rightarrow \W_{twc}$ que associa a cada reta $l\in \mathbb{G}(2,4)$ a cúbica reversa degenerada\- obtida como o quadrado do ideal definindo $l$.  Cabe ressaltar que $\mathbb{G}$ é disjunta do centro de explosão $Z$ (vide diagrama \ref{diagrama_twc}), assim para uma eventual explosão de $\mathbb{W}_{twc}$ ao longo de $\mathbb{G}$ podemos nos restringir a $\mathbb{X}\subset \mathbb{G}(3, \mathcal{F}_2)$ a variedade de rede de quádricas do tipo determinantal. {As contas locais mostram que $\mathbb{G}$ é o lugar de indeterminação do mapa  $\mathbb{X}\dashrightarrow  \mathbb{G}(6,\mathcal{F}_4)$ definido por $(q_1,q_2,q_3) \dashrightarrow (q_1,q_2,q_3)^2$.} 

Por exemplo, os pontos fixos do tipo \ref{pf_4}: $\langle x_i^2,
x_ix_j, x_j^2\rangle$, são mapeados em\break  $\langle x_i^4, x_i^3x_j,
x_i^2x_j^2, x_ix_j^3, x_j^4 \rangle$. Observe que
temos somente 5 monômios de grau 4, onde o
esperado é $6 = \binom{4+3}{3} -(9\times 4 -
7)$), indicando uma queda de posto do mapa em questão. Assim, para obtermos de fato um fibrado no caso $d=4$ precisamos resolver a indeterminação do mapa, o que será efetivado na sequência.

{Sejam $\mathcal{C}$ o subfibrado tautológico de posto 3 da grassmanniana de redes de quádricas $\mathbb{G}(3,\mathcal{F}_2)$ e $\mathcal{D}$ o subfibrado tautológico de posto 6 da variedade grassamanniana de subespaços 6-dimensional do espaço de formas quárticas $\mathbb{G}(6,\mathcal{F}_4)$. 
\bprop \label{prop_mu_twc} Seja $\mu: \mathcal{S}_2(\mathcal{C})|_\mathbb{X}\rightarrow \mathcal{D}$ o mapa natural induzido por multiplicação. O posto de $\mu$ fora de $\mathbb{G}$ é 6 (seis) e $\mathbb{G}$ é o esquema de zeros de $\stackrel{6}\wedge\mu$.
\eprop
\bdem Denote por $\mathfrak{Z}$ o esquema de zeros em questão. Temos que a inclusão $\mathbb{G}\subset \mathfrak{Z}$ é imediata. Seja $o = (x_0^2, x_0x_1, x_1^2)$ um representante da órbita fechada de $\mathbb{G}$ dado como o ideal gerado pelos menores $2\times 2$ da matriz $$u = \left(\begin{array}{cc}
x_0 & 0 \\ 
x_1 & x_0 \\ 
0 & x_1
\end{array}\right).$$
Mostraremos a igualdade sobre os espaços tangentes $\mathcal{T}_o\mathbb{G}=\mathcal{T}_o\mathfrak{Z}$. Isto implica que $\mathbb{G}$ é a única componente de $\mathfrak{Z}$ através de $o$. Como $\mathfrak{Z}$ é $GL_4$-invariante, qualquer componente precisa conter a única órbita fechada e, deste modo, $\mathbb{G}=\mathfrak{Z}$ como afirmado. Observe que é suficiente verificarmos que $$\dim \mathcal{T}_o\mathbb{G}\geq \dim\mathcal{T}_o\mathfrak{Z}.$$
Assim, precisamos apenas mostrar que o ideal de $\mathfrak{Z}$ em uma vizinhança do ponto $o$ contém $\mbox{codim } \mathbb{G}$ elementos com termos lineares independentes.
Para o cálculo de $\mathcal{T}_o\mathfrak{Z}$ precisamos escolher coordenadas. Para isso, c.f. Vainsencher \citep{Vainsencher_1987}, Proposição 1.2, pág. 83, temos que o mapa de fibrados sobre $\mathbb{G}(3,\mathcal{F}_2)$ induzido por multiplicação $\lambda: \mathcal{F}_2 \otimes \mathcal{F}_1 \rightarrow \mathcal{F}_3$ tem posto menor do que ou igual a 10.\\
\indent Sejam $\mathcal{U}=\mathbb{A}^{21}$ a vizinhança afim padrão de $o$ com funções coordenadas $a_{ij}$, $i=1,2,3$ e $j=1,\cdots, 7$, o que nos fornece 
$$\begin{array}{ccc}
V_1  & = & x_0^2+a_{1,1}x_0x_2+a_{1,2}x_0x_3+a_{1,3}x_1x_2+a_{1,4}x_1x_3+a_{1,5}x_2^2+a_{1,6}x_2x_3+a_{1,7}x_3^2 \\ 
V_2 & = & x_0x_1+ a_{2,1}x_0x_2+a_{2,2}x_0x_3+a_{2,3}x_1x_2+a_{2,4}x_1x_3+a_{2,5}x_2^2+a_{2,6}x_2x_3+a_{2,7}x_3^2 \\ 
V_3 & = & x_1^2+a_{3,1}x_0x_2+a_{3,2}x_0x_3+a_{3,3}x_1x_2+a_{3,4}x_1x_3+a_{3,5}x_2^2+a_{3,6}x_2x_3+a_{3,7}x_3^2
\end{array}$$

Agora, impondo a condição de posto explicitada acima (vide Apêndice \ref{codigos_singular_twc} para as contas locais feitas no Singular \cite{Singular}), obtemos
 
\begin{center}\begin{tabular}{ccp{14cm}}
$V'_1$ &=& $x_0^2+a_{1,1}x_0x_2+a_{1,3}x_1x_2+(-a_{1,3}a_{2,1}+a_{1,1}a_{2,3}-a
     _{2,3}^2+a_{1,3}a_{3,3})x_2^2+a_{1,2}x_0x_3+a_{1,4}x_1x_3+(-a_{1,4}
     a_{2,1}-a_{1,3}a_{2,2}+a_{1,2}a_{2,3}+a_{1,1}a_{2,4}-2a_{2,3}a_{2,4}+
     a_{1,4}a_{3,3}+a_{1,3}a_{3,4})x_2x_3+(-a_{1,4}a_{2,2}+a_{1,2}a_{2,4
     }-a_{2,4}^2+a_{1,4}a_{3,4})x_3^2$\\     
$V'_2$ &=& $x_0x_1+a_{2,1}x_0x_2+a_{2,3}x_1x_2+(a_{2,1}a_{2,3}-a_{1,3}a_{3,1
      })x_2^2+a_{2,2}x_0x_3+a_{2,4}x_1x_3+(a_{2,2}a_{2,3}+a_{2,1}a_{2,4
      }-a_{1,4}a_{3,1}-a_{1,3}a_{3,2})x_2x_3+(a_{2,2}a_{2,4}-a_{1,4}a_{3,
      2})x_3^2$\\      
$V'_3$ &=& $x_1^2+a_{3,1}x_0x_2+a_{3,3}x_1x_2+(-a_{2,1}^2+a_{1,1}a_{3,1}-a_{2,3
      }a_{3,1}+a_{2,1}a_{3,3})x_2^2+a_{3,2}x_0x_3+a_{3,4}x_1x_3+(-2a_{2
      ,1}a_{2,2}+a_{1,2}a_{3,1}-a_{2,4}a_{3,1}+a_{1,1}a_{3,2}-a_{2,3}a_{3,
      2}+a_{2,2}a_{3,3}+a_{2,1}a_{3,4})x_2x_3+(-a_{2,2}^2+a_{1,2}a_{3,2}-a
      _{2,4}a_{3,2}+a_{2,2}a_{3,4})x_3^2$
\end{tabular}\end{center}
o que nos fornece uma base para a restrição (trivial) do fibrado tautológico $\mathcal{C}$ à vizinhança de $(x_0^2, x_0x_1, x_1^2)$. Note que aparecem somente 12 (= $\dim W_{twc}$) indeterminadas: $a_{1,1}, a_{1,2},a_{1,3},a_{1,4},a_{2,1},a_{2,2},a_{2,3},a_{2,4}
,a_{3,1},a_{3,2},a_{3,3},a_{3,4}.$ Denotamos por 
$$\cl I = \langle a_{1,1},a_{1,2},a_{1,3},a_{1,4},a_{2,1},a_{2,2},a_{2,3},a_{2,4}
,a_{3,1},a_{3,2},a_{3,3},a_{3,4} \rangle$$ o ideal gerado por essas indeterminadas.

Consideremos na sequência uma matriz de representação local de $\mu$ adequada para o cálculo de $\stackrel 6\wedge\mu$. Para isso, escolhendo bases adequadas para $\mathcal{S}_2(\cl C)$ e $\mathcal{F}_4$ e efetuando operações elementares nas linhas módulo $\cl I^2$, obtemos uma representação de $\mu$ na forma 
\begin{center}$\left(\begin{array}{cc}
J' & A' \\ 
0 & C'
\end{array}\right)_{6\times 35},$ \end{center}
onde $J'$ é uma matriz triangular superior $5\times 5$ com entradas na diagonal principal iguais a 1.\\
\indent A matriz $C'_{1\times 35}$ contém as seguintes entradas módulo quadrados, em que entre parênteses indicamos os monômios correspondentes a cada coluna:
\begin{eqnarray}\label{monomios_fibras_twc} & -a_{3,1}(x_0^3x_2),-a_{3,2}(x_0^3x_3), (2a_{2,1}-a_{3,3})(x_0^2x_1x_2)\nonumber\\
& (2a_{2,2}-a_{3,4})(x_0^2x_1x_3),(-a_{1,1}+2a_{2,3})(x_0x_1^2x_2), \\
& (-a_{1,2}+2a_{2,4})(x_0x_1^2x_3)
-a_{1,3}(x_1^3x_2),
-a_{1,4}(x_1^3x_3)\nonumber\end{eqnarray}
\indent Isto mostra que $\mathcal{T}_o\mathfrak{Z}$ é de dimensão no máximo $12 - 8 = 4 = \dim \mathbb{G}$.
\edem
\bobs Denotemos por $\mathcal{F}_i^{(L)^j}$ o espaço de formas de grau $i$ que se anulam sobre $(L)^j$, onde $(L)$ indica o ideal definindo a reta $L$. Os monômios aparecendo em (\ref{monomios_fibras_twc}) formam uma base para o espaço de formas quárticas $\mathcal{F}_3^{(x_0,x_1)^3}\cdot \mathcal{F}_1$ módulo o espaço $\mathcal{F}_4^{(x_0,x_1)^4}$. O seu significado geométrico será esclarecido na Proposição \ref{fibra_blowup_twc}.
\eobs
\indent A próxima proposição traz à luz o fibrado de formas quárticas. Para uma rede de quádricas $\pi = (q_1,q_2, q_3) \in \mathbb{X}\setminus \mathbb{G}$, o espaço de formas quárticas $\mathcal{D}_\pi = \pi^2$ é de posto correto 6. Explodir $\mathbb{X}$ ao longo de $\mathbb{G}$ permite-nos estender a família $(\mathcal{D}_\pi)_{\pi \in \mathbb{X}\setminus \mathbb{G}}$ sobre uma compactificação  suave $ \widehat{\mathbb{X}}$ de $\mathbb{X}\setminus \mathbb{G}$
}

{\bprop\label{prop_mu_twc_bwup}
Considere o diagrama de explosão de $\mathbb{X}$ ao longo de $\mathbb{G}$
\beq 
\begin{gathered}
\xymatrix{
 &&\widehat{\mathbb{X}} \ar[d] &\supset& \widehat{\mathbb{E}} \ar[d] \\
 \mathbb{G}(3,\mathcal{F}_2)&\supset& \mathbb{X} &\supset & \mathbb{G}
}
\end{gathered}
\eeq
Então $\widehat{\mathbb{X}}$ mergulha em $\mathbb{X}\times \mathbb{G}(6,\mathcal{F}_4)$ de tal forma que o pullback $\mathcal{D}$ do subfibrado tautológico de posto 6 de $\mathcal{F}_4$ contém a imagem de $\mu_{\widehat{\mathbb{X}}}$.
\eprop
}
{\bdem
O mapa racional $\mathbb{X}\dashrightarrow \mathbb{G}(6,\mathcal{F}_4)$ definido por $\mu$, em vista da Proposição \ref{prop_mu_twc}, estende a um morfismo $\widehat{\mathbb{X}} \rightarrow \mathbb{G}(6,\mathcal{F}_4)$. Assim, temos que seu gráfico produz o mergulho desejado.
\edem
}

{A seguir descrevemos as fibras do fibrado $\mathcal{Q}$ de formas quárticas sobre o divisor excepcional $\widehat{\mathbb{E}}$.
}

\bprop\label{fibra_blowup_twc} Sejam $L = (H,H')$ uma reta, $\mathcal{F}_i^{(L)^i}$ o espaço de formas de grau $i$ que se anulam sobre $(H,H')^i$. Dado $\widehat{e}\in \widehat{\mathbb{E}}$ na fibra sobre $e = (H^2, H\cdot H',H'^2)\in \mathbb{G}$ temos que $$\mathcal{Q}_{\widehat{e}} =  \mathcal{F}_4^{(L)^4}+\langle Q \rangle,$$
para algum $Q\in \mathbb{P}\left(\frac{\mathcal{F}_3^{(L)^3}\cdot \mathcal{F}}{\mathcal{F}_4^{(L)^4}}\right)$.
\eprop
\bdem
Seja $e = (x_0^2, x_0x_1, x_1^2)$ um representante da órbita de $\mathbb{G}$. Ponha $L = (x_0,x_1)$. Temos que a imagem de $\mu_e$ é igual a $\mathcal{F}_4^{(L)^4}$. Agora, qualquer $Q\in \mathcal{F}_3^{(L)^3}\cdot \mathcal{F}$ pode ser escrito na forma $x_0^2Q_1-2x_0x_1Q_2+x_1^2Q_3$ para $Q_i\in \mathcal{F}_2$ adequadas. Forme a família a um parâmetro de redes de quádricas $$\gamma_t = \langle x_0^2+tQ_3, x_0x_1+tQ_2, x_1^2+tQ_1\rangle, t\in \mathbb{A}^1.$$

Temos que a imagem de $\mu_{\gamma_t}$ é um espaço de formas quárticas que contém $t(x_0^2Q_1-2x_0x_1Q_2+x_1^2Q_3)$.
Assim, ela contém $Q$ para todo $t$ não nulo. Deste modo, se $Q\neq 0$ módulo $\mathcal{F}_4^{(L)^4}$ segue que o espaço de formas quárticas $\mathcal{F}_4^{(L)^4}+\langle Q\rangle$ é igual a segunda coordenada de algum ponto em $\widehat{\mathbb{E}}$ estando sobre $e$. Isto mostra que a imagem do mergulho natural 

$$\begin{array}{ccc}
\mathbb{P}\left(\frac{\mathcal{F}_3^{(L)^3}\cdot \mathcal{F}}{\mathcal{F}_4^{(L)^4}}\right) & \hookrightarrow & \mathbb{G}(6,\mathcal{F}_4) \\ 
Q & \longmapsto & \mathcal{F}_4^{(L)^4} + \langle Q \rangle
\end{array} $$
está contido na imagem da fibra de $\widehat{\mathbb{E}}$ sobre $e$. E como as dimensões são as mesmas, segue a fórmula desejada.
\edem

\bobs \label{obs_fibra_twc}
Mantendo em mente a inclusão $\widehat{\mathbb{X}} \subset \mathbb{X}\times \mathbb{G}(6,\mathcal{F}_4)$ um ponto em $\widehat{\mathbb{X}}$ pode ser representado por um par $(\pi, \omega)$ consistindo de uma rede de quádricas do tipo determinantal e um espaço de quárticas, o qual, pela Proposição \ref{fibra_blowup_twc}, para $\pi = (L)^2\in \mathbb{G}$  é da forma $\mathcal{F}_4^{(L)^4}+\langle Q \rangle$ para uma única $Q\in \mathbb{P}\left(\frac{\mathcal{F}_3^{(L)^3}\cdot \mathcal{F}}{\mathcal{F}_4^{(L)^4}}\right)$. Portanto, um ponto $e'\in \widehat{\mathbb{E}}$ pode ser escrito como $\langle (L)^2, Q\rangle$.
\eobs

Passamos agora a descrever as contribuições dos pontos fixos no divisor excepcional $\widehat{\bb E}$ da explosão de $\bb X$ ao longo de $\bb G$. Sobre cada ponto fixo $f\in \bb G$, a fibra $\widehat{\bb E}_f$ é o espaço projetivo $\bb P (\cl N_{_f\bb G/\bb X})$. Assim, calculando os espaços normais sobre cada um dos 6 pontos fixos do tipo \ref{pf_4}, por exemplo, sobre $f = \langle x_0^2, x_0x_1, x_1^2\rangle$, obtemos: 

$$\begin{array}{lcl}\cl N_{_f\bb G/\bb X} &=& \left(\frac{x_0x_2}{x_1^2}+ \frac{2x_2}{x_1}+ \frac{2x_2}{x_0}+ \frac{x_1x_2}{x_0^2}+\frac{x_0x_3}{x_1^2}+ \frac{2x_3}{x_1}+ \frac{2x_3}{x_0}+ \frac{x_1x_3}{x_0^2}\right)-\left(\frac{x_2}{x_0}+\frac{x_2}{x_1}+\frac{x_3}{x_0}+\frac{x_3}{x_1}\right)\\
&=& \frac{x_0x_2}{x_1^2}+ \frac{x_2}{x_1}+ \frac{x_2}{x_0}+ \frac{x_1x_2}{x_0^2}+ \frac{x_0x_3}{x_1^2}+      \frac{x_3}{x_1}+ \frac{x_3}{x_0}+ \frac{x_1x_3}{x_0^2}\end{array}.$$

Note que a decomposição do espaço normal apresenta somente caracteres distintos, ao todo 8 (= $\dim \cl N_{_f\bb G/\bb X}$). Assim, em $\bb P (\cl N_{_f\bb G/\bb X})$ teremos 8 pontos fixos sobre $f$. A mesma verificação pode ser realizada sobre os demais 5 pontos  fixos em $\bb G$, totalizando $48 \, (= 6\times 8)$ pontos fixos sobre o divisor excepcional $\widehat{\bb E}$.

Para obtermos a contribuição do denominador em $\dfrac{c_{12}^{f'}(-\mathcal{E}_d)}{c_{12}^{f'}(\mathcal{T}_{f'}\widehat{\bb X})}$ sobre um ponto fixo $f'$ do divisor excepcional na fibra sobre $f\in \bb G$, precisamos estudar o espaço tangente $\cl T_{f'}\widehat{\bb X}$. Por exemplo, sobre o ponto fixo $f = \langle x_0^2, x_0x_1, x_1^2\rangle \in \bb G$ temos 8 pontos fixos, um para cada auto-espaço que aparece na decomposição do fibrado normal $\cl N_{_f\bb G/\bb X}$. Tomando para $f'$ o ponto fixo correspondente ao auto-espaço com caracter $\frac{x_0x_2}{x_1^2}$, obtemos:

$$\begin{array}{lcl}
\cl T_{f'} \widehat{\bb X} &=& \left(\frac{x_0x_2}{x_1^2}\right) + \left(\frac{x_2}{x_0}+\frac{x_2}{x_1}+\frac{x_3}{x_0}+\frac{x_3}{x_1}\right) \\
&+& \left(\frac{x_2}{x_1}+ \frac{x_2}{x_0}+ \frac{x_1x_2}{x_0^2}+ \frac{x_0x_3}{x_1^2}+      \frac{x_3}{x_1}+ \frac{x_3}{x_0}+ \frac{x_1x_3}{x_0^2}\right)\otimes \left(\frac{x_0x_2}{x_1^2}\right)^\vee\\
&=& \frac{x_1}{x_0}+\frac{ x_1^2}{x_0^2}+\frac{ x_1^3}{x_0^3}+\frac{x_0x_2}{x_1^2}+\frac{ x_2}{x_1}+\frac{x_2}{x_0}+\frac{x_3}{x_2}+\frac{x_1x_3}{x_0x_2}+\frac{x_1^2x_3}{x_0^2x_2}+\frac{x_1^3x_3}{x_0^3x_2}+\frac{ x_3}{x_1}+\frac{x_3}{x_0}
\end{array}.$$
E para a contribuição numérica do denominador em $\dfrac{c_{12}^{f'}(-\mathcal{E}_d)}{c_{12}^{f'}(\mathcal{T}_{f'}\widehat{\bb X})}$, obtemos o produto dos pesos da representação acima. De modo análogo, obtemos a contribuição numérica do denominador para todos os 48 pontos fixos sobre o divisor excepcional $\widehat{\bb E}$. Para os demais pontos fixos dos tipos \ref{pf_1} - \ref{pf_3} e \ref{pf_5_1} - \ref{pf_5_4} que  estão fora do centro de explosão $\mathbb{G}$ as contribuições do denominador em (\ref{deg_twc_bott}) podem ser calculadas pelas fórmulas expressas em (\ref{rep_1_4_twc}) e (\ref{rep_5_twc}).

 Por outro lado, pela descrição da Proposição \ref{fibra_blowup_twc}, temos que $\widehat{E}_{f}$ é a projetivização do espaço gerado pelas formas quárticas $x_0^3x_2, x_0^2x_1x_2, x_0x_1^2x_2, x_1^3x_2, x_0^3x_3, x_0^2x_1x_3, x_0x_1^2x_3, x_1^3x_3$. E por comparação direta da auto-decomposição dos espaços tangentes segue que
 $$ {x_0^3x_2\leftrightarrow \frac{x_0x_2}{x_1^2}}, {x_0^2x_1x_2\leftrightarrow \frac{x_2}{x_1}}, 
{x_0x_1^2x_2\leftrightarrow \frac{x_2}{x_0}}, 
{x_1^3x_2\leftrightarrow \frac{x_1x_2}{x_0^2}},$$
$$ 
{x_0^3x_3\leftrightarrow \frac{x_0x_3}{x_1^2}},
{x_0^2x_1x_3\leftrightarrow \frac{x_3}{x_1}}, 
{x_0x_1^2x_3\leftrightarrow \frac{x_3}{x_0}}, 
{x_1^3x_3\leftrightarrow  \frac{x_1x_3}{x_0^2}}$$
 
 Daí segue que podemos identificar cada um dos pontos fixos em $\bb W'$ provenientes de $f = \langle x_i^2, x_ix_j, x_j^2\rangle \in \bb G$, como:
\beq \label{pt_fixo_blowup_2_twc} f' = \langle x_i^4, x_i^3x_j, x_i^2x_j^2, x_ix_j^3, x_j^4,[Q]\rangle, \eeq
onde $[Q]$ indica a classe de uma das formas de grau 4:
\beq x_i^3x_k, x_i^2x_jx_k, x_ix_j^2x_k, x_j^3x_k, x_i^3x_l, x_i^2x_jx_l, x_ix_j^2x_l, x_j^3x_l, \eeq
que corresponde aos geradores monomiais da fibra do divisor excepcional (vide Proposição \ref{fibra_blowup_twc}).

\bobs Cada um dos ideais da forma (\ref{pt_fixo_blowup_2_twc}) possui polinômio de Hilbert igual a $9t-7$, bem como os ideais obtidos como quadrado dos ideais do tipo \ref{pf_1} - \ref{pf_3} e \ref{pf_5_1} - \ref{pf_5_4}. \eobs

Em resumo,  indicamos na Tabela \ref{tabela_pf_twc} cada um dos tipos de classe de isomorfismo de pontos fixos  em $\W_{twc}{}$ e os correspondentes pontos em $\bb W'$, bem como o número de pontos fixos em cada classe.
\begin{table}[h!]
\centering
\noindent\begin{scriptsize}
\begin{tabular}{|c|c|c|c|}
\hline 
Tipo & Pontos fixos em $\bb W_{twc}$ & Pontos fixos em $\W'_{twc}$ & $\#$ pontos \\ 
\hline \hline
(1) & $\langle x_0x_1, x_1x_2, x_2x_3\rangle$ & $\langle x_2^2x_3^2,x_1x_2^2x_3,x_0x_1x_2x_3,x_1^2x_2^2,x_0x_1^2x_2,x_0^2x_1^2\rangle$ & 12 \\ 
\hline 
(2) & $\langle x_0x_1, x_1x_2, x_0x_2\rangle$ & $\langle x_1^2x_2^2,x_0x_1x_2^2,x_0^2x_2^2,x_0x_1^2x_2,x_0^2x_1x_2,x_0^2x_1^2\rangle$ & 4 \\ 
\hline 
(3) & $\langle x_0x_1, x_2^2, x_0x_2\rangle$ & $\langle x_2^4,x_0x_2^3,x_0x_1x_2^2,x_0^2x_2^2,x_0^2x_1x_2,x_0^2x_1^2\rangle$ & 24 \\ 
\hline 
(4) & $\langle x_i^2, x_ix_j, x_j^2\rangle$ & $\langle x_i^4, x_i^3x_j, x_i^2x_j^2, x_ix_j^3, x_j^4, [Q] \rangle$ & 48 ($ = 6\times 8$) \\ 
\hline 
(5.1) & $\langle x_0^2, x_0x_1, x_0x_2, x_1x_2x_3\rangle$ &$\langle x_0^2x_2^2,x_0^2x_1x_2,x_0^3x_2,x_0^2x_1^2,x_0^3x_1,x_0^4,x_0x_1x_2^2x_3,x_0x_1^2x_2x_3,x_1^2x_2^2x_3^2\rangle$ & 12 \\ 
\hline 
(5.2) & $\langle x_0^2, x_0x_1, x_0x_2, x_1x_2^2\rangle$ &$\langle x_0^2x_2^2,x_0^2x_1x_2,x_0^3x_2,x_0^2x_1^2,x_0^3x_1,x_0^4,x_0x_1x_2^3,x_0x_1^2x_2^2,x_1^2x_2^4 \rangle$ & 24 \\ 
\hline 
(5.3) & $\langle x_0^2, x_0x_1, x_0x_2, x_2^2x_3\rangle$ &$\langle x_0^2x_2^2,x_0^2x_1x_2,x_0^3x_2,x_0^2x_1^2,x_0^3x_1,x_0^4,x_0x_2^3x_3,x_0x_1x_2^2x_3,x_2^4x_3^2 \rangle$ & 24 \\ 
\hline 
(5.4) & $\langle x_0^2, x_0x_1, x_0x_2, x_2^3\rangle$ &$\langle x_0^2x_2^2,x_0^2x_1x_2,x_0^3x_2,x_0^2x_1^2,x_0^3x_1,x_0^4,x_0x_2^4,x_0x_1x_2^3,x_2^6 \rangle$ & 24 \\ 
\hline 
\multicolumn{3}{|l|}{Total de pontos fixos} & 172\\
\hline
\end{tabular} 
\end{scriptsize}
\caption{Pontos fixos em $\W_{twc}{}$}\label{tabela_pf_twc}
\end{table}

Em referência a Tabela \ref{tabela_pf_twc}, os pontos fixos em $\bb W'$ do tipo (1) - (4) são 4-regular no sentido de Castelnuovo-Mumford, ao passo  que os pontos fixos do tipo (5.1) - (5.4) são 6-regular e por semi-continuidade de cohomologia (Hartshorne \cite{Hartshorne_1977}, Teorema 12.8) segue que qualquer feixe ideal correspondendo a um ponto em $\W'$ é 6-regular. A seguir veremos como essa informação se traduz  na obtenção da representação  das fibras de $\mathcal{E}_d$.

Para $d=4$ temos que o posto do fibrado $\mathcal{E}_4$ é dado por $rk (\mathcal{E}_4) = \binom{4+3}{3}-(9\times 4-7) = 6$ e assim, para obtermos a representação da fibra de  $\mathcal{E}_4$ sobre cada um dos pontos fixos da Tabela \ref{tabela_pf_twc}, basta tomarmos os 6 monômios de grau 4 que aparecem no ideal de cada um dos pontos fixos. Por exemplo, tomando $f= \langle x_0x_1, x_1x_2, x_2x_3 \rangle$ a representação da fibra $\mathcal{E}_{4_f}$ é dada por
$$\mathcal{E}_{4_f} = x_2^2x_3^2+x_1x_2^2x_3+x_0x_1x_2x_3+x_1^2x_2^2+x_0x_1^2x_2+x_0^2x_1^2,$$
onde os monômios indicam os geradores do espaço de formas de grau 4 com gradiente nulo sobre $f$ e $x_0^{\alpha_0}x_1^{\alpha_1}x_2^{\alpha_2}x_3^{\alpha_3}$ indica o $\mathbb{T}$-espaço com pesos $\alpha_0w_0+\alpha_1w_1+\alpha_2w_2+\alpha_3w_3$.

Na Figura \ref{Q_twc} vemos um exemplo de uma superfície quártica singular ao longo de uma cúbica reversa.
\begin{figure}[ht]
\centering\ifpdf
\includegraphics[scale=0.2]{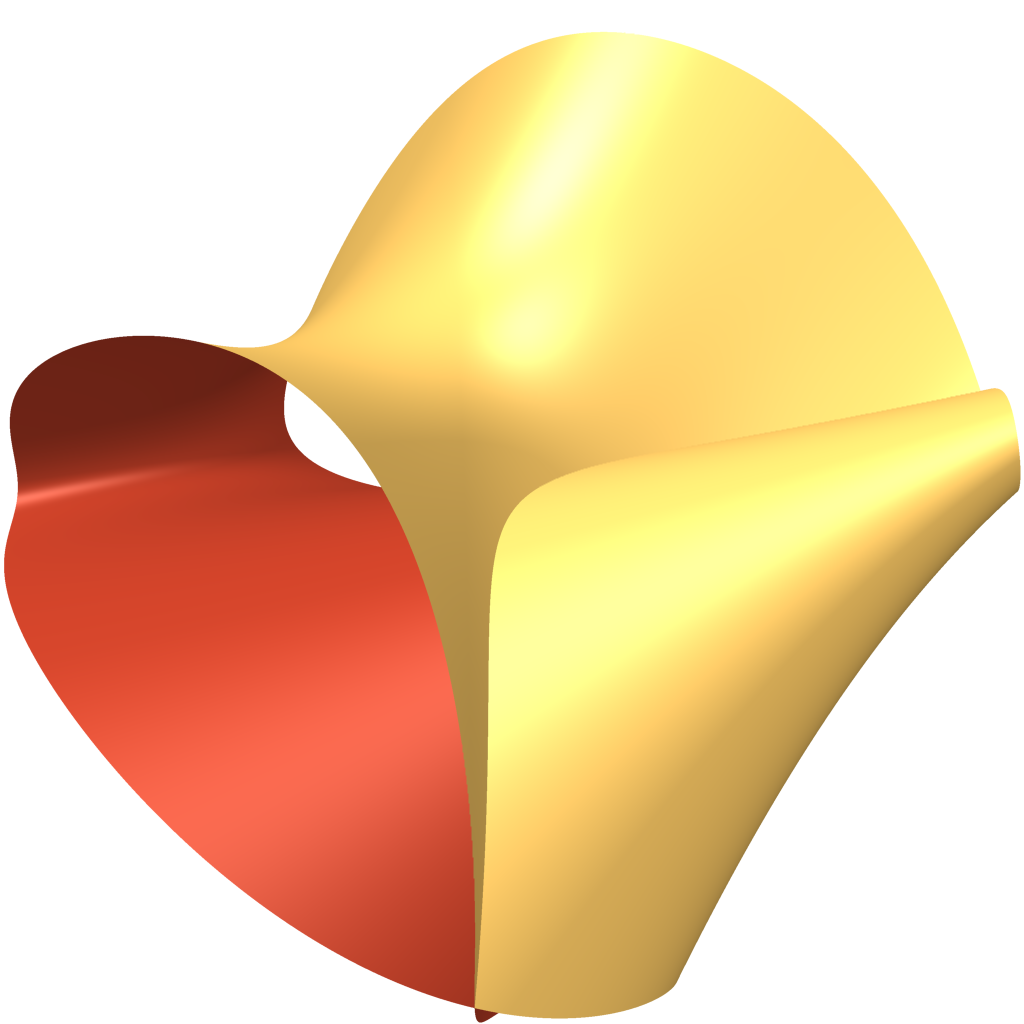}\fi
\caption{$-8 y^4+16 x y^2 z+8 y^3 z-8 x^2 z^2-8 x y z^2-8 y^2 z^2+2 x z^3+2 y z^3+4 z^4-8 x y^2+2 y^3+8 x^2 z+10 x y z-2 y^2 z-2 x z^2-8 y z^2-6 x^2+2 x y+4 y^2=0 $}
\label{Q_twc}
\end{figure}

Para $d=5$ temos que o posto do fibrado $\mathcal{E}_5$ é dado por $rk (\mathcal{E}_5) = \binom{5+3}{3}-(9\times 5-7) = 18$ e assim, a representação da fibra de $\mathcal{E}_5$ sobre cada um dos pontos fixos da Tabela \ref{tabela_pf_twc} deverá conter 18 monômios de grau 5 que são os geradores do espaço de formas de grau 5 com gradiente anulando-se sobre $f$. Esses 18 monômios são obtidos da seguinte forma: no caso dos pontos fixos 4-regulares ((1)-(3) e (4)) basta completarmos o grau, i.e, tomarmos os monômios de grau 5 que aparecem no produto $\mathcal{E}_{4_f}\otimes(x_0+x_1+x_2+x_3)$. Por exemplo, no caso de\break $f= \langle x_0x_1, x_1x_2,x_2x_3\rangle$, temos que
$$
\ba {ccl} \mathcal{E}_{5_f} &=& x_2^2x_3^3+x_2^3x_3^2+x_1x_2^2x_3^2+x_0x_2^2x_3^2+x_0x_1x_2x_3^2+x_1x_2^3x_3\\\na9 
&+&x_1^2x_2^2x_3+x_0x_1x_2^2x_3+x_0x_1^2x_2x_3+x_0^2x_1x_2x_3+x_0^2x_1^2x_3\\\na9 
&+&x_1^2x_2^3+x_1^3x_2^2+x_0x_1^2x_2^2+x_0x_1^3x_2+x_0^2x_1^2x_2+x_0^2x_1^3+x_0^3x_1^2.
\ea $$
Já para os pontos fixos do tipo (5.1) - (5.4) que são 6-regular, a representação da fibra para $\mathcal{E}_5$ é obtida tomando inicialmente os monômios de grau 5 que já estão presentes no ideal (2 ao todo), e posteriormente completar os 18 monômios com aqueles provenientes do produto $\mathcal{E}_{4_f}\otimes(x_0+x_1+x_2+x_3)$ . 
Por exemplo, para $$f^2 = \langle x_0^2x_2^2,x_0^2x_1x_2,x_0^3x_2,x_0^2x_1^2,x_0^3x_1,x_0^4,x_0x_1x_2^2x_3,x_0x_1^2x_2x_3,x_1^2x_2^2x_3^2\rangle$$ devemos tomar inicialmente os monômios 
$x_0x_1x_2^2x_3,x_0x_1^2x_2x_3$ e posteriormente completar com os 16 monômios 
$$\ba cx_0^2x_2^2x_3,x_0^2x_1x_2x_3,x_0^3x_2x_3,x_0^2x_1^2x_3,x_0^3x_1x_3,x_0^4x_3,x_0^2x_2^3,\\\na9
x_0^2x_1x_2^2,x_0^3x_2^2,x_0^2x_1^2x_2,x_0^3x_1x_2,x_0^4x_2,x_0^2x_1^3,x_0^3x_1^2,x_0^4x_1,x_0^5,\ea $$  
resultantes do produto $\mathcal{E}_{4_f}\otimes(x_0+x_1+x_2+x_3)$. Portanto, para este exemplo, temos a representação:

$$ \ba {ccl}
\mathcal{E}_{5_f} &=& x_0^2x_2^2x_3+x_0^2x_1x_2x_3+x_0^3x_2x_3+x_0^2x_1^2x_3+x_0^3x_1x_3+x_0^4x_3\\\na9
&+&x_0^2x_2^3+x_0^2x_1x_2^2+x_0^3x_2^2+x_0^2x_1^2x_2+x_0^3x_1x_2+x_0^4x_2+x_0^2x_1^3\\\na9
&+&x_0^3x_1^2+x_0^4x_1+x_0^5+x_0x_1x_2^2x_3+x_0x_1^2x_2x_3.
\ea$$

Para $d=6$ temos que o posto do fibrado é igual $rk (\mathcal{E}_6) = \binom{6+3}{3}-(9\times 6-7) = 37$, e isto nos diz que a representação da fibra de $\mathcal{E}_6$ sobre cada um dos pontos fixos da Tabela \ref{tabela_pf_twc} contém 37 monômios de grau 6, os quais são produzidos da seguinte forma: para os pontos do tipo (1) - (3) e (4) basta completarmos o grau, ou seja, tomarmos os monômios de grau 6 que aparecem no produto $\mathcal{E}_{5_f} \otimes (x_0+x_1+x_2+x_3)$. Já para os pontos do tipo (5.1) - (5.4) procedemos como no caso  $d=5$, a saber, inicialmente tomamos o monômio de grau 6 presente no ideal do ponto e posteriormente coletamos os monômios de grau 6 (ao todo 36) presentes no produto $\mathcal{E}_{5_f} \otimes (x_0+x_1+x_2+x_3)$. 

A partir da regularidade ($d>6$), para obtermos a representação da fibra de $\mathcal{E}_d$ sobre o ponto fixo $f$, basta tomarmos os monômios de grau $d$ (ao todo $\binom{d+3}{3}-(9d-7)$) presentes no produto $\mathcal{E}_{6_f} * S_{(d-6)}$, onde estamos indicando por $S_i$ a soma de todos os monômios de grau $i$ nas variáveis $x_0, x_1, x_2, x_3$.

Reunidos os dados necessários para o cálculo de $\deg \Sigma(\W_{twc}{},d)$ para cada valor de $d\gg 0$, e sabendo que o grau da família de superfícies de grau $d$ singulares ao longo de alguma cúbica reversa (variável) é polinomial em $d$ (vide Proposição \ref{grau_polinomio}) de grau menor do que ou igual a $36  (= 3\times 12)$, segue que para obtermos o polinômio em questão basta calcularmos $\deg \Sigma(\W_{twc}{},d)$, utilizando a fórmula de resíduos de Bott, para 37 valores distintos de $d$ e em seguida utilizar interpolação polinomial. O polinômio que nos dá o grau de $\Sigma(\W_{twc}{},d)$ é explicitado a seguir:

\medskip
{\small
\beq \label{grau_twc}
\ba c
\deg \Sigma(\W_{twc}{},d) = \frac{1095687}{50462720}d^{24}\!-\!\frac{19230291}{18022400}d^{23}+\frac{24114591}{985600}d^{22}\\\na9
\!-\! \frac{3932462817}{11468800}d^{21}+

\frac{73665592101}{22937600}d^{20}\!-\!\frac{23321377833}{1146880}d^{19}+\frac{4087404048523}{51609600}d^{18}\\\na9
\!-\! \frac{205245946577}{2457600}d^{17}
\!-\!\frac{79029321809671}{68812800}d^{16}
+\!\frac{2854774357217311}{309657600}d^{15}\\\na9
\!-\! \frac{6688891988137}{143360}d^{14}
\!+\!\frac{895445339622112187}{3406233600}d^{13}\!-\!\frac{4177328126526143027}{2270822400}d^{12}\\\na9
+\frac{1134029525022301939}{94617600}d^{11}\!-\!\frac{29052565860084958379}{464486400}d^{10}
+\frac{1100107099486708819}{4300800}d^9\\\na9
\!-\!\frac{31950097995158831119}{38707200}d^8\!+\!
\frac{365421773568911927}{172800}d^7\!-\!
\frac{8318629615873057099}{1935360}d^6
\\\na9
\!+\!\frac{615395937691427021}{89600}d^5\!-\!\frac{337777058982513508747}{39916800}d^4+
\frac{5167781409451915223}{665280}d^3\\\na9
\!-\!\frac{693707469384158233}{138600}d^2+\frac{466431399017887}{231}d-\mbox{\footnotesize
383398629664}.
\ea
\eeq
}

Observe que o grau do polinômio em (\ref{grau_twc}) é igual a $2\times \dim(\W_{twc}) = 2\times 12$.

No Apêndice \ref{codigos_macaulay2_twc} o leitor pode consultar os procedimentos/funções escritos no Macaulay2, \cite{Macaulay2}, utilizados para obtenção do $\deg\Sigma(\W_{twc}{},d)$ em (\ref{grau_twc}) .

\section{Hipersuperfícies singulares ao longo de uma superfície cúbica regrada em \texorpdfstring{$\P4$}{Lg}\label{sec_QD_rc}}

Uma superfície cúbica regrada em $\P4$ é uma
superfície de grau 3, com polinômio de Hilbert
dado por  $(3/2)t^2+(5/2)t+1$, formada por uma
rede de quádricas do tipo determinantal, isto é,
gerada  pelos menores $2\times 2$ de uma matriz
$3\times 2$ de formas lineares.
É imagem da explosão de \p2 em um ponto,
  mergulhada pelo sistema de cônicas com o ponto
  base (vide Beauville \cite{Beauville_1996}, Proposição IV.1 (iii), pág. 40).
Denotemos por $\W_{rc}$ a família de superfícies cúbicas regrada em $\P4$. A construção que Vainsencher \& Xavier \cite{Vainsencher_Xavier_02} utilizaram para obter uma compactificação suave do espaço de cúbicas reversas em $\P3$ também se aplica para o caso da família de subesquemas definido por redes de quádricas do tipo determinantal, e.g., superfícies cúbicas regradas em $\P4$ ou variedades de Segre em $\P5$. Em particular, segue da construção, que a família $\W_{rc}$ tem dimensão 18, enquanto a família das variedades de Segre em $\P5$, denotada por $\W_{sg}$, tem dimensão 24.

Nesta seção obtemos o $\deg\Sigma({\W{}_{rc}},d)$ da família de hipersuperfícies em $\P4$ que são singulares em algum membro $W \in \W_{rc}$. Para isso, da Proposição \ref{prop_fibrado} e do Teorema \ref{teorema_injetividade_generica}, segue que para $d\gg 0$ existe um fibrado $\mathcal{E}_d$ tal que $\widetilde{\Sigma}(\W{}_{rc},d) = \mathbb{P}(\mathcal{E}_d)$ e o grau 
\beq\label{grau_rc_1} \deg \Sigma(\W{}_{rc},d) = \int \segre(18, \mathcal{E}_d)\cap [\W'{}_{rc}]= \sum\limits_f \dfrac{c^f_{18}(-\mathcal{E}_d)}{c^f_{18}(\mathcal{T}\W'_{rc})},\eeq
onde a última igualdade segue da fórmula de resíduos de Bott.

\subsection{O espaço de parâmetros\label{espaco_parametros_rc}}
Como no caso de redes de quádricas do tipo determinantal em $\P3$, cúbicas reversas, temos que o espaço de parâmetros para $\W_{rc}$ é obtido como a explosão da rede de quádricas determinantais $\mathbb{X}\subset \mathbb{G}(3,  \mathcal{F}_2)$ ao longo da subvariedade $Z$ consistindo nas redes de quádricas com um hiperplano fixado e uma reta imersa nesse hiperplano, i.e, são da forma $L_0\cdot \langle L_0, L_1, L_2\rangle$, com $L_0=0$ a equação do hiperplano e $L_0=L_1=L_2=0$ a equação da reta imersa. O diagrama de explosão \ref{diagrama_blowup_rc} dá um resumo do exposto acima.
\beq \label{diagrama_blowup_rc}
\begin{gathered}
\xymatrix{
\W_{rc} &=& \widetilde{\mathbb{X}} \ar[d] &\supset& \widetilde{\mathbb{E}} \ar[d] \\
\mathbb{G}(3,\mathcal{F}_2) & \supset & \mathbb{X} &\supset & Z
}
\end{gathered}
\eeq

\subsection{Aplicando Bott\label{aplicando_bott_rc}}

Aqui, como nos demais casos tratados no texto, tomamos uma ação do Toro $\mathbb{T}= \mathbb{C}^*$ agindo diagonalmente sobre $(\mathbb{C}^5)^\vee$ via $t\circ x_i := t^{w_i}x_i$, com pesos apropriados, digamos:
\beq \label{pesos_rc} w_0=11, w_1 = 17, w_2 = 32, w_3 = 55, w_4=95. \eeq 

A ação de $T$ induz ações naturais sobre
$\mathbb{G}(3,\mathcal{F}_2)$,
$\mathbb{X}\subset\mathbb{G}(3,\mathcal{F}_2)$ e
$\widetilde{\mathbb{X}}$ de tal forma que os
pontos fixos sejam isolados. Na sequência fazemos
uma descrição dos pontos fixos dessa ação, bem como as contribuições numéricas em (\ref{grau_rc_1}).

\subsubsection{Pontos fixos em \texorpdfstring{$\mathbb{X}$}{}\label{subsecao_pt_1_5}}

Temos, como no caso de cúbicas reversas e a menos de permutação das variáveis de um dado $\P4$, que os pontos fixos (isolados) da ação induzida sobre $\mathbb{X}$ são projetivamente equivalentes {àqueles listados em \ref{pf_1} - \ref{pf_5} e em cada classe de isomorfismo temos, respectivamente,  60, 10, 60, 10 e 30 pontos fixos.}

De acordo com a descrição do centro de explosão $Z$ (vide (\ref{diagrama_blowup_rc})), temos que os pontos fixos do tipo \ref{pf_5} (30 pontos) estão sobre $Z$ e assim, devemos considerar a contribuição proveniente do divisor excepcional sobre estes pontos. Os demais 140 pontos fixos, aqueles oriundos dos tipos \ref{pf_1} - \ref{pf_4}, estão fora do centro de explosão $Z$ e assim podemos calcular a contribuição em (\ref{grau_rc_1}) diretamente sobre $\mathbb{X}$, ou seja, sem influência do divisor excepcional. Portanto, dado um ponto fixo $f$ dos tipos \ref{pf_1} - \ref{pf_4}, temos 
\beq\label{pontos_tipo_1-4-rc} f \Rightarrow \dfrac{c^f_{18}(-\mathcal{E}_d)}{c^f_{18}(\mathcal{T}\W_{rc})}, \eeq 
onde o denominador em (\ref{pontos_tipo_1-4-rc}) é o produto dos pesos da representação de $\mathcal{T}_f \mathbb{X}$. 

Para a determinação da representação de $\mathcal{T}_f \mathbb{X}$ sobre cada um dos pontos fixos dos tipos \ref{pf_1} - \ref{pf_5}, utilizamos a receita presente em Vainsencher \cite{Vainsencher_1987} e remetemos as contas para o Apêndice \ref{ap_cubica_regrada_P4}, onde utilizamos o software Macaulay2. Os resultados obtidos para cada uma das classes de isomorfismo \ref{pf_1} - \ref{pf_5} são os seguintes:

\begin{small}
\mathleft{
\beq \label{tg_1_rc}
\begin{aligned}
&\mathcal{T}_{\langle x_0x_1, x_1x_2, x_2x_3\rangle}\mathbb{X}= (x_0x_1x_2x_3)^\vee \otimes (x_0^3x_1+x_0^2x_1x_2+x_0x_1x_2^2+x_0^2x_1x_3\\
&+x_0x_1^2x_3+x_0^2x_2x_3+x_1^2x_2x_3+x_0x_2^2x_3+x_0x_1x_3^2+x_0x_2x_3^2+x_1x_2x_3^2\\
&+x_2x_3^3+x_0^2x_1x_4+x_0x_1x_2x_4+x_0x_1x_3x_4+x_0x_2x_3x_4+x_1x_2x_3x_4\\
&+ x_2x_3^2x_4)
\end{aligned}\eeq}%

{\setlength{\mathindent}{0cm}

\beq \label{tg_2_rc}
\begin{aligned}&\mathcal{T}_{\langle x_0x_1, x_1x_2, x_0x_2\rangle}\mathbb{X} =(x_0x_1x_2)^\vee\otimes (x_0^2x_1+x_0x_1^2+x_0^2x_2+x_1^2x_2+x_0x_2^2\\
&+x_1x_2^2+2x_0x_1x_3+2x_0x_2x_3+2x_1x_2x_3+2x_0x_1x_4+2x_0x_2x_4+2x_1x_2x_4)
\end{aligned}\eeq}%

{\setlength{\mathindent}{0cm}
\beq \label{tg_3_rc}
\begin{aligned}
&\mathcal{T}_{\langle x_0x_1, x_2^2, x_0x_2\rangle}\mathbb{X} =(x_0x_1x_2^2)^\vee\otimes(x_0x_1^3+x_0^2x_1x_2+x_0x_1^2x_2+x_0^2x_2^2+x_1^2x_2^2\\
&+x_1x_2^3+x_0x_1^2x_3+2x_0x_1x_2x_3+x_0x_2^2x_3+x_1x_2^2x_3+x_2^3x_3+x_0x_1^2x_4\\
&+2x_0x_1x_2x_4+x_0x_2^2x_4+x_1x_2^2x_4+x_2^3x_4)
\end{aligned}\eeq}%

{\setlength{\mathindent}{0cm}
\beq 
\begin{aligned}
&\mathcal{T}_{\langle x_0^2, x_0x_1, x_1^2\rangle}\mathbb{X} =(x_0^2x_1^2)^\vee\otimes (x_0^3x_2+2x_0^2x_1x_2+2x_0x_1^2x_2+x_1^3x_2+x_0^3x_3\\
&+2x_0^2x_1x_3+2x_0x_1^2x_3+x_1^3x_3+x_0^3x_4 +2x_0^2x_1x_4+2x_0x_1^2x_4+x_1^3x_4)
\end{aligned}
\eeq}%

{\setlength{\mathindent}{0cm}
\beq \label{tg_f_rc}
\begin{aligned}
&\mathcal{T}_{\langle x_0^2, x_0x_1, x_0x_2\rangle}\mathbb{X} =(x_0x_1x_2)^\vee \otimes (x_1^3+2x_1^2x_2+2x_1x_2^2+x_2^3+x_0x_1x_3+x_1^2x_3\\
&+x_0x_2x_3+2x_1x_2x_3+x_2^2x_3+x_0x_1x_4+x_1^2x_4+x_0x_2x_4+2x_1x_2x_4+x_2^2x_4)
\end{aligned}
\eeq}%
\end{small}
\subsubsection{Pontos fixos em \texorpdfstring{$\widetilde{\mathbb{X}}$}{}}

De modo análogo ao caso de cúbicas reversas (seção \ref{sec_QD_twc}), temos que os pontos fixos do tipo \ref{pf_5} dão origem, a menos de equivalência projetiva, aos pontos fixos no divisor excepcional  $\widetilde{\mathbb{E}}$ listados {em \ref{pf_5_1} à \ref{pf_5_4}, sendo a quantidade de pontos fixos em cada classe de equivalência, respectivamente, 60, 60, 120 e 60.} Além disso, a fibra do divisor excepcional $\widetilde{\bb E}$ sobre cada ponto fixo $f = \langle x_i^2, x_ix_j, x_ix_k\rangle$ é descrito geometricamente como o espaço projetivo de superfícies cúbicas no hiperplano $x_i=0$ singulares no ponto $(0:0:0:1:1)$, ou seja, $\widetilde{\bb E}_f$ é a projetivização do espaço gerado pelas formas cúbicas 
\beq \label{cubicas_rc}x_j^3, x_j^2x_k, x_jx_k^2, x_k^3, x_j^2x_l, x_jx_kx_l, x_k^2x_l,
     x_j^2x_m, x_jx_kx_m, x_k^2x_m, \eeq
     onde $x_i,x_j,x_k,x_l,x_m$ indicam as funções coordenadas de $\P4$.

Dessa maneira, cada ponto fixo sobre o divisor excepcional $\widetilde{\bb E}$ toma a forma  
\beq \label{ponto_fixo_rc_5_1_4}f' = \langle x_i^2, x_ix_j, x_ix_k, [c]\rangle,\eeq
onde $[c]$ indica a classe de uma das formas cúbicas aparecendo em (\ref{cubicas_rc}).
     
Para os 300 pontos fixos $f'$ dos tipos \ref{pf_5_1} - \ref{pf_5_4} temos a contribuição imediata 
\beq\label{contribuicao_1_blowup_rc} f' \Rightarrow \dfrac{c^{f'}_{18}(-\mathcal{E}_d)}{c^{f'}_{18}(\mathcal{T}\widetilde{\mathbb{X}})},\eeq
onde sobre o denominador devemos considerar também a contribuição dada pelo divisor excepcional $\widetilde{\mathbb{E}}$. Em virtude do divisor excepcional ser dado como o fibrado projetivo $\mathbb{P}(\mathcal{N}_{Z/\mathbb{X}})$ e sobre $f = \langle x_i^2, x_ix_j, x_ix_k \rangle$ termos a sequência de representações 
$$0 \rightarrow \mathcal{T}_f Z \rightarrow \mathcal{T}_f \mathbb{X} \rightarrow \mathcal{N}_{_f Z/\mathbb{X}}\rightarrow 0,$$
segue das fórmulas para os espaços tangentes $\mathcal{T}_f Z$ e $\mathcal{T}_f \mathbb{X}$ que 
\beq \label{equacao_normal_rc}\ba {lcl } \mathcal{N}_{_f Z/\mathbb{X}} &=& \mathcal{T}_f \mathbb{X} - \mathcal{T}_f Z\\\na9
&=& (x_ix_jx_k)^\vee \otimes (x_j^3+2x_j^2x_k+2x_jx_k^2+x_k^3+x_ix_jx_l+x_j^2x_l\\\na9
&+&x_ix_kx_l+2x_jx_kx_l+x_k^2x_l+x_ix_jx_m+x_j^2x_m+x_ix_kx_m\\\na9 
&+& 2x_jx_kx_m+x_k^2x_m) - ((x_i)^\vee \otimes(\mathcal{F}_1-x_i)\\\na9 
&+&(x_j+x_k)^\vee\otimes((\mathcal{F}_1-x_i)-(x_j+x_k)))\\\na9
&=& (x_ix_jx_k)^\vee\otimes (x_j^3+x_j^2x_k+x_jx_k^2+x_k^3+x_j^2x_l\\\na9
&+& x_jx_kx_l+x_k^2x_l+x_j^2x_m+x_jx_kx_m+x_k^2x_m).
\ea \eeq

Note que a decomposição do espaço normal apresenta todos os caracteres distintos, ao todo $10 = \dim \mathcal{N}_{_f Z/\mathbb{X}}$. Assim, em $\bb P (\mathcal{N}_{_f Z/\mathbb{X}})$ teremos 10 pontos fixos sobre $f = \langle x_i^2, x_ix_j, x_ix_k\rangle$, totalizando os $300 = (30 \times 10)$ pontos fixos mencionados no parágrafo anterior.

Passamos agora ao estudo do espaço tangente $\cl T_{f'}\widetilde{\bb X}$ sobre um ponto fixo do divisor excepcional na fibra sobre $f\in \bb X$, o qual é dado pela decomposição análoga à (\ref{decomposicao_tangente}) $$\cl T_{f'}\widetilde{\bb X} = \cl L_{f'}\oplus\cl T_f \bb X \oplus \cl T_{[\cl L_{f'}]}\bb P(\cl N_{_fZ/\bb X}).$$

Por exemplo, sobre o ponto fixo $f= \langle x_0^2,x_0x_1, x_0x_2\rangle$ segue de (\ref{equacao_normal_rc}) que 
$$\cl N_{_f Z/\bb X} = \frac{x_1^2}{x_0x_2}+\frac{x_1}{x_0}+\frac{x_2}{x_0}+\frac{x_2^2}{x_0x_1}+\frac{ x_1x_3}{x_0x_2}+ \frac{x_3}{x_0}+\frac{x_2x_3}{x_0x_1}+\frac{x_1x_4}{x_0x_2}+\frac{x_4}{x_0}+\frac{x_2x_4}{x_0x_1}$$
e tomando para $f'$ o ponto correspondente ao auto-espaço com caracter $\frac{x_1^2}{x_0x_2}$, obtemos
$$\begin{array}{lcl} \cl T_{f'}\widetilde{\bb X} &=& \left(\frac{x_1^2}{x_0x_2}\right)+\left(\frac{x_1}{x_0}+\frac{x_2}{x_0}+\frac{ x_3}{x_2}+\frac{ x_3}{x_1}+\frac{ x_3}{x_0}+\frac{x_4}{x_2}+\frac{ x_4}{x_1}+\frac{x_4}{x_0}\right)\\
&+&\left( \frac{x_1}{x_0}+\frac{x_2}{x_0}+\frac{x_2^2}{x_0x_1}+\frac{ x_1x_3}{x_0x_2}+ \frac{x_3}{x_0}+\frac{x_2x_3}{x_0x_1}+\frac{x_1x_4}{x_0x_2}+\frac{x_4}{x_0}+\frac{x_2x_4}{x_0x_1}\right)\otimes \left( \frac{x_1^2}{x_0x_2}\right)^\vee\\
&=& \frac{x_3}{x_2}+\frac{x_4}{x_2}+\frac{x_3}{x_1}+\frac{x_4}{x_1}+\frac{x_1^2}{x_0x_2}+\frac{2x_1}{x_0}+\frac{x_1x_3}{x_0x_2}+\frac{x_1x_4}{x_0x_2}+\frac{2x_2}{x_0}+\frac{2x_3}{x_0}+\frac{2x_4}{x_0}+\frac{x_2^2}{x_0x_1}+\frac{x_2x_3}{x_0x_1}+\frac{x_2x_4}{x_0x_1}.
\end{array}
 $$ 
 
Por comparação direta da autodecomposição dos espaços tangentes, obtemos

$$x_1^3 \leftrightarrow \frac{x_1^2}{x_0x_2}, x_1^2x_2 \leftrightarrow \frac{x_1}{x_0}, x_1x_2^2 \leftrightarrow \frac{x_2}{x_0} , x_2^3 \leftrightarrow \frac{x_2^2}{x_0x_1}, x_1^2x_3 \leftrightarrow \frac{ x_1x_3}{x_0x_2}, $$ 
$$x_1x_2x_3 \leftrightarrow \frac{x_3}{x_0}, x_2^2x_3 \leftrightarrow \frac{x_2x_3}{x_0x_1},
     x_1^2x_4 \leftrightarrow \frac{x_1x_4}{x_0x_2}, x_1x_2x_4 \leftrightarrow \frac{x_4}{x_0}, x_2^2x_4 \leftrightarrow \frac{x_2x_4}{x_0x_1}.$$

E para a contribuição do denominador em (\ref{contribuicao_1_blowup_rc}) obtemos o produto dos pesos da representação de $\mathcal{T}_{f'}\widetilde{\mathbb{X}}$. De modo análogo, procedemos sobre cada um dos 300 pontos fixos em $\widetilde{\bb X}$, bem como sobre os demais pontos fora do centro de explosão $Z$. 

\subsubsection{Fibras de \texorpdfstring{$\mathcal{E}_d$}{}\label{fibras_ed_rc}}
Seja  a família formada pelos  subesquemas de
$\P4$ definido por $\mathcal{I}^2_W$ para algum $W
\in \W_{rc}$. Essa família não é plana, e.g., os
pontos fixos do tipo \ref{pf_4}: $\langle x_0^2, x_0x_1,
x_1^2\rangle$ são membros legítimos de $\W_{rc}$
(polinômio de Hilbert $p_{\W_{rc}}(t) =
(3/2)t^2+(5/2)t+1$), mas seu quadrado tem
polinômio de Hilbert $5t^2 - 5t
+5$ que é diferente do esperado, a saber:
$p_{\W'_{rc}} = (9/2)t^2-(5/2)t+2$. A
subvariedade de $\W_{rc}$ que contém os pontos
fixos do tipo \ref{pf_4}, a qual vamos denotar por
$\mathbb{G}$, é obtida como imagem da
grassmanniana $\mathbb{G}(3,5)$ via o mapa $\rho:
\mathbb{G}(3,5) \rightarrow \P4$ que associa a
cada $l\in \mathbb{G}(3,5)$ a superfície cúbica
regrada limite obtida como quadrado do ideal definindo $l$. 

{Análogo} ao caso de cúbicas reversas (Proposição \ref{prop_mu_twc}) temos que a explosão de $\mathbb{X}$ - rede de quádricas do tipo determinantal - ao longo de $\mathbb{G}$, que vamos denotar por $\widehat{\mathbb{X}}$, mergulha em $\mathbb{X}\times \mathbb{G}(6,\mathcal{F}_4)$ e a fibra do fibrado de formas quárticas sobre o divisor excepcional $\widehat{\mathbb{E}}$ é descrita como na Proposição \ref{fibra_blowup_twc}.

Na sequência vamos descrever as contribuições dos pontos fixos no divisor excepcional $\widehat{\bb E}$. Sabemos que a fibra $\widehat{\bb E}_f$ sobre cada ponto fixo $f\in \bb G$ é o espaço projetivo $\bb P (\cl N_{_f \bb G/\bb X})$, e assim sobre cada um dos 10 pontos fixos em $\bb G$, por exemplo $f=\langle x_0^2, x_0x_1, x_1^2\rangle$ obtemos:

$$\begin{array}{lcl}\cl N_{_f\bb G/\bb X} &=& \left({\frac {x_{{0}}x_{{2}}}{{x_{{1}}}^{2}}}+{\frac {x_{{0}}x_{{3}}}{{x_{{1
}}}^{2}}}+{\frac {x_{{0}}x_{{4}}}{{x_{{1}}}^{2}}}+{\frac {2x_{{2}}}{
x_{{1}}}}+{\frac {2x_{{3}}}{x_{{1}}}}+{\frac {2x_{{4}}}{x_{{1}}}}+
{\frac {2x_{{2}}}{x_{{0}}}}+{\frac {2x_{{3}}}{x_{{0}}}}+{\frac 
{2x_{{4}}}{x_{{0}}}}+{\frac {x_{{1}}x_{{2}}}{{x_{{0}}}^{2}}}+{\frac {x_
{{1}}x_{{3}}}{{x_{{0}}}^{2}}}+{\frac {x_{{1}}x_{{4}}}{{x_{{0}}}^{2}}}\right) \\
&-& \left({\frac {x_{{2}}}{x_{{0}}}}+{\frac {x_{{2}}}{x_{{1}}}}+{\frac {x_{{3}}
}{x_{{0}}}}+{\frac {x_{{3}}}{x_{{1}}}}+{\frac {x_{{4}}}{x_{{0}}}}+{
\frac {x_{{4}}}{x_{{1}}}}\right)\\
&=& {\frac {x_{{0}}x_{{2}}}{{x_{{1}}}^{2}}}+{\frac {x_{{0}}x_{{3}}}{{x_{{1
}}}^{2}}}+{\frac {x_{{0}}x_{{4}}}{{x_{{1}}}^{2}}}+{\frac {x_{{2}}}{x_{
{1}}}}+{\frac {x_{{3}}}{x_{{1}}}}+{\frac {x_{{4}}}{x_{{1}}}}+{\frac {x
_{{2}}}{x_{{0}}}}+{\frac {x_{{3}}}{x_{{0}}}}+{\frac {x_{{4}}}{x_{{0}}}
}+{\frac {x_{{1}}x_{{2}}}{{x_{{0}}}^{2}}}+{\frac {x_{{1}}x_{{3}}}{{x_{
{0}}}^{2}}}+{\frac {x_{{1}}x_{{4}}}{{x_{{0}}}^{2}}}
\end{array}.
$$

Note que temos 12 auto-espaços na decomposição do fibrado normal e todos com caracteres distintos. Logo, sobre o ponto fixo $f$ temos 12 pontos fixos isolados no divisor excepcional, um para cada auto-espaço na decomposição acima. De modo análogo, verifica-se que sobre cada um dos 10 pontos fixos em $\bb G$ a decomposição do fibrado normal sempre apresenta caracteres distintos (12 ao todo) e dessa forma temos um total de $120 = (10\times 12)$ pontos fixos isolados no divisor excepcional provenientes dos pontos do tipo \ref{pf_4}.

Tomamos, por exemplo, o ponto fixo $f'$ correspondente ao auto-espaço com caracter $\frac{x_0x_2}{x_1^2}$, o que nos fornece:

$$\begin{array}{lcl}\cl T_{f'}\widehat{\bb X} &=& \left(\frac{x_0x_2}{x_1^2}\right) + \left( {\frac {x_{{2}}}{x_{{0}}}}+{\frac {x_{{2}}}{x_{{1}}}}+{\frac {x_{{3}}
}{x_{{0}}}}+{\frac {x_{{3}}}{x_{{1}}}}+{\frac {x_{{4}}}{x_{{0}}}}+{
\frac {x_{{4}}}{x_{{1}}}}\right) \\
&+& \left({\frac {x_{{0}}x_{{3}}}{{x_{{1}}}^{2}}}+{\frac {x_{{0}}x_{{4}}}{{x_{{1
}}}^{2}}}+{\frac {x_{{2}}}{x_{{1}}}}+{\frac {x_{{3}}}{x_{{1}}}}+{
\frac {x_{{4}}}{x_{{1}}}}+{\frac {x_{{2}}}{x_{{0}}}}+{\frac {x_{{3}}}{
x_{{0}}}}+{\frac {x_{{4}}}{x_{{0}}}}+{\frac {x_{{1}}x_{{2}}}{{x_{{0}}}
^{2}}}+{\frac {x_{{1}}x_{{3}}}{{x_{{0}}}^{2}}}+{\frac {x_{{1}}x_{{4}}
}{{x_{{0}}}^{2}}}
\right)\otimes \left(\frac{x_0x_2}{x_1^2}\right)^\vee\\
&=& {\frac {x_{{0}}x_{{2}}}{{x_{{1}}}^{2}}}+{\frac {x_{{2}}}{x_{{0}}}}+{
\frac {x_{{2}}}{x_{{1}}}}+{\frac {x_{{3}}}{x_{{0}}}}+{\frac {x_{{3}}}{
x_{{1}}}}+{\frac {x_{{4}}}{x_{{0}}}}+{\frac {x_{{4}}}{x_{{1}}}}+{
\frac {x_{{3}}}{x_{{2}}}}+{\frac {x_{{4}}}{x_{{2}}}}+{\frac {x_{{1}}}{
x_{{0}}}}+{\frac {x_{{1}}x_{{3}}}{x_{{0}}x_{{2}}}}+{\frac {x_{{1}}x_{{
4}}}{x_{{0}}x_{{2}}}}+{\frac {{x_{{1}}}^{2}}{{x_{{0}}}^{2}}}\\\
&+&{\frac {{
x_{{1}}}^{2}x_{{3}}}{{x_{{0}}}^{2}x_{{2}}}}+{\frac {{x_{{1}}}^{2}x_{{4
}}}{{x_{{0}}}^{2}x_{{2}}}}+{\frac {{x_{{1}}}^{3}}{{x_{{0}}}^{3}}}+{
\frac {{x_{{1}}}^{3}x_{{3}}}{{x_{{0}}}^{3}x_{{2}}}}+{\frac {{x_{{1}}}^
{3}x_{{4}}}{{x_{{0}}}^{3}x_{{2}}}}
\end{array}$$
e para a contribuição numérica do denominador em (\ref{grau_rc_1}) obtemos o produto dos pesos da representação acima. De forma análoga, obtemos a contribuição numérica sobre todos os 120 pontos fixos sobre o divisor excepcional $\widehat{\bb E}$.

Ao efetuarmos  comparações diretas da auto-decomposição dos espaços tangentes, obtemos as seguintes relações 
$$\begin{array}{c}
x_0^3x_2\leftrightarrow\frac{ x_0x_2}{x_1^2},
x_0^2x_1x_2\leftrightarrow\frac{ x_2}{x_1}, 
x_0x_1^2x_2\leftrightarrow\frac{x_2}{x_0}, 
x_1^3x_2\leftrightarrow\frac{ x_1x_2}{x_0^2}, 
x_0^3x_3\leftrightarrow\frac{ x_0x_3}{x_1^2},
x_0^2x_1x_3\leftrightarrow\frac{ x_3}{x_1}, \\
x_0x_1^2x_3\leftrightarrow\frac{ x_3}{x_0}, 
x_1^3x_3\leftrightarrow\frac{x_1x_3}{x_0^2}, 
x_0^3x_4\leftrightarrow\frac{ x_0x_4}{x_1^2}, 
x_0^2x_1x_4\leftrightarrow\frac{ x_4}{x_1},
x_0x_1^2x_4\leftrightarrow\frac{ x_4}{x_0}, 
x_1^3x_4\leftrightarrow\frac{ x_1x_4}{x_0^2}.\end{array}$$

Como no caso de cúbicas reversas, os pontos fixos do tipo \ref{pf_4}: $\langle x_i^2,
x_ix_j, x_j^2\rangle$, dão origem em $\W'$ aos
pontos da forma $\langle x_i^4, x_i^3x_j,
x_i^2x_j^2, x_ix_j^3, x_j^4\rangle$. Note que
temos somente 5 monômios de grau 4, onde o
esperado é $6  = \binom{4+4}{4} - [(9/2)\times 4^2
  - (5/2)\times 4 + 2]$. Em vista de observação análoga a \ref{obs_fibra_twc}, temos que cada um desses pontos
fixos produzem os 12 pontos fixos da forma:

\beq \label{pt_fixo_blowup_2_rc} f' = \langle x_i^4, x_i^3x_j, x_i^2x_j^2, x_ix_j^3, x_j^4,[Q]\rangle, \eeq
onde $[Q]$ indica a classe de uma das formas de grau 4: 
\beq\ba{c}  x_i^3x_k, x_i^2x_jx_k, x_ix_j^2x_k, x_j^3x_k, x_i^3x_l, x_i^2x_jx_l, \\\na9
x_ix_j^2x_l, x_j^3x_l, x_i^3x_m, x_i^2x_jx_m, x_ix_j^2x_m, x_j^3x_m \ea \eeq
que corresponde aos geradores da fibra do divisor excepcional que é dada como o espaço $\mathbb{P}\left(\frac{\mathcal{F}_3^{(x_i,x_j)^3}\cdot \mathcal{F}}{\mathcal{F}_4^{(x_i,x_j)^4}}\right)$.

\bobs Todos os pontos fixos em $\W'$ do tipo (\ref{pt_fixo_blowup_2_rc}) possuem polinômio de Hilbert $p_{\W'{}_{rc}} = (9/2)t^2-(5/2)t+2$, bem como aqueles provenientes dos tipos \ref{pf_1} - \ref{pf_3} e \ref{pf_5_1} - \ref{pf_5_4}. \eobs

Indicamos na Tabela \ref{tabela_ptfixos_rc} cada uma das classes de isomorfismos de pontos fixos em $\W{}_{rc}$ e os correspondentes pontos em $\W'_{rc}$, bem como o número de pontos fixos em cada classe.
\begin{table}[!h]
\centering
\begin{small}
\begin{tabular}{|c|c|c|c|}
\hline 
Tipo & Pontos fixos em $\W_{rc}$ & Pontos fixos em $\W'_{rc}$ & $\#$ pontos \\ 
\hline 
(1) & $\langle x_ix_j, x_jx_k, x_kx_l\rangle$ & $\langle x_i^2x_j^2,x_ix_j^2x_k,x_ix_jx_kx_l,x_j^2x_k^2,x_jx_k^2x_l,x_k^2x_l^2\rangle$ & 60 \\ 
\hline 
(2) & $\langle x_ix_j, x_jx_k, x_ix_k\rangle$ & $\langle x_i^2x_j^2,x_ix_j^2x_k,x_i^2x_jx_k,x_j^2x_k^2,x_ix_jx_k^2,x_i^2x_k^2 \rangle$ & 10 \\ 
\hline 
(3) & $\langle x_ix_j, x_jx_k, x_ix_k \rangle$ & $\langle x_i^2x_j^2,x_ix_jx_k^2,x_i^2x_jx_k,x_k^4,x_ix_k^3,x_i^2x_k^2 \rangle$ & 60 \\ 
\hline 
(4) & $\langle x_i^2, x_ix_j, x_j^2 \rangle$ & $\langle x_i^4, x_i^3x_j, x_i^2x_j^2, x_ix_j^3, x_j^4,[Q] \rangle$ & 120 (= $10\times 12$) \\ 
\hline
(5) & $\langle x_i^2, x_ix_j, x_ix_k, c \rangle$ & $\langle x_i^2\cdot\langle x_i^2,x_ix_j,x_ix_k,c,x_j^2,x_jx_k,x_k^2\rangle,x_ix_j
     c,x_ix_kc,c^2\rangle$ & 300 ($=30\times 10$) \\ 
\hline
\multicolumn{3}{|l|}{Total de pontos fixos} & 550\\
\hline
\end{tabular}\end{small}
\caption{Pontos fixos em $\W{}_{rc}$}\label{tabela_ptfixos_rc} 
\end{table}

Em referência à Tabela \ref{tabela_ptfixos_rc}, temos que os pontos fixos em $\W'$ dos tipos (1) - (4) são 4-regular, enquanto os do tipo (5) são 6-regular. Como no caso de cúbicas reversas, segue que qualquer feixe ideal correspondendo a um ponto em $\W{}_{rc}$ é 6-regular no sentido de Castelnuovo-Mumford. 

O procedimento de obtenção da representação das fibras de $\mathcal{E}_d$ é análogo ao apresentado na seção \ref{sec_QD_twc}, onde, a partir da regularidade 6, obtemos a representação da fibra de $\mathcal{E}_d$ sobre o ponto fixo $f$ tomando os monômios de grau $d$ (ao todo $\binom{d+4}{4}- [(9/2)d^2-(5/2)d+2]$) presentes no produto $\mathcal{E}_{6_f}*\mathcal{S}_{d-6}$. 

Reunido as informações necessárias para o cálculo do $\deg \Sigma(\W{}_{rc},d)$ via fórmula de resíduos de Bott e sabendo que o grau anterior é polinomial em $d$ (vide Proposição \ref{grau_polinomio}) de grau no máximo $72 ( = 4 \times 18)$, segue que é suficiente encontrarmos $\deg \Sigma(\W{}_{rc},d)$ para 73 valores diferentes de $d$ e, posteriormente, interpolar os resultados para obter o polinômio em questão. O polinômio que nos dá o grau de $\Sigma(\W{}_{rc},d)$ é explicitado a seguir:

\begin{scriptsize}
\beq\label{grau_cubica_regrada_P4}
\ba c
\deg \Sigma(\W{}_{rc},d) = \frac{1089331}{2820745970948505600}d^{54}-
\frac{4609327}{138135296519700480}d^{53}\\\na9
+\frac{17053361977}{12432176686773043200}d^{52}-\frac{44006738257}{1243217668677304320}d^{51}+\frac{43540862009}{68559797904998400}d^{50}\\\na9
-\frac{6776065867607}{822717574859980800}d^{49}+\frac{25203282464989}{329087029943992320}d^{48}-\frac{95461703632727}{205679393714995200}d^{47}\\\na9
+\frac{3121945759267787}{3290870299439923200}d^{46}+\frac{13975371538743871}{987261089831976960}d^{45}
-\frac{1762263793046822003}{9872610898319769600}d^{44}\\\na9
+\frac{1571373547792223293}{1645435149719961600}d^{43}-\frac{18657333817850689}{21095322432307200}d^{42}-\frac{21162893089184824063}{822717574859980800}d^{41}\\\na9
+\frac{8817237395388371983}{42070785078067200}d^{40}-\frac{7285835577039579827299}{7404458173739827200}d^{39}+\frac{18439965173115436460101}{2278294822689177600}d^{38}\\\na9
-\frac{30625726302752154570146789}{251751577907154124800}d^{37}
+\frac{286671605346783151488709819}{201401262325723299840}d^{36}\\\na9
-\frac{5957731889573498708183240461}{503503155814308249600}d^{35}
+\frac{946219385360559194318492423}{13078004047124889600}d^{34}\\\na9
-\frac{28843644632003758667785804741}{88853498084877926400}d^{33}
+\frac{3586612308873070845414316631}{3702229086869913600}d^{32}\\\na9
-\frac{2772990057804229211772760003}{3173339217317068800}d^{31}
-\frac{173239617944054456458227898277}{17770699616975585280}d^{30}\\\na9
+\frac{3107360934070968268891455300733}{44426749042438963200}d^{29}
-\frac{1302777164405876523072798778669}{4936305449159884800}d^{28}\\\na9
+\frac{2175543494720246680252051667789}{3748506950455787520}d^{27}
-\frac{15324266643945858395023213928441}{88853498084877926400}d^{26}\\\na9
-\frac{583723723691983350730395768869707}{133280247127316889600}d^{25}
+\frac{295008612506350533900867771909281}{14808916347479654400}d^{24}\\\na9
-\frac{72882298518045984492971696381249}{1514548262810419200}d^{23}
+\frac{3179423312365559691881647284007591}{59235665389918617600}d^{22}\\\na9
+\frac{65074915758634148942372090942475703}{799681482763901337600}d^{21}
-\frac{17650658027740832446748837419090939}{33566877054287216640}d^{20}\\\na9
+\frac{43155219287681067897344483362302109}{35402565643193548800}d^{19}
-\frac{30042531700267289895379997718912521}{22377918036191477760}d^{18}\\\na9
-\frac{749894075579299475576086383836784223}{906305680465754849280}d^{17}
+\frac{1152884114126290978903651885817821}{176296623184281600}d^{16}\\\na9
-\frac{679247544279215190070362388445065693}{49980092672743833600}d^{15}
+\frac{27244209645180356835326895182601977}{1851114543434956800}d^{14}\\\na9
-\frac{14180655522525890878698424573977769}{8330015445457305600}d^{13}
-\frac{17786673868531949329900173945074227}{694167953788108800}d^{12}\\\na9
+\frac{8140256480874854682039834827204717}{148750275811737600}d^{11}
-\frac{15847193428252892198587722393037621}{231389317929369600}d^{10}\\\na9
+\frac{51203085967146132778275681925029671}{851933397830860800}d^{9}
-\frac{415833099791358148948760413114949}{10846374277939200}d^{8}\\\na9
+\frac{190922280640278098795730933090799}{10846374277939200}d^{7}
-\frac{47833769039838754264953305641}{8608233553920}d^{6}\\\na9
+\frac{2764737243980163013076109790463}{2560949482291200}d^{5}
-\frac{1553358364438869321892260077}{17784371404800}d^{4}\\\na9
-\frac{1981299728200259795937983}{242514155520}d^{3}
+\frac{15743878343562160667}{7623616}d^{2}\\\na9
-\frac{655521591855018725}{7351344}d+4625512425
\ea
\eeq
\end{scriptsize}

Observe que o grau em (\ref{grau_cubica_regrada_P4}) é $54 =  (2+1)\times 18$.

No Apêndice \ref{ap_cubica_regrada_P4} encontra-se a disposição do leitor os procedimentos e funções escritos no Macaulay2 \cite{Macaulay2} e Maple \cite{Maple_2015}   utilizados para obter o $\deg \Sigma(\W{}_{rc},d)$ em (\ref{grau_cubica_regrada_P4}).

\section{Hipersuperfícies singulares ao longo de uma 3-variedade de Segre em \texorpdfstring{$\mathbb{P}^5$}{P5}\label{sec_QD_sg}}

A variedade de Segre em $\P5$ é uma 3-variedade obtida pelo mergulho de Segre $\P1\times \P2 \rightarrow \P5$. Seu polinômio de Hilbert é dado por $(1/2)t^3+2t^2+(5/2)t +1$. É bem conhecido (cf. Harris \cite{Harris_1992}) que a variedade de Segre é do tipo determinantal, definida por uma rede de quádricas obtidas pelos menores $2\times 2$ de uma matriz $2\times 3$ de formas lineares. Denotemos por $\W_{sg}$ a família das variedades de Segre em $\P5$, a qual tem dimensão 24 (cf.  Vainsencher \& Xavier \cite{Vainsencher_Xavier_02}). 

O propósito desta seção é a obtenção do $\deg (\W{}_{sg},d)$ da família de hipersuperfícies de grau $d$ em $\P5$ que são singulares em algum membro (variável) $W \in \W_{sg}$. De forma análoga à seção \ref{sec_QD_twc}, para $d \gg 0$ existe um fibrado $\mathcal{E}_d$ tal que $\widetilde{\Sigma}(\W{}_{sg}, d) = \mathbb{P}(\mathcal{E}_d)$. {O $\deg\Sigma(\W{}_{sg}, d)$ é obtido pela fórmula explicitada em (\ref{grau_geral}), a saber 
\beq \label{f_grau_v_segre}\deg\Sigma(\W{}_{sg} = \int \segre(24, \mathcal{E}_d)\cap [\W'{}_{rc}]= \sum\limits_f \dfrac{c^f_{24}(-\mathcal{E}_d)}{c^f_{24}(\mathcal{T}\W'_{rc})},\eeq
onde a última igualdade segue da fórmula de resíduos de Bott. }

\subsection{O espaço de parâmteros\label{espaco_parametros_sg}}

Analogamente ao caso de superfícies cúbicas
regradas em $\P4$ (vide seção \ref{sec_QD_rc}),
temos que $\W_{sg}$ é obtido como a explosão da
variedade das redes de quádricas determinantais $\mathbb{X}\subset \mathbb{G}(3,\mathcal{F}_2)$ ao longo da subvariedade $Z$, consistindo nas redes de quádricas da forma $L_0\cdot\langle L_0, L_1, L_2\rangle$, com $L_i$ sendo formas lineares. No diagrama (\ref{diagrama_blowup_segre}) temos um resumo do exposto:
\beq \label{diagrama_blowup_segre}
\begin{gathered}
\xymatrix{
\W_{sg} &=& \widetilde{\mathbb{X}} \ar[d] &\supset& \widetilde{\mathbb{E}} \ar[d] \\
\mathbb{G}(3,\mathcal{F}_2) & \supset & \mathbb{X} &\supset & Z
}
\end{gathered}
\eeq

\subsection{Aplicando Bott\label{aplicando_bott_sg}} 
Novamente, tomamos o toro $\mathbb{T} = \mathbb{C}^*$ agindo diagonalmente sobre $(\mathbb{C}^6)^\vee$ via $t\circ x_i = t^{w_i}x_i$, com pesos apropriados, digamos:
\beq w_0=11, w_1=17, w_2 = 32, w_3 = 55, w_4 = 95, w_5 = 160. \eeq

Temos ações naturais induzidas sobre $\mathbb{G}(3, \mathcal{F}_2)$, $\mathbb{X}\subset \mathbb{G}(3, \mathcal{F}_2)$  e $\widetilde{\mathbb{X}}$ de forma que os pontos fixos sejam isolados. 
{A descrição dos pontos fixos e as contribuições numéricas necessárias a aplicação da fórmula de resíduos de Bott seguem de forma totalmente análoga à seção \ref{sec_QD_rc}. Por exemplo, os  pontos fixos da ação induzida sobre $\mathbb{X}$ são projetivamente equivalentes aqueles listados em \ref{pf_1} - \ref{pf_5}, alterando somente a quantidade de pontos fixos em cada classe de equivalência, a saber: 180, 20, 120, 15 e 60. Para os 335 pontos, oriundos dos tipos \ref{pf_1} - \ref{pf_4}, que estão fora do centro de explosão o cálculo das contribuições na fórmula de resíduos de Bott é feito diretamente sobre $\mathbb{X}$ (vide subseção \ref{subsecao_pt_1_5}), com os pontos fixos do tipo \ref{pf_4} tratados como na subseção \ref{fibras_ed_rc}. A representação do $\mathcal{T}_f\mathbb{X}$ em cada ponto fixo dos tipos \ref{pf_1} - \ref{pf_5} segue o mesmo procedimento do Apêndice \ref{ap_cubica_regrada_P4}, trocando-se o valor de $n=4$ para $n=5$ no código da seção \ref{codigo-macaulay2_tangente_P4}\ . Já os pontos fixos do tipo \ref{pf_5}: $\langle x_i^2, x_ix_j, x_ix_k\rangle$, dão origem aos pontos fixos no divisor excepcional $\widetilde{\mathbb{E}}$, da forma:
\beq\label{pontos_fixos_tipo_5_sg} f' = \langle x_i^2, x_ix_j, x_ix_k, [c]\rangle, \eeq
onde $[c]$ indica a classe de uma das seguintes formas cúbicas:
$$x_j^3, x_j^2x_k, x_jx_k^2, x_k^3, x_j^2x_l, x_jx_kx_l, x_k^2x_l,x_j^2x_m, x_jx_kx_m, x_k^2x_m, x_j^2x_n, x_jx_kx_n, x_k^2x_n,$$ correspondendo aos geradores monomiais da fibra do divisor excepcional.
}

Para os pontos fixos $f'$ temos a contribuição imediata
\beq\label{contribuicao_1_blowup_sg} f' \Rightarrow \dfrac{c^{f'}_{24}(-\mathcal{E}_d)}{c^{f'}_{24}(\mathcal{T}\widetilde{\mathbb{X}})},\eeq
em que sobre o denominador devemos considerar a contribuição proveniente do divisor excepcional $\widetilde{\mathbb{E}}$. Sobre $f = \langle x_i^2, x_ix_j, x_ix_k \rangle$ temos a sequência de representações 

$$0 \rightarrow \mathcal{T}_f Z \rightarrow \mathcal{T}_f \mathbb{X} \rightarrow \mathcal{N}_{_f Z/\mathbb{X}}\rightarrow 0,$$

e das fórmulas para os espaços tangentes $\mathcal{T}_f Z$ e $\mathcal{T}_f \mathbb{X}$ segue que

{\setlength{\mathindent}{0cm}
\beq \ba {lcl } \mathcal{N}_{_f Z/\mathbb{X}} &=& \mathcal{T}_f \mathbb{X} - \mathcal{T}_f Z\\\na9
&=& (x_ix_jx_k)^\vee \otimes (x_j^3+2x_j^2x_k+2x_jx_k^2+x_k^3+x_ix_jx_l+x_j^2x_l+x_ix_kx_l\\\na9
&+&2x_jx_kx_l+x_k^2x_l+x_ix_jx_m+x_j^2x_m+x_ix_kx_m+ 2x_jx_kx_m+x_k^2x_m\\\na9
&+&x_ix_jx_n+x_j^2x_n+x_ix_kx_n+2x_jx_kx_n+x_k^2x_n- ((x_i)^\vee \otimes(\mathcal{F}_1-x_i)\\\na9
&+&(x_j+x_k)^\vee\otimes((\mathcal{F}_1-x_i)-(x_j+x_k)))\\\na9
&=& (x_ix_jx_k)^\vee \otimes (x_j^3+x_j^2x_k+x_jx_k^2+x_k^3+x_j^2x_l+x_jx_kx_l+x_k^2x_l\\\na9 
&+&x_j^2x_m+x_jx_kx_m+x_k^2x_m+x_j^2x_n+x_jx_kx_n+x_k^2x_n).
\ea \eeq

Como no caso de cúbicas reversas e superfícies cúbicas regradas, temos que por comparação direta da auto-decomposição dos espaços tangentes segue que 
$$c \leftrightarrow \langle c\rangle = \dfrac{c}{x_ix_jx_k},$$
onde $\langle c\rangle$ indica um caracter específico  da decomposição em auto-espaços do fibrado normal. 

Para a contribuição do denominador em (\ref{contribuicao_1_blowup_sg}) obtemos o produto dos pesos da representação
\beq\label{tg_5_sg}\ba {ccl} 
\mathcal{T}_{f'}\widetilde{\mathbb{X}} &=& \mathcal{T}_f\mathbb{Z} + \mathcal{T}_{[c]}\mathbb{P}(\mathcal{N}_{\mathbb{Z}/\mathbb{X}})+\langle c\rangle\\\na9
&=& (x_i)^\vee\otimes(\mathcal{F}_1-x_i)+(x_j+x_k)^\vee\otimes((\mathcal{F}_1-x_i)-(x_j+x_k))\\\na9
&+&(\mathcal{N}_{_f\mathbb{Z}/\mathbb{X}} - \langle c\rangle)\otimes \langle c\rangle^\vee + \langle c\rangle.
\ea \eeq

\subsubsection{Fibras de \texorpdfstring{$\mathcal{E}_d$}{Ed}\label{fibras_ed_segre}}

Análogo à seção \ref{sec_QD_rc}\ , a família formada pelos  subesquemas de $\P5$ definido por $\mathcal{I}^2_W$ para algum $W \in \W_{sg}$ não é plana, e.g., os pontos fixos do tipo \ref{pf_4}\ : $\langle x_0^2, x_0x_1, x_1^2\rangle$ são membros legítimos de $\W_{sg}$ (polinômio de Hilbert $p_{\W_{sg}}(t) = (1/2)t^3+2t^2+(5/2)t+1$), mas seu quadrado tem polinômio de Hilbert $(5/3)t^3+(10/3)t+1$ que é diferente do esperado, a saber: $p^2_{\W'_{sg}} = (3/2)t^3+t^2+(3/2)t+2$. A subvariedade de $\W_{sg}$ que contém os pontos fixos do tipo \ref{pf_4}, a qual vamos denotar por $\mathbb{G}$, é obtida como imagem da grassmanniana $\mathbb{G}(4,6)$ via o mapa $\rho: \mathbb{G}(4,6) \rightarrow \P5$ que associa a cada $l\in \mathbb{G}(4,6)$ a variedade de Segre limite obtida como quadrado do ideal definindo $l$. 

Análogo à Proposição \ref{prop_mu_twc} temos que a explosão de $\mathbb{X}$ - rede de quádricas do tipo determinantal - ao longo de $\mathbb{G}$, que vamos denotar por $\widehat{\mathbb{X}}$, mergulha em $\mathbb{X}\times \mathbb{G}(6,\mathcal{F}_4)$ e a fibra do fibrado de formas quárticas sobre o divisor excepcional $\widehat{\mathbb{E}}$ é descrita como na Proposição \ref{fibra_blowup_twc}.

A contribuição dos pontos fixos no divisor excepcional $\widehat{\bb E}$, provenientes dos pontos fixos da forma $f = \langle x_i^2, x_ix_j, x_j^2\rangle$, segue os mesmos princípios como no caso de cúbicas reversas e superfície cúbica regrada. Primeiramente, a fibra $\widehat{\bb E}_f$ é o espaço projetivo $\bb P (\cl N_{_f\bb G/ \bb X})$, em que a decomposição em auto-espaços de $\cl N_{_f\bb G/ \bb X}$ é dada por

\break 

{\setlength{\mathindent}{0cm}
\beq\ba {ccl}  
\mathcal{N}_{_f\mathbb{G}/\bb X} &=& \mathcal{T}_f \W - \mathcal{T}_f\mathbb{G}\\\na9
&=& (x_i^2x_j^2)^\vee\otimes (x_i^3x_k+2x_i^2x_jx_k+2x_ix_j^2x_k+x_j^3x_k+x_i^3x_l+2x_i^2x_jx_l\\\na9
&+& 2x_ix_j^2x_l+x_j^3x_l+x_i^3x_m+2x_i^2x_jx_m+2x_ix_j^2x_m+x_j^3x_m + x_i^3x_n\\\na9
&+&2x_i^2x_jx_n+2x_ix_j^2x_n+x_j^3x_n)- (x_i+x_j)^\vee\otimes(\mathcal{F}_1-(x_i+x_j))\\\na9
&=& (x_i^2x_j^2)^\vee (x_i^3x_k+x_i^2x_jx_k+x_ix_j^2x_k+x_j^3x_k+x_i^3x_l+x_i^2x_jx_l+x_ix_j^2x_l\\\na9
&+& x_j^3x_l+x_i^3x_m+x_i^2x_jx_m+x_ix_j^2x_m+x_j^3x_m+x_i^3x_n+x_i^2x_jx_n+x_ix_j^2x_n\\\na9
&+&x_j^3x_n).
\ea \eeq}

Note que os caracteres referentes à decomposição em auto-espaços de $\mathcal{N}_{_f\mathbb{G}/\bb X}$ são todos distintos (16 ao todo), donde seque que sobre cada ponto fixo $f\in \bb G$, temos 16 pontos fixos isolados sobre o divisor excepcional $\widehat{\bb E}$. Daí segue que temos ao todo $240 (= 15\times 16)$ pontos fixos sobre $\widehat{\bb E}$. Sobre cada um desses 240 pontos fixos o espaço tangente $\cl T_{f'}\widehat{\bb X}$ apresenta a seguinte decomposição em auto-espaços:

\beq\label{tangente_blowup_2_v_segre}\ba {ccl} 
\mathcal{T}_{f'}\widehat{\bb X} &=& \mathcal{T}_f\mathbb{G} + \mathcal{T}_{[q]}\mathbb{P}(\mathcal{N}_{\mathbb{G}/\bb X})+\langle q\rangle\\\na9
&=& (x_i+x_j)^\vee\otimes(\mathcal{F}_1 - (x_i+x_j))+(\mathcal{N}_{_f\mathbb{G}/\bb X} - \langle q\rangle)\otimes \langle q\rangle^\vee + \langle q\rangle.
\ea \eeq
Além disso, da comparação direta da decomposição dos espaços tangentes, temos as seguinte relações
$$q \leftrightarrow \langle q\rangle = \dfrac{q}{x_i^2x_j^2},$$
onde $\langle q\rangle$ indica um caracter na decomposição do fibrado normal e $q$ é a correspondente forma de grau 4 relacionada a esse caracter dentre as listadas abaixo:

\beq\ba c x_i^3x_k, x_i^2x_jx_k, x_ix_j^2x_k, x_j^3x_k, x_i^3x_l,x_i^2x_jx_l, x_ix_j^2x_l,x_j^3x_l,\\ x_i^3x_m, x_i^2x_jx_m, x_ix_j^2x_m, x_j^3x_m, x_i^3x_n, x_i^2x_jx_n, x_ix_j^2x_n,x_j^3x_n
\ea\eeq

A contribuição numérica do denominador em (\ref{f_grau_v_segre}) referente aos pontos fixos no divisor excepcional $\widehat{\bb E}$ é obtida pelo produto dos pesos da decomposição em (\ref{tangente_blowup_2_v_segre}).

Indicamos na Tabela \ref{tabela_ptfixos_sg} cada uma das classes de isomorfismos de pontos fixos em $\W{}_{sg}$ e os correspondentes pontos em $\W'{sg}$, bem como o número de pontos em cada classe.

\begin{table}[!h]
\centering
\begin{small}
\begin{tabular}{|c|c|c|c|}
\hline 
Tipo & Pontos fixos em $\W_{sg}$ & Pontos fixos em $\W'_{sg}$ & $\#$ pontos \\ 
\hline 
(1) & $\langle x_ix_j, x_jx_k, x_kx_l\rangle$ & $\langle x_i^2x_j^2,x_ix_j^2x_k,x_ix_jx_kx_l,x_j^2x_k^2,x_jx_k^2x_l,x_k^2x_l^2\rangle$ & 180 \\ 
\hline 
(2) & $\langle x_ix_j, x_jx_k, x_ix_k\rangle$ & $\langle x_i^2x_j^2,x_ix_j^2x_k,x_i^2x_jx_k,x_j^2x_k^2,x_ix_jx_k^2,x_i^2x_k^2 \rangle$ & 20 \\ 
\hline 
(3) & $\langle x_ix_j, x_jx_k, x_ix_k \rangle$ & $\langle x_i^2x_j^2,x_ix_jx_k^2,x_i^2x_jx_k,x_k^4,x_ix_k^3,x_i^2x_k^2 \rangle$ & 120 \\ 
\hline 
(4) & $\langle x_i^2, x_ix_j, x_j^2 \rangle$ & $\langle x_i^4, x_i^3x_j, x_i^2x_j^2, x_ix_j^3, x_j^4,[q] \rangle$ & 240 (= $15\times 16$) \\ 
\hline
(5) & $\langle x_i^2, x_ix_j, x_ix_k, c \rangle$ & $\langle x_i^2\cdot\langle x_i^2,x_ix_j,x_ix_k,c,x_j^2,x_jx_k,x_k^2\rangle,x_ix_j
     c,x_ix_kc,c^2\rangle$ & 780 ($=60\times 13$) \\ 
\hline
\multicolumn{3}{|l|}{Total de pontos fixos} & 1340\\
\hline
\end{tabular}\end{small}
\caption{Pontos fixos em $\W{}_{sg}$}\label{tabela_ptfixos_sg} 
\end{table}

\bobs Todos os pontos fixos em $\W'_{sg}$ possuem polinômio de Hilbert $p_{\W'{}_{sg}} = (3/2)t^3+t^2+(3/2)t+2 $.\eobs

Em referência a Tabela \ref{tabela_ptfixos_sg}, temos que os pontos fixos em $\W'_{sg}$ dos tipos (1) - (4) são 4-regular, enquanto os do tipo (5) são 6-regular. Como no caso de cúbicas reversas, segue que qualquer feixe ideal correspondendo a um ponto em $\W'_{sg}$ é 6-regular no sentido de Castelnuovo-Mumford. 

O procedimento de obtenção da representação das fibras de $\mathcal{E}_d$ é análogo ao apresentado na seção \ref{sec_QD_twc}, onde, a partir da regularidade 6, obtemos a representação da fibra de $\mathcal{E}_d$ sobre o ponto fixo $f$ tomando os monômios de grau $d$ (ao todo $\binom{d+5}{5}- [(3/2)\times d^3+d^2+(3/2)\times d+2]$) presentes no produto $\mathcal{E}_{6_f}*\mathcal{S}_{d-6}$. 

Temos as informações necessárias para o cálculo do $\deg \Sigma(\W{}_{sg},d)$ via fórmula de resíduos de Bott e sabemos que o grau anterior é polinomial em $d$ (Vide Proposição \ref{grau_polinomio}) de grau no máximo $120 ( = 5 \times 24)$, donde segue que é suficiente encontrarmos $\deg \Sigma(\W{}_{sg},d)$ para 121 valores diferentes de $d$ e, posteriormente, interpolar os resultados para obter o polinômio em questão. No entanto, por questões computacionais, conseguimos calcular até o momento o grau para $4 \leq d \leq 28$, cujos valores explicitamos na Tabela \ref{graus_segre}. 

\begin{table}[H]
\centering
\noindent\begin{tiny}
\begin{tabular}{|c|c|}
\hline 
d & grau\\
\hline
4& 4985292672535\\
 \hline
5& 38085453623924002125608\\
 \hline 
6& 75285508677103874434199729447346\\
 \hline
7& 6919928722801305898152558631141006297978\\
 \hline
8& 42181954432466686484802366327946036350563667373\\
 \hline
9& 30538531184782134440883223805188165885850765266730973\\
 \hline
10& 4224340951726565859342587822879909669270072209918091111509\\
 \hline
11& 158437528281133532734337703310993668084277908103801228619349318\\
 \hline
12& 2080035353059957499641534559924163791462457116358313751435919907641\\
 \hline
13& 11549735996636189943619254985547139290129087463355134074887299468381440\\
 \hline
14& 31296770227603270473657644859463859788303319257226489697655766935282861144\\
 \hline
15& 46218251138854455896028288030807107836206397262026919058025989004860345865068\\
 \hline
16& 40573178025017053248163455791995253138333248830219749681901524680514920694647875\\
 \hline
17& 22696403460389782282918120220096612693066990902486735463037695748458355012102065130\\
 \hline
18& 8560094850432050145388608162764331545974912158826771912534187363304630242378140685505\\
 \hline
19& 2280218446179281906894436399299532691147069188695294809825377606946754403932028306244123\\
 \hline
20& 445913122370782785268625533245649250274532741301118606978525517483582671680154337345798650\\
 \hline
21& 66136044830890785552763166513088475675562647217232322960605533153943919181299528743949231995\\
 \hline
22& 7648060182749239379957328222725038044389468341441118678038359708033154622298562461760431031987\\
 \hline
23& 706122123807470790783755440277242773510506336035275026531817959540511125994738199269690002855831\\
 \hline
24& 53126049393404266440928946834127714486755547747259961078332514818104185788066175129628674092418346\\
 \hline
25& 3315561352388199144671538442416320830215174679718026171794913021436371907104234446224006732128647329\\
 \hline
26& 174334857471395667347731728239322112964231210603282356701591598514084204408291171080634141000772703155\\
 \hline
27& 7829482987143513944990986949407455476367377747552625701320278751035638067417139596314040948040344857400\\
 \hline
28& 303991364820542511002698414336553281396075120749252336213971319871871164262548779281153647072907136671375\\
 \hline
\end{tabular} 
\end{tiny}
\caption{
 {(-:UAU:-)}$\deg \Sigma(\W{}_{sg},d)$}\label{graus_segre}
\end{table}

No Apêndice \ref{ap_segre} encontram-se a disposição do leitor os procedimentos e funções escritos no Macaulay2 \cite{Macaulay2} e Maple \cite{Maple_2015}   utilizados para obter o $\deg \Sigma(\W{}_{sg},d)$ em (\ref{graus_segre}).

\chapter{Superfícies singulares ao longo de quárticas elípticas\label{cap_eqc}}
Uma curva quártica elíptica é a interseção
completa de um (único) feixe de superfícies
quádricas. Avritzer \& Vainsencher 
\cite{Vainsencher_Avritzer_92}, \cite{Avritzer_Vainsencher_1999} obtiveram uma
descrição explícita da componente $\W_{eqc}$ de
quárticas elípticas do esquema de Hilbert
$\Hilb_{4t}(\mathbb{P}^3)$, a qual tem dimensão
16. Essa descrição foi utilizada em Ellingsrud \& Str{\o}mme 
\cite{Ellingsrud_Stromme_Bott_96} para enumerar
curvas na Calabi-Yau 3-folds e em
Cukierman, Lopez \& Vainsencher \cite{Cukierman_Lopez_Vainsencher_14} para
enumerar superfícies contendo uma curva quártica
elíptica.
G.\,Gotzmann \cite{Gotzmann_2008}
mostrou que ${\Hilb}_{4t}(\mathbb{P}^3)$ é formado
por duas componentes irredutíveis.

Para aplicarmos a fórmula de resíduos de Bott no cálculo do grau de $\Sigma(\W_{eqc}{},d) \subset \mathbb{P}^{N_d}$, a família de superfícies de grau $d$ que contêm algum membro $W \in \W_{eqc}$ em seu lugar singular,  utilizamos a 
descrição dos pontos fixos e respectivos tangentes explicitados em Araújo \cite{Araujo_2009}. Aqui aparecem mais pontos fixos por conta da planificação da família formada pelos subesquemas de $\P3$ definido por $\mathcal{I}^2_W$ para algum $W \in \W_{eqc}$ (vide Seção \ref{sub_fibra_eqc}). 

Como nos casos anteriores temos que  
\beq\label{grau_eqc_1} \deg \Sigma(\W_{eqc}{},d) = \int \segre(16, \mathcal{E}_d)\cap [\W'_{eqc}{}]= \sum\limits_f \dfrac{c^f_{16}(-\mathcal{E}_d)}{c^f_{16}(\mathcal{T}\W'_{eqc})}.\eeq

Na Figura \ref{eqc_cruz_malta} encontramos um exemplo de superfície de grau 6 singular ao longo da quártica elíptica $\langle x_0^2+x_1^2-x_3^2,x_1^2+x_2^2-x_3^2\rangle$ (redutível) no aberto afim $x_3=1$. Ao passo que na Figura \ref{Exemplo_Eqc} vemos um outro exemplo de uma superfície de grau 6 cujo lugar singular contém a quártica elíptica  $\langle x_0^2+x_1^2+x_2^2-100x_3^2, x_0^2+4x_1^2+90x_2x_3 \rangle$ (agora não redutível) no aberto afim $x_3=1$.

\begin{figure}[!h]
\centering\ifpdf
\includegraphics[scale=0.15]{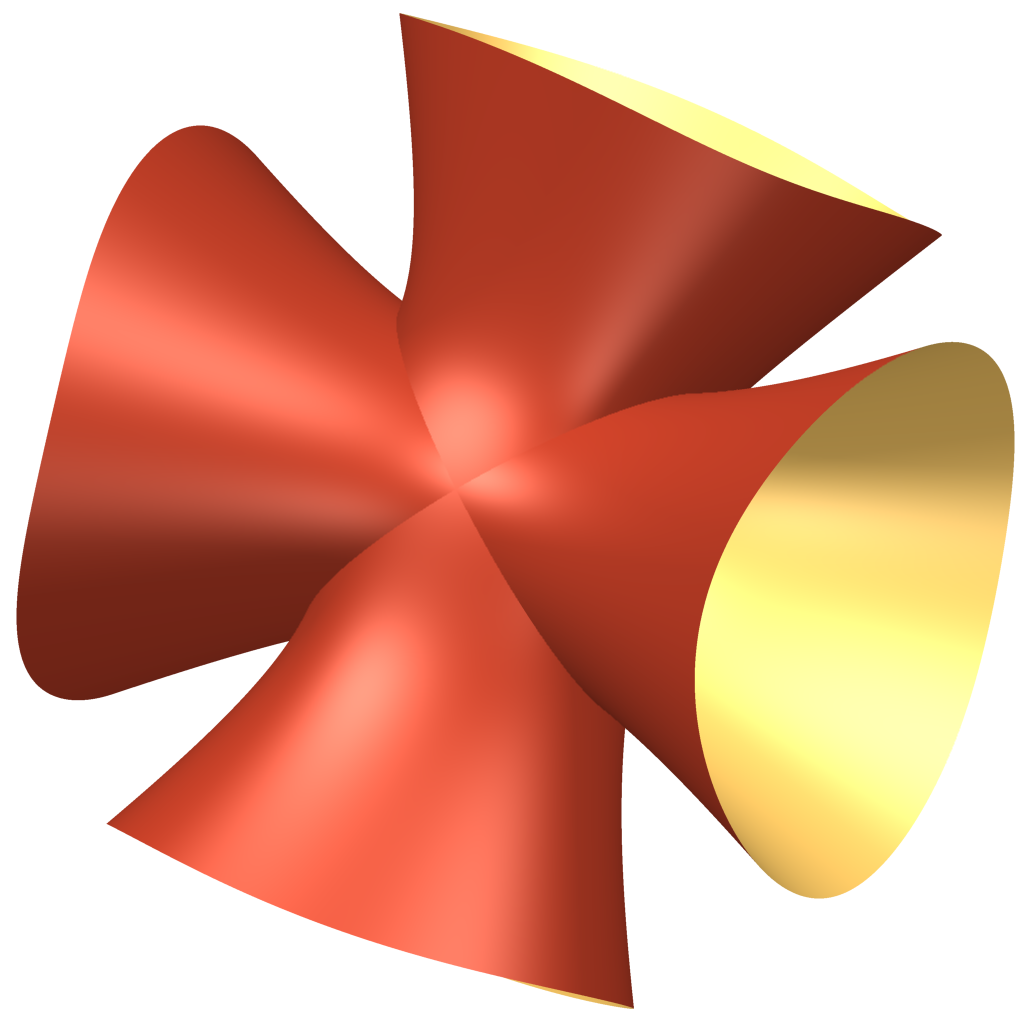}\fi
\caption{$-16x_0^6-6x_0^5x_1+42x_0^4x_1^2-9x_0^3x_1^3+63x_0^2x_1^4-x_0x_1^5+19x_1^6-18x_0^5x_2+4x_0^4x_1x_2-14x_0^3x_1^2x_2+4x_0^2x_1^3      x_2-3x_0x_1^4x_2+42x_0^4x_2^2+3x_0^3x_1x_2^2+13x_0^2x_1^2x_2^2+7x_0x_1^3x_2^2-6x_1^4x_2^2+22x_0^3x_2^3-4x_0^2x_1x_2^3+8x_0      x_1^2x_2^3-4x_1^3x_2^3-2x_0^2x_2^4+2x_0x_1x_2^4+2x_1^2x_2^4-7x_0x_2^5-5x_2^6-4x_0^5-2x_0^4x_1-11x_0^3x_1^2+3x_0^2x_1^3-2x_0x_1^4+
x_1^5-13x_0^4x_2-21x_0^2x_1^2x_2-6x_1^4x_2-3x_0^3x_2^2+7x_0^2x_1x_2^2+7x_0x_1^2x_2^2-x_1^3x_2^2+5x_0^2x_2^3+9x_1^2x_2^3+5x_0x_2^4-4x_1x_2^4+2x_2^5-16x_0^4+9x_0^3x_1-66x_0^2x_1^2+2x_0x_1^3-29x_1^4+14x_0^3x_2-4x_0^2x_1x_2+6x_0x_1^2x_2-5x_0^2x_2^2
-7x_0x_1x_2^2+8x_1^2x_2^2-8x_0x_2^3+4x_1x_2^3-8x_2^4+11x_0^3-3x_0^2x_1+4x_0x_1^2-2x_1^3+21x_0^2x_2+12x_1^2x_2-7x_0x_2^2+x_1x_2^2-9x_2^3+3x_0^2
-x_0x_1+x_1^2-3x_0x_2-2x_2^2-2x_0+x_1-6x_2+9$}
\label{eqc_cruz_malta}
\end{figure}

\begin{figure}[!h]
\centering\ifpdf
\includegraphics[scale=0.15]{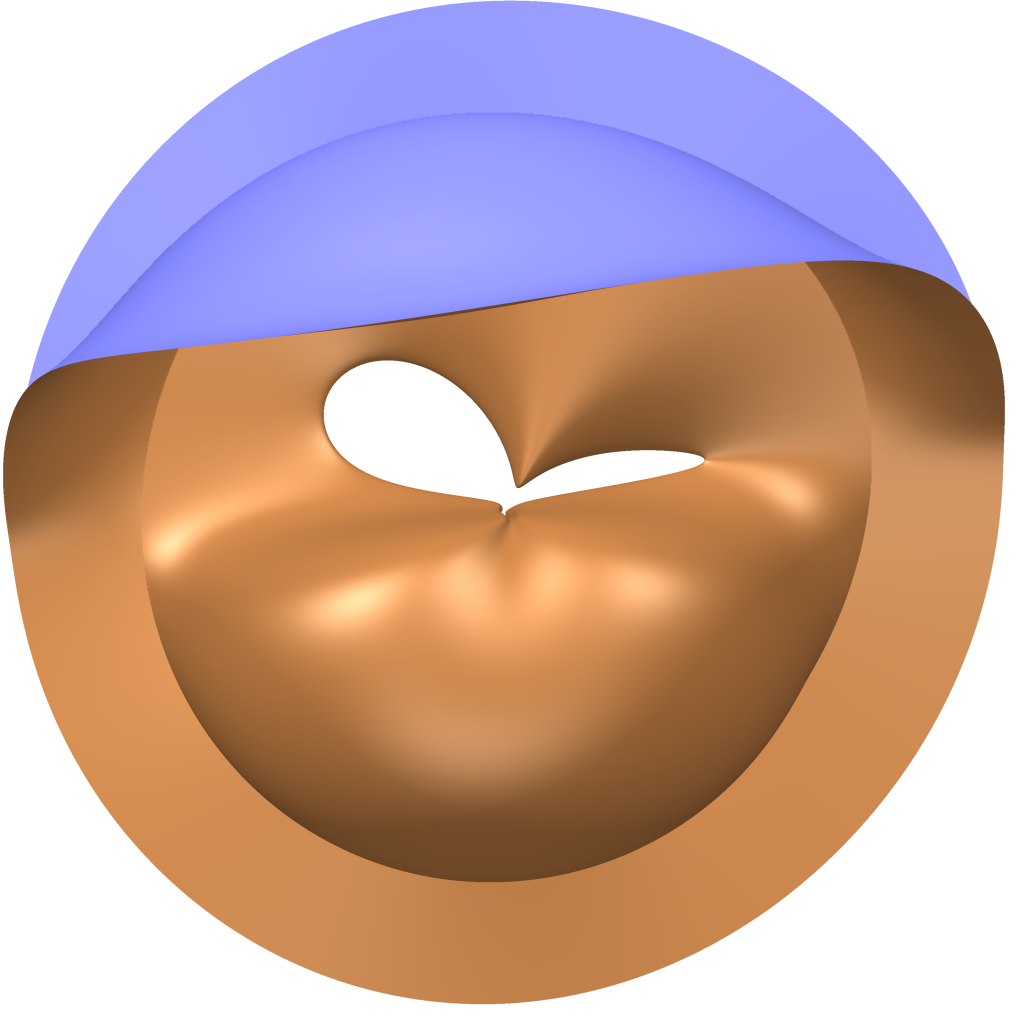}\fi
\caption{$(87017/364500000)x_0^6-(137/270000)x_0^5x_1-(57667/12150000)x_0^4x_1^
     2-(29/33750)x_0^3x_1^3-(975509/24300000)x_0^2x_1^4-(437/270000)x_0x_
     1^5-(25419073/364500000)x_1^6+(18335489/2025000)x_0^5x_2+(41726021/
     8100000)x_0^4x_1x_2+(15068581/810000)x_0^3x_1^2x_2+(17729639/1620000
     )x_0^2x_1^3x_2+(10671493/1012500)x_0x_1^4x_2+(12869111/2025000)x_1^
     5x_2+(413143379/81000000)x_0^4x_2^2-(2/1875)x_0^3x_1x_2^2+(425307677
     /40500000)x_0^2x_1^2x_2^2-(1/15000)x_0x_1^3x_2^2+(109296941/20250000
     )x_1^4x_2^2+(8074997/450000)x_0^3x_2^3+(9079999/900000)x_0^2x_1x_2^
     3+(1999997/112500)x_0x_1^2x_2^3+(2329999/225000)x_1^3x_2^3+(2840603/
     281250)x_0^2x_2^4-(7/10000)x_0x_1x_2^4+(46798403/4500000)x_1^2x_2^4
     +9x_0x_2^5+5x_1x_2^5+5x_2^6+(1341817/135000)x_0^5+(1145053/202500)x
     _0^4x_1+(1336819/27000)x_0^3x_1^2+(1157101/40500)x_0^2x_1^3+(5268037/
     135000)x_0x_1^4+(19285417/810000)x_1^5+(22737131/2025000)x_0^4x_2+(7/
     1500)x_0^3x_1x_2+(22473893/405000)x_0^2x_1^2x_2-(107/1500)x_0x_1^3
     x_2+(174916543/4050000)x_1^4x_2+(46543/1875)x_0^3x_2^2+(744013/30000
     )x_0^2x_1x_2^2+(1714349/15000)x_0x_1^2x_2^2+(2256043/30000)x_1^3x_
     2^2+(1513837/75000)x_0^2x_2^3+(3/100)x_0x_1x_2^3+(8016881/150000)x_1
     ^2x_2^3-5x_0x_2^4+8x_1x_2^4+9x_2^5+(198721/32400)x_0^4+(8/75)x_0^3
     x_1+(19985293/405000)x_0^2x_1^2+(1/150)x_0x_1^3+(16075727/162000)x_1
     ^4-(1018634/1125)x_0^3x_2-(4419467/9000)x_0^2x_1x_2-(4148009/4500)x_
     0x_1^2x_2-(2088529/4500)x_1^3x_2+8x_0^2x_2^2-x_0x_1x_2^2+7x_1^2x
     _2^2-8x_0x_2^3+6x_1x_2^3+x_2^4-(450049/450)x_0^3-(124997/225)x_0^2x
     _1-(900017/225)x_0x_1^2-(999949/450)x_1^3-5x_0^2x_2-3x_0x_1x_2+9x
     _1^2x_2+5x_0x_2^2+4x_1x_2^2+x_2^3-x_0^2-7x_0x_1-4x_1^2+x_0x_2+3x
     _1x_2+2x_2^2+6x_0+x_1-7x_2+7$}
\label{Exemplo_Eqc}
\end{figure}

\section{O espaço de parâmetros\label{espaco_parametros_eqc}}
A seguir fazemos um resumo/transcrição da construção de $\W_{eqc}$ como explicado em Avritzer \& Vainsencher  \cite{Vainsencher_Avritzer_92}, explicitada no diagrama (\ref{diagrama_eqc}) e que aparece em \cite{Cukierman_Lopez_Vainsencher_14}.
Ponha $\mathbb{X}=\mathbb{G}(2,\mathcal{F}_2)$, a grassmanniana de feixes de quádricas em $\P3$. 
\beq \label{diagrama_eqc}
\begin{gathered}
\xymatrix{ 
\mathbb{G}(19,\mathcal{F}_4) \supset \W_{eqc} &=& \widehat{\mathbb{X}}\ar[d] & & & & \widehat{\mathbb{E}}\ar[d]\ar@{_{(}->}[llll]\\
\mathbb{G}(8,\mathcal{F}_3){\times \mathbb{X}}
&\supset & \widetilde{\mathbb{X}}\ar[d] &\supset & \widetilde{\mathbb{E}}\ar[d] & \supset & \widetilde{Y}\ar[d]\\
\mathbb{G}(2,\mathcal{F}_2) &=& \mathbb{X} &\supset & Z & \supset & Y
}
\end{gathered}
\eeq

\noindent onde

\noindent$\left\{\begin{array}{l}
Z\cong \check{\mathbb{P}}^3 \times \mathbb{G}(2, \mathcal{F}_1) \mbox{ consiste em feixes com um plano fixado}; \\ 
Y\cong \{(p,l)|p\supset l\} = \mbox{ órbita fechada de } Z; \\ 
\widetilde{Y}\rightarrow Y = \P2\mbox{-fibrado de divisores de grau 2 sobre a reta variável } l \subset p; \\
\widetilde{\mathbb{X}} = \mbox{ explosão de } \mathbb{X} \mbox{ ao longo de } Z \\ 
\widehat{\mathbb{X}} =  \mbox{ explosão de } \widetilde{\mathbb{X}} \mbox{ ao longo de } \widetilde{Y}\\ 
\end{array}\right.$
\begin{center}
$Z = \left\{
\resizebox{!}{1.1cm}{
\begin{tikzpicture}
\draw (0,0) -- (1,1);
\draw (0,0) -- (2,0);
\draw (1,1) -- (3,1);
\draw (2,0) -- (3,1);
\draw (1.5, 0.5) -- (1.5, 1.5);
\end{tikzpicture}}\right\}
 \supset Y = \left\{
\resizebox{!}{1.1cm}{\begin{tikzpicture}
\draw (0,0) -- (1,1);
\draw (0,0) -- (2,0);
\draw (1,1) -- (3,1);
\draw (2,0) -- (3,1);
\draw (1, 0) -- (2, 1);
\draw (1.5,1.5);
\end{tikzpicture}}\right\}
\leftarrow \widetilde{Y}=
\left\{
\resizebox{!}{1.1cm}{\begin{tikzpicture}
\draw (0,0) -- (1,1);
\draw (0,0) -- (2,0);
\draw (1,1) -- (3,1);
\draw (2,0) -- (3,1);
\draw (1, 0) -- (2, 1);
\draw (1.5,1.5);
\draw (1.2,0.2) node[black]{$\star$} (1.8,0.8) node[black]{$\star$};
\end{tikzpicture}}\right\}$
\end{center}

Seja $$\mathcal{A}\subset \mathcal{F}_2\times
\mathbb{X}$$ o subfibrado tautológico de posto 2
sobre nossa grassmanniana de feixes de
quádricas. Existe um mapa natural de fibrados
vetoriais sobre $\mathbb{X}$ induzido por
multiplicação, $$\mu_3:
\mathcal{A}{\bigotimes} \mathcal{F}_1
\rightarrow \mathcal{F}_3 \times \mathbb{X},$$ com
posto genérico 8. Temos uma queda de posto
precisamente sobre $Z$. Daí temos induzido um mapa
racional $\kappa:\mathbb{X}\dashrightarrow
\mathbb{G}(8,\mathcal{F}_3)$. Explodindo
$\mathbb{X}$ ao longo de $Z$, encontramos o fecho
$\widetilde{\mathbb{X}}\subset
\mathbb{G}(8,\mathcal{F}_3)\times \mathbb{X}$ do
gráfico de $\kappa$. Similarmente, sobre
$\widetilde{\mathbb{X}}$ temos um subfibrado
$\mathcal{B}\subset \mathcal{F}_3\times
\widetilde{\mathbb{X}}$ de posto 8 e um mapa de
multipĺicação
$$\mu_4:\mathcal{B}{\bigotimes} \mathcal{F}_1 \rightarrow \mathcal{F}_4 \times \widetilde{\mathbb{X}}$$ com posto genérico 19. O esquema de zeros de $\stackrel{19}{\wedge} \mu_4$ é igual a $\widetilde{Y}$. De fato, pode-se verificar que cada fibra de $\mathcal{B}$ é um sistema linear de cúbicas que
\bi
\item ou tem local de base igual a uma curva com polinômio de Hilbert $p_{\W_{eqc}}(t) = 4t$
\item ou é da forma $p\cdot \mathcal{F}^{**}_2$, significando o sistema linear com componente fixada um plano $p$, e $\mathcal{F}^{**}_2$ denotando um espaço de quádricas de dimensão 8 que define um subesquema de $p$ de dimensão 0 e grau 2.  
\ei

O divisor excepcional $\widehat{\mathbb{E}}$ é um $\P8$-fibrado sobre $\widetilde{Y}$. A fibra de $\widehat{\mathbb{E}}$ sobre $(p,y_1+y_2)\in \widetilde{Y}$ é o sistema de curvas quárticas no plano $p$ que são singulares no "doublet"\, $y_1+y_2$. Precisamente, se $x_0, \dots, x_3$ denotam as coordenadas homogêneas sobre $\P3$, assumindo $p:=x_0$,  $l = \langle x_0,x_1\rangle$, um típico "doublet"\, tem ideal homogêneo da forma $\langle x_0, x_1, f(x_2,x_3)\rangle$, com $\deg f = 2$. Nosso sistema de quárticas planas está no ideal $\langle x_1, f\rangle^2 = \langle x_1^2, x_1f, f^2\rangle$. Dada uma quártica não-nula $g$ neste ideal, podemos formar o ideal $J = \langle x_0^2, x_0x_1, x_0f, g\rangle$. Verifica-se que $J$ contém precisamente 19 quárticas independentes e seu polinômio de Hilbert é correto. De fato, qualquer tal ideal é 4-regular (no sentido de Castelnuovo-Mumford).

\section{Aplicando Bott\label{aplicando_bott_eqc}}

Novamente, tomamos o toro  $\mathbb{T}=\mathbb{C}^*$ agindo diagonalmente sobre $(\mathbb{C}^4)^\vee$ via $t\circ x_i := t^{w_i}x_i$, com pesos apropriados (ver Apêndice \ref{ap_quartica_eliptica}), digamos:
\beq \label{pesos_eqc} w_0=55, w_1 = 95, w_2 = 160, w_3 = 267. \eeq 

A ação de $T$ induz ações naturais sobre $\mathbb{X}=\mathbb{G}(2,\mathcal{F}_2)$, $\widetilde{\mathbb{X}} \mbox{ e } \widehat{\mathbb{X}}$. Na sequência damos uma descrição dos pontos fixos dessa ação, bem como as contribuições necessárias à aplicação da fórmula de Bott.

\tocless\subsection{Pontos fixos em \texorpdfstring{$\mathbb{X}$}{}}

Para uma escolha adequada dos pesos $w_0, \dots,
w_3$, o mergulho de Plücker mostra que os  pontos
fixos são isolados, dados pelos feixes de monômios da
forma 
\beq 
\langle x_ix_j, x_rx_s\rangle, 0 \leq i,j,,r,s \leq 3, x_ix_j \neq x_rx_s.
\eeq
Isto  nos dá um total de $\binom{\binom{3+2}{2}}{2} = \binom{10}{2} = 45$ pontos fixos.

É fácil observar que qualquer um dos pontos fixos da ação de $\mathbb{T}$ sobre $\mathbb{X}$ é projetivamente equivalente a um dos 5 tipos:

\begin{multicols}{3}
\begin{enumerate}[label=(\arabic*)]
\item\label{pf_eqc_1} $\langle x_0^2,x_1^2 \rangle$
\item \label{pf_eqc_2}$\langle x_0^2, x_1x_2\rangle$
\item \label{pf_eqc_3}$\langle x_0x_1, x_2x_3\rangle$
\item \label{pf_eqc_4}$\langle x_0^2, x_0x_1\rangle$
\item \label{pf_eqc_5}$\langle x_0x_1, x_0x_2 \rangle$
\end{enumerate}
\end{multicols}

\bobs Existem vários pontos fixos sobre cada classe de isomorfismo acima. De fato, é fácil verificar que permutando as variáveis de um dado $\P3$ o número de pontos fixos do tipo \ref{pf_eqc_1} ao \ref{pf_eqc_5} são, respectivamente, 6, 12, 3, 12 e 12. \eobs

Os pontos fixos que correspondem a um feixe com
componente fixa $x_i$, isto é, são da forma
$\langle x_ix_j, x_i x_k\rangle$ pertencem ao
primeiro centro de explosão $Z$ (vide
\ref{diagrama_eqc}). Deste modo, o primeiro centro
de explosão $Z$ contém $24 (=
(3+1)\cdot\binom{3+1} {2})$ pontos fixos, aqueles
provenientes dos tipos \ref{pf_eqc_4} e \ref{pf_eqc_5}. Os outros 21
pontos fixos,  provenientes dos tipos \ref{pf_eqc_1},\  \ref{pf_eqc_2} \ e \ \ref{pf_eqc_3}  estão fora do centro de explosão
$Z$, bem como do segundo 
centro de explosão $\widetilde{Y}$  (vide \ref{diagrama_eqc}). Dessa forma, para estes 21 pontos fixos a contribuição é calculada sobre $\mathbb{X}$, ou seja, os divisores excepcionais não contribuem nestes pontos fixos. 

Portanto, temos contribuições imediatas para (\ref{grau_eqc_1} )
$$f = \langle x_ix_j, x_rx_s\rangle \Rightarrow \frac{c_{16}^{f}(-\mathcal{E}_d)}{c_{16}^{f}(\mathcal{T}\mathbb{X})},$$
onde o denominador $c_{16}^{f}(\mathcal{T}\mathbb{X})$ é o produto dos pesos da representação
\beq \label{rep_tang_eqc} \ba {lcl}\mathcal{T}_f\mathbb{X} &=& (\mathcal{F}_2/\mathcal{A}_{|_f})\otimes \mathcal{A}^{\vee}_{|_f} \\ &=& \left(\left(\sum\limits_{0\leq i\leq j \leq 2} t^{w_i+w_j}\right)/\langle x_ix_j, x_rx_s\rangle \right) \otimes \langle x_ix_j, x_rx_s\rangle^\vee.\ea\eeq

Por exemplo, tomando o ponto fixo do tipo (1), $f = \langle x_0^2, x_1^2\rangle$, temos que a representação em (\ref{rep_tang_eqc}) exprime-se da seguinte forma:

$\ba {lcl} \mathcal{T}_f\mathbb{X} &=& \langle x_0x_1, x_0x_2, x_0x_3, x_1x_2, x_1x_3, x_2^2, x_2x_3, x_3^2\rangle\otimes \langle x_0^2, x_1^2\rangle^\vee\\ &=& \frac{x_0x_1}{x_0^2}+\frac{x_0x_2}{x_0^2}+\frac{x_0x_3}{x_0^2} +\frac{x_1x_2}{x_0^2}+\frac{x_1x_3}{x_0^2}+\frac{x_2^2}{x_0^2}+\frac{x_2x_3}{x_0^2}+\frac{x_3^2}{x_0^2}  + \frac{x_0x_1}{x_1^2}\\ &+&\frac{x_0x_2}{x_1^2}+\frac{x_0x_3}{x_1^2} +\frac{x_1x_2}{x_1^2}+\frac{x_1x_3}{x_1^2}+\frac{x_2^2}{x_1^2}+\frac{x_2x_3}{x_1^2}+\frac{x_3^2}{x_1^2}\ea$, 

\noindent onde $\frac{x_rx_s}{x_i^2}$ indica o $\mathbb{T}$-espaço com peso $w_r+w_s - 2w_i$. 
\bobs A contribuição  do numerador em (\ref{grau_eqc_1}) para cada ponto fixo será descrita na seção \ref{sub_fibra_eqc}.  \eobs 

\tocless\subsection{Pontos fixos em \texorpdfstring{$\widetilde{\mathbb{X}}$}{}}

Consideremos a fibra de $\widetilde{\mathbb{E}}$ sobre $f = \langle x_ix_j, x_ix_k\rangle$, ponto fixo em $Z$. De acordo com Avritzer \& Vainsencher \cite{Avritzer_Vainsencher_1999} esta fibra é dada por 
\beq \label{fibra_divisor_1}\mathbb{P}^{3^2-1} \cong \mathbb{P}\left(\mathcal{F}_3^{\langle x_j,x_k\rangle}/x_i\cdot \mathcal{F}_2^{\langle x_j,x_k\rangle}\right)\eeq
isto é, o espaço projetivo de formas cúbicas que anulam-se ao longo da reta $l = \langle x_j, x_k\rangle$ módulo $x_i$ vezes as formas quadráticas que anulam-se sobre a mesma reta $l$. Daí segue que os pontos fixos do tipo (4) dão origem aos pontos fixos no divisor excepcional $\widetilde{\mathbb{E}}$ da forma:
\beq\label{pontos_fixos_4_1 - 4_9} \langle x_i\cdot\langle x_i, x_j\rangle, [c]\rangle, \eeq
onde $[c]$ indica a classe de uma forma cúbica descrita em (\ref{fibra_divisor_1}), ou seja, $c$ é uma das formas cúbicas $x_ix_k^2, x_ix_kx_l, x_ix_l^2, x_j^3, x_j^2x_k, x_j^2x_l, x_jx_k^2, x_jx_kx_l, x_jx_l^2$.

De modo análogo, os pontos fixos do tipo (5) dão origem aos pontos fixos no divisor excepcional $\widetilde{\mathbb{E}}$ da forma:
\beq \label{pontos_fixos_5_1} \langle x_i\cdot\langle x_j, x_k\rangle, [c]\rangle, \eeq
onde $[c]$ indica a classe de uma forma cúbica descrita em (\ref{fibra_divisor_1}), ou seja, $c$ é uma das formas cúbicas $x_j^2x_k, x_jx_k^2, x_jx_kx_l, x_k^3, x_k^2x_l, x_kx_l^2, x_j^3, x_j^2x_l, x_jx_l^2$.

Dessa forma, cada ponto fixo dos tipos (4) e (5) produzem 9 pontos fixos sobre o divisor excepcional $\widetilde{\mathbb{E}}$. No entanto, temos que desconsiderar os pontos fixos que estão no segundo centro de explosão $\widetilde{Y}$ (vide (\ref{diagrama_eqc})). Da descrição de $\widetilde{Y}$ temos que os únicos pontos fixos sobre este último são aqueles do tipo 
\beq \label{pontos_fixos_4_1 - 4_3}\langle x_i\cdot\langle x_i, x_j, q\rangle\rangle, \eeq
onde $q$ é uma das formas quadráticas: $x_k^2, x_kx_l$ ou $x_l^2$, ou seja, $q$ é uma forma quadrática que não se anula na reta $\langle x_i, x_j\rangle$. Portanto, os pontos fixos dos tipos (4) e (5) dão origem, a menos de permutação das variáveis, a 9+9-3 = 15 pontos fixos fora do segundo centro de explosão $\widetilde{Y}$. Para estes, a contribuição para (\ref{grau_eqc_1}) é calculada sobre $\widetilde{\mathbb{X}}$ e o primeiro divisor excepcional $\widetilde{\mathbb{E}}$ contribui nestes pontos fixos.

Para cada ponto fixo $f' = \langle x_i\cdot \langle x_j,x_k\rangle, [c]\rangle$, onde $[c]$ indica a classe da forma cúbica, temos a seguinte contribuição: 
\beq\label{contribuicao_1_blowup} f' = \langle x_i\cdot \langle x_j,x_k\rangle, [c]\rangle\Rightarrow\frac{c_{16}^{f'}(-\mathcal{E}_d)}{c_{16}^{f'}(\mathcal{T}\widetilde{\mathbb{X}})}.\eeq

Assim, para obtermos as contribuições do divisor
excepcional $\widetilde{\mathbb{E}}$  e do fibrado
tangente $\mathcal{T}_{f'}\widetilde{\mathbb{X}}$,
precisamos lembrar que o divisor excepcional
$\widetilde{\mathbb{E}}\rightarrow Z$ é o\lb
fibrado projetivo $\mathbb{P}(\mathcal{N}_{Z/\mathbb{X}}).$ Agora, sobre $f = \langle x_ix_j, x_ix_k\rangle$, obtemos a sequência de representações $$0 \rightarrow \mathcal{T}_f Z \rightarrow \mathcal{T}_f \mathbb{X} \rightarrow \mathcal{N}_{_f Z/\mathbb{X}} \rightarrow 0$$
e a partir das fórmulas para o espaço tangente da grassmanniana $\mathcal{T}_f \mathbb{X}$ e $\mathcal{T}_f Z$, escritos como decomposição de auto-espaços, obtemos:
\begin{footnotesize}
\beq \ba {lcl} 
    \mathcal{N}_{_f Z/\mathbb{X}} & = & \mathcal{T}_f \mathbb{X}-\mathcal{T}_f Z \\
    &=& (x_ix_j+x_ix_k)^\vee \otimes (\mathcal{F}_2 - (x_ix_j +x_ix_k)) \\
    &-& \left(\langle x_i\rangle^\vee \otimes (\mathcal{F}_1 - x_i) + (x_j+x_k)^\vee \otimes (\mathcal{F}_1 -(x_j+x_k))\right)\\
    &=& (x_ix_j+x_ix_k)^\vee \otimes (\mathcal{F}_2 - (x_ix_j +x_ix_k)) \\
    &-&(\langle x_i\rangle^\vee \otimes (\mathcal{F}_1 - x_i) + (x_ix_j+x_ix_k)^\vee \otimes (x_i\mathcal{F}_1 -(x_ix_j+x_ix_k)))\\
    &=& (x_ix_j+x_ix_k)^\vee \otimes (\mathcal{F}_2 - x_i \mathcal{F}_1) - \langle x_i\rangle^\vee \otimes (\mathcal{F}_1 - x_i)\\
    &=& \langle x_ix_jx_k\rangle^\vee \otimes (x_k \mathcal{F}_2 + x_j \mathcal{F}_2 - x_i(x_j+x_k)\mathcal{F}_1 - x_jx_k \mathcal{F}_1 + x_ix_jx_k)\\
    &=& \langle x_ix_jx_k\rangle^\vee \otimes \left(\mathcal{F}_3^{\langle x_j,x_k\rangle}/ x_i\cdot \mathcal{F}_2^{\langle xj,x_k\rangle}\right)
  \ea 
\eeq
\end{footnotesize}

Por comparação direta da descrição geométrica e dos pesos obtidos pela decomposição em auto-espaços do fibrado normal, obtemos

$$ c \leftrightarrow \langle c\rangle = \frac{c}{x_ix_j x_k},$$
onde $\langle c\rangle$ indica o caracter da decomposição em auto-espaços do fibrado normal correspondente à cúbica $c$. E para a contribuição do denominador em (\ref{contribuicao_1_blowup}) obtemos o produto dos pesos da representação
\begin{footnotesize}
{\setlength{\mathindent}{0cm}
\beq \label{rep_tang_eqc_2}
\begin{gathered}
\ba {lcl} 
\mathcal{T}_{f'} \widetilde{\mathbb{X}} &=& \mathcal{T}_f Z + \mathcal{T}_{[c]} \mathbb{P}(\mathcal{N}_{Z/X})+\langle c\rangle \\
&=& \langle x_i \rangle^\vee \otimes (\mathcal{F}_1 - x_i) + (x_j+x_k)^\vee \otimes (\mathcal{F}_1 - (x_j+x_k))\\
&+& \left((x_ix_j+x_ix_k)^\vee \otimes (\mathcal{F}_2-x_i\mathcal{F}_1) - \langle x_i \rangle^\vee \otimes (\mathcal{F}_1 -x_i)-\langle c\rangle\right)\otimes \langle c\rangle^\vee\\
 &+& \langle c \rangle
\ea
\end{gathered}
\eeq}
\end{footnotesize}

\tocless\subsection{Pontos fixos em \texorpdfstring{$\widehat{\mathbb{X}}$}{}}

De modo análogo, consideremos a fibra de $\widehat{\mathbb{E}}$, divisor excepcional da explosão de $\widetilde{\mathbb{X}}$ ao longo de $\widetilde{Y}$, sobre o ponto fixo $x_i\cdot\langle x_i,x_j, q\rangle$ (vide (\ref{pontos_fixos_4_1 - 4_3})) em $\widetilde{Y}$, a qual, cf. \cite{Avritzer_Vainsencher_1999} , é dada por 
\beq\label{fibra_divisor_2} \P8 \cong \mathbb{P}\left(\mathcal{F}_4^{\langle x_i, x_j, q\rangle^2}/x_i\cdot \mathcal{F}_3^{\langle x_i, x_j, q\rangle}\right),\eeq
isto é, o espaço projetivo das formas de grau 4 contidas no ideal $\langle x_i,x_j, q\rangle^2$ módulo $x_i$ vezes as formas cúbicas contidas no ideal $\langle x_i, x_j, q\rangle$.

Agora, cada ponto fixo sobre o segundo centro de explosão $\widetilde{Y}$ (vide (\ref{pontos_fixos_4_1 - 4_3})), produz 9 pontos fixos sobre o divisor excepcional $\widehat{\mathbb{E}}$, a saber aqueles da forma 
\beq \label{pontos_4.1.1-4.3.9} \langle x_i\cdot\langle x_i, x_j, q\rangle, [q_4]\rangle, \eeq
onde $[q_4]$ indica a classe de uma forma de grau 4 descrita por (\ref{fibra_divisor_2}). Portanto, a menos de permutação das variáveis, temos um total de $27 (= 3 \times 9)$ pontos fixos sobre o divisor excepcional $\widehat{\mathbb{E}}$. Sobre estes pontos a contribuição dos pesos é calculada sobre $\widehat{\mathbb{X}}$ e os divisores excepcionais também contribuem para este cálculo. Assim, a partir de (\ref{grau_eqc_1}), temos que 
\beq\label{contribuicao_blowup_2} f^{''} = (x_i\cdot\langle x_i, x_j, q\rangle, [q_4]) \Rightarrow \frac{c_{16}^{f^{''}}(-\mathcal{E}_d)}{c_{16}^{f^{''}}(\mathcal{T}\widehat{\mathbb{X}})}.\eeq

Sabemos que o divisor excepcional $\widehat{\mathbb{E}} \rightarrow \widetilde{Y}$ é dado como a projetivização do fibrado normal, $\mathbb{P}\left(\mathcal{N}_{\widetilde{Y}/\widetilde{\mathbb{X}}}\right).$ Agora, sobre $f^{'} = ( x_i\cdot\langle x_i, x_j, q\rangle)$, obtemos a sequência de representações $$0 \rightarrow \mathcal{T}_{f'} \widetilde{Y} \rightarrow \mathcal{T}_{f'} \widetilde{\mathbb{X}} \rightarrow \mathcal{N}_{_{f'} \widetilde{Y}/\widetilde{\mathbb{X}}} \rightarrow 0$$
e a partir das fórmulas para os espaços tangentes $\mathcal{T}_{f'} \widetilde{Y}$ e $\mathcal{T}_{f'} \widetilde{\mathbb{X}}$, escritos como decomposição de auto-espaços, obtemos:
\begin{footnotesize}
\beq \ba {lcl} 
\mathcal{N}_{_{f'} \widetilde{Y}/\widetilde{\mathbb{X}}} & = & \mathcal{T}_{f'} \widetilde{\mathbb{X}} - \mathcal{T}_{f'} \widetilde{Y}= \left(\mathcal{T}_f Z + \mathcal{T}_{[c]} \mathbb{P}(\mathcal{N}_{Z/X})+\langle c\rangle\right) -\mathcal{T}_{f'} \widetilde{Y}\\
&=& \langle x_i\rangle^\vee \otimes (\mathcal{F}_1 - x_i) + (x_i+x_j)^\vee \otimes (\mathcal{F}_1-(x_i+x_j)) \\
&+& \langle x_iq\rangle^\vee \otimes (x_i\mathcal{F}_2 + x_j\mathcal{F}_2 - x_ix_j\mathcal{F}_1 - x_i(x_i\mathcal{F}_1 +x_j\mathcal{F}_1 - x_ix_j) - x_iq) \\
&+& \langle\frac{x_iq}{x_i^2x_j}\rangle - \langle x_i\rangle^\vee \otimes (\mathcal{F}_1-x_i) + \langle x_j\rangle^\vee \otimes(\mathcal{F}_1-(x_i+x_j))\\
&-& \langle q\rangle^\vee \otimes (\mathcal{F}_2 - x_i\mathcal{F}_1- x_j\mathcal{F}_1+x_ix_j\mathcal{F}_1-q)\\
&=& \langle x_ix_jq\rangle^\vee \otimes \left(\mathcal{F}_4^{\langle x_i,x_j,q\rangle^2}/x_i\cdot \mathcal{F}_3^{\langle x_i,x_j,q\rangle^2}\right).
\ea \eeq
\end{footnotesize}

Por comparação direta da descrição geométrica e dos pesos obtidos pela decomposição em auto-espaços do fibrado normal, obtemos

$$ q_4 \leftrightarrow \langle q_4\rangle = \frac{q_4}{x_ix_j q},$$
e para a contribuição do denominador em (\ref{contribuicao_blowup_2}) obtemos o produto dos pesos da representação
\beq \label{rep_tang_eqc_3}
\mathcal{T}_{f''} \widehat{\mathbb{X}}= \mathcal{T}_{f'} \widetilde{Y} + \mathcal{T}_{[q_4]}\mathbb{P}\left(\mathcal{N}_{\widetilde{Y}/\widetilde{\mathbb{X}}}\right) + \langle q_4\rangle.
\eeq

\tocless\subsection{Fibras de \texorpdfstring{$\mathcal{E}_d$}{Ed}\label{sub_fibra_eqc}}
Para determinar a contribuição  do numerador em (\ref{grau_eqc_1}), precisamos conhecer a representação da fibra de $\mathcal{E}_d$ sobre cada um dos pontos fixos. Como no caso de redes de quádricas do tipo determinantal (vide Capítulo \ref{cap_QD}) , a família formada pelos subesquemas de $\P3$ definido por $\mathcal{I}^2_W$ para algum $W \in \W_{eqc}$ não é plana. De fato, os pontos fixos do tipo 
\beq \label{pontos_4_4-4_6}\langle x_i\cdot \langle x_i, x_j\rangle, [c]\rangle, \eeq 
onde $[c]$ indica a classe de uma forma cúbica do tipo $x_j^2\cdot(\mathcal{F}_1\big/{\langle x_i\rangle})$, são membros legítimos de $\W_{eqc}$ (polinômio de Hilbert $p_{\W_{eqc}}(t) = 4t$), mas seu quadrado tem polinômio de Hilbert "ruim", a saber $13t-20$, enquanto o esperado é $12t-16$. 

{As contas locais (vide Apêndice \ref{codigos_singular_eqc}) mostram que a variedade $\mathbb{M}$ descrita no Lema \ref{lema_M} é o local onde $\mathbb{W}'$ deixa de ser plana. Sejam $\check{\mathbb{P}}^3$ o espaço projetivo
  dual com sequência tautológica 
  \beq \mathcal{L} =
  \mathcal{O}_{\check{\mathbb{P}}^3}(-1)
  \rightarrowtail \mathcal{F}_1 \twoheadrightarrow
 {\mathcal{F}_1}
\big/
    {\mathcal{L}},
   \eeq
$\mathbb{G}(2,\mathcal{F}_1)$ a grassmanniana de
  retas em $\P3$ com sequência
  tautológica 
  \beq\mathcal{S}\rightarrowtail
  \mathcal{F}_1 \twoheadrightarrow \mathcal{Q},
  \eeq
  com $rk(\mathcal{S}) = 2$ e por fim, escrevemos
  $\mathbb{X}=\mathbb{G}(2,S_2(\mathcal{F}_1))$
  para a variedade grassmanniana de feixes de
  quádricas em $\P3$
  e
  \beq \label{subfibrado_A} \mathcal{A}\rightarrowtail
  S_2(\mathcal{F}_1) \twoheadrightarrow
  {S_2(\mathcal{F}_1)}\big/{\mathcal{A}}
  \eeq
a sequência tautológica de $\mathbb{X}$, onde o $rk(\mathcal{A}) = 2$.
}
{
\bobs Em referência ao diagrama \ref{diagrama_eqc}
temos que um ponto $\widetilde{e}\in
\widetilde{\mathbb{E}}$ pode ser escrito como uma
trinca $\langle H,L,C\rangle$, onde $\langle
H,L\rangle\in Z= \check{\mathbb{P}}^3\times
\mathbb{G}(2,\mathcal{F}_1)$ e $C\in
\mathbb{P}\left(
       {S_3\mathcal{F}_1^L}\big/
       {H\cdot S_2\mathcal{F}_1^L}\right)$ (cf.  Avritzer $\&$ Vainsencher \cite{Avritzer_Vainsencher_1999}, Corolário 1.6, pág. 2999).
\eobs}

{
\blem \label{lema_M}
Seja $\mathbb{M}=
\mathbb{P}\left(({\mathcal{S}}\big/{\mathcal{L}})^2\otimes
{\mathcal{F}_1}\big/{\mathcal{L}}\right)
\times_{\mathbb{X}} Y$ o $\P2$-fibrado sobre $Y$
que parametriza as trincas $\langle H,L,C\rangle$,
onde $H$ denota um plano, $L = \langle H,
H'\rangle$ uma reta no plano $H$ e $C$ uma cúbica
da forma $(H')^2\cdot (H'')$ com $H'' \in
\mathcal{F}_1$ um representante de
$\overline{H''}\in
{\mathcal{F}_1}\big/{\mathcal{L}_{|_H}}$. Então,
temos um mergulho $$\mathbb{M}=
\mathbb{P}\left(({\mathcal{S}}\big/{\mathcal{L}})^2\otimes
{\mathcal{F}_1}\big/{\mathcal{L}}\right)
\hookrightarrow
\mathbb{P}\left({\mathcal{R}_3}\big/{\mathcal{L}
 { \otimes} \mathcal{R}_2}\right)\cong \widetilde{\mathbb{E}},$$
onde $\mathcal{R}_i = ker(S_i(\mathcal{F}_1) \longrightarrow S_i(\mathcal{Q}))$.
\elem
\bdem 
Temos um mapa natural injetivo localmente split de
fibrados sobre $Y$,
$({\mathcal{S}}\big/{\mathcal{L}})^2\otimes
{\mathcal{F}_1}\big/{\mathcal{L}} \longrightarrow
{\mathcal{R}_3}\big/{\mathcal{L}\otimes
  \mathcal{R}_2}$ definido a seguir. Na fibra
sobre um par $H\supset L = \langle H,H'\rangle$
seja $\overline{H''} \in
{\mathcal{F}_1}\big/{\mathcal{L}}|_H$ representado
por $H''\in \mathcal{F} _1$. Então, $(H')^2\cdot
H''$ é uma forma cúbica anulando-se sobre $L$ e
mapeamos $(H')^2\cdot H''$ na classe $\overline{(H')^2\cdot H''}$.
\edem}

{Note que $\mathbb{M}$ é disjunta do segundo centro de explosão $\widetilde{Y}$ (vide diagrama \ref{diagrama_eqc}) e $\dim (\mathbb{M}) = 7 (= 5+ 2)$. Além disso, em $\mathbb{M}$ temos uma única órbita fechada sob a ação de $GL_4$, representada por $\widetilde{o}_1 = \langle x_0^2, x_0x_1, x_1^3\rangle$ (cf. \cite{Avritzer_Vainsencher_1999}, Lema 1.8, pág. 3000).}

Sejam $\mathcal{A}$ o subfibrado de posto 2 descrito em \ref{subfibrado_A}, $\mathcal{B}$ o subfibrado de posto 8 da
variedade grassmanniana
$\mathbb{G}(8,\mathcal{F}_3)$ (vide
\ref{diagrama_eqc}) e $\mathcal{D}$ o subfibrado
tautológico de posto 12 da variedade grassmanniana
$\mathbb{G}(12,\mathcal{F}_5)$. Temos um mapa natural de fibrados vetoriais sobre $\widetilde{\bb X}$ induzido por multiplicação com posto genérico 12 (vide Proposição \ref{prop_nu_w_eqc}),  onde a queda de posto ocorre precisamente sobre $\bb M$.

\bprop \label{prop_nu_w_eqc} Seja $\nu: \cl A\bigotimes\mathcal{B}\rightarrow \mathcal{D}$ o mapa definido por multiplicação. O posto de $\nu$ fora de $\mathbb{M}$ é 12 e $\mathbb{M}$ é o esquema de zeros de $\stackrel{12}\wedge \nu$.
\eprop
\bdem
Denote por $\mathfrak{Z}$ o esquema de zeros do enunciado. Temos que $\mathbb{M}\subseteq \mathfrak{Z}$. Como  na prova da Proposição \ref{prop_mu_twc}, resta mostrarmos que o espaço tangente de $\mathbb{M}$ e $\mathfrak{Z}$ coincidem na órbita fechada $\widetilde{o}_1 = \langle x_0^2, x_0x_1, x_1^3\rangle$. Para isso, precisamos computar o ideal de $\mathfrak{Z}$ módulo quadrados em uma vizinhança dessa órbita e mostrar que ele contém o número correto de elementos com termos lineares independentes, a saber, a $\mbox{codim} \mathbb{M}$ em $\mathbb{W}_{eqc}$: 16 - 7 = 9.

Para o cálculo de $\mathcal{T}_{\widetilde{o}_1}\mathfrak{Z}$ precisamos introduzir coordenadas. Tome $U \cong \mathbb{A}^{16}$ a vizinhança padrão afim de $\langle x_0^2, x_0x_1\rangle$ em $\mathbb{X}$ com funções coordenadas $a_{i,j}$ $i=1,2$, $j=1, \cdots, 8$. A restrição do fibrado tautológico $\mathcal{A}$ de $\mathbb{G}(2,\mathcal{F}_2)$ à $U$ é trivial com base

$V_1 = x_0^2+a_{1,1}x_0x_2+a_{1,2}x_0x_3+a_{1,3}x_1^2+a_{1,4}x_1x_2+a_{1,5}x_1x_3+a_{1,6}x_2^2+
a_{1,7}x_2x_3+a_{1,8}x_3^2$

$V_2 = x_0x_1+a_{2,1}x_0x_2+a_{2,2}x_0x_3+a_{2,3}x_1^2+a_{2,4}x_1x_2+a_{2,5}x_1x_3+a_{2,6}x_2^2+
a_{2,7}x_2x_3+a_{2,8}x_3^2$.

Tomamos $b_i, 0\leq i \leq 8$, obtidos como na Proposição 1.2, pág. 2996 (c.f Avritzer \& Vainsencher \cite{Avritzer_Vainsencher_1999}) ou conforme o Apêndice \ref{codigos_singular_eqc}, a saber:

\begin{eqnarray}\label{eq_eqc_cubicas}
& b_1=-a_{1,3}(x_{1}^3), b_2=-a_{1,4}(x_{1}^2x_{2}),b_3=-a_{1,5}(x_{1}^2x_{3})\nonumber \\
& b_4=-a_{1,6}(x_{1}x_{2}^2), b_5=-a_{1,7}(x_{1}x_{2}x_{3}),
b_6=-a_{1,8}(x_{1}x_{3}^2) \\
& b_7=a_{2,6}(x_{0}x_{2}^2),
b_8=a_{2,7}(x_{0}x_{2}x_{3}),b_9=a_{2,8}(x_{0}x_{3}^2)\nonumber
\end{eqnarray}
onde os monômios entre parênteses formam uma base para o espaço $\mathbb{P}(\mathcal{F}_3^L/H\mathcal{F}_2^L)$ que dá a descrição do divisor excepcional $\widetilde{\mathbb{E}}$ (vide \ref{diagrama_fibrado_eqc}), em que $L = \langle x_0,x_1\rangle$ e $H = \langle x_0\rangle$.

Considere o ideal $\mathcal{I}=\langle a_{1,1}, a_{1,2}, a_{1,4}, \cdots, a_{1,8}, a_{2,1}, \cdots, a_{2,8}\rangle^2+\langle a_{1,3} \rangle^2$. Sejam $c_1, c_2, \cdots,$ $c_9$ coordenadas homogêneas para $\mathbb{P}^8$. Então, $U\times_{\mathbb{X}} \widetilde{\mathbb{X}}$  é a subvariedade de $U\times \P8$ definida pelas equações $b_ic_j=b_jc_i$, $1\leq i, j \leq 9$. Seja $U'\subset U\times_{\mathbb{X}} \widetilde{\mathbb{X}}$ o subconjunto aberto afim definido fazendo $c_1=1$. Esta escolha é realizada com o intuito de por o ponto $\widetilde{o}_1$ na origem das cartas coordenadas. Daí segue que coordenadas locais para $U'$ em torno de $\widetilde{o}_1$ são dadas por
$$b_1, c_2, \cdots, c_8, a_{1,1}, a_{1,2}, a_{2,1}, \cdots, a_{2,5}.$$ 

Nessas condições temos que $b_1$ é a equação do divisor excepcional e $b_i=b_1c_i$, $1\leq i \leq 9$, módulo o ideal $\mathcal{I}$. Resolvendo essas relações para os termos lineares (\ref{eq_eqc_cubicas}) e substituindo na matriz de representação de $\mu_3|_U$, onde $\mu_3:\mathcal{A}\bigotimes \mathcal{F}_1 \rightarrow \mathcal{F}_3\times \mathbb{X}$ é o mapa de fibrados induzidos por multiplicação (vide \ref{diagrama_eqc}), obtemos uma representação local para $\mu_3|_{\widetilde{\mathbb{X}}}$. Dividindo a última linha por $b_1$ as 8 linhas da nova matriz produzem seções linearmente independentes $w_i (1 \leq i \leq 8)$ de $\mathcal{F}_3|_{U'}$ que geram o feixe ideal $\mathcal{B}|_{U'}$. As 7 primeiras seções geram o mesmo submódulo que $V_1x_0, V_1x_1, V_1x_2, V_1x_3, V_2x_1, V_2x_2, V_2x_3$. A última delas, escrevendo somente os termos lineares, é dada por:

$w_8 = -x_1^3 - c_1x_0x_2^2-c_2x_0x_2x_3-c_3x_0x_3^2+(-c_4-a_{2,1})x_1^2x_2+(-c_5-a_{2,2})x_1^2x_3-c_6x_1x_2^2-c_7x_1x_2x_3 - c_8x_1x_3^2$.

Escolhendo uma base adequada para $\mathcal{F}_5$ e efetuando operações elementares nas linhas módulo o quadrado do ideal gerado pelas funções coordenadas de $U'$, obtemos uma representação de $\nu$ na forma $$\left(\begin{array}{cc}
J' & A' \\ 
0 & R \\
* & *
\end{array} \right)$$
onde $J'$ é uma matriz bloco triangular superior de tamanho $11\times 11$ e cujas entradas na diagonal principal são todas iguais a 1. Aqui $R$ é uma matriz linha cujas entradas módulo quadrado são as seguintes:

\beq \label{geradores_divisor_exp_eqc}\begin{array}{ccc}
-c_{1}(x_{0}^3x_{2}^2)& -c_{2}(x_{0}^3x_{2}x_{3}) &
-c_{3}(x_{0}^3x_{3}^2)\\
-c_{6}(x_{0}^2x_{1}x_{2}^2) &-c_{7}(x_{0}^2x_{1}x_{2}x_{3})&
-c_{8}(x_{0}^2x_{1}x_{3}^2)\\
(-a_{1,1}+2a_{2,4})(x_{0}x_{1}^3x_{2})&
(-a_{1,2}+2a_{2,5})(x_{0}x_{1}^3x_{3})&
-a_{1,3}(x_{1}^5)
\end{array} \eeq
onde entre parênteses indicamos os monômios correspondentes a cada coluna módulo quadrados.

Portanto, vemos que o número de elementos com termos lineares independentes no ideal de $\mathfrak{Z}$ em torno de $\tilde{o}_1$ é de fato $9=\mbox{codim}_{\mathbb{W}_{eqc}}\mathbb{M}$, como era para ser mostrado.
\edem

Na próxima proposição entra em cena o fibrado de formas quínticas. Para $\pi \in \widetilde{\mathbb{X}}\setminus \mathbb{M}$, o espaço de formas quínticas obtido como imagem do mapa $\mu$ é de posto correto 12. Explodir $\widetilde{\mathbb{X}}$ ao longo de $\mathbb{M}$ permite-nos estender a família $(\mathcal{Q}_\pi)_{\pi \in \widetilde{\mathbb{X}}\setminus \mathbb{M}}$ sobre uma compactificação  suave $ \widetilde{\mathbb{X}}' $ de $\widetilde{\mathbb{X}}\setminus \mathbb{M}$.

\bprop 
Considere o diagrama de explosão de $\widetilde{\mathbb{X}}$ ao longo de $\mathbb{M}$
 
\beq \label{diagrama_fibrado_eqc}
\begin{gathered}
\xymatrix{
\widetilde{\mathbb{\bb X}}' \ar[d]
&\supset& \widetilde{\mathbb{M}} \ar[d]
\\
 \widetilde{\bb X} &\supset & \mathbb{M}
}
\end{gathered}
\eeq
Então, temos que $\widetilde{\mathbb{X}}'$ mergulha em $\widetilde{\mathbb{X}}\times \mathbb{G}(12,\mathcal{F}_5)$ de tal forma que o pullback $\mathcal{Q}$ do subfibrado tautológico de posto 12 de $\mathcal{F}_5$ contém a imagem de $\nu|_{\widetilde{\mathbb{X}}'}$.
\eprop
\bdem
Em vista da proposição anterior, temos que o mapa racional $\widetilde{\mathbb{X}}\dashrightarrow \mathbb{G}(12,\mathcal{F}_5)$ definido por $\nu$ estende a um morfismo $\widetilde{\mathbb{X}}'\rightarrow \mathbb{G}(12,\mathcal{F}_5)$. Assim, seu gráfico produz o mergulho desejado.
\edem

Passamos agora a análise das contribuições dos pontos fixos sobre o divisor excepcional $\widetilde{\bb M}=\bb P(\cl N_{{\bb M}/\widetilde{\bb X}})$. Temos que sobre cada ponto fixo $f\in {\bb M}$ a fibra $\widetilde{\bb M}_f$ é o espaço projetivo $\bb P(\cl N_{_f{\bb M}/\widetilde{\bb X}}) = \cl T_f\widetilde{\bb X}/\cl T_f\bb M$. Por exemplo, sobre $f=\langle x_0^2, x_0x_1, x_1^3\rangle$ temos:

\beq\label{normal_eqc}
\begin{array}{lcl}
\cl N_{_f\bb M/\widetilde{\bb X}} &=& \cl T_f \widetilde{\bb X} - \cl T_f \bb M
 = \left( \frac{x_1}{x_0}+ \frac{x_1^2}{x_0^2}+
\frac{2x_2}{x_1}+ \frac{2x_2}{x_0}+\frac{x_0x_2^2}{x_1^3}+ \frac{x_2^2}{x_1^2}+ \frac{2x_3}{x_1}+\frac{2x_3}{x_0}\right.\\
&+&\frac{x_0x_2x_3}{x_1^3}+ \left. \frac{x_2x_3}{x_1^2}+\frac{x_0x_3^2}{x_1^3}+\frac{x_3^2}{x_1^2}\right)-\left( \frac{x_1}{x_0}+\frac{2x_2}{x_1}+\frac{x_2}{x_0}+\frac{2x_3}{x_1}+\frac{x_3}{x_0}\right)\\
&=& \frac{x_1^2}{x_0^2}+\frac{ x_2}{x_0}+\frac{ x_0x_2^2}{x_1^3}+\frac{ x_2^2}{x_1^2}+\frac{ x_3}{x_0}+\frac{x_0x_2x_3}{x_1^3}+\frac{x_2x_3}{x_1^2}+\frac{x_0x_3^2}{x_1^3}+\frac{
      x_3^2}{x_1^2} \\
 &=& \frac{1}{x_0^2x_1^3}\left(x_0x_1^3x_2+ x_0x_1^3x_3, x_0^2x_1x_2^2+x_0^2x_1x_2x_3+ x_0^2x_1x_3^2+x_0^3x_2^2\right.\\&+& \left. x_0^3x_2x_3 + x_0^3x_3^2+x_1^5\right)
\end{array}
\eeq

Note que os monômios de grau 5 em (\ref{normal_eqc}) são exatamente os mesmos em (\ref{geradores_divisor_exp_eqc}). De forma análoga, dado um ponto fixo $f=\langle x_i^2, x_ix_j, c\rangle \in \bb M$, temos que $\cl N_{_f \bb M /\widetilde{\bb X}}$ é dado por 

\beq\label{normal_blowup_3_eqc} \ba {lcl}
\mathcal{N}_{_{f}\mathbb{M}/\widetilde{\bb X}} &=& (x_i^2c)^\vee \otimes(x_ix_kc + x_ix_lc+x_i^2x_jx_k^2+x_i^2x_jx_kx_l\\
 &+& x_i^2x_jx_l^2+x_i^3x_k^2+x_i^3x_kx_l+x_i^3x_l^2+(c^2/x_j))
\ea \eeq
o que nos diz que a fibra $\widetilde{\bb M}_f$ é o espaço projetivo gerado pelas formas de grau 5 
\beq \label{geradores_normal_grau_5} x_ix_kc , x_ix_lc, x_i^2x_jx_k^2, x_i^2x_jx_kx_l, x_i^2x_jx_l^2, x_i^3x_k^2, x_i^3x_kx_l, x_i^3x_l^2,(c^2/x_j).
\eeq

Observe que a decomposição em auto-espaços de $\cl N_{_f\bb M /\widetilde{\bb X}}$ (vide \ref{normal_blowup_3_eqc}) apresenta todos os caracteres distintos. Dessa forma, a $\bb T$-ação induzida sobre $\widetilde{\bb X}'$ apresenta pontos fixos isolados, sendo ao todo $27 (= 3\times 9)$ pontos a menos de permutação das variáveis.  

Agora, seja $f'$ um ponto fixo do divisor excepcional  $\widetilde{\bb M}$ na fibra sobre $f\in \bb M$, para o qual temos que descrever o espaço tangente $\cl T_{f'} \widetilde{\bb X}'$. Por exemplo, dado $f=\langle x_0^2, x_0x_1, x_1^3\rangle$  já obtemos a decomposição do normal $\cl N_{_f\bb M/\widetilde{\bb X}}$ (vide \ref{normal_eqc})  e tomando $f'$  o ponto correspondente ao auto-espaço com caracter $\frac{x_1^2}{x_0^2}$, obtemos

\beq \begin{array}{lcl} \cl T_{f'}\widetilde{\bb X}' &=& \cl L_{f'}\oplus \cl T_f\bb M \oplus \cl T_{[\cl L_{f'}]}\bb P(\cl N_f \bb M/ \widetilde{\bb X})\\

&=& \left( \frac{x_1^2}{x_0^2}\right) + \left( \frac{x_1}{x_0}+\frac{2x_2}{x_1}+\frac{x_2}{x_0}+\frac{2x_3}{x_1}+\frac{x_3}{x_0}\right)\\
&+& \left(\frac{ x_2}{x_0}+\frac{ x_0x_2^2}{x_1^3}+\frac{ x_2^2}{x_1^2}+\frac{ x_3}{x_0}+\frac{x_0x_2x_3}{x_1^3}+\frac{x_2x_3}{x_1^2}+\frac{x_0x_3^2}{x_1^3}+\frac{
      x_3^2}{x_1^2}\right)\otimes \left(\frac{x_1^2}{x_0^2}\right)^\vee\\
&=& \frac{x_1}{x_0}+\frac{ x_1^2}{x_0^2}+\frac{ x_0x_2}{x_1^2}+\frac{ 2x_2}{x_1}+\frac{ x_2}{x_0}+\frac{
      x_0^3x_2^2}{x_1^5}+\frac{ x_0^2x_2^2}{x_1^4}+\frac{ x_0x_3}{x_1^2}+\frac{ 2x_3}{x_1}+\frac{
      x_3}{x_0}+\frac{ x_0^3x_2x_3}{x_1^5}\\
&+& \frac{ x_0^2x_2x_3}{x_1^4}+\frac{
      x_0^3x_3^2}{x_1^5}+\frac{ x_0^2x_3^2}{x_1^4}.
\end{array}
\eeq

E por comparação da auto-decomposição do espaço tangente e da descrição da fibra do divisor excepcional, segue as seguintes identificações:
\beq \begin{array}{c} x_1^5 \leftrightarrow\frac{x_1^2}{x_0^2}, \quad x_0x_1^3x_2\leftrightarrow \frac{ x_2}{x_0}, \quad x_0^3x_2^2\leftrightarrow \frac{ x_0x_2^2}{x_1^3}, \quad x_0^2x_1x_2^2 \leftrightarrow \frac{ x_2^2}{x_1^2}, \quad x_0x_1^3x_3 \leftrightarrow\frac{ x_3}{x_0},\\
x_0^3x_2x_3 \leftrightarrow\frac{x_0x_2x_3}{x_1^3}, \quad x_0^2x_1x_2x_3\leftrightarrow\frac{x_2x_3}{x_1^2}, \quad x_0^3x_3^2 \leftrightarrow\frac{x_0x_3^2}{x_1^3}, \quad x_0^2x_1x_3^2 \leftrightarrow \frac{
      x_3^2}{x_1^2}\end{array}
\eeq

De modo mais geral, dado $f=\langle x_i^2, x_ix_j, c\rangle$ ponto fixo em $\bb M$, e $q_5$ uma das formas quínticas em (\ref{geradores_normal_grau_5}), temos as seguintes identificações: $$q_5 = \langle q_5\rangle = \frac{q_5}{x_i^2c},$$
onde $\langle q_5\rangle$ indica o carácter correspondente à quíntica $q_5$ que aparece na decomposição em auto-espaços do normal (vide \ref{normal_blowup_3_eqc}) e, assim, a decomposição do tangente no ponto fixo $f'$ referente ao carácter $\langle q_5\rangle$ é dada por:

\beq\ba {ccl} 
\mathcal{T}_{f'}\widetilde{\bb X}' &=& \cl L_{f'}\oplus \cl T_f\bb M \oplus \cl T_{[\cl L_{f'}]}\bb P(\cl N_f \bb M/ \widetilde{\bb X})\\\na9
&=& \langle q_5\rangle + 
((x_i)^\vee \otimes (\mathcal{F}_1-x_i) + (x_j)^\vee\otimes(\mathcal{F}_1 - (x_i+x_j))\\\na9 
&+&(c)^\vee\otimes(x_j^2\otimes(\mathcal{F}_1 - x_i)-c))+(\mathcal{N}_{_{f''}\mathbb{M}'/\W'} - \langle q_5\rangle)\otimes \langle q_5\rangle^\vee.
\ea \eeq

E para a contribuição numérica do denominador na fórmula de resíduos de Bott obtemos o produto dos pesos da representação acima.

Na Tabela \ref{tabela_pf_eqc} indicamos cada um dos tipos de classes de isomorfismo de pontos fixos em $\W_{eqc}{}$ e os correspondentes pontos em $\W'_{eqc}$, bem como o número de pontos em cada classe:

\begin{table}[!h]
\centering
\begin{footnotesize}
	\begin{tabular}{|c|c|c|c|}
	\hline 
	Tipo & Pontos fixos em $\W_{eqc}$ & Pontos fixos em $\W'_{eqc}$ & $\#$ pontos \\ 
	\hline 
	(1) &$\langle x_i^2, x_j^2\rangle$ & $\langle x_i^4,x_i^2x_j^2,x_j^4\rangle$ & 6 \\ 
	\hline 
	(2)& $\langle x_i^2, x_jx_k\rangle$ & $\langle x_i^4,x_i^2x_jx_k,x_j^2x_k^2 \rangle$ & 12 \\ 
	\hline 
	(3)& $\langle x_ix_j, x_kx_l\rangle$  & $\langle x_i^2x_j^2,x_ix_jx_kx_l,x_k^2x_l^2 \rangle$ & 3 \\ 
	\hline 
	(4) & $\langle x_i\cdot\langle x_i,x_j\rangle, c\rangle$ \tablefootnote{Aqui $c$ é uma das 3 formas cúbicas: $x_jx_k^2, x_jx_kx_l, x_jx_l^2$}  & $\langle x_i^4, x_i^3x_j x_i^2c, x_i^2x_j^2, x_ix_jc, c^2 \rangle$& $36( = 3\times 12)$ \\ 
	\hline 
	(5) & $\langle x_i\cdot\langle x_j,x_k\rangle, c\rangle$ \tablefootnote{Aqui $c$ indica uma das formas cúbicas: $x_j^2x_k, x_jx_k^2, x_jx_kx_l, x_k^3, x_k^2x_l, x_kx_l^2, x_j^3, x_j^2x_l, x_jx_l^2$} & $\langle x_i^2x_j^2, x_i^2x_jx_k, x_ix_jc, x_i^2x_k^2, x_ix_kc, c^2\rangle$ & 108 ($= 12\times 9)$ \\ 
	\hline 
	(6) & $\langle x_i\langle x_i,x_j,q\rangle,q_4\rangle$ \tablefootnote{$q$ e $q_4$ são descritas em (\ref{pontos_fixos_4_1 - 4_3}) e (\ref{fibra_divisor_2}), respectivamente}  & $\langle x_i^2\cdot \langle x_i^2,x_ix_j,x_iq,q_4,x_j^2,x_jq,q^2\rangle, x_ix_jq_4,
     x_iqq_4,q_4^2  \rangle$  & 324 ($=12\times 3 \times 9$) \\ 
	\hline 
	 (7)& $\langle x_i\cdot \langle x_i, x_j\rangle, c \rangle$ &  $\langle x_i^4, x_i^3x_j, x_i^2c, x_i^2x_j^2, x_ix_jc, c^2,[q_5]\rangle$ \tablefootnote{$c$ é uma formas cúbicas: $x_j^3, x_j^2x_k, x_j^2x_l$. E $q_5$ é uma das formas quínticas: $\{ x_ix_kc, x_ix_lc, x_i^2x_jx_l^2, x_i^3x_k^2, x_i^3x_kx_l, x_i^3x_l^2, (c/x_j)c\}$ } &  324 ($=12\times 3 \times 9$)\\ 
	\hline 
	\multicolumn{3}{|l|}{Total de pontos fixos} & 813\\
\hline
	\end{tabular} 
	\end{footnotesize}
	\caption{Pontos fixos em $\W_{eqc}{}$}\label{tabela_pf_eqc}
\end{table}

Em referência a Tabela \ref{tabela_pf_eqc}, temos que todos os ideais em $W'_{eqc}$ possuem polinômio de Hilbert $12t-16$. Além disso, os pontos fixos em $W'_{eqc}$ do tipo (1) - (3) são 5-regular no sentido de Castelnuovo-Mumford, ao passo  que os pontos fixos dos tipos (4) e  (5) são 6-regular, os pontos do tipo (6) são 8-regular e os pontos fixos do tipo (7) apresentam regularidade 5 ou 6. E por um argumento de semi-continuidade de cohomologia segue que qualquer feixe ideal correspondendo a um ponto em $\W'_{eqc}$ é 8-regular. 

O procedimento de obtenção da representação das fibras de $\mathcal{E}_d$ é análogo ao discutido na seção \ref{sec_QD_twc}, em que a partir da regularidade 8, obtemos a representação da fibra de $\mathcal{E}_d$ sobre o ponto fixo $f$ tomando os monômios de grau $d$ (ao todo $\binom{d+3}{3}-(12d-16)$) presentes no produto $\mathcal{E}_{8_f} * S_{(d-8)}$.  

Reunidas as informações necessárias para o cálculo do $\deg \Sigma(\W_{eqc}{},d)$ via fórmula de resíduos de Bott e sabendo que o grau anterior é polinomial em $d$ (vide Proposição \ref{grau_polinomio}) de grau no máximo $48 ( = 3\times 16)$, segue que é suficiente encontrarmos $\deg \Sigma(\W_{eqc}{},d)$ para 49 valores diferentes de $d$ e, posteriormente, interpolar os resultados para obter o polinômio em questão. O polinômio que nos dá o grau de $\Sigma(\W_{eqc}{},d)$ é explicitado a seguir:

{
\beq \label{grau_singular_eqc}
\ba c
\deg \Sigma_{\W_{eqc}{},d} =
\frac{77991978249}{47023181004800}d^{32}-\frac{142130943}{922746880}d^{31}\\\na9
+\frac{8109239447979}{1175579525120}d^{30}-\frac{4150267051797}{20992491520}d^{29}+\frac{47676232841150619}{11755795251200}d^{28}\\\na9
-\frac{6615027446596551}{104962457600}d^{27}+\frac{128385059997089001}{167939932160}d^{26}\\\na9
-\frac{103459871906659801}{14129561600}d^{25}+\frac{893796960041917863271}{16277254963200}d^{24}\\\na9
-\frac{312845973151702414313}{1017328435200}d^{23}+\frac{4312587609200253695639}{4069313740800}d^{22}\\\na9
+\frac{6155781582234103357}{7266631680}d^{21}-\frac{1105621403101024328482787}{24415882444800}d^{20}\\\na9
+\frac{2134617904050477326290337}{5410337587200}d^{19}-\frac{1027704290752048951537337771}{476109707673600}d^{18}\\\na9
+\frac{1568309607110425883232529237}{223176425472000}d^{17}+\frac{399314335681097660200615893191}{57133164920832000}d^{16}\\\na9
-\frac{127911974311612787565094357769}{396758089728000}d^{15}+\frac{729760755266942589134714032019}{238054853836800}d^{14}\\\na9
-\frac{18285322486683264514566399967249}{892705701888000}d^{13}+\frac{15050777906503580350914982390277}{137339338752000}d^{12}\\\na9
-\frac{8362721204990643447960751421719}{17167417344000}d^{11}+\frac{178565283439979930078484872809}{98099527680}d^{10}\\\na9
-\frac{2731787128737717049736180171243}{476872704000}d^9+\frac{1125598445944774654288515801691861}{74392141824000}d^8\\\na9
-\frac{58025484355390407710374488759691}{1743565824000}d^7+\frac{16796039461040747482814365174429}{278970531840}d^6\\\na9
-\frac{8521350244073783951990040324653}{96864768000}d^5+\frac{599422208545470260381592707347}{5930496000}d^4\\\na9
-\frac{796327032680715287225577370219}{9081072000}d^3+\frac{434272227079029305979707333}{8072064}d^2\\\na9
-\frac{14906420412807524159489839}{720720}d+3713124778880030320\\
\ea
\eeq}
Observe mais uma vez que o grau do polinômio em (\ref{grau_singular_eqc}) é  igual a $(1+1)\times \dim(\W_{eqc}) = 2\times 16$. 

No Apêndice \ref{ap_quartica_eliptica} o leitor pode consultar os procedimentos/funções escritos no Macaulay2, \cite{Macaulay2}, utilizados para obtenção do $\deg\Sigma(\W_{eqc}{},d)$ em (\ref{grau_singular_eqc}) .

\newpage

\chapter*{Apêndices}

Na sequência disponibilizamos os códigos escritos em Macaulay2 \cite{Macaulay2}, Singular \cite{Singular} e Maple \cite{Maple_2015} necessários aos cálculos enumerativos presentes na tese. 

Os códigos do Macaulay2 podem ser executados através do copiar/colar em uma seção disponível em \url{http://habanero.math.cornell.edu:3690/}. Já para os códigos do Singular, o leitor interessado pode abrir uma seção no SageMathCell disponível em \url{https://sagecell.sagemath.org/}.

\appendix

\chapter{Código: Cálculo do \texorpdfstring{$\deg\Sigma(\W_{(2,4)}{},d)$}{} via Schubert2\label{ap_line_schubert2}}
Código para calcular o grau da família de superfícies de grau $d$ singulares ao longo de uma reta usando o pacote Schubert2 do Macaulay2 \cite{Macaulay2}. 
\begin{verbatim}
-------------------------------------------------------------------------
loadPackage "Schubert2";
B = base(d);
G = flagBundle({2,2},B);
(S,Q) = bundles G; 
Ed = symmetricPower(2,S)*symmetricPower(d-2,S++Q) 
     - S * exteriorPower(2,S)*symmetricPower(d-3,S++Q);
deg = integral segre(4,Ed)
factor deg
-------------------------------------------------------------------------
\end{verbatim}

\chapter{Código: Cálculo do \texorpdfstring{$\deg\Sigma(W_m{},d)$}{} via Schubert2\label{ap_curva_m_schubert2}}

Código para calcular o grau da família de superfícies de grau $d$ singulares ao longo de uma curva plana de grau $m$ usando o pacote Schubert2 do Macaulay2 \cite{Macaulay2}. 
\begin{verbatim}
-------------------------------------------------------------------------
loadPackage "Schubert2";

X = base(d);

P3dual = projectiveSpace(3,X);

--Q = polinômios lineares sobre o plano variável
(S,Q) = bundles(P3dual);

m = 2; -- Grau da curva plana

--Sistema de curvas planas de grau m
Wm = projectiveBundle'(dual symmetricPower(m,Q)); 

DIM = dim Wm;

PiOCm = (d1,m1) -> (
Oc1 = symmetricPower(d1,Q)-symmetricPower(d1-m1,Q)**OO_Wm(-1); 
return(Oc1)
);

NdCmPQd = OO_Wm(-1)**PiOCm(d-m,m);

NPQd = OO_P3dual(-1)**PiOCm(d-1,m);

NCmP3 = NdCmPQd ++ NPQd;

PiOsingCm = NCmP3++PiOCm(d,m);

Grau = factor (integral(chern(DIM, PiOsingCm)))
-------------------------------------------------------------------------
\end{verbatim}

\chapter{Funções e procedimentos gerais do Macaulay2 e Maple\label{ap_rotinas_gerais}} 
Este capítulo contém os procedimentos e funções, alguns dos quais foram adaptadas de Meurer \cite{Meurer_96} para a linguagem do software Macaulay2 \cite{Macaulay2} e, basicamente os mesmos procedimentos escritos para a linguagem do Maple \cite{Maple_2015},  que são utilizadas para o cálculo do grau de zero-ciclos descritos em função de classes de Chern de fibrados equivariantes. Aqui, como em todo o texto, temos $T = \mathbb{C}^{*}$ agindo diagonalmente em $\mathbb{P}^n$ , de modo que as coordenadas homogêneas $x_0, \dots , x_n$ são auto-vetores com pesos $w_0 , \dots , w_n$.

\tocless\section{Funções do Macaulay2\label{funcoes_Macaulay2}}
As funções a seguir devem ser carregadas antes dos cálculos específicos dos casos considerados nos Capítulos \ref{cap_pklinear}, \ref{cap_eqc} e seção \ref{sec_QD_twc}, cujos códigos e rotinas estão explicitados nos Apêndices \ref{ap_curva_plana}, \ref{ap_twc} e \ref{ap_quartica_eliptica}. 
\begin{verbatim}
-------------------------------------------------------------------------
--S(d) Retorna a representação da d-ésima potência simétrica nas
--variáveis x_0, ..., x_n. S(d) Necessita do valor n = dimensão de P^n
--Exemplo; dado n = 3, temos que S(2) retorna: 
--x_0^2+x_0*x_1+...+x_1*x_3+x_2^2+x_2*x_3+x_3^2
S = (d) -> (return(sum((ideal(basis(d,r)))_*)));
-------------------------------------------------------------------------
\end{verbatim}
\begin{verbatim}
-------------------------------------------------------------------------
--dualRep retorna a representação dual
dualRep = (p) -> sub(p,dualString);
-------------------------------------------------------------------------
\end{verbatim}

\newpage
\begin{verbatim}
-------------------------------------------------------------------------
--tGrass determina o espaço Tangente à Grassmaniana G=G(r,s) no
--ponto p, T_p G: onde Sr = Symm(r,F) = S(r)
tGrass =(Sr, p) -> (Sr - p)* dualRep(p);
-------------------------------------------------------------------------
\end{verbatim}
\begin{verbatim}
-------------------------------------------------------------------------
--dotProd Efetua o "produto interno"
dotProd =(L1,L2) ->(p = 0;
         if (class L1 === Ideal) then(Laux1 = L1_*;);
         if (class L1 === List) or (class L1 === Sequence) or 
         (class L1 === Array) or (class L1 === Set)  then(
           Laux1 = toList(L1););
         if (class L2 === Ideal) then(Laux2 = L2_*;);
         if (class L2 === List) or (class L2 === Sequence) or 
         (class L2 === Array) or (class L2 === Set)  then(
           Laux2 = toList(L2););
         lm = min(length Laux1, length Laux2);
         for i to (lm -1) do(p = p + Laux1_i * Laux2_i;);
         return(p);
);
-------------------------------------------------------------------------
\end{verbatim}
\begin{verbatim}
-------------------------------------------------------------------------
--pesos retorna uma lista com os pesos da representação de um
--fibrado. Ex.: (2*x_0^3 + x_1 *x_2)/(x_3*x_4) é transformado
--na lista {3*x_0-x_3-x_4,3*x_0-x_3-x_4,x_1+x_2-x_3 - x_4}
pesos = ipesos -> (
       if (isConstant (ipesos)) then(return({0});)
       else(
         lipesos1 ={}; lpts = {};
         if (class ring ipesos === PolynomialRing) then(
           lvspt = (ideal (vars ring ipesos))_*; 
           lipesos2 = terms ipesos;
           for i in lipesos2 do(
             lc = sub(leadCoefficient (i),ZZ);
             lm = leadMonomial(i);
             if (lc != 1) then(
               for j to (lc-1) do(
                 lipesos1 = append (lipesos1, lm););)
             else(lipesos1 = append (lipesos1, lm);););
           for i in lipesos1 do(
             lexp = (exponents i)_0; 
             lpts = append (lpts, dotProd(lvspt, lexp)););
         return(lpts););
         if (class ring ipesos === FractionField) then(
           nipesos = numerator(ipesos);
           dipesos = denominator (ipesos);
           lptsn = pesos(nipesos);
           lptsd = sum(pesos(dipesos));
           lpts = for i in lptsn list i - lptsd;
           return(lpts);)
         else(
           return("Tipo de entrada não suportado");););
);
-------------------------------------------------------------------------
\end{verbatim}
\begin{verbatim}
-------------------------------------------------------------------------
--subPesos substitui a lista de opções lpesos (já definida), 
--Ex.: lpesos = {x_0=>0,x_1=>2,x_2=>7,x_3=>11,x_4=>19,x_5=>37} 
--em uma lista, obtida por exemplo pela função pesos.
--Ex.: subPesos(pesos(tg))
subPesos = ilpesos -> (
         return(for i in ilpesos list sub(i,lpesos));
);
-------------------------------------------------------------------------
\end{verbatim}
\begin{verbatim}
-------------------------------------------------------------------------
--cherns determina as classes de Chern da lista de pesos H.
--Necessita que DIM e o anel rt = QQ[t] estejam definidos.
cherns = (H) -> (lchern = {}; sigma = 1;
   for i in H do(
     sigma = sub(sub((1+i*t)*sigma, rt/ideal(t^(DIM+1))),rt););
   return(for i from 1 to DIM 
           list leadCoefficient(part(i,sigma)));
);
-------------------------------------------------------------------------
\end{verbatim}
\begin{verbatim}
-------------------------------------------------------------------------
--topChern determina a classe Chern top da lista de pesos H.
--assuma rank = dimensão
topChern = (H) -> (return(product (H)););
-------------------------------------------------------------------------
\end{verbatim}
\begin{verbatim}
-------------------------------------------------------------------------
-- bott recebe uma lista {"P1:pt fixo, P2: Tangente,P3:fibra"}, com 
--todos os pontos fixos, tangentes e fibra e o valor d. 
-- Retorna uma lista {d, grau}, onde grau = deg(Sigma(W,d))
bott =(LFTaux, d) -> (
        LFTaux1 = LFTaux;
        dg = d;  resp = 0; 
        for i to (length (LFTaux1)-1) do (
          fb = S(dg) -LFTaux1_i_3; 
          tg = LFTaux1_i_2; 
          pfb = pesos(fb);
          ptg = pesos(tg);
          num = (cherns(subPesos(pfb)))_(DIM-1);
          den = topChern(subPesos(ptg));
          resp = resp+num/den; 
        );
        return {dg, resp};
);
-------------------------------------------------------------------------
\end{verbatim}
\begin{verbatim}
-------------------------------------------------------------------------
--polynomialInterpolation encontra uma função polinomial, pelo método de 
--Lagrange, ajustando uma lista L={{x_i,y_i=f(x_i)}}.
--Necessita de um anel de polinômio em uma variável R = QQ[d].
polynomialInterpolation =(L, R) -> (
       if (numgens R != 1) then 
         error "Número de variáveis excede 1";
       pl:= 0;  n = #L; t:= first gens R;
       for j from 0 to (n-1) do (
         Ld:= 0..(j-1)| (j+1)..(n-1);
         p := for k in Ld list ((d-L_k_0)/(L_j_0 - L_k_0));
         pl = pl + L_j_1 * product(p); 
       );
       return(pl);
);
-------------------------------------------------------------------------
\end{verbatim}

\tocless\section{Procedimentos do Maple\label{funcoes_Maple}}
Os procedimentos a seguir basicamente efetuam as mesmas tarefas que as funções listadas na seção \ref{funcoes_Macaulay2}, sendo que para algumas tarefas elas são mais eficientes do ponto de vista computacional por utilizarem séries geradoras, como no caso do procedimento \verb+S(d)+.
\begin{verbatim}
#########################################################################
#S(d) Retorna a representação da d-ésima potência simétrica nas 
#variáveis x[0], ..., x[N] . S(d) Necessita do valor N = dimensão de P^N
#Exemplo; dado N = 3, temos que S(2) retorna: 
#x[0]^2+x[0]*x[1]+...+x[1]*x[3]+x[2]^2+x[2]*x[3]+x[3]^2
S := proc (d) local genf, i; 
     genf := expand(convert(series(1/(product('-t*x[i]+1', i = 0 .. N)), 
     t, d+1), polynom)); 
     return sort(coeff(genf, t, d))
end:
#########################################################################
\end{verbatim}
\begin{verbatim}
#########################################################################
# dualrep(C) retorna a representação dual
dualrep := proc (C) local i; 
     sort(expand(subs(seq(x[i] = 1/x[i], i = 0 .. N), C))) end:
#########################################################################
\end{verbatim}
\begin{verbatim}
#########################################################################
# simpleweights(A) retorna a soma de uma representação 1-dimensional, 
#uma para cada peso ocorrendo em A. 
# Por exemplo, simpleweights(S(j)*S(k)) irá retornar S(j+k).
simpleweights := proc (p) local mon; 
                 expand(p); 
                 if p = 0 then RETURN(0) fi: 
                 coeffs(%, variables, 'mon'); 
                 convert([mon], `+`) 
end:
#########################################################################
\end{verbatim}
\begin{verbatim}
#########################################################################
#Seja C^* -> T o subgrupo a um parâmetro dado por pesos w[0],w[1],...,w[N]. 
#Desse modo o caracter x[i] restringe a  t^(w[i]), i=0..N. 
# Precisamos encontrar todos os pesos em cada ponto fixo
weightlist := proc (f) local t, cof, mon, res, i,L; 
   L := [seq(x[i] = t^w[i], i = 0 .. N)]:
   cof := [coeffs(f, variables, mon)]; 
   mon := subs(L, [mon]); 
   res := NULL; 
   for i to nops(mon) do 
     res := res, `$`(subs(t = 1, diff(mon[i], t)), cof[i]) od: 
   res 
end:
#########################################################################
\end{verbatim}
\begin{verbatim}
#########################################################################
# subweights substitui os pesos de listweights in f
subweights := proc (f) local apf, saida,L;
              L := [seq(w[i], i = 0 .. N)]: 
              apf := unapply(f, L): 
              saida := apf(op(listweights)) 
end:
#########################################################################
\end{verbatim}
\begin{verbatim}
#########################################################################
#TGrass retorna o tangente da Grassmanniana no ponto p
TGrass := proc (Sr, p) expand((Sr-p)*dualrep(p)): end:
#########################################################################
\end{verbatim}
\begin{verbatim}
#########################################################################
# Blow_up de X ao longo de B. 
# TX = Tangente de X e TB = Tangente de B.
# Retorna NBX o normal relativo, TL direções normais, 
# ETT espaço tangente total
blowUp := proc (TX, TB) local i, NBX, dNBX, L1; 
  L1 := []; 
  NBX := expand(TX-TB); 
  dNBX := denom(NBX); 
  for i in op(NBX) do 
    L1 := [op(L1), ["P2:NBX,P3:TL,P4:ETT",NBX,i,expand(TB+i+(NBX-i)/i)]]:
  od: 
  return (L1):
end:
#########################################################################
\end{verbatim}
\begin{verbatim}
#########################################################################
# partdegree retorna a parte de grau deg do polinômio pol
partdegree := proc (pol, deg) local f, part, i; 
    f := expand(pol); 
    part := 0; 
    for i to nops(f) do 
      if degree(op(i, f)) = deg then part := part+op(i, f) fi:
    od:
    part 
end:
#########################################################################
\end{verbatim}
\begin{verbatim}
#########################################################################
# um retorna o coeficiente do monômio
um := proc (f) local Lvar, L1, apf, saida; 
      Lvar := [seq(x[i], i = 0 .. N)]; 
      L1 := [seq(1, i = 0 .. N)]; 
      apf := unapply(f, Lvar); 
      saida := apf(op(L1)) 
end:
#########################################################################
\end{verbatim}
\begin{verbatim}
#########################################################################
# topchern retorna a classe de Chern top dimensional
topchern := proc (tg) local lpesos, ppesos, k; 
  lpesos := [weightlist(tg)]: 
  ppesos := subweights(product(lpesos[k],k=1 .. nops(lpesos))): 
  ppesos 
end:
#########################################################################
\end{verbatim}
\begin{verbatim}
#########################################################################
# cherns retorna as classes de Chern do fibrado até a dimensão DIM
cherns := proc (tg) local t, c, i1, u, m, z, wtg, ws, lpesos; 
     c := 1; 
     wtg := expand(tg); 
     if whattype(wtg) = `+` or whattype(wtg) = list then 
       for i1 to nops(wtg) do 
         u := op(i1, wtg): 
         m := um(u): 
         z := u/m: 
         lpesos := [weightlist(z)]: 
         ws := subweights(add(lpesos[k],k=1..nops(lpesos))): 
         c := mtaylor(c*(1+t*ws)^m, t, DIM+1) 
       od:
     else 
       u := wtg: 
       m := um(u): 
       z := u/m: 
       lpesos := [weightlist(z)]: 
       ws := subweights(add(lpesos[k], k=1..nops(lpesos))): 
       c := mtaylor(c*(1+t*ws)^m, t, DIM+1) 
     fi:
     c := collect(c, t): 
     [seq(coeff(c, t, i1), i1 = 1 .. DIM)] 
end:
#########################################################################
\end{verbatim}
\begin{verbatim}
#########################################################################
# Bott recebe uma lista [["P2:pt fixo,P3:fibra,P4:Tangente", pf fixo, 
#fibra, tangente]], com todos os pontos fixos, fibra, tangentes e o grau d. 
#Retorna uma lista [d, grau], onde grau = deg(Sigma(W,d))
Bott := proc (LFTaux, d) 
        local LFTaux1, resp, fp, tg, num, den, dg, i; 
        LFTaux1 := LFTaux: 
        dg := d: 
        resp := 0: 
        for i to nops(LFTaux1) do 
          fp := S(dg)-LFTaux1[i][3]: 
          tg := LFTaux1[i][4]: 
          num := cherns(fp)[DIM]: 
          den := topchern(tg): 
          resp := resp+num/den: 
        od:
        return ([dg, resp]):
end:
#########################################################################
\end{verbatim}

\chapter{Código: Cálculo do \texorpdfstring{$\deg\Sigma(\W_{(k+1,n+1)}{},d)$}{}\label{ap_pk_linear}}
Neste apêndice disponibilizamos o código para executar o cálculo do grau da subvariedade $\Sigma(W_{(k+1,n+1)}{},d)\subset \mathbb{P}^{N_d}$, hipersuperfícies em $\P{n}$ de grau $d$ singulares ao longo de algum $\P{k}$-linear(variável), usando o software Maple \cite{Maple_2015}. Para execução dos cálculos é necessário a inclusão dos procedimentos da seção \ref{funcoes_Maple}, bem como os procedimentos adicionais listados abaixo.

\vspace{0.5cm}
\noindent\verb+#########################################################################+

\noindent\verb+--Inserir os procedimentos da seção+ \ref{funcoes_Maple} \verb+aqui+
\begin{verbatim}
#########################################################################
#Procedimentos adicionais
#list é uma lista não vazia, 1 <= ini <=nops(list)
#retorna os elementos da list de ini até nops(list)
sublist := proc (list, ini) 
  [seq(list[i], i = ini .. nops(list))]:
end proc:

#elt é qualquer elemento válido, e o argumento "list" é uma lista de 
#listas. #elt é adicionado ao início de cada elemento (o qual é uma 
#lista) no argumento "list"
addElt := proc (elt, list) 
local templist, i; 
templist := []:
 for i to nops(list) do 
   templist := [op(templist), [elt, op(list[i])]] 
 od:
 templist 
end proc:

#L é uma lista não vazia, 1 <= i <= nops(L)
#retorna uma lista de todos os possíveis subconjuntos 
#(sublistas) de) L com i elementos
listSubsets := proc (L, i) 
 local n, j, fl, temp: # fl - final list
 n := nops(L):
 if i = 1 then #base case
  fl := []:
  for j to n do 
   fl := [op(fl), [L[j]]]:
  end do:
 else 
  fl := []:
  for j to n-i+1 do 
    temp := listSubsets(sublist(L, j+1), i-1):
    temp := addElt(L[j], temp):
    fl := [op(fl), op(temp)]:
  end do:
 end if:
 fl:
end proc:

#PFTpklinear retorna os pontos fixos de G(kl+1,nl+1) 
PFTpklinear := proc (kl, nl)
 local lpaux, lp, var, i, j, LPFTaux:
 var := [seq(x[i], i = 0 .. nl)]:
 lpaux := listSubsets(var, nl-kl):
 lp := []:
 for i in lpaux do 
   lp := [op(lp), add(i)]:
 od:
 LPFTaux := []:
 for j in lp do 
   LPFTaux := [op(LPFTaux), ["P2: pt fixo, P3: Tangente",j,TGrass(S(1),j)]]:
 od:
 return LPFTaux:
end proc:
#########################################################################
\end{verbatim}
\begin{verbatim}
#########################################################################
#Cálculos iniciam aqui
#Escolha N e k como desejar (K < N)
N := 4: #Dimensão de P^N
variables := [seq(x[i], i = 0 .. N), seq(1/x[i], i = 0 .. N)]:
K := 2: #Dimensão de P^K
DIM := (K+1)*(N-K): #Dimensão da Grassmanniana G(K+1,N+1) 
M := (K+1)*DIM+4: #Limite para o grau

#Deve-se adequar os pesos em listweights de acordo com o valor de N. 
#São necessários pelo menos N+1 pesos distintos. A lista abaixo serve
#para os casos de N<=8. Caso ocorra algum erro na execução, como 
#division_by_zero, deve-se alterar os pesos 
listweights := [4, 11, 17, 32, 55, 95, 160, 267, 441]:

PFT := PFTpklinear(K, N): #Lista com os pontos fixos/tangentes

listresp := []:
for dg from 3 to M do 
 PFTF := []:
 for i in PFT do 
  PFTF := [op(PFTF), ["P2: pt fixo, P3: fibra, P4: Tangente",i[2], 
  simpleweights(i[2]^2*S(dg-2)), i[3]]]:
 od:
 listresp := [op(listresp), Bott(PFTF, dg)]:
 print(listresp[dg-2]):
od:
#Cálculo do polinômio gerador
with(CurveFitting):
list1 := [seq(listresp[i][1], i = 1 .. M-2)]:
list2 := [seq(listresp[i][2], i = 1 .. M-2)]:
Grau := sort(PolynomialInterpolation(list1, list2, d)):
Grau := factor(Grau); 
convert(Grau, string);
#########################################################################
\end{verbatim}

\chapter{Código: Cálculo do \texorpdfstring{$\deg \Sigma(\W_m{},d)$}{}\label{ap_curva_plana}}
Nesse apêndice disponibilizamos o código para efetuar o cálculo do grau da subvariedade $\Sigma(W_m{},d)\subset \mathbb{P}^{N_d}$, superfícies em $\P3$ de grau $d$ singulares ao longo de alguma curva plana de grau $m>1$, usando o software Macaulay2 \cite{Macaulay2}. Para execução dos cálculos é necessário a inclusão das funções da seção \ref{funcoes_Macaulay2}. 

\vspace{0.5cm}
\noindent\verb+--Inserir as funções da seção+ \ref{funcoes_Macaulay2} \verb+aqui+
\vspace{-0.3cm}
\begin{verbatim}-------------------------------------------------------------------------
--Cálculos iniciam aqui
n = 3; --dimensão de P^3
r = QQ[x_0..x_n]; 
dualString = for i from 0 to n list(x_i => 1/x_i);
m = 3; --Grau da curva (m>=2)
--Dimensão da Variedade das curvas planas de grau m em P3
DIM = 3+(binomial(m+2,2) - 1) 
rg = 2*m --regularidade 
M = 3*DIM + rg; 

--Obs.: Caso apareça a mensagem:  error: division by zero no cálculo 
--usando Bott os pesos abaixo devem ser alterados
lpesos = {x_0 => 11, x_1 => 17, x_2 => 32, x_3 => 55}; 
rt = QQ[t]; 

--Obs.: As fibras aqui são para o caso d = rg. Devido a regularidade para 
--d>rg basta completar o grau da fibra do caso d = rg multiplicando por 
--formas lineares, etc
PFTF = {}; 
for i to n do(--xx formas de grau m que não contém x_i
   xx = (trim((ideal(terms(S(1) - x_i)))^m))_*; 
   for j in xx do(
     idd = sum(((ideal(x_i, j))^2)_*);
     fbaux = 0;
     for k from 2 to rg do(
       fbaux = fbaux + part(k,idd)*S(rg-k););
       fb = sum((trim(ideal(terms(fbaux))))_*);
   PFTF = append(PFTF, {"P1: pt fixo, P2: Tangente, P3: fibra", x_i + j, 
         tGrass(S(1),x_i)+tGrass(sum(xx),j),fb});
   );
);

--Cálculo do grau usando a Fórmula de Bott 
listresp = {};
for d from rg to M do(
   listresp = append(listresp, bott(PFTF,d));
   PFTF = for j to (#PFTF -1) list{PFTF_j_0,PFTF_j_1, PFTF_j_2,
    sum((trim(ideal( terms(PFTF_j_3 * S(1)))))_*)};
);

listresp --lista com os dados calculados

--Polinômio que dá o grau de \Sigma(W_m,d)
Grau = factor polynomialInterpolation(listresp,QQ[d])
Grau = toString Grau
-------------------------------------------------------------------------
\end{verbatim}

\chapter{Código: Cálculo do \texorpdfstring{$\deg\Sigma(\W_{twc}{},d)$}{}\label{ap_twc}}
Neste apêndice disponibilizamos o código para o cálculo do grau da subvariedade $\Sigma(\W_{twc}{},d)\subset \mathbb{P}^{N_d}$, superfícies em $\P3$ de grau $d$ singulares ao longo de alguma cúbica reversa. Inicialmente, na Seção \ref{codigos_singular_twc} utilizamos o Singular \cite{Singular} munido da caixa
de ferramentas "myprocs" , procedimentos/funções criadas por Vainsencher (vide Seção \ref{codigos_Israel}) , para as
contas locais que revelam o local de indeterminação do mapa racional $\W_{twc}\dashrightarrow \bb G(6, \cl F 4 )$
definido por $\mu$ (vide Proposição \ref{prop_mu_twc_bwup}). Na Seção \ref{codigos_singular_twc} efetuamos a explosão de $\W_{twc}$
ao longo de $\bb G$ e verificamos que isto resolve a indeterminação. Já na Seção \ref{codigos_macaulay2_twc} utilizamos 
o software Macaulay2 \cite{Macaulay2} para gerar todos os pontos fixos, tangentes e a representação da fibra de $\cl E_d$, necessárias a aplicação da fórmula de resíduos de Bott, bem como o cálculo das contribuições numéricas. Além dos procedimentos/funções da seção \ref{funcoes_Macaulay2} é necessário a inclusão das 
funções expand e blowUp listadas em \ref{codigos_macaulay2_twc}.

\tocless\section{myprocs por Vainsencher\label{codigos_Israel}}

\vspace{0.5cm}



\tocless\section{Códigos do Singular\label{codigos_singular_twc}}

\noindent\verb+//Copiar/colar códigos do Apêndice+ \ref{codigos_Israel}
\begin{verbatim}

LIB "primdec.lib";
proj(3);
ring r=0,(a(1..3)(1..7)),dp;r=r+P3;setring r;imapall(P3);
def xx2,xx3,xx4,xx5,xx6=xx^2,xx^3,xx^4,xx^5,xx^6 ;
def aa = ideal(a(1..3)(1..7));
def I = ideal(x(0)^2, x(0)*x(1), x(1)^2);
def quads = I;

//Vizinhança afim U
for(int i=1;i<=3;i++){
   quads[i]=quads[i]+dotprod([a(i)(1..7)],omit(xx2,I));};

//Mapa de multiplicação \lambda
def cubs = quads*xx;
def m=transpose(coeffs(cubs,xx3,xxx));
ran(m);//12
def l=rowsnopivo(m);m=l[1];l[2];l[3];
def m0=submat(m,l[2],1..ncols(m));
m=submat(m,omit(1..nrows(m),l[2]),1..ncols(m));
nrows(m);//10
putsolvform(ideal(m0));
def Y=_;

//Impondo a condição de posto <=10
def twcs = dosubs(Y,quads);

hilbsp(aa,twcs);//1+3*t

//indeterminadas (12)
def mtwcs=transpose(coeffs(twcs,xx2,xxx));
pol2id(indets(mtwcs)); aa=_;
//_[1]=a(1)(1)
//_[2]=a(1)(2)
//_[3]=a(1)(3)
//_[4]=a(1)(4)
//_[5]=a(2)(1)
//_[6]=a(2)(2)
//_[7]=a(2)(3)
//_[8]=a(2)(4)
//_[9]=a(3)(1)
//_[10]=a(3)(2)
//_[11]=a(3)(3)
//_[12]=a(3)(4)

def vs = aa;

//Mapa de multiplicação \mu 
def twcs2 = twcs^2;
hilbsp(aa,twcs2);//-7+9*t
std(origin(vs,twcs2));  def twcs20=_;
def m=transpose(coeffs(twcs2,xx4,xxx));
ran(m);//6
def l=rowsnopivo(m);m=l[1];l[2];l[3];
def m0=submat(m,l[2],1..ncols(m));
m=submat(m,omit(1..nrows(m),l[2]),1..ncols(m));
nrows(m);ran(origin(aa,m));//5

//Geradores da fibra do divisor excepcional 
def im02 = reduce(ideal(m0),std(aa^2));
for(int i=1;i<=35;i++){
  if(im02[i]<>0){see(im02[i],xx4[i]);};};

//-a(3)(1),x(0)^3*x(2)
//-a(3)(2),x(0)^3*x(3)
//2*a(2)(1)-a(3)(3),x(0)^2*x(1)*x(2)
//2*a(2)(2)-a(3)(4),x(0)^2*x(1)*x(3)
//-a(1)(1)+2*a(2)(3),x(0)*x(1)^2*x(2)
//-a(1)(2)+2*a(2)(4),x(0)*x(1)^2*x(3)
//-a(1)(3),x(1)^3*x(2)
//-a(1)(4),x(1)^3*x(3)

//Família a 1-parâmetro
def TW = twcs[1]*twcs[3] - twcs[2]^2;
for(int i1=1; i1<=12; i1++){
origin(omit(aa,aa[i1]),TW);
def quarts = origin(omit(aa,aa[i1]),twcs2);
hilbsp(aa,quarts); //9t-7
origin(aa[i1],sat(quarts,aa[i1])[1]);
hilbsp(_);}; //9t - 7

//Colocando na forma localmente gráfico
putsolvform(ideal(m0));
def W =_;

//Local de indeterminação do mapa \mu 
def Indet = dosubs(W,twcs2);
list G = primdecGTZ(Indet);
G;

//Explosão de W_{twc} em G
ring r1=0,d(1..size(W)-1),dp;r1=r1+r;setring r1;imapall(r);
l=m0,m,vs;
int i2;

i2++; //Escolha do divisor excepcional
def exc2=W[i2];def W1=omit(W,exc2);
W1=seq("W1[i]-exc2*d(i)",1,size(W1));
def l2=dosubs(W1,l);
m0,m,vs=l2[1..size(l2)];vs=pol2id(indets(vs));
mysat(m0,m);
twcs2=ideal(m*transpose(xx4));

if(ran(origin(vs,m0))<>0){
twcs2=twcs2+ideal(m0*transpose(xx4));};
hilbsp(vs,twcs2);

def twcs21 = std(origin(vs,twcs2));
see("As 6 quarticas = ", twcs21);

//Nova quártica
def QN = std(reduce(twcs21,std(twcs20)));

\end{verbatim}

\tocless\section{Códigos do Macaulay2\label{codigos_macaulay2_twc}}
\vspace{0.5cm}
\noindent\verb+--Inserir as funções da seção+ \ref{funcoes_Macaulay2} \verb+aqui+
\begin{verbatim}
-------------------------------------------------------------------------
--expand recebe um polinômio com coeficientes inteiros/ retorna uma lista
--com os termos expandindos do polinômio
--Ex.: expand(2*x_0) retorna {x_0,x_0}
expand = (p) -> ( LT = {};
             (M,C) = coefficients(p);
             for i from 0 to (numColumns(M)-1) do(
               c = sub(C_0_i,ZZ);
               if (c >= 1) then(
                 for j from 1 to c do (
                    LT = append(LT, M_i_0););)
               else (
                   for j from 1 to abs(c) do (
                     LT = append(LT, - M_i_0);););); 
             return(LT););
-------------------------------------------------------------------------
\end{verbatim}
\begin{verbatim}
-------------------------------------------------------------------------
--blowUp efetua a explosão de X ao longo de B. 
--Entrada: TX = Tangente de X e TB = Tangente de B
--Retorna NBX o normal relativo, TL direções normais, 
--ETT espaço tangente total
blowUp = (TX, TB) -> ( 
   NBX = TX - TB; 
   nNBX = expand(numerator NBX);
   dNBX = denominator NBX;
   TL = for i in nNBX list i/dNBX; 
   return for i in TL 
            list {"P1:NBX,P2:TL,P3:ETT",NBX,i,TB+i+(NBX-i)/i};
);
-------------------------------------------------------------------------
\end{verbatim}
\begin{verbatim}
-------------------------------------------------------------------------
--Cálculos iniciam aqui
n = 3; -- Dimensão de P3
r = QQ[x_0..x_n]; 
dualString = for i from 0 to n list(x_i => 1/x_i);
DIM = 12; --Dimensão de W_twc
m = 40; -- 4 + 3*12
lpesos = {x_0=>11, x_1 => 17, x_2 => 32, x_3=>55}; 
rt = QQ[t];

--Pontos fixos em W_twc-  Tipos (1) - (4) e (5.1) - (5.4)
PF_1 = ideal(x_0*x_1,x_1*x_2,x_2*x_3);
PF_2 = ideal(x_0*x_1,x_0*x_2,x_1*x_2);
PF_3 = ideal(x_0*x_1,x_0*x_2,x_2^2);
PF_4 = ideal(x_0^2,x_0*x_1,x_1^2);
PF_51 = ideal(x_0^2,x_0*x_1,x_0*x_2,x_1*x_2*x_3);
PF_52 =  ideal(x_0^2,x_0*x_1,x_0*x_2,x_1*x_2^2);
PF_53 = ideal(x_0^2,x_0*x_1,x_0*x_2,x_2^2*x_3);
PF_54 = ideal(x_0^2,x_0*x_1,x_0*x_2,x_2^3);

--Tangentes dos Pontos fixos em W_twc -  
--Tipos (1) - (4) e (5.1) - (5.4)

TG_1 = (x_0^3*x_1+x_0^2*x_1*x_2+x_0*x_1*x_2^2+x_0^2*x_1*x_3+x_0*x_1^2*x_3+
x_0^2*x_2*x_3+x_1^2*x_2*x_3+x_0*x_2^2*x_3+x_0*x_1*x_3^2+x_0*x_2*x_3^2+
x_1*x_2*x_3^2+x_2*x_3^3)/(x_0*x_1*x_2*x_3);

TG_2 = (x_0^2*x_1+x_0*x_1^2+x_0^2*x_2+x_1^2*x_2+x_0*x_2^2+x_1*x_2^2+
2*x_0*x_1*x_3+2*x_0*x_2*x_3+2*x_1*x_2*x_3)/(x_0*x_1*x_2);

TG_3 = (x_0*x_1^3+x_0^2*x_1*x_2+x_0*x_1^2*x_2+x_0^2*x_2^2+x_1^2*x_2^2+
x_1*x_2^3+x_0*x_1^2*x_3+2*x_0*x_1*x_2*x_3+x_0*x_2^2*x_3+x_1*x_2^2*x_3+
x_2^3*x_3)/(x_0*x_1*x_2^2);

TG_4 = (x_0^3*x_2+2*x_0^2*x_1*x_2+2*x_0*x_1^2*x_2+x_1^3*x_2+x_0^3*x_3+
2*x_0^2*x_1*x_3+2*x_0*x_1^2*x_3+x_1^3*x_3)/(x_0^2*x_1^2);

TG_51 = (x_0*x_1^3+x_0*x_1^2*x_2+x_0*x_1*x_2^2+x_0*x_2^3+x_0*x_1^2*x_3+
x_1^2*x_2*x_3+x_0*x_2^2*x_3+x_1*x_2^2*x_3+x_0*x_1*x_3^2+x_0*x_2*x_3^2+
2*x_1*x_2*x_3^2)/
(x_0*x_1*x_2*x_3);

TG_52 = (x_0*x_1^3+x_0*x_1^2*x_2+x_1^2*x_2^2+x_0*x_2^3+2*x_1*x_2^3+
x_0*x_1^2*x_3+2*x_0*x_1*x_2*x_3+2*x_0*x_2^2*x_3+x_1*x_2^2*x_3)/
(x_0*x_1*x_2^2);

TG_53 = (x_0*x_1^4+x_0*x_1^3*x_2+x_0*x_1^2*x_2^2+x_0*x_1*x_2^3+
x_0*x_1^3*x_3+x_0*x_1^2*x_2*x_3+x_1^2*x_2^2*x_3+x_1*x_2^3*x_3+
x_0*x_1*x_2*x_3^2+x_0*x_2^2*x_3^2+x_1*x_2^2*x_3^2+x_2^3*x_3^2)/
(x_0*x_1*x_2^2*x_3);

TG_54 = (x_0*x_1^4+x_0*x_1^3*x_2+x_0*x_1^2*x_2^2+x_1^2*x_2^3+x_1*x_2^4+
x_2^5+x_0*x_1^3*x_3+x_0*x_1^2*x_2*x_3+2*x_0*x_1*x_2^2*x_3+x_0*x_2^3*x_3+
x_1*x_2^3*x_3)/(x_0*x_1*x_2^3);

--Lista com os pontos fixos e tangentes (a menos de permutação das 
--variáveis)
LPFT = {};
for i from 1 to 4 do(
 LPFT=append(LPFT, {"P1: pt fixo, P2: tangente", PF_i, TG_i});
);

for i from 51 to 54 do(
 LPFT=append(LPFT, {"P1: pt fixo, P2: tangente", PF_i, TG_i});
);

--Explosão sobre os pontos fixos tipo (4)
TB = tGrass(S(1),x_0+x_1); -- tangente do centro de explosão
TX = TG_4; --tangente do ambiente
LB4 = blowUp(TX,TB);

--Pontos fixos/tangentes em W'_twc do tipo (4.1) - (4.8)
LPFT4 = for i in LB4 list {"P1: pt fixo, P2: tangente", 
 trim(PF_4^2 +sub(i_2*denominator(i_1),r)), i_3};

--Lista com todos os pontos fixos/tangentes em W'_twc(a menos
-- de permutação das variáveis)
LPFTA = {};
for i from 0 to 7 do(
   if (i != 3) then (
     LPFTA = append(LPFTA, {"P1: pt fixo, P2: tangente", 
      trim (LPFT_i_1^2), LPFT_i_2});
   )
   else(
    for j in LPFT4 do(
     LPFTA=append(LPFTA,{"P1: pt fixo, P2: tangente",j_1,j_2});
    );
   );
);

--Geração de todos os pontos fixos/tangentes (permmutação da variáveis)
--Lista de Permutações
ListP = permutations toList (x_0..x_n);
ListOP = for i in ListP 
         list {x_0 => i_0, x_1 => i_1, x_2 => i_2, x_3 => i_3};

for i to (length LPFTA - 1) do (
PFT_(i+1) = for j in ListOP list {"P1: pt fixo, P2: tangente",
  sum((sub(sub(LPFTA_i_1,r), j))_*),
sub(numerator LPFTA_i_2, j)/sub(denominator LPFTA_i_2, j)};);

LPFTT = {};
for i from 1 to 15 do(
   for j in PFT_i do(
     LPFTT = append(LPFTT, {j_0, j_1, j_2});
   );
);
LPFTT = unique LPFTT;
length LPFTT -- 172 pontos fixos
            
-- Geração das fibras de E_d
-- Cálculo do grau de \Sigma(W_twc,d) para d = 4, 5 ou 6
--Obs. 1: 
--Para salvar os dados calculados no arquivo Grau_TWC.txt 
--criado/salvo na pasta /home/.../Singular_Surface_TWC
--utilizamos os comandos a seguir:
--arquivo_saida = openOutAppend "/home/.../Grau_TWC.txt";
--arquivo_saida << (result) <<endl<<close;

-- Geração das fibras e graus para os casos 4 <= d <= 6
listresp = {};
for d from 4 to 6 do(
  LPFTF = for i in LPFTT 
     list {"P1: pt fixo, P2: tangente, P3: fibra", i_1, i_2, 
     sum((trim(ideal(terms(part(4,i_1)*S(d-4)+
  part(5,i_1)*S(d-5)+part(6,i_1)*S(d-6)))))_*) };
  result = bott(LPFTF,d); print(result);
  listresp = append(listresp, result);
  arquivo_saida = openOutAppend "/home/weversson/Dropbox/
DOUTORADO/Contas_Macaulay_Singular_Maple/
Singular_Surface_TWC/Grau_TWC.txt";
  arquivo_saida << (result) <<endl<<close;
);

-- Cálculo do grau para d > 6 (= regularidade)
--Obs.: Se quiser calcular o grau para um valor específico, sem rodar 
--todos os casos, basta iniciar d com o valor desejado e trocar m pelo
--mesmo valor no código abaixo. 
--Por exemplo: para d = 20 -> for d from 20 to 20 do

for d from 7 to m do (
  LPFTFr = for i in LPFTF 
  list {"P1: pt fixo, P2: tangente, P3: fibra", i_1,i_2, 
  sum(((trim(ideal(terms (i_3 * S(d-6))))))_*)};
  result = bott(LPFTFr,d); print(result);
  listresp = append(listresp, result);
  arquivo_saida = openOutAppend "/home/weversson/Dropbox/
  DOUTORADO/Contas_Macaulay_Singular_Maple/
  Singular_Surface_TWC/Grau_TWC.txt";
  arquivo_saida << (result) <<endl<<close;
);

--Polinômio que dá o grau de \Sigma(W_twc,d)
Grau = factor polynomialInterpolation(listresp,QQ[d])

--Salvando o polinômio
Grau = toString Grau
arquivo_saida = openOutAppend "/home/weversson/Dropbox/
DOUTORADO/Contas_Macaulay_Singular_Maple/
Singular_Surface_TWC/Grau_TWC.txt";
arquivo_saida << (Grau) <<endl<<close;
-------------------------------------------------------------------------
\end{verbatim}

\chapter{Código: Cálculo do \texorpdfstring{$\deg\Sigma(\W_{rc}{},d)$}{}\label{ap_cubica_regrada_P4}}

Inicialmente utilizamos o Macaulay2 para o cálculo dos tangentes dos pontos fixos antes da primeira explosão (vide \ref{diagrama_blowup_rc}), onde adaptamos o código em Vainsencher \cite{Script_Israel}. Posteriormente, por questões técnicas computacionais (tempo/memória), utilizamos o Maple para gerar as explosões necessárias, bem como todos os pontos fixos/tangentes a partir dos pontos fixos/tangentes iniciais (vide (\ref{tg_1_rc}) - (\ref{tg_f_rc}) ) e fazer as substituições de pesos para o cálculo dos graus utilizando a fórmula de resíduos de Bott.

\tocless\section{Códigos do Macaulay2\label{codigo-macaulay2_tangente_P4}}
As funções abaixo são para o cálculo dos respectivos tangentes dos pontos fixos do tipo \ref{pf_1} a \ref{pf_5} (vide (\ref{aplicando_bott_twc})).
\begin{verbatim}
--Procedimentos/funções
-------------------------------------------------------------------------
S = (d) -> (return(sum((ideal(basis(d,r)))_*)));
-------------------------------------------------------------------------
---------------------------------------------------------------
--retorna a representação dual
dualRep = (p) -> sub(p,dualString);
-------------------------------------------------------------------------
-------------------------------------------------------------------------
-- Determina o espaço Tangente à Grassmaniana G=G(r,s) no ponto
-- p, T_p G: onde Sr = Symm(r,F) = S(r)
tGrass =(Sr, p) -> (Sr - p)* dualRep(p);
-------------------------------------------------------------------------
\end{verbatim}
\begin{verbatim}
-------------------------------------------------------------------------
--omit omite do(a) (ideal, lista,sequence, conjunto) A  o elemento 
--(ou lista/sequencia de elementos) B
omit = (A,B) -> ( Laux2 = {};
        if (class B === List) or (class B === Sequence) then (
            Baux = B;)
         else (
            Baux = {B};);
         if (class A === Ideal) then (
           Laux1 = A_*;
           for i in Laux1 do(
             if (not member(i,Baux)) then (
               Laux2 = append (Laux2, i);););
           saida  = ideal(Laux2););
         if (class A === List) then(
           Laux1 = A;
           for i in Laux1 do(
             if (not member(i,Baux)) then (
               Laux2 = append (Laux2, i);););
           saida = Laux2;);
        if (class A === Sequence) then(
           Laux1 = A;
           for i in Laux1 do(
             if (not member(i,Baux)) then (
               Laux2 = append (Laux2, i);););
           saida = toSequence (Laux2););  
         if (class A === Set) then(
           Laux1 = A;
           Laux2 = Laux1 - set(Baux);
           saida = Laux2;);
         return(saida);
);
-------------------------------------------------------------------------
\end{verbatim}
\begin{verbatim}
-------------------------------------------------------------------------
--origin recebe v (ideal com variáveis/ou uma soma de variáveis de um 
--anel A) e p (um polinômio, um ideal ou lista de polinômios do mesmo
--anel A) e retorna p zerando as variáveis que aparecem em v
origin = (v,p,rA) -> (
            use(rA);
            laux = {};
            lvv= support(v); --lista com variáveis de v
            lov = for i in lvv list sub(i,rA)=> 0; 
            qo = p;
            if (class qo === List) then( 
             laux = for i in qo list sub(i,lov);
             return(laux);)
            else(
              return(sub(qo,lov)););
);
-------------------------------------------------------------------------
\end{verbatim}
\begin{verbatim}
-------------------------------------------------------------------------
--wt2frac transforma um polinômio g em uma fração, onde:
--numerador = produto dos monômios com coeficientes positivos, estes 
--transformados em  potência (Ex.: 2x_0 -> x_0^2)
--denominador= produto dos monômios com coeficientes negativos, estes 
--transformados em potência (Ex.:-2x_0 - x_1-> x_0^2*x_1)  
wt2frac = (g) -> (
           ng=1; dg=1;
           lw2f:= terms(g);
           for i in lw2f do(
             lc := sub(leadCoefficient(i),ZZ);
             lm :=  leadMonomial(i);
             if ( lc < 0) then(
               dg= dg*(lm)^(-lc);)
             else(
               ng = ng*(lm)^lc;););
             return(ng/dg);
);
-------------------------------------------------------------------------
\end{verbatim}
\begin{verbatim}
-------------------------------------------------------------------------
wt = ip -> (
 if (not(class ip === Sequence)) then(
   return("The function wt needs at least 2 arguments, the 
           latter being a ring of polynomials");)
 else(
   if (not (class last(ip) === PolynomialRing)) then(
    return("The function wt needs at least 2 arguments, the 
           latter being a ring of polynomials");)
   else(
     br = last(ip); -- Base Ring
     use(br);
     if (length ip == 3) then(
       s = length ip; p := sub(ip_0,br); q = sub(ip_1, br);)
     else(
       p = sub(ip_0,br);
       s= 1;);););
 if ((p == 0) or (class p === ZZ)) then 
   return(0)
 else ( 
   if (class remember === Symbol) then (
     remember = sub(matrix{{x_0..x_n},{x_0..x_n}},br);
   ); 
   known = sub(trim(ideal(remember^{0})),br);
   vb = sub(ideal(generators(br,CoefficientRing=>QQ)),br); 
   p = sub(leadMonomial (p),br);
   vs = sub(ideal(support(p)),br);
   vsmodknown := sub(trim(sub(vs,br/known)),br);
   if (vsmodknown == 0) then(
     if ( sum(degree p) == 1) then (
   for j from 0 to (numcols remember -1) do (
     if (p == sub(remember_(0,j),br)) then ( sj = j; break;););
       return(remember_(1,sj));)
     else(
       iv = (exponents(p))_0;
       vbn = 0;
       for j from 0 to (length(iv) -1) do(
         if (iv_j != 0) then(
           vbn = vbn + iv_j * wt(vb_j,br);););
   return(vbn););)
 else(
   if (s == 1) then(
     return("--Think! need ideal");)
   else( 
    for j from 0 to (numgens q -1) do(
      vsmodknown = trim sub(sub(ideal(support(vsmodknown*q_j)),
                        br/known),br);
      qj = sub(origin(vsmodknown, q_j,br),br);
      lqj = terms(qj);
      if (sub(ideal(support(qj)),br/known) !=0) then(
        return("bad news: -(");)
      else(
       qi = sub(q_j - qj,br);
       lqi = terms(qi);
       for i in lqi do(
         newv =sub(origin(known,sum(support(i)),br),br); 
         if (size newv == 1) then (
           wtnewv = - wt(sub(i,br)//newv,br)+wt(lqj_0,br);
           mnewv = matrix{{sub(newv,br)},{sub(wtnewv,br)}};
           known = trim(ideal(remember^{0}));
           remember = sub(sub(remember,br)|mnewv,br);
          ););););
        if (sub(p, br/known) ==0) then (
          return(wt(p,br));)
        else(
          return("work harder"););
  ););););
-------------------------------------------------------------------------
\end{verbatim}
\begin{verbatim}
-------------------------------------------------------------------------
--Cálculos dos tangentes iniciam aqui
n = 4; --Dimensão de P^n. n= 3,4 ou 5
r = QQ[x_0..x_n]; 
dualString = for i from 0 to n list(x_i => 1/x_i);
n1 = binomial(2+n,n) - 3;
n2 = 2*3*(n+1);
r1 = QQ[x_0..x_n,a_(1,1)..a_(3,n1),b_1..b_n2];

LVPn = toList(x_0..x_n); -- Lista das variáveis de Pn;

IVPn = ideal(LVPn); --ideal das variáveis de Pn

--Pontos Fixos Tipos 1 - 5
PF1 = ideal(x_0*x_1, x_1*x_2,x_2*x_3); 
PF2 = ideal(x_0*x_1, x_1*x_2, x_0*x_2);
PF3 = ideal(x_0*x_1, x_2, x_0*x_2);
PF4 = ideal(x_0^2, x_0*x_1, x_1^2);
PF5 = ideal(x_0^2, x_0*x_1, x_0*x_2);

--Matrizes que geram os pontos fixos
mPF1 = matrix{{x_1, x_3, 0},{0, x_0, x_2}};
mPF2 = matrix{{x_0, 0, -x_2},{0, x_1, x_2}};
mPF3 = matrix{{x_0, x_2, 0},{0, x_1, x_2}};
mPF4 = matrix{{x_0, x_1, 0},{0, x_0, x_1}};
mPF5 = matrix{{x_0, 0, -x_2},{0, x_0, x_1}};

LmPF = {mPF1, mPF2, mPF3, mPF4,  mPF5}; -- Lista das matrizes

--Lista com os pontos fixos Tipo 1 - 5
LPF = for i in LmPF list minors(2, i); 

--Precisamos da imagem do mapa tangente
bb = ideal(b_1..b_n2);

--Pertubação linear de cada matriz
Lpm = {}; --Lista das matrizes após pertubação
for i in LmPF do(
 m0 = mutableMatrix i;
 cont = 0;
 for j to (numrows m0 -1) do (
   for k to (numcols m0 -1) do (
      m0_(j,k)=m0_(j,k)+dotprod(bb_cont .. bb_(cont+n), LVPn);
      cont = cont + (n+1);
    );
 );
 Lpm = append(Lpm, matrix(m0));
);

--Redução dos menores 2x2 de cada matriz pertubada módulo quadrados (bb^2)

bb2 = trim bb^2;

Lpm2 = for i in Lpm list sub(sub(minors(2,i),r1/bb2),r1);

--pertubação dos pontos fixos(redes de quádricas) por quádricas
--( n1 = binomial(2+n,n) - 3)

PLPF = for i in LPF list
 ideal (for j from 0 to (numgens i -1) list 
    i_j + dotprod(a_(j+1,1)..a_(j+1,n1), omit(IVPn^2,i_*)));

--Novo anel para agrupar os coeficientes em relação às variáveis x_i
s = QQ[b_1..b_n2][x_0..x_n];

LPFs = for i in LPF list sub(i,s);

Lpm3 = {};
for i to (length(Lpm2) -1) do(
    u = Lpm2_i;
    q0s = LPFs_i;
    for j from 0 to 2 do
        uaux_j = sub(u_j, s);
    for j from 0 to 2 do(
        iv = omit(0..2,j);
        cq = coefficient(q0s_j, uaux_j);
        uaux_j = sub(sub(uaux_j*(2-cq),r1/bb2),s);
        for k from 0 to 1 do(
            l = iv_k;
            cq = coefficient(q0s_j,uaux_l);
            uaux_l = sub(sub(uaux_l - cq*uaux_j,r1/bb2),s);
        );
     );
     Lpm3 = append(Lpm3,ideal(uaux_0,uaux_1,uaux_2));
);

--Teste de limpeza (opcional)
      u = Lpm3_0; q0s = LPFs_0;
      for i from 0 to 2 do
        for j from 0 to 2 do
           print(i,j, coefficient(q0s_j, u_i));

Lu1 = {};
for i to (length(PLPF)-1) do(
   Lu1 = append (Lu1, ideal(gens(sub(PLPF_i,r1)) - 
         gens(sub(Lpm3_i,r1))));
);

--Novo anel para agrupar os coeficientes de u1 em relação às 
--variáveis x_i
s2 = QQ[b_1..b_n2, a_(1,1)..a_(3,n1)][x_0..x_n];

Lu1 = for i in Lu1 list sub(i,s2);

--Coletando os coeficientes e agrupando em um mesmo ideal
Lu1c = {};
for i in Lu1 do (
    (M1,C1) = coefficients (i_0);
    (M2,C2) = coefficients (i_1);
    (M3,C3) = coefficients (i_2);
    Lu1c = append(Lu1c, ideal(C1,C2,C3));
);

--Anel sem as variáveis x_i
s3 = QQ[b_1..b_n2, a_(1,1)..a_(3,n1)];

Lu1c = for i in Lu1c list sub(i, s3);

--Eliminar as variáveis b_i
v = toList(b_1..b_n2);
Lu1c = for i in Lu1c list eliminate(v,i);

--Teste do nº de equações para o espaço tangente do detlocus
for i in Lu1c list numgens i 

--Cálculo dos tangentes: tgs
--Obs.: -wt é devido ao conormal

s4 = QQ[a_(1,1)..a_(3,n1), x_0..x_n];

LPFT = {}; --Lista dos pontos fixos e respectivos tangentes
for i to (length(PLPF) - 1) do(
  remember = symbol remember;
  wt(a_(1,1), PLPF_i, s4);
  tg = tGrass(S(2),sub(sum((LPF_i)_*),r)); 
  tgs = tg;
  for j in ((Lu1c_i)_*) do(
    tgs = tgs - sub(wt2frac(-wt(j,s4)),frac(r)););
  LPFT = append(LPFT, {"P1: pt fixo, P2: Tangente",LPF_i,tgs});
);

--Lista com os pontos fixos (1 a 5) e respectivos tangentes
LPFT
-------------------------------------------------------------------------
\end{verbatim}

\tocless\section{Códigos do Maple\label{codigo_Maple_P4_P5}}
Nesta seção tomamos os pontos fixos (1 a 5) e seus respectivos tangentes obtidos com o Macaulay2 na seção \ref{codigo-macaulay2_tangente_P4} (vide (\ref{tg_1_rc}) - (\ref{tg_f_rc})) e efetuamos as permutações das variáveis para gerar todos os pontos fixos e tangentes, as explosões necessárias e a substituição dos pesos para o cálculo do $\deg \Sigma(W_{rc}{},d)$ utilizando a fórmula de resíduos de Bott. Para a execução dos cálculos específicos é necessário a inclusão dos procedimentos da seção \ref{funcoes_Maple}.

\vspace{0.5cm}
\noindent\verb+#########################################################################+

\noindent\verb+# Inserir aqui os procedimentos da seção+ \ref{funcoes_Maple}
\begin{verbatim}
#########################################################################
# Os cálculos iniciam aqui

N := 4: M := 60: 
variables := [seq(x[i], i = 0 .. N), seq(1/x[i], i = 0 .. N)]:
DIM := 18: 
listweights := [11, 17, 32, 55, 95]:

# Pontos fixos e respectivos tangentes a menos de permutação 
# das variáveis

PF[1] := x[0]*x[1]+x[1]*x[2]+x[2]*x[3]:
PF[2] := x[0]*x[1]+x[0]*x[2]+x[1]*x[2]:
PF[3] := x[0]*x[1]+x[0]*x[2]+x[2]^2:
PF[4] := x[0]^2+x[0]*x[1]+x[1]^2:
PF[5] := x[0]^2+x[0]*x[1]+x[0]*x[2]:

TG[1] := expand((x[0]^3*x[1]+x[0]^2*x[1]*x[2]+x[0]^2*x[1]*x[3]
+x[0]^2*x[1]*x[4]+x[0]^2*x[2]*x[3]+x[0]*x[1]^2*x[3]
+x[0]*x[1]*x[2]^2+x[0]*x[1]*x[2]*x[4]+x[0]*x[1]*x[3]^2
+x[0]*x[1]*x[3]*x[4]+x[0]*x[2]^2*x[3]+x[0]*x[2]*x[3]^2
+x[0]*x[2]*x[3]*x[4]+x[1]^2*x[2]*x[3]+x[1]*x[2]*x[3]^2
+x[1]*x[2]*x[3]*x[4]+x[2]*x[3]^3+x[2]*x[3]^2*x[4])
/(x[0]*x[1]*x[2]*x[3])):

TG[2] := expand((x[0]^2*x[1]+x[0]^2*x[2]+x[0]*x[1]^2
+2*x[0]*x[1]*x[3]+2*x[0]*x[1]*x[4]+x[0]*x[2]^2
+2*x[0]*x[2]*x[3]+2*x[0]*x[2]*x[4]+x[1]^2*x[2]
+x[1]*x[2]^2+2*x[1]*x[2]*x[3]+2*x[1]*x[2]*x[4])
/(x[0]*x[1]*x[2])):

TG[3] := expand((x[0]^2*x[1]*x[2]+x[0]^2*x[2]^2+x[0]*x[1]^3
+x[0]*x[1]^2*x[2]+x[0]*x[1]^2*x[3]+x[0]*x[1]^2*x[4]
+2*x[0]*x[1]*x[2]*x[3]+2*x[0]*x[1]*x[2]*x[4]
+x[0]*x[2]^2*x[3]+x[0]*x[2]^2*x[4]+x[1]^2*x[2]^2+x[1]*x[2]^3
+x[1]*x[2]^2*x[3]+x[1]*x[2]^2*x[4]+x[2]^3*x[3]+x[2]^3*x[4])
/(x[0]*x[1]*x[2]^2)):

TG[4] := expand((x[0]^3*x[2]+x[0]^3*x[3]+x[0]^3*x[4]
+2*x[0]^2*x[1]*x[2]+2*x[0]^2*x[1]*x[3]+2*x[0]^2*x[1]*x[4]
+2*x[0]*x[1]^2*x[2]+2*x[0]*x[1]^2*x[3]+2*x[0]*x[1]^2*x[4]
+x[1]^3*x[2]+x[1]^3*x[3]+x[1]^3*x[4])/(x[0]^2*x[1]^2)):

TG[5] := expand((x[0]*x[1]*x[3]+x[0]*x[1]*x[4]+x[0]*x[2]*x[3]
+x[0]*x[2]*x[4]+x[1]^3+2*x[1]^2*x[2]+x[1]^2*x[3]+x[1]^2*x[4]
+2*x[1]*x[2]^2+2*x[1]*x[2]*x[3]+2*x[1]*x[2]*x[4]+x[2]^3
+x[2]^2*x[3]+x[2]^2*x[4])/(x[0]*x[1]*x[2])):

# Lista com todos os pontos fixos e respectivos tangentes 1 - 4
LPFT := []:
for i to 4 do 
  LPFT:=[op(LPFT), ["P2: pt fixo, P3: Tangente",PF[i],TG[i]]]: 
od:

# Primeira explosão (Pontos tipo 5)

# Tangente do centro de explosão
tgc := TGrass(S(1),x[0])+TGrass(x[1]+x[2]+x[3]+x[4],x[1]+x[2]):

# Tangente do ambiente
tgs := TG[5]:

LBPF5 := blowUp(tgs, tgc): 

# Lista de pontos fixos a menos de permutação das variáveis
# Pontos tipo 5 após explosão (LPFT5)
LP := [5, 6, 8, 9]:
LPFT5 := []:
for i in LP do 
  LPFT5 := [op(LPFT5), ["P2: pt fixo, P3: Tangente", 
  expand(PF[5]+LBPF5[i][3]*denom(LBPF5[i][2])), LBPF5[i][4]]]:
od:

#Junção das listas de ponto fixo
LPFT := [op(LPFT), op(LPFT5)]:

# Efetuando a permutação das variáveis para gerar todos os 
# pontos fixos

with(combinat): with(ListTools):

p1 := permute([x[0], x[1], x[2], x[3], x[4]]):
LF := []:
for i to nops(LPFT) do 
  LF := [op(LF), [LPFT[i][1],
      unapply(LPFT[i][2],x[0],x[1],x[2],x[3],x[4]), 
      unapply(LPFT[i][3], x[0], x[1], x[2], x[3], x[4])]]:
od:

LPT := []:
for i to nops(LPFT) do 
  Laux := []:
  for j to nops(p1) do 
    Laux := [op(Laux), [LF[1][1], LF[i][2](op(p1[j])), 
             LF[i][3](op(p1[j]))]]:
  od: 
  LPT := [op(LPT), [MakeUnique(Laux)]]:
od:

# (opcional) Total de pontos fixos antes da 2ª explosão #440
npf := 0:
for i to nops(LPT) do 
  npf := npf+nops(LPT[i][1]):
od:
print(npf):

#Gerando as fibras e tangentes

# Caso d=4
# Casos i = 1 a 3
for i to 3 do 
  LFT[i] := []:
  for j to nops(LPT[i][1]) do 
    LFT[i] := 
    [op(LFT[i]), ["P2: pt fixo, P3: fibra, P4: Tangente", 
    LPT[i][1][j][2], expand(simpleweights(LPT[i][1][j][2]^2)),
    LPT[i][1][j][3]]]:
 od:
od:

# Caso i = 4 (2ª explosão)
LFT[4] := []:
for j in LPT[4][1] do 
   tgc4 := TGrass(S(1), add(indets(j[2]))):
   tgs4 := j[3]: 
   LB4 := blowUp(tgs4, tgc4):
   for k in LB4 do 
     LFT[4] := 
     [op(LFT[4]),["P2: pt fixo, P3: fibra, P4: Tangente",j[2],
     expand(simpleweights(j[2]^2+k[3]*denom(k[2]))), k[4]]]:
   od:
od:

#Casos i= 5 a 8
for i from 5 to 8 do 
  LFT[i] := []:
  for j to nops(LPT[i][1]) do 
    LFT[i] := 
    [op(LFT[i]), ["P2: pt fixo, P3: fibra, P4: Tangente", 
    LPT[i][1][j][2], 
    partdegree(expand(simpleweights(LPT[i][1][j][2]^2)), 4),
    LPT[i][1][j][3]]]:
  od:
od:

#Grau para d=4 
LFTatual := []:
for i to 8 do 
   for j in LFT[i] do 
     LFTatual := [op(LFTatual), [j[1], j[2], j[3], j[4]]]:
   od:
od:
listresp := []:
listresp := [op(listresp), Bott(LFTatual, 4)]:
print(listresp):
save listresp, "Grau_Maple_P4.txt":

#Casos d=5 e 6
for d from 5 to 6 do 
  for i to 4 do 
    LFTd[i] := []:
    for j in LFT[i] do 
      LFTd[i] := [op(LFTd[i]), 
      ["P2: pt fixo, P3: fibra, P4: Tangente", j[2], 
      simpleweights(j[3]*S(1)), j[4]]]:
    od: 
    LFT[i] := LFTd[i]:
  od: 
  for i from 5 to 8 do 
    LFTd[i] := []: 
    for j in LFT[i] do 
     LFTd[i] := 
     [op(LFTd[i]),["P2: pt fixo, P3: fibra, P4: Tangente",j[2],
     simpleweights(j[3]*S(1)+partdegree(j[2]^2, d)),j[4]]]:
    od: 
    LFT[i] := LFTd[i]:
  od:
  # Grau para d = 5 e 6
  LFTatual := []:
  for i to 8 do 
    for j in LFT[i] do 
      LFTatual := [op(LFTatual), [j[1], j[2], j[3], j[4]]]:
    od:
  od:
  listresp := [op(listresp), Bott(LFTatual, d)]:
  print(listresp):
  save listresp, "Grau_Maple_P4.txt":
od:

# Casos d>= 7 (>regularidade = 6)
# Obs.: Não tentar rodar muitos casos de cada vez, 
# pois ocupa uma quantidade razoável de memória
for d from 7 to 10 do 
   for i to 8 do 
      LFTd[i] := []:
      for j in LFT[i] do 
        LFTd[i] := [op(LFTd[i]), [j[1], j[2], 
        simpleweights(j[3]*S(d-6)), j[4]]]:
      od:
    od:
    LFTatual := []:
    for i to 8 do 
      for j in LFTd[i] do 
        LFTatual := [op(LFTatual), [j[1], j[2], j[3], j[4]]]:
      od:
    od: 
    listresp := [op(listresp), Bott(LFTatual, d)]:
    print(listresp):
    save listresp, "Grau_Maple_P4.txt":
od:
#########################################################################
\end{verbatim}
\chapter{Código: Cálculo do \texorpdfstring{$\deg\Sigma(\W_{sg}{},d)$}{}\label{ap_segre}}
Inicialmente, utilizamos o Macaulay2 para o cálculo dos tangentes dos pontos fixos antes da primeira explosão (vide \ref{diagrama_blowup_segre}), onde utilizamos o código da seção \ref{codigo-macaulay2_tangente_P4}, trocando $n=4$ por $n=5$. Posteriormente, utilizamos o Maple para gerar as explosões necessárias, bem como todos os pontos fixos/tangentes a partir dos pontos fixos/tangentes iniciais e fazer as substituições de pesos para o cálculo dos graus utilizando a fórmula de resíduos de Bott. Para isso, é necessário os procedimentos da seção \ref{funcoes_Maple} e os códigos listados a seguir.

\vspace{0.5cm}
\noindent\verb+#########################################################################+

\noindent\verb+# Inserir aqui os procedimentos da seção +\ref{funcoes_Maple}
\begin{verbatim}
#########################################################################
# Os cálculos iniciam aqui

N := 5: M := 19: #M precisa ir pelo menos até (5)*24 + 6 = 126 
variables := [seq(x[i], i = 0 .. N), seq(1/x[i], i = 0 .. N)]:
DIM := 24: 
listweights := [11, 17, 32, 55, 95,160]:

# Pontos fixos e respectivos tangentes a menos de permutação 
# das variáveis
PF[1] := x[0]*x[1]+x[1]*x[2]+x[2]*x[3]:
PF[2] := x[0]*x[1]+x[0]*x[2]+x[1]*x[2]:
PF[3] := x[0]*x[1]+x[0]*x[2]+x[2]^2:
PF[4] := x[0]^2+x[0]*x[1]+x[1]^2:
PF[5] := x[0]^2+x[0]*x[1]+x[0]*x[2]:

TG[1] := expand((x[0]^3*x[1]+x[0]^2*x[1]*x[2]+x[0]^2*x[1]*x[3]
+x[0]^2*x[1]*x[4]+x[0]^2*x[1]*x[5]+x[0]^2*x[2]*x[3]+
x[0]*x[1]^2*x[3]+x[0]*x[1]*x[2]^2+x[0]*x[1]*x[2]*x[4]
+x[0]*x[1]*x[2]*x[5]+x[0]*x[1]*x[3]^2+x[0]*x[1]*x[3]*x[4]+
x[0]*x[1]*x[3]*x[5]+x[0]*x[2]^2*x[3]+x[0]*x[2]*x[3]^2+
x[0]*x[2]*x[3]*x[4]+x[0]*x[2]*x[3]*x[5]+x[1]^2*x[2]*x[3]+
x[1]*x[2]*x[3]^2+x[1]*x[2]*x[3]*x[4]+x[1]*x[2]*x[3]*x[5]+
x[2]*x[3]^3+x[2]*x[3]^2*x[4]+x[2]*x[3]^2*x[5])/
(x[0]*x[1]*x[2]*x[3])):

TG[2] := expand((x[0]^2*x[1]+x[0]^2*x[2]+x[0]*x[1]^2+
2*x[0]*x[1]*x[3]+2*x[0]*x[1]*x[4]+2*x[0]*x[1]*x[5]+
x[0]*x[2]^2+2*x[0]*x[2]*x[3]+2*x[0]*x[2]*x[4]+2*x[0]*x[2]*x[5]
+x[1]^2*x[2]+x[1]*x[2]^2+2*x[1]*x[2]*x[3]+2*x[1]*x[2]*x[4]+
2*x[1]*x[2]*x[5])/(x[0]*x[1]*x[2])):

TG[3] := expand((x[0]^2*x[1]*x[2]+x[0]^2*x[2]^2+x[0]*x[1]^3+
x[0]*x[1]^2*x[2]+x[0]*x[1]^2*x[3]+x[0]*x[1]^2*x[4+
x[0]*x[1]^2*x[5]+2*x[0]*x[1]*x[2]*x[3]+2*x[0]*x[1]*x[2]*x[4]
+2*x[0]*x[1]*x[2]*x[5]+x[0]*x[2]^2*x[3]+x[0]*x[2]^2*x[4]
+x[0]*x[2]^2*x[5]+x[1]^2*x[2]^2+x[1]*x[2]^3+x[1]*x[2]^2*x[3]
+x[1]*x[2]^2*x[4]+x[1]*x[2]^2*x[5]+x[2]^3*x[3]+x[2]^3*x[4]
+x[2]^3*x[5])/(x[0]*x[1]*x[2]^2)):

TG[4] := expand((x[0]^3*x[2]+x[0]^3*x[3]+x[0]^3*x[4]
+x[0]^3*x[5]+2*x[0]^2*x[1]*x[2]+2*x[0]^2*x[1]*x[3]
+2*x[0]^2*x[1]*x[4]+2*x[0]^2*x[1]*x[5]+2*x[0]*x[1]^2*x[2]
+2*x[0]*x[1]^2*x[3]+2*x[0]*x[1]^2*x[4]+2*x[0]*x[1]^2*x[5]
+x[1]^3*x[2]+x[1]^3*x[3]+x[1]^3*x[4]+x[1]^3*x[5])/
(x[0]^2*x[1]^2)):

TG[5] := expand((x[0]*x[1]*x[3]+x[0]*x[1]*x[4]+x[0]*x[1]*x[5]
+x[0]*x[2]*x[3]+x[0]*x[2]*x[4]+x[0]*x[2]*x[5]+x[1]^3
+2*x[1]^2*x[2]+x[1]^2*x[3]+x[1]^2*x[4]+x[1]^2*x[5]
+2*x[1]*x[2]^2+2*x[1]*x[2]*x[3]+2*x[1]*x[2]*x[4]+
2*x[1]*x[2]*x[5]+x[2]^3+x[2]^2*x[3]+x[2]^2*x[4]+
x[2]^2*x[5])/(x[0]*x[1]*x[2])):

# Lista com todos os pontos fixos e respectivos tangentes 1-4
LPFT := []:
for i to 4 do LPFT := [op(LPFT), ["P2: pt fixo, P3: Tangente", 
 PF[i], TG[i]]]: od:

# Primeira explosão (Pontos tipo 5)
# Tangente do centro de explosão
tgc := TGrass(S(1),x[0])+TGrass(x[1]+x[2]+x[3]+x[4]+x[5],
x[1]+x[2]):

# Tangente do ambiente
tgs := TG[5]:

LBPF5 := blowUp(tgs, tgc): 

# Lista de pontos fixos a menos de permutação das variáveis
# Pontos tipo 5 após explosão (LPFT5)
# LP Lista com as posições dos pontos fixos após explosão 
# que geram os demais pontos a menos de permutação
LP := []: 
for i to nops(LBPF5) do 
  if (LBPF5[i][3] = x[2]/x[0]) then LP := [op(LP), i]: fi:
od:
for i to nops(LBPF5) do 
  if (LBPF5[i][3]=x[2]^2/(x[0]*x[1])) then LP:=[op(LP), i]:fi:
od:
for i to nops(LBPF5) do 
  if (LBPF5[i][3] = x[3]/x[0]) then LP := [op(LP), i]:fi:
od:
for i to nops(LBPF5) do 
  if(LBPF5[i][3]=x[2]*x[3]/(x[0]*x[1])) then LP:=[op(LP),i]:fi:
od:

LPFT5 := []:
for i in LP do 
  LPFT5 := [op(LPFT5), ["P2: pt fixo, P3: Tangente", 
  expand(PF[5]+LBPF5[i][3]*denom(LBPF5[i][2])), LBPF5[i][4]]]:
od:

#Junção das listas de ponto fixo
LPFT := [op(LPFT), op(LPFT5)]:

# Efetuando a permutação das variáveis para gerar todos os pontos fixos

with(combinat): with(ListTools):

p1 := permute([x[0], x[1], x[2], x[3], x[4],x[5]]):
LF := []:
for i to nops(LPFT) do 
  LF := [op(LF), [LPFT[i][1], 
  unapply(LPFT[i][2], x[0], x[1], x[2], x[3], x[4],x[5]), 
  unapply(LPFT[i][3], x[0], x[1], x[2], x[3], x[4],x[5])]]:
od:

LPT := []:
for i to nops(LPFT) do 
  Laux := []:
  for j to nops(p1) do 
    Laux := [op(Laux), [LF[1][1], LF[i][2](op(p1[j])), 
     LF[i][3](op(p1[j]))]]:
  od: 
  LPT := [op(LPT), [MakeUnique(Laux)]]:
od:

# (opcional) Total de pontos fixos antes da 2ª explosão #1115
npf := 0:
for i to nops(LPT) do 
  npf := npf+nops(LPT[i][1]):
od:
print(npf):

#Gerando as fibras e tangentes

# Caso d=4
# Casos i = 1 a 3
for i to 3 do 
  LFT[i] := []:
  for j to nops(LPT[i][1]) do 
    LFT[i] := [op(LFT[i]), 
    ["P2: pt fixo, P3: fibra, P4: Tangente", LPT[i][1][j][2], 
  expand(simpleweights(LPT[i][1][j][2]^2)), LPT[i][1][j][3]]]:
 od:
od:

# Caso i = 4 (2ª explosão)
LFT[4] := []:
for j in LPT[4][1] do 
   tgc4 := TGrass(S(1), add(indets(j[2]))):
   tgs4 := j[3]: 
   LB4 := blowUp(tgs4, tgc4):
   for k in LB4 do 
     LFT[4] := [op(LFT[4]), 
     ["P2: pt fixo, P3: fibra, P4: Tangente", j[2], 
     expand(simpleweights(j[2]^2+k[3]*denom(k[2]))), k[4]]]:
   od:
od:

#Casos i= 5 a 8
for i from 5 to 8 do 
  LFT[i] := []:
  for j to nops(LPT[i][1]) do 
    LFT[i] := [op(LFT[i]), 
    ["P2: pt fixo, P3: fibra, P4: Tangente", LPT[i][1][j][2], 
    partdegree(expand(simpleweights(LPT[i][1][j][2]^2)), 4), 
    LPT[i][1][j][3]]]:
  od:
od:

#Grau para d=4 
LFTatual := []:
for i to 8 do 
   for j in LFT[i] do 
     LFTatual := [op(LFTatual), [j[1], j[2], j[3], j[4]]]:
   od:
od:
listresp := []:
listresp := [op(listresp), Bott(LFTatual, 4)]:
print(listresp):
save listresp, "Grau_Maple_P5.txt":


#Casos d= 5 e 6 - ao final sai com as fibras e tangentes para d = 6
for d from 5 to 6 do 
  for i to 4 do 
    LFTd[i] := []:
    for j in LFT[i] do 
      LFTd[i] := [op(LFTd[i]), 
      ["P2: pt fixo, P3: fibra, P4: Tangente", j[2], 
      simpleweights(j[3]*S(1)), j[4]]]:
    od: 
    LFT[i] := LFTd[i]:
  od: 
  for i from 5 to 8 do 
    LFTd[i] := []: 
    for j in LFT[i] do 
      LFTd[i] := [op(LFTd[i]), 
      ["P2: pt fixo, P3: fibra, P4: Tangente", j[2], 
      simpleweights(j[3]*S(1)+partdegree(j[2]^2, d)), j[4]]]:
    od: 
    LFT[i] := LFTd[i]:
  od:
  # Grau para d = 5 e 6
  LFTatual := []:
  for i to 8 do 
    for j in LFT[i] do 
      LFTatual := [op(LFTatual), [j[1], j[2], j[3], j[4]]]:
    od:
  od:
  listresp := [op(listresp), Bott(LFTatual, d)]:
  print(listresp):
  save listresp, "Grau_Maple_P5.txt":
od:

# Casos d>= 7 (>regularidade = 6)
# Obs.: Não tentar rodar muitos casos de cada vez, pois ocupa uma 
# quantidade razoável de memória
for d from 19 to 19 do 
   for i to 8 do 
      LFTd[i] := []:
      for j in LFT[i] do 
        LFTd[i] := [op(LFTd[i]), [j[1], j[2], 
        simpleweights(j[3]*S(d-6)), j[4]]]:
      od:
    od:
    LFTatual := []:
    for i to 8 do 
      for j in LFTd[i] do 
        LFTatual := [op(LFTatual), [j[1], j[2], j[3], j[4]]]:
      od:
    od: 
    listresp := [op(listresp), Bott(LFTatual, d)]:
    print(listresp):
    save listresp, "Grau_Maple_P5.txt":
od:
#########################################################################
\end{verbatim}  

\chapter{Código: Cálculo do \texorpdfstring{$\deg\Sigma(\W_{eqc}{},d)$}{}\label{ap_quartica_eliptica}}

Neste apêndice disponibilizamos o código para o cálculo do grau da subvariedade $\Sigma(W_{eqc}{},d)\subset \mathbb{P}^{N_d}$, superfícies em $\P3$ de grau $d$ singulares ao longo de alguma curva quártica elíptica. Inicialmente, na Seção \ref{codigos_singular_eqc} utilizamos o Singular \cite{Singular} munido da caixa de ferramentas "myprocs"\,, procedimentos/funções criadas por Vainsencher (vide Seção \ref{codigos_Israel}), para as contas locais que revelam o local de indeterminação do mapa racional $\W_{eqc} \dashrightarrow \bb G (12,\cl F _5)$ definido por $\nu$ (vide Proposição \ref{prop_nu_w_eqc}). Ainda na Seção \ref{codigos_singular_eqc} efetuamos a explosão de $\W_{eqc}$ ao longo de $\bb M$ e verificamos que isto resolve a indeterminação. Posteriormente, na Seção \ref{codigos_macaulay2_eqc} utilizamos o software Macaulay2 \cite{Macaulay2} para gerar todos os pontos fixos, tangentes e as contribuições numéricas destes para o cálculo efetivo do $\deg\Sigma(\W_{eqc},d)$.
Além das funções da seção \ref{funcoes_Macaulay2} é necessário a inclusão das 
funções expand e blowUp listadas no Apêndice \ref{ap_twc}.

\tocless\section{Códigos do Singular\label{codigos_singular_eqc}}
\vspace{0.5cm}
\noindent\verb+//Inserir os códigos myprocs disponível na seção+  \ref{codigos_Israel}

\begin{verbatim}
LIB "primdec.lib";
proj(3);
ring r=0,(a(1..2)(1..8)),dp;r=r+P3;setring r;imapall(P3);
def xx2,xx3,xx4,xx5,xx6=xx^2,xx^3,xx^4,xx^5,xx^6 ;
def aa = ideal(a(1..2)(1..8));
def qe = ideal(x(0)^2,x(0)*x(1));

//Vizinhança afim padrão do feixe de quádricas <x(0)^2, x(0)*x(1)>
for(int i=1;i<=2;i++){
  qe[i]=qe[i]+dotprod([a(i)(1..8)],omit(xx2,ideal(x(0)^2,x(0)*x(1))));}

//Fibrado de cúbicas (posto genérico 8)
def cubs = qe*xx;
def m=transpose(coeffs(cubs,xx3,xxx));
ran(m);
def l=rowsnopivo(m);m=l[1];l[2];l[3];
def m0=submat(m,l[2],1..ncols(m));
m=submat(m,omit(1..nrows(m),l[2]),1..ncols(m));
nrows(m);
ran(origin(aa,m));//7
std(origin(aa,cubs)); def cubs0=_;

//Geradores da fibra divisor excepcional da primeira explosão
def im02 = reduce(ideal(m0),std(aa^2));
def b = pol2id(dotprod(im02,xx3)); b;
//b[1]=-a(1)(3)*x(1)^3
//b[2]=-a(1)(4)*x(1)^2*x(2)
//b[3]=-a(1)(5)*x(1)^2*x(3)
//b[4]=-a(1)(6)*x(1)*x(2)^2
//b[5]=-a(1)(7)*x(1)*x(2)*x(3)
//b[6]=-a(1)(8)*x(1)*x(3)^2
//b[7]=a(2)(6)*x(0)*x(2)^2
//b[8]=a(2)(7)*x(0)*x(2)*x(3)
//b[9]=a(2)(8)*x(0)*x(3)^2

//Coloca o ideal(m0) na forma localmente gráfico
putsolvform(ideal(m0));
def Y=_;

//Z = Local de indeterminação do mapa k: X --> G(8,F_3)
def Indet = dosubs(Y,qe);
list Z = primdecGTZ(Indet);
Z;

//Explosão de X ao longo de Z
ring r1=0,c(1..size(Y)-1),dp;
r1=r1+r;setring r1;imapall(r);


//Divisor excepcional de interesse Y[4] 
//(o que produz a órbita fechada - <x_0^2, x_0x_1, x_1^3>)
int i1=4; 

def exc1=Y[i1];def Y1=omit(Y,exc1);
Y1=seq("Y1[i]-exc1*c(i)",1,size(Y1));
def vs=aa;

l=qe,m,m0,vs;l=dosubs(Y1,l);

qe,m,m0,vs=l[1..size(l)];

mysat(m0,m);

nrows(m); 
vs=pol2id(indets(vs));vs;
ran(origin(vs,m));
cubs=ideal(m*transpose(xx3));

if(ran(origin(vs,m0))<>0){
cubs=cubs+ideal(m0*transpose(xx3));};
hilbsp(vs,cubs);

//seções linearmente independentes w_i
def w = reduce(cubs,std(vs^2)); w;

//As 8 cúbicas
std(origin(vs,cubs)); def cubs1 = _; 

//A nova cúbica gerada após explosão
std(reduce(cubs1,std(cubs0)));

//Fibrado de quínticas (posto genérico = 12)
def quits = cubs*qe;
std(origin(vs,quits));  def quits0=_; //11 quínticas
hilbsp(vs,quits);
def m=transpose(coeffs(quits,xx5,xxx));
ran(m);
def l=rowsnopivo(m);m=l[1];l[2];l[3];
def m0=submat(m,l[2],1..ncols(m));
m=submat(m,omit(1..nrows(m),l[2]),1..ncols(m));
nrows(m);
ran(origin(vs,m)); ran(origin(vs,m0));
putsolvform(ideal(m0));
def W=_;

//M = (Local de indeterminação do mapa W_eqc --> G(12,F_5))
def T = ideal(cubs[8]);
def TW = dosubs(W,T);
M = myfactor(TW);

//Geradores da fibra do divisor excepcional
def im02 = reduce(ideal(m0),std(vs^2));
list L5 = xx5[(1..56)],xx5[(1..56)];
for(int i=1; i<=size(L5); i++){
  if(im02[i]<>0){
    see(im02[i],L5[i]);};};
    
//Explosão de W_{eqc} ao longo de M
ring r2=0,d(1..size(W)-1),dp;r2=r2+r1;setring r2;imapall(r1);
l=m0,m,vs; 
int i2;

//Variando as vizinhanças
i2++; 
def exc2=W[i2];def W1=omit(W,exc2);
W1=seq("W1[i]-exc2*d(i)",1,size(W1));
def l2=dosubs(W1,l);
m0,m,vs=l2[1..size(l2)];vs=pol2id(indets(vs));vs;
mysat(m0,m);
nrows(m);
quits=ideal(m*transpose(xx5)); ran(origin(vs,m0));

if(ran(origin(vs,m0))<>0){
quits=quits+ideal(m0*transpose(xx5));}

//As 12 quínticas
def quits1 =  std(origin(vs,quits)); quits1;

//Nova quíntica
std(reduce(quits1,std(quits0))); def QN =_;

hilbsp(vs,quits); //12t-16 (Polinômio de Hilbert correto)

quits1+x(1)^6; hilbsp(_);//12t-16
\end{verbatim}

\tocless\section{Códigos do Macaulay2\label{codigos_macaulay2_eqc}}

\vspace{0.5cm}
\noindent\verb+--Inserir as funções da seção+ \ref{funcoes_Macaulay2} \verb+aqui+
\vspace{0.5cm}

\noindent\verb+--Inserir as funções expand e blowUp do Apêndice+ \ref{ap_twc} \verb+aqui+
\begin{verbatim}
-------------------------------------------------------------------------
--Cálculos iniciam aqui
n = 3; --dimensão de Pn
r = QQ[x_0..x_n]; 
dualString = for i from 0 to n list(x_i => 1/x_i);
DIM = 16;--Dimensão da Variedade das curvas Quárticas Elípticas
m = 56; -- = 3*16+8. 
-- Obs.1: Na verdade é suficiente calcular \Sigma(W^2_eqc,d) 
-- para 33 (= 2*DIM +1) valores diferentes de d, o que implica
-- tomar m = 32 + 6 = 38

--Pontos fixos tipo 1 - 5
--Tipo 1-3 estão fora do 1º centro de explosão Z=P3^* x G(2,4)
PF_1 = ideal(x_0^2, x_1^2);
PF_2 = ideal(x_0^2, x_1*x_2);
PF_3 = ideal(x_0*x_1, x_2*x_3);

--Tipo 4-5 estão no primeiro centro de explosão
PF_4 = ideal(x_0^2, x_0*x_1);
PF_5 = ideal(x_0*x_1, x_0*x_2);

--Lista com os Pontos fixos/Tangentes do tipo 1 - 3
PFT = for i from 1 to 3 
list {"P1:pt fixo,P2:Tangente",PF_i,tGrass(S(2),sum((PF_i)_*))}; 

--1ª Explosão (pontos tipo 4 e 5)
for i from 4 to 5 do(
    plane = gcd(PF_i_*); 
    line = sub(sum(PF_i_*)/plane, ring plane);
    tg = tGrass(S(2),sum(PF_i_*)); 
    tgb = tGrass(S(1),plane)+tGrass(S(1),line); 
    LB1_i = blowUp(tg,tgb);
);

--pontos fixos/tangentes após 1ª explosão
for i from 4 to 5 do(
    PFT1_i = for  j in LB1_i 
      list {"P1:pt fixo,P2:Tangente", 
            PF_i+sub((j_2)*(denominator j_1), ring PF_i), j_3};
);

--Pontos fixos no 2º centro de Explosão
--h = 4t (Polinômio de Hilbert de uma curva quártica elíptica)
h = hilbertPolynomial(PF_1, Projective=>false); 
PFT2 = {}; --Lista com os pontos fixos/tangentes
for i from 4 to 5 do(
  for j in PFT1_i do(
      if (hilbertPolynomial(j_1,Projective=> false) != h) then(
        PFT2 = append(PFT2, j);)
      else(
        PFT = append(PFT, j););
  );
);

--2ª Explosão
for i to (length PFT2 -1) do(
  id2p = sum((PFT2_i_1)_*);
  q_1 = part(3,id2p);
  q_2 = part(2, id2p);
  L = gcd((PFT2_i_1)_*); 
  L1 = sum(support(q_2))-L;
  q0_1 = sub(q_1/L, ring L);
  q0_2 = sub(S(2), {L=>0, L1=>0});
  tg = PFT2_i_2; 
  tgb = tGrass(S(1),L)+tGrass(S(1)-L,L1)+tGrass(q0_2, q0_1);
  LB2_i = blowUp(tg,tgb);
);

--pontos fixos/tangentes após 2ª explosão
for i to (length PFT2 -1)  do(
   PFT3_i = for  j in LB2_i 
     list {"P1:pt fixo,P2:Tangente", PFT2_i_1+
       sub((j_2)*(denominator j_1), ring PFT2_i_1), j_3};
);

--Todos os pontos fixos/tangentes a menos de permutação
PFT = PFT|PFT3_0|PFT3_1|PFT3_2;

--Pontos fixos no 3º Centro de Explosão
--h2  = Polinômio de Hilbert do ideal(pt fixo)^2 = 12t - 16
h2 = hilbertPolynomial((PFT_0_1)^2, Projective=>false); 
PFTA ={}; 
PFTB = {}; --Pontos fixos no 3º centro de explosão
for i in PFT do(
  if (hilbertPolynomial((i_1)^2, Projective=>false) ==h2) then(
     PFTA = append(PFTA, i);)
  else(
     PFTB = append(PFTB,i));
);

--3ª Explosão (correção do polinômio de Hilbert de I^2)
for i to (length PFTB -1) do(
   id2p = sum((PFTB_i_1)_*);
   plane = gcd(terms part(2, id2p));
   line = sub(part(2,id2p)/plane,ring plane);
   pcl = line - plane;
   fs = S(1)-plane;
   SE1b = pcl^2*fs;
   u1 = part(3, id2p);
   tg = PFTB_i_2; 
   tgc = tGrass(S(1),plane)+tGrass(fs,pcl)+tGrass(SE1b,u1);
   LB3_i = blowUp(tg,tgc);
);

--Pontos fixos/tangentes em W'_eqc
--Coletando pontos após explosão
for i to (length PFTB -1) do(
  LF_i =  for j in LB3_i list {"P1: pt fixo, P2: Tangente", 
(PFTB_i_1)^2+sub((j_2)*(denominator j_1), ring PFTB_i_1), j_3};
);

LF1 = LF_0|LF_1|LF_2; -- Junção das listas

LF2 = for i in PFTA 
                list{"P1: pt fixo, P2: Tangente",(i_1)^2,i_2};

LF3 = LF1|LF2;

-- Transforma os ideais em um polinômio com as somas dos 
-- geradores (representação dos pts fixos e fibra)
LF3s = for i in LF3 list {i_0, sum(i_1_*), i_2};

--Aplicando as permutações para gerar todos os pts fixos 
--Lista de Permutações
ListP = permutations toList (x_0..x_n);
ListOP = for i in ListP 
        list {x_0 => i_0, x_1 => i_1, x_2 => i_2, x_3 => i_3};

for i to (length LF3s -1) do(
LFTF_i=unique for j in ListOP list {LF3s_i_0,sub(LF3s_i_1, j),
  sub(numerator LF3s_i_2, j)/sub(denominator LF3s_i_2, j)};
);

LPFT = {};
for i to (length LF3s - 1) do(
  LPFT = LPFT|LFTF_i;
);

LPFT = unique LPFT; --Lista com todos os pts fixos/tangentes
print("Total de pontos fixos = ",length LPFT); --813

--Aplicando a Fórmula de Resíduos de Bott

--obs.2: 
--para obtenção do polinômio gerador dos graus de 
--\Sigma(W_eqc,d) devemos tomar d>=6

--Obs. 3: 
--Para salvar os dados calculados no arquivo Grau_EQC.txt 
--criado/salvo na pasta /home/.../Singular_Surface_EQC
--utilizamos os comandos a seguir
--arquivo_saida = openOutAppend "/home/.../Grau_EQC.txt";
--arquivo_saida << (result) <<endl<<close;

-- Geração das fibras e graus para os casos 6 <= d <= 8
listresp = {};
for d from 6 to 8 do(

  LPFTF = for i in LPFT 
  list {"P1: pt fixo, P2: tangente, P3: fibra",i_1, i_2, 
sum((trim(ideal(terms(part(4,i_1)*S(d-4)+part(5,i_1)*S(d-5)+
part(6,i_1)*S(d-6)+part(7,i_1)*S(d-7)+
part(8,i_1)*S(d-8)))))_*)};
  result = bott(LPFTF,d); print(result);
  arquivo_saida = openOutAppend "/home/ufvjm/Dropbox/DOUTORADO/
  Contas_Macaulay_Singular_Maple/
  Singular_Surface_EQC/Grau_EQC.txt";
  arquivo_saida << (result) <<endl<<close;
  listresp = append(listresp,result);
);

--Casos d>8 (=regularidade)
--Obs.4:
--Se quiser calcular o grau para um valor específico
--sem rodar todos os casos, basta iniciar d com o valor 
--desejado e trocar m pelo mesmo valor no código abaixo.
--Por exemplo: para d = 20 -> for d from 20 to 20 do

for d from 9 to m do(
  LPFTFr = for i in LPFTF 
  list {"P1: pt fixo, P2: tangente, P3: fibra", i_1,i_2, 
  sum(((trim(ideal(terms (i_3 * S(d-8))))))_*)};
  result = bott(LPFTFr,d); print(result);
  listresp = append(listresp, result);
  arquivo_saida = openOutAppend "/home/ufvjm/Dropbox/
  DOUTORADO/Contas_Macaulay_Singular_Maple/
  Singular_Surface_EQC/Grau_EQC.txt";
  arquivo_saida << (result) <<endl<<close;
);

--Polinômio que dá o grau de \Sigma(W_eqc,d)
Grau = factor polynomialInterpolation(listresp, QQ[d])

--Salvando o polinômio
Grau = toString Grau
arquivo_saida = openOutAppend "/home/ufvjm/Dropbox/DOUTORADO/
Contas_Macaulay_Singular_Maple/Singular_Surface_EQC/
Grau_EQC.txt";
arquivo_saida << (result) <<endl<<close;
-------------------------------------------------------------------------
\end{verbatim}

\bibliography{mybib} 

\begin{thebibliography}{36}
\providecommand{\natexlab}[1]{#1}
\providecommand{\url}[1]{\texttt{#1}}
\expandafter\ifx\csname urlstyle\endcsname\relax
  \providecommand{\doi}[1]{doi: #1}\else
  \providecommand{\doi}{doi: \begingroup \urlstyle{rm}\Url}\fi

\bibitem[Altman and Kleiman(1977)]{Altman_Kleiman_1977}
A.~B. Altman and S.~L. Kleiman.
\newblock Foundations of the theory of {F}ano schemes.
\newblock \emph{Compositio Mathematica}, 34\penalty0 (1):\penalty0 3--47, 1977.
\newblock URL \url{http://www.numdam.org/item?id=CM_1977__34_1_3_0}.

\bibitem[Araujo(2009)]{Araujo_2009}
A.~L.~M. Araujo.
\newblock \emph{Aplica\c{c}\~oes da F\'ormula de {B}ott \`{a} Geometria
  Enumerativa}.
\newblock PhD thesis, Universidade Federal de Minas Gerais, 2009.
\newblock URL
  \url{http://www.mat.ufmg.br/intranet-atual/pgmat/TesesDissertacoes/uploaded/Tese019.pdf}.

\bibitem[Avritzer and Vainsencher(1999)]{Avritzer_Vainsencher_1999}
D.~Avritzer and I.~Vainsencher.
\newblock The hilbert scheme component of the intersection of two quadrics.
\newblock \emph{Communications in Algebra}, 27\penalty0 (6):\penalty0
  2995--3008, 1999.
\newblock \doi{10.1080/00927879908826606}.
\newblock URL \url{https://doi.org/10.1080/00927879908826606}.

\bibitem[Beauville(1996)]{Beauville_1996}
A.~Beauville.
\newblock \emph{Complex algebraic surfaces}.
\newblock London Mathematical Society student texts 34. Cambridge University
  Press, 2nd ed edition, 1996.
\newblock ISBN 9780521495103,0521495105,0521498422.

\bibitem[Coray and Vainsencher(1986)]{Coray1986}
D.~F. Coray and I.~Vainsencher.
\newblock Enumerative formulae for ruled cubic surfaces and rational quintic
  curves.
\newblock \emph{Commentarii Mathematici Helvetici}, 61\penalty0 (1):\penalty0
  501--518, Dec 1986.
\newblock ISSN 1420-8946.
\newblock \doi{10.1007/BF02621930}.
\newblock URL \url{https://doi.org/10.1007/BF02621930}.

\bibitem[Cuadrado(2010)]{Ferrer_2010}
V.~F. Cuadrado.
\newblock \emph{Enumerative Aspects of Holomorphic Foliations}.
\newblock PhD thesis, Universidade Federal de Minas Gerais, 2010.

\bibitem[Cukierman et~al.(October 2014)Cukierman, Lopez, and
  Vainsencher]{Cukierman_Lopez_Vainsencher_14}
F.~Cukierman, A.~Lopez, and I.~Vainsencher.
\newblock Enumeration of surfaces containing an elliptic quartic curve.
\newblock \emph{Proceedings of the American Mathematical Society}, 142\penalty0
  (10):\penalty0 3305--3313, October 2014.
\newblock URL \url{https://doi.org/10.1090/S0002-9939-2014-11998-8}.

\bibitem[{Dao} et~al.(2017){Dao}, {De Stefani}, {Grifo}, {Huneke}, and
  {N{\'u}{\~n}ez-Betancourt}]{2017arXiv170803010D}
H.~{Dao}, A.~{De Stefani}, E.~{Grifo}, C.~{Huneke}, and
  L.~{N{\'u}{\~n}ez-Betancourt}.
\newblock {Symbolic powers of ideals}.
\newblock \emph{ArXiv e-prints}, Aug. 2017.
\newblock URL \url{https://arxiv.org/abs/1708.03010}.

\bibitem[Decker et~al.(2018)Decker, Greuel, Pfister, and
  Sch\"onemann]{Singular}
W.~Decker, G.-M. Greuel, G.~Pfister, and H.~Sch\"onemann.
\newblock {\sc Singular} {4-1-1} --- {A} computer algebra system for polynomial
  computations.
\newblock \url{http://www.singular.uni-kl.de}, 2018.

\bibitem[Eisenbud(1995)]{Eisenbud_95}
D.~Eisenbud.
\newblock \emph{Commutative Algebra with a View Toward Algebraic Geometry},
  volume 150 of \emph{Graduate Texts in Mathematics}.
\newblock Springer New York, reprint of 1st edition edition, 1995.

\bibitem[Ellingsrud and Str{\o}mme(1995)]{Ellingsrud_Stromme_95}
G.~Ellingsrud and S.~Str{\o}mme.
\newblock The number of twisted cubic curves on the general quintic threefold.
\newblock \emph{Math. Scand.}, 76:\penalty0 5--34, 1995.
\newblock URL \url{http://www.jstor.org/stable/24491856}.

\bibitem[Ellingsrud and Str{\o}mme(1996)]{Ellingsrud_Stromme_Bott_96}
G.~Ellingsrud and S.~Str{\o}mme.
\newblock Bott's formula and enumerative geometry.
\newblock \emph{J. Amer. Math. Soc.}, 9\penalty0 (1):\penalty0 175--193, 1996.
\newblock URL \url{https://doi.org/10.1090/S0894-0347-96-00189-0}.

\bibitem[Ellingsrud and Stromme(Jul., 1989)]{Ellingsrud_Stromme_89}
G.~Ellingsrud and S.~A. Stromme.
\newblock On the {C}how ring of a geometric quotient.
\newblock \emph{Annals of Mathematics, Second Series}, 130\penalty0
  (1):\penalty0 159--187, Jul., 1989.
\newblock URL \url{http://www.jstor.org/stable/1971479}.

\bibitem[Ellingsrud et~al.(1987)Ellingsrud, Piene, and
  Str{\o}mme]{Ellingsrud_Piene_Stromme_87}
G.~Ellingsrud, R.~Piene, and S.~A. Str{\o}mme.
\newblock \emph{On the Variety of Nets of Quadrics Defining Twisted Cubic
  Curves, Space curves (Rocca Di Papa, 1985)}, volume 1266 of \emph{Lecture
  Notes in Mathematics}, pages 84--96.
\newblock Springer, 1987.
\newblock URL \url{https://doi.org/10.1007/BFb0078179}.

\bibitem[Fulton(1998)]{Fulton_1998}
W.~Fulton.
\newblock \emph{Intersection Theory}.
\newblock Springer, 2ª edition, 1998.

\bibitem[{Gotzmann}(2008)]{Gotzmann_2008}
G.~{Gotzmann}.
\newblock {The irreducible components of Hilb\^{}$\{$4n$\}$(P\^{}3)}.
\newblock \emph{ArXiv e-prints}, Nov. 2008.
\newblock URL \url{https://arxiv.org/abs/0811.3160}.

\bibitem[Grayson and Stillman()]{Macaulay2}
D.~Grayson and M.~E. Stillman.
\newblock Macaulay2 version 1.9.2, a software system for research in algebraic
  geometry.
\newblock available at https://faculty.math.illinois.edu/Macaulay2/.

\bibitem[Harris(1992)]{Harris_1992}
J.~Harris.
\newblock \emph{Algebraic {G}eometry: {A} {F}irst {C}ourse}.
\newblock Number 133 in Graduate Texts in Mathematics. Springer-Verlag, 1992.

\bibitem[Hartshorne(1977)]{Hartshorne_1977}
R.~Hartshorne.
\newblock \emph{Algebraic Geometry}.
\newblock Graduate Texts in Mathematics 52. Springer, 1977.

\bibitem[Hochster(1973)]{Hochster_73}
M.~Hochster.
\newblock Criteria for equality of ordinary and symbolic powers of primes.
\newblock \emph{Mathematische Zeitschrift | Mathematische Zeitschrift},
  133:\penalty0 53--65, 1973.
\newblock URL \url{http://eudml.org/doc/171940}.

\bibitem[Kazarian(2003)]{Kazarian_03}
M.~E. Kazarian.
\newblock Multisingularities, cobordisms, and enumerative geometry.
\newblock \emph{Russian Mathematical Surveys}, 58\penalty0 (4):\penalty0 665,
  2003.
\newblock URL \url{http://stacks.iop.org/0036-0279/58/i=4/a=R02}.

\bibitem[Kool et~al.(2011)Kool, Shende, and Thomas]{kool2011}
M.~Kool, V.~Shende, and R.~P. Thomas.
\newblock A short proof of the göttsche conjecture.
\newblock \emph{Geom. Topol.}, 15\penalty0 (1):\penalty0 397--406, 2011.
\newblock \doi{10.2140/gt.2011.15.397}.
\newblock URL \url{https://doi.org/10.2140/gt.2011.15.397}.

\bibitem[Liu(2000)]{liu2000}
A.-K. Liu.
\newblock Family blowup formula, admissible graphs and the enumeration of
  singular curves, i.
\newblock \emph{J. Differential Geom.}, 56\penalty0 (3):\penalty0 381--579, 11
  2000.
\newblock \doi{10.4310/jdg/1090347696}.
\newblock URL \url{https://doi.org/10.4310/jdg/1090347696}.

\bibitem[Maple()]{Maple_2015}
Maple.
\newblock Maplesoft, a division of {W}aterloo maple inc., {W}aterloo,
  {O}ntario.
\newblock Version 2015.

\bibitem[Meireles and Vainsencher(2001)]{Meireles_Vainsencher_2001}
A.~L. Meireles and I.~Vainsencher.
\newblock Equivariant intersection theory and {B}ott’s residue formula -
  {XVI} {E}scola de {Á}lgebra – {P}art 1.
\newblock \emph{Matemática Contemporânea}, 20:\penalty0 1--70, 2001.
\newblock URL \url{http://mc.sbm.org.br/docs/mc/pdf/20/a1.pdf}.

\bibitem[Meurer(1996)]{Meurer_96}
P.~Meurer.
\newblock The number of rational quartics on calabi-yau hypersurfaces in
  weighted projective space p(2,$1^4$).
\newblock \emph{Mathematica Scandinavica}, 78\penalty0 (1):\penalty0 63--83,
  1996.
\newblock URL \url{http://www.jstor.org/stable/24492817}.

\bibitem[Mumford(1966)]{Mumford_1966}
D.~Mumford.
\newblock \emph{Lectures on curves on an algebraic surface, With a section by
  G. M. Bergman}.
\newblock Number~59 in Annals of Mathematics Studies. Princeton University
  Press, Princeton New Jersey, 1966.

\bibitem[Piene and Schlessinger(Aug.,1985)]{Piene_Schlessinger_1985}
R.~Piene and M.~Schlessinger.
\newblock On the {H}ilbert scheme compactification of the space of twisted
  cubics.
\newblock \emph{American Journal of Mathematics}, 107\penalty0 (4):\penalty0
  761--774, Aug.,1985.
\newblock URL \url{http://www.jstor.org/stable/2374355}.

\bibitem[Rennemo(2017)]{Rennemo_2017}
J.~V. Rennemo.
\newblock Universal polynomials for tautological integrals on hilbert schemes.
\newblock \emph{Geometry $\&$ Topology}, 21\penalty0 (1):\penalty0 253 -- 314,
  2017.
\newblock URL \url{https://doi.org/10.2140/gt.2017.21.253}.

\bibitem[Tzeng(2012)]{tzeng2012}
Y.-J. Tzeng.
\newblock A proof of the göttsche-yau-zaslow formula.
\newblock \emph{J. Differential Geom.}, 90\penalty0 (3):\penalty0 439--472, 03
  2012.
\newblock \doi{10.4310/jdg/1335273391}.
\newblock URL \url{https://doi.org/10.4310/jdg/1335273391}.

\bibitem[Vainsencher()]{Script_Israel}
I.~Vainsencher.
\newblock Computer algebra scripts.
\newblock URL \url{http://www.mat.ufmg.br/~israel/Projetos/degNL/}.

\bibitem[Vainsencher(1987)]{Vainsencher_1987}
I.~Vainsencher.
\newblock A note on the {H}ilbert scheme of twisted cubics.
\newblock \emph{Bol. Soc. Bras. Mat}, 18\penalty0 (1):\penalty0 81--89, 1987.
\newblock URL \url{https://doi.org/10.1007/BF02584832}.

\bibitem[Vainsencher(2003)]{Vainsencher_2013}
I.~Vainsencher.
\newblock Hypersurfaces with up to six double points.
\newblock \emph{Communications in Algebra}, 31\penalty0 (8):\penalty0
  4107--4129, 2003.
\newblock URL \url{https;//doi.org/10.1081/AGB-120022456}.

\bibitem[Vainsencher(2015)]{Vainsencher_2015}
I.~Vainsencher.
\newblock Foliations singular along a curve.
\newblock \emph{Trans. London Math. Soc}, 2\penalty0 (1):\penalty0 80--92, July
  2015.
\newblock URL \url{https://doi.org/10.1112/tlms/tlv004}.

\bibitem[Vainsencher and Avritzer(1992)]{Vainsencher_Avritzer_92}
I.~Vainsencher and D.~Avritzer.
\newblock Compactifying the space of elliptic quartic curves.
\newblock \emph{Complex Projective Geometry(Trieste, 1989/Bergen, 1989), London
  Math. Soc. Lecture Note Ser.}, 179, Cambridge Univ. Press,
  Cambridge:\penalty0 47--58, 1992.
\newblock URL \url{https://doi.org/10.1017/CBO9780511662652.005}.

\bibitem[Vainsencher and Xavier(2002)]{Vainsencher_Xavier_02}
I.~Vainsencher and F.~Xavier.
\newblock A compactification of the space of twisted cubics.
\newblock \emph{Mathematica Scandinavica}, 91\penalty0 (2):\penalty0 221--243,
  2002.
\newblock URL \url{http://dx.doi.org/10.7146/math.scand.a-14387}.

\end{thebibliography}

\label{LastPage}

\end{document}